%% file: DISSER.TEX
\documentclass[12pt]{report} %{amsart}

\usepackage{amssymb}
\usepackage{epsf}
\usepackage{amsmath,amsfonts}

\usepackage[cp866]{inputenc}
\usepackage[russian]{babel}

\def\topfraction{1}
\def\bottomfraction{1}
\begin{document}

\normalsize
\setcounter{page}{0}
\input{DEF.TEX}
%\input{title.tex}
\input{Title-eng.tex}
\tableofcontents %\newpage
\input{INTRO.TEX}%\newpage
\input{CHAPTER1.TEX} %\newpage
\input{CHAPTER2.TEX} %\newpage
 \input{CHAPTER3.TEX} %\newpage
 \input{CHAPTER4.TEX} %\newpage
 \input{CHAPTER5.TEX} %\newpage
 \input{CONCL.TEX}  %\newpage
\normalsize
    %\input{appendix.tex}
\input{Refsfull.tex}

\end{document}

%% file: DEF.TEX
\newcommand{\rb}[1]{\raisebox{1.5ex}[-1.5ex]{#1}}
\renewcommand{\topfraction}{1}
\renewcommand{\bottomfraction}{1}
\newcommand{\be}{\begin{equation}}
\newcommand{\ee}{\end{equation}}
\newcommand{\nn}{\nonumber}

\newtheorem{theorem}{Теорема}
\newtheorem{proposition}[theorem]{Предложение}
\newtheorem{lemma}[theorem]{Лемма}
\newtheorem{corollary}[theorem]{Следствие}

\newtheorem{definition}{Определение}
\newtheorem{conjecture}{Гипотеза}

\newtheorem{remark}{Замечание}

\newcommand{\pqinf}{(p;p)_{\infty}(q;q)_{\infty}}
\newcommand\fac[2]{\mbox{$\frac{#1}{#2}$}}
\newcommand\s{\sigma}
\newcommand{\N}{\mathbb N}
\newcommand{\Q}{\mathbb Q}
\renewcommand{\C}{\mathbb C}
\newcommand{\R}{\mathbb R}
\newcommand{\Z}{\mathbb Z}
\newcommand{\F}{\mathbb F}
\newcommand{\T}{\mathbb T}
\renewcommand{\L}{\mathbb L}
\renewcommand{\P}{\mathcal P}
\newcommand{\la}{\langle}
\newcommand{\ra}{\rangle}
\newcommand{\ba}{\begin{eqnarray}}
\newcommand{\ea}{\end{eqnarray}}
\newcommand{\baa}{\begin{eqnarray*}}
\newcommand{\eaa}{\end{eqnarray*}}
\newcommand{\bb}{\begin{thebibliography}}
\newcommand{\eb}{\end{thebibliography}}
\newcommand{\ci}[1]{\cite{#1}}
\newcommand{\bi}[1]{\bibitem{#1}}
\newcommand{\lab}[1]{\label{#1}}
\newcommand{\re}[1]{(\ref{#1})}
\def\stackreb#1#2{\ \mathrel{\mathop{#1}\limits_{#2}}}

\newcommand{\tb}   {{\ifmmode \tan\beta     \else $\tan\beta$      \fi}}
\def\bea{\begin{eqnarray}}
\def\eea{\end{eqnarray}}

\renewcommand{\chaptername}{Глава}
\renewcommand{\bibname}{Список литературы}
\renewcommand{\contentsname}{Оглавление} %{Содержание}
%\renewcommand{\figurename}{Рис.}
%\renewcommand{\tablename}{Таблица}

%% file: Title-eng.tex
%
% Титульная страница
%
\thispagestyle{empty}

\begin{center}
ОБЪЕДИНЕННЫЙ ИНСТИТУТ ЯДЕРНЫХ ИССЛЕДОВАНИЙ
\end{center}

%\vspace{2mm}

\hspace{9cm}  На правах рукописи

%\vspace{5mm}

\begin{center} \large
СПИРИДОНОВ Вячеслав Павлович
\end{center}

\vspace{2mm}

\begin{center} \Large
Эллиптические гипергеометрические функции
\end{center}

\vspace{2mm}

\begin{center} \large
Специальность: 01.01.03 --- математическая физика
\end{center}

% \vspace{5mm}

\begin{center} \large
Диссертация на соискание ученой степени \\[-3mm]
доктора физико-математических наук
\end{center}

%\vspace{5mm}
%\vfill

\begin{center} \large
Дубна, 2004 год
\end{center}

\centerline{---------------------------------------------------------------------------------------{} }

\begin{center}
JOINT INSTITUTE FOR NUCLEAR RESEARCH
\end{center}

%\vspace{5mm}

\hspace{9cm}  Copyright by the author

\vspace{2mm}

\begin{center} \large
SPIRIDONOV Vyacheslav Pavlovich
\end{center}

\vspace{2mm}

\begin{center} \Large
Elliptic hypergeometric functions%
\footnote{Submitted in September 2004 and defended in April 2005 at the Saint-Petersburg
branch of the V.A. Steklov mathematical institute of the Russian Academy of Sciences.
Noticed missprints and simple technical omissions have been corrected.}
\end{center}

\vspace{2mm}

\begin{center} \large
Speciality: 01.01.03 --- mathematical physics
\end{center}

%\vspace{2mm}

\begin{center} \large
Doctor of Physical and Mathematical Sciences dissertation \\[-3mm]
(Habilitation thesis)
\end{center}

%\vspace{.5cm}
%\vfill

\begin{center} \large
Dubna, 2004
\end{center}

\clearpage

%% file: INTRO.TEX
%
% Introduction
%
% % % % % % % % % % % % % % % % % % % % % % % % % % % % % %

%
% =========================================================
%
\chapter*{Введение}
\addcontentsline{toc}{chapter}{Введение}

Специальные функции играют важную прикладную роль в теоретической и
математической физике. Особенно интенсивно они используются в спектральных задачах
для дифференциальных и конечно-разностных операторов в гильбертовом
пространстве. С практической точки зрения такие задачи считаются
точно решенными если они сводятся к уравнениям, решаемым или в элементарных
функциях или в терминах классифицированных специальных функций.
При этом существенную роль играют теоретико-групповые методы, позволяющие
выделить системы для которых основные физические величины вычисляются
в замкнутом виде. Расширение круга таких универсальных моделей
является актуальной задачей математической физики.

Существует много справочников и учебников по специальным функциям, например,
\cite{abr-ste:handbook,aar:special,vil,erd:higher,grad,whi-wat:course}.
Однако, ни один из них не содержит списка формальных свойств, которыми
должна обладать  функция для того чтобы быть ``специальной".
Обычно обсуждается какой-либо класс функций специфического вида
(гипергеометрические, автоморфные, и т.д.). Р. Аски предложил называть
специальной любую функцию если она настолько полезна, что
получает какое-либо собственное имя. Другая существенная, но не столь
универсальная как полезность, характеристика таких функций основана
на их асимптотическом поведении. Конкретнее, для специальных функций
естественно ожидать, что известное локальное поведение функции
позволяет вывести асимптотику на бесконечности, то есть проблема пересвязки
асимптотик должна быть решаемой. Такой подход к специальным функциям
характерен для специалистов, работающих над функциями типа Пенлеве
и общими изомонодромными деформациями \cite{kit}.

По мнению автора, эти два определения опираются на второстепенные характеристики
специальных функций. Необходимо уже иметь в руках функции для того,
чтобы начать определять их свойства. Если иметь целью поиск новых
специальных функций, то тогда необходимо найти определение, предоставляющее
более широкий набор технических средств для работы. В этом отношении
необходимо подчеркнуть, что даже термин ``классические специальные функции"
оказывается не таким уж устойчивым. Например, к настоящему  времени
признано, что семейство классических ортогональных полиномов одной переменной
включает в себя не только полиномы Якоби и их упрощения, но и существенно более
сложные полиномы Аски-Вильсона, открытые всего два десятилетия назад
\cite{ask-wil:some}, и всю иерархию их предельных случаев \cite{koe-swa:askey}.

Теория групп и связанных с ними алгебр предоставляет достаточно богатый набор
средств для построения новых функций, но, к сожалению, теория их представлений
зачастую приводит к интерпретации функций, уже определенных каким-либо другим
способом. Тем не менее, подход через группы симметрий кажется центральным в
теории специальных функций. В частности, все основные ``старые" специальные
функции появляются из разделения переменных в очень простых (и, таким образом,
универсальных и полезных) уравнениях \cite{Mi3}. Исследования автора в этом
направлении были привязаны к следующему рабочему определению: {\it специальными
функциями являются функции, связанные с автомодельными редукциями цепочек
спектральных преобразований для линейных задач на собственные значения.}

Это определение связывает специальные функции с фиксированными точками
различных непрерывных и дискретных преобразований симметрии для линейных
спектральных проблем. Оно хорошо работает только для специальных
функций одной независимой переменной (которые могут зависеть от бесконечного
числа параметров) и даже для них оно не претендует на покрытие всех возможных
случаев. С одной стороны, это определение происходит из теории полностью
интегрируемых систем, для которых поиск автомодельных решений нелинейных
эволюционных уравнений является стандартной задачей  \cite{AS,sol:book}.
С другой стороны, многие примеры функций выводимых таким образом
показывают, что центральное место в этом подходе занимают
{\it соотношения сопряжения (или связи)}---линейные или нелинейные уравнения связывающие
специальные функции при различных значениях их параметров.

Схематически, данный эвристический подход к поиску таких ``спектральных"
специальных функций состоит из следующих шагов:

1. Необходимо взять линейную спектральную задачу, определяемую
дифференциальными, конечно-разностными или интегральными уравнениями.

2. Необходимо построить другое линейно независимое уравнение по
переменным, не обязательно входящим в первое уравнение, на пространстве
решений первого уравнения.

3. Необходимо разрешить условие совместности  взятых линейных уравнений
и вывести нелинейные соотношения для функций, входящих как
свободные коэффициенты в эти уравнения. Если (условно) второе уравнение
дифференциальное, то возникают непрерывные потоки описываемые уравнениями
типа Кортвега-де Вриза, Кадом\-це\-ва-Пет\-виа\-шви\-ли, Тоды и т.д. Если же второе
уравнение конечно-разностное, то возникают цепочки уравнений с дискретным
временем типа Тоды, Вольтерра и прочие, схожие по смыслу с цепочками
преобразований Дарбу, меняющими спектральные данные дискретным образом.

4. Необходимо проанализировать дискретные и непрерывные симметрии полученных
нелинейных уравнений с помощью Лиевских теоретико-групповых методов, то есть
найти непрерывные и дискретные преобразования, отображающие пространство
решений нелинейных уравнений на себя.

5. Необходимо построить {\it автомодельные} решения полученных нелинейных
уравнений, которые инвариантны относительно определенных допустимых преобразований
симметрии. В результате автомодельных редукций возникают конечные множества нелинейных
дифференциальных, дифференциально-разностных, двумерных разностных и т.д.
уравнений, решения которых определяют ``нелинейные" специальные функции.
Решения самих начальных линейных уравнений с коэффициентами, связанными
с указанными автомодельными функциями, определяют ``линейные" специальные
функции.

Последние два шага носят эвристический характер, поскольку, несмотря на
существенный прогресс в общей теории автомодельных редукций (см., например,
\cite{LW,Mae}), полностью регулярные методы построения автомодельных
решений еще не построены. Например, редукции, использованные в работах
\cite{SZ0,SZ1} для построения рекуррентных соотношений ассоциированных полиномов
Аски-Вильсона и в работе \cite{spi-zhe:spectral} при открытии новых биортогональных
рациональных функций, еще не нашли теоретико-групповой
интерпретации. Эти редукции описаны в первой главе настоящей диссертации.
Перечисление проявлений понятия автомодельности в различных математических структурах,
включая спектральные задачи, дано, например, в сборнике статей \cite{PS}.

Другим важным составляющим элементом теории специальных функций,
не указанным в приведенной схеме, является теория трансцендентности.
Известно, что функции Пенлеве трансцендентны над дифференциальными полями,
построенными с помощью конечного числа расширений Пикара-Вессио поля
рациональных функций. При решении дифференциальных (разностных или любых
других) уравнений необходимо в конце концов определить какому
дифференциальному (конечно-разностному) полю принадлежит полученное решение.
Например, эти решения могут принадлежать дифференциальному полю, над которым определено
начальное дифференциальное уравнение. До настоящего времени остается
открытой проблема интерпретации автомодельных решений цепочек спектральных
преобразований с точки зрения дифференциальной (или разностной) теории Галуа.

Спектральные задачи типа Штурма-Лиувилля имеют много приложений в физике.
Квантовая механика и теория солитонов во многом основаны на спектральном
анализе оператора Шредингера.
Метод факторизации был предложен Шредингером как удобный
формализм для нахождения спектров некоторых операторов в
квантовой механике \cite{Sc1}. Инфельд переформулировал
задачу поиска гамильтонианов, собственные значения которых легко
находятся с помощью этого метода, как проблему поиска
решений факторизационной цепочки \cite{Inf,IH}. Относительно недавно
была обнаружена глубокая связь между этим формализмом и
суперсимметрией и ряд исследований автора в этом направлении
был посвящен более детальному изучению этой связи \cite{RS,AIS,spi:deformed}.

Реализация факторизационных операторов дифференциальными операторами первого
порядка соответствует преобразованиям Дарбу для линейного дифференциального
уравнения второго порядка, рассматривавшихся еще в XIX веке.
В современной теории  интегрируемых
систем различные версии этого подхода фигурируют под названиями
преобразований Ла\-пла\-са, Дар\-бу, Бэк\-лун\-да, одевания, и т.д.
\cite{AS}. В теории ортогональных полиномов истоком
этого метода служит теория ядерных полиномов Крис\-тоф\-фе\-ля и т.д.
В теории специальных функций такие преобразования соответствуют
соотношениям сопряжения. Важные результаты при исследовании этого метода
были получены Бурхналлом и Чонди \cite{BC1}, включая некоторые элементы
операторного подхода. Строгий математический анализ некоторых аспектов
преобразований Дарбу и метода факторизации дан в исследованиях \cite{Crum,D,F,Kr,Schm}.
Отметим, что в этом формализме специальные функции возникают как функции,
связанные с автомодельными редукциями факторизационных цепочек согласно
приведенной выше схеме.

Основное содержание настоящей диссертации составляют результаты исследований,
проводившихся автором в течение последних пяти лет. В ней содержится описание
основных положений теории принципиально нового класса специальных функций
математической физики --- эллиптических гипергеометрических функций.
Впервые такие объекты возникли в рамках квантового метода обратной задачи
рассеяния \cite{bik,ks,stf,tf} в качестве эллиптических решений уравнения Янга-Бакстера
\cite{bax:partition,abf:eight,djkmo:exactly1,djkmo:exactly2,djmo:fusion},
которые, как это было
продемонстрировано Френкелем и Тураевым в \cite{fre-tur:elliptic}, выражаются
через эллиптическое обобщение специального $q$-гипергеометрического ряда $_{10}\varphi_9$.
В настоящей диссертации не рассматриваются соответствующие результаты,
а описывается независимый подход к этим функциям, зародившийся в рамках
указанной выше схемы генерирования специальных функций  \cite{spi-zhe:spectral}
и детально разработанный в работах автора. Обычные и $q$-гипергеометрические
ряды и интегралы нашли очень большое число приложений в различных физических теориях.
Поэтому изучение специальных функций, обобщающих их на качественно новый
уровень, представляет большой интерес как с чисто математической
точки зрения так и в перспективе практических применений. Диссертация
является в некотором смысле дополнительной к фундаментальной
работе Милна \cite{mil:infinite}, отражающей недавний прогресс в
классической теории тета-функций Якоби.

Основными целями диссертации являются:
1) построение общей теории рядов и интегралов гипергеометрического
типа, связанных c эллиптическими функциями и тета-функциями Якоби,
и классификация эллиптических бета-интегралов;
2) изучение семейства непрерывных биортогональных функций одной
переменной, выражающихся в виде произведения двух обрывающихся совершенно
уравновешенных $_{12}E_{11}$  эллиптических гипергеометрических
рядов со специальным выбором параметров и обобщающих полиномы
Аски-Вильсона \cite{ask-wil:some} и биортогональные рациональные функции
Рахмана \cite{rah:integral} и Вильсона \cite{wil:hypergeometric,wil:orthogonal}.

\underline{В первой главе} диссертации дается краткий обзор автомодельных
редукций нелинейных цепочек с дискретным временем, связанных
с различными спектральных задачами. Она начинается с описания метода
факторизации для уравнения Шредингера и автомодельных потенциалов, для которых
дискретный спектр гамильтониана состоит из конечного
числа геометрических прогрессий. Такие спектры генерируются специфическими
полиномиальными квантовыми алгебрами, включающими в себя алгебру
$q$-гармонического осциллятора и $q$-аналог $su(1,1)$ алгебры. Когерентные состояния
таких систем описываются дифференциальными уравнениями с запаздывающим
аргументом. Сами же потенциалы связаны с бесконечно-солитонными
системами специального типа, цепочками Изинга, случайными матрицами,
двумерным Кулоновским газом и т.д. Полученные результаты не связаны напрямую
с эллиптическими гипергеометрическими функциями, но такие приложения специальных
функций носят универсальный характер и, поэтому, они приведены для полноты описания.

Так же в этой главе кратко излагается модификация метода факторизации для
ко\-неч\-но-раз\-но\-ст\-но\-го
оператора Шредингера или трехчленного рекуррентного соотношения для
ор\-то\-го\-наль\-ных полиномов. Дискретные аналоги преобразований Дарбу
определяются преобразованием Кристоффеля от обычных к ядерным полиномам и
преобразованием Геронимуса, обратным к нему \cite{Ger1}. Обобщение
факторизационного подхода
на этот случай рассмотрено впервые в статье \cite{Mi1}. Та же самая техника была
открыта и в работах по численным расчетам собственных значений матриц.
Такие численные алгоритмы как $LR, QR, g$-алгоритм и т.д. представляют собой
различные модификации цепочек дискретных преобразований Дарбу.

Ключевые новые результаты первой главы изложены в параграфах 1.3--1.5.
В параграфе 1.3 выводится цепочка спектральных преобразований ($R_{II}$-цепочка)
для обобщенной задачи на собственные значения для двух матриц Якоби. В следующем
параграфе анализируются симметрии этой цепочки и выводится широкий класс ее
автомодельных решений, выражающихся через рациональные,
элементарные и эллиптические функции. Эти решения оказываются связанными с
совершенно уравновешенным 2-сбалансированным гипергеометрическим рядом
$_9F_8$, его $q$-гипергеометрическим аналогом $_{10}\varphi_9$
и эллиптическим гипергеометрическим рядом $_{12}E_{11}$.
В последнем параграфе выводится нелинейная цепочка спектральных
преобразований для обобщенной спектральной задачи, порождающей
симметричные $R_{II}$-полиномы.

\underline{Вторая глава} посвящена общей теории тета-гипергеометрических рядов.
Сначала вводятся формальные ряды гипергеометрического типа, для которых отношение
последовательных членов ряда есть мероморфная функция номера члена ряда, обладающая
квазипериодичностью характерной для тета-функций Якоби (т.е. она имеет экспоненциальные
множители квазипериодичности). Это определение представляет собой
обобщение старых идей Похгаммера и Хорна (см., например, обзор \cite{ggr:general})
на случай функций с двумя независимыми
квазипериодами. Ряды, ограниченные с одной стороны, обозначаются как $_rE_s$,
а двусторонние ряды как $_rG_s$. В определенном пределе, они сводятся к
хорошо известным $q$-гипергеометрическим рядам $_r\varphi_s$ и $_r\psi_s$,
соответственно.

Эллиптическими функциями называются мероморфные двояко-периодические функции.
Они играют фундаментальную роль в математике и одним из главных результатов
диссертации является введение понятия общих эллиптических гипергеометрических
рядов и интегралов.
Формальные ряды $\sum_n c_n$ называются эллиптическими гипергеометрическими
рядами, если отношение $h(n)=c_{n+1}/c_n$ равно ограничению некоторой
эллиптической функции $h(y),\, y\in \C,$ на дискретную
решетку $y\in \N, \Z$ или $\Z_N$. Это приводит к ряду ограничений
на параметры $_rE_s$ и $_rG_s$ рядов, которые называются условиями балансировки.
В конечном счете они оказываются связанными со старыми хорошо известными
условиями балансировки для обычных и $q$-гипергеометрических рядов.
Эта спецификация тета-рядов носит фундаментальный характер и указывает
путь для введения дальнейших структурных элементов. Функция $h(y)$ обладает
конечным набором параметров описывающих положения нулей и полюсов, по которым она
квазипериодична. Естественно потребовать, чтобы $h(y)$ была двояко-периодична по этим
переменным. Это приводит к ограничениям на параметры, известным как условия
вполне уравновешенности рядов гипергеометрического типа. Далее вводится условие
совершенной уравновешенности для эллиптических гипергеометрических рядов, также
носящее достаточно естественный характер и связанное с удвоением аргумента тета-функций.

Оказывается, что при специальном значении степенного аргумента в
обрывающемся совершенно уравновешенном сбалансированном $_{10}E_9$ ряде,
он сворачивается в замкнутое выражение равное отношению произведений
тета-функций Якоби с явно заданными аргументами. Эта формула суммирования была впервые
выведена Френкелем и Тураевым в работе \cite{fre-tur:elliptic} и при $p\to 0$
она сводится к сумме Джексона для обрывающегося совершенно уравновешенного
сбалансированного $_8\varphi_7$ ряда \cite{gas-rah:basic}. С помощью этой формулы
суммирования, в
параграфе 2.3 строится эллиптическая цепочка Бэйли, позволяющая найти бесконечные
последовательности нетривиальных тождеств для эллиптических гипергеометрических рядов.
Она обобщает хорошо известные результаты по цепочкам преобразований Бэйли для
$q$-гипергеометрических рядов \cite{aab:bailey,and:bailey,and-ber:bailey,bre:some,war:50},
ключевая итеративная природа которых была обнаружена Эндрюсом
\cite{and:multiple} (см. также \cite{pau:identities}).
Сама оригинальная работа Бэйли \cite{bai:some,bai:identities} была написана с
целью прояснения общего универсального механизма в доказательствах знаменитых тождеств
Роджерса-Рамануджана.

Многократные эллиптические гипергеометрические ряды были рассмотрены
впервые Варнааром в работе \cite{war:summation}. В параграфе 2.4
дается общее определение многократных эллиптических гипергеометрических
рядов и приводится несколько типов обобщений формулы суммирования Френкеля-Тураева
для многократных рядов на корневых системах $A_n$ и $C_n$.

Последний параграф этой главы посвящен обсуждению возможности обобщения эллиптических
результатов на Римановы поверхности произвольного рода. При этом удается
получить довольно простую формулу суммирования для однократного многопараметрического
телескопирующегося ряда, построенного из многомерных тета-функций Римана.
Этот результат обобщает один из результатов Варнаара, связанный с
эллиптическим аналогом многопараметрической формулы суммирования Макдональда.

\underline{В третьей главе} излагаются результаты исследований автора по
тета-гипер\-гео\-ме\-три\-чес\-ким интегралам,
обобщающим обычные и $q$-гипергеометрические
интегралы на эллиптический уровень. Согласно их общему
определению, данному в параграфе 3.1,
необходимо взять мероморфную функцию $\Delta(y_1,\ldots,$ $y_n)$
и сконструировать многократные интегралы
$I_n=\int_{\cal D} \; \Delta(y_1,\ldots,y_n)\: dy_1\cdots dy_n$
с некоторым многомерным циклом ${\cal D}\in {\bf C}^n$. При этом требуется, чтобы
ядро $\Delta(y_1,\ldots,$ $y_n)$ удовлетворяло системе разностных уравнений
первого порядка с коэффициентами, обладающими свойствами аналогичными
квазипериодичности тета-функций Якоби. В параграфе 3.1 выводится
общий вид таких интегралов в одномерном
случае, что приводит к эллиптическому и тета-аналогам функции Мейера.
При построении таких объектов используются эллиптические гамма-функции.

Стандартная эллиптическая гамма-функция с двумя базисными переменными
$q$ и $p$, удовлетворяющими ограничениям $|q|,|p|<1$, как независимый
объект теории специальных функций, была рассмотрена Рюйсенаарсом в
недавней работе \cite{rui:first}
с удовлетворительным с современной точки зрения анализом ее свойств.
Впервые функции такого типа появились в работе Джексона 1905 г.
\cite{jac:basic}. В неявном виде эта функция возникала также во
многих исследованиях по решаемым моделям статистической механики, начиная со
статьи Бахтера 1972 г. по восьмивершинной модели \cite{bax:partition}.
Некоторые дополнительные важные свойства этой функции были установлены
Фельдером и Варченко в работе \cite{fel-var:elliptic}. При этом точных формул
интегрирования, связанных с этой обобщенной гамма-функцией, ни в одном
из упомянутых исследований не предлагалось. Другой тип эллиптической
гамма-функции, у которой один из базисных параметров может лежать
на единичной окружности, $|q|=1$, предложен автором диссертации \cite{spi:integrals}.
В тригонометрическом пределе $p\to 0$ эта функция вырождается в известную
обобщенную гамма-функцию, рассмотренную  Шинтани \cite{shi:kronecker}
и названную функцией двойного синуса \cite{kur:multiple}.
Все эти гамма-функции выражаются в виде различных комбинаций
многократной гамма-функции Барнса \cite{bar:theory,bar:multiple}, играющей
таким образом основополагающую роль в теории функций гипергеометрического типа.

В параграфе 3.2 получен один из самых важных результатов диссертации, состоящий
в вычислении одномерного эллиптического бета-интеграла, выражающегося
через стандартную эллиптическую гамма-функцию $\Gamma(z;q,p)$,
хорошо определенную при $|p|,|q|<1$.
С его помощью, в параграфе 3.4 построено два интегральных
аналога цепочек Бэйли, которые приводят к двоичному дереву тождеств для
эллиптических гипергеометрических интегралов. Этот результат описывает самый
первый пример интегрального аналога леммы Бэйли.

В параграфах 3.5, 3.6 и 3.7 выводятся $n$-кратные эллиптические бета-интегралы
на корневых системах $A_n$ и $C_n$. Всего получено семь различных точных
формул интегрирования, разбитых на три типа. Тип I характеризуется наличием
$2n+3$ свободных параметров и доказываются с помощью разностных
уравнений. Тип II выводится из типа I с помощью составных многомерных
интегралов путем разных порядков интегрирования в них и такие интегралы содержат
меньше чем $2n+3$ параметра. Интегралы типа III вычисляются как нетривиальные
детерминанты от одномерного эллиптического бета-интеграла. Часть из этих интегралов
обобщает многократные $q$-бета интегралы Густафсона
\cite{gus:generalization,gus:some1,gus:some2,gus:some3} и Густафсона-Ракха
\cite{gus-rak:beta}. Наиболее важный из них представляет собой эллиптический аналог
интеграла Сельберга \cite{sel:bemerkninger}.

\underline{В четвертой главе} диссертации конструируется система биортогональных
функций одной переменной, для которых одномерный эллиптический бета-интеграл
определяет меру в соотношениях биортогональности.
Начинается она с описания соотношений сопряжения для эллиптической гипергеометрической
функции, заданной контурным интегралом от определенной комбинации эллиптических
гамма-функций. Аналогичные соотношения справедливы
для обрывающегося совершенно уравновешенного
$_{12}E_{11}$ эллиптического гипергеометрического ряда с аргументом $x=-1$.
Эти соотношения эквивалентны уравнениям пары Лакса
для обобщенной спектральной задачи первой главы.

В параграфе 4.2 вводятся функции $R_n(z;q,p),\: n\in\N,$
выражающиеся через совершенно уравновешенный $_{12}E_{11}$ ряд
со специальной параметризацией. С помощью соотношений сопряжения для этих
функций выводится трехчленное рекуррентное соотношение по номеру $n$
и разностное уравнение второго порядка по $z$.

Параграф 4.3 содержит еще один ключевой результат диссертации. Для функций
$R_{nk}(z)\equiv R_n(z;q,p)R_k(z;p,q)$ и их партнеров из дуального
пространства $T_{nk}(z)\equiv T_n(z;q,p)$ $\times T_k(z;p,q)$, где $T_n(z;q,p)$
есть последовательность функций аналогичных $R_n(z;q,p)$,
доказано соотношение двухиндексной биортогональности. По всей видимости,
такие соотношения, характерные для функций двух независимых переменных,
еще не встречались в теории специальных функций одной переменной.
Кроме того,
функции $R_{nk}(z)$ и $T_{nk}(z)$ не имеют предела $p\to 0$ при $k\neq 0$,
то есть в общем случае они существуют только на эллиптическом уровне.

В параграфе 4.4 выводится интегральное представление для симметричного
произведения двух обрывающихся совершенно уравновешенных $_{12}E_{11}$ рядов.
В следующем параграфе выводится обрывающаяся цепная дробь, связанная
с трехчленным рекуррентным соотношением для функций $R_n(z;q,p)$ и
описываемая аналогичным $_{12}E_{11}$ рядом.
Эта дробь включает в себя цепную дробь Рамануджана \cite{Ramanu,Berndt} и ее
обобщение, найденное Массоном \cite{Masson}, которые выражаются
через обычную гипергеометрическую функцию $_9F_8$.
Также она содержит и сопутствующие $q$-гипергеометрические цепные дроби Ватсона
\cite{Watson} и Гупты-Массона \cite{gup-mas:contiguous}, выражающиеся через
$_{10}\varphi_9$ ряд.
В параграфе 4.6 кратко описывается связь разностного уравнения для
$R_n(z;q,p)$ с одночастичным сектором обобщенной конечно-разностной
эллиптической модели типа Калоджеро-Сазерленда-Мозера, предложенной
в работах \cite{die:integrability,die:difference}.

\underline{В пятой главе} диссертации теория тета и эллипти\-чес\-ких
гипер\-гео\-ме\-три\-чес\-ких функций обобщена на случай когда один из
базисных параметров, например $q$, лежит на единичной окружности, $|q|=1$,
но при этом $|p|<1$.

В разделе 5.1 вводится модифицированная эллиптическая гамма-функция
$G(u;\boldsymbol{\omega})$, параметризуемая тремя комплексными переменными
 $\omega_{1,2,3}$. Она удовлетворяет трем разностным уравнениям
первого порядка, одно из которых совпадает с ключевым уравнением для
$\Gamma(z;q,p)$ при $q=e^{2\pi i\frac{\omega_1}{\omega_2}}$,
$p=e^{2\pi i\frac{\omega_3}{\omega_2}} $ и $z=e^{2\pi i\frac{u}{\omega_2}}$.
Эта функция хорошо определена при $|p|<1, |q|=1$, в отличие от $\Gamma(z;q,p)$.
При тригонометрическом вырождении возникает функция двойного синуса,
связанная с концепцией модулярного дубля Фаддеева
\cite{fad:discrete,fad:modular,fkv:strongly}
(в этой серии работ эта функция называется некомпактным квантовым дилогарифмом).

В параграфе 5.2 выводятся модифицированные аналоги основных  эллиптических
бета-интегралов третьей главы, которые хорошо определены и
в режиме $|p|<1, |q|=1$.
В следующем параграфе рассматривается предел $\mbox{Im}(\omega_3/\omega_2)\to+\infty$
взятым так, что одновременно с $p \to 0$ параметр
$r=e^{2\pi i\frac{\omega_1}{\omega_3}}$ также стремится к нулю,
$r \to 0$. Это приводит к
некомпактным $q$-бета-интегралам по бесконечным контурам, выражающимся через
функцию двойного синуса. Эти интегралы доказываются
независимым образом, так как взятый предел носит формальный характер.
В разделе 5.3.1 представлен интеграл такого типа, связанный с модифицированным
эллиптическим бета-интегралом типа II для системы корней $C_n$.
В разделе 5.3.2 рассмотрены соответствующие партнеры для эллиптических бета-интегралов
типа I на системах корней $A_n$ и $C_n$. В параграфе 5.4 кратко описываются одномерные
биортогональные функции с мерой, определяемой модифицированным
одномерным эллиптическим бета-интегралом, и их $q$-гипергеометрическое вырождение.

\underline{В заключении} перечислены основные результаты, выносимые
на защиту.

\bigskip
%\clearpage

Результаты диссертации докладывались и обсуждались на семинарах
Лаборатории теоретической физики им. Н.Н. Боголюбова ОИЯИ,
в отделах дифференциальных уравнений и математической физики
Математическом института им. В.А. Стеклова РАН (г. Москва),
в университетах в городах Монреаль (Канада), Париж (Франция),
Сантьяго и Талька (Чили), Луван-ла-Нев (Бельгия),
Институте Математики им. Макса Планка (г. Бонн, Германия).
Также результаты диссертации были представлены в докладах на следующих
международных конференциях и совещаниях: ``Special functions" (Гон Конг, 1999),
школе передовых исследований НАТО ``Special functions-2000: Current perspective
and future directions" (Аризона, США, 2000), ``Fifth International Conference on
Difference Equations and Applications" (Темуко, Чили, 2000), ``Integrable Structures
of Exactly Solvable Two-Dimensional Models of Quantum Field Theory" (Киев, Украина,
2000), ``Algebra and Theory of Numbers" (Талька, Чили, 2000),
``Asymptotic Combinatorics with Applications to Mathematical Physics"
(Санкт Петербург, 2001), ``Special Functions in the Digital Age"
(Институт Прикладной Математики, Миннеаполис, США, 2002),
``Foundations of Computational Mathematics" (Миннесота, США, 2002),
``Number Theory and Combinatorics in Physics"
(Флорида, США, 2003), ``Jack, Hall-Littlewood, and Macdonald Polynomials"
(Эдинбург, 2003), ``Special functions, orthogonal polynomials, quantum groups
and related topics" (Бексбах, Германия, 2003), ``Классические и квантовые
интегрируемые системы" (Протвино, 2002, 2003;  Дубна, 2004), международном семинаре
``Quarks-2004" (Псков, 2004).
По теме диссертации написана 21 статья, цитированные как
\cite{die-spi:elliptic}--\cite{die-spi:unit}, \cite{zhe-sip:umn},
\cite{spi:solitons}--\cite{spi:short} и
\cite{spi-zhe:spectral}--\cite{spi-zhe:theory}.

%% file: CHAPTER1.TEX
\chapter[Нелинейные цепочки с дискретным временем и их автомодельные решения]
{Нелинейные цепочки с дискретным временем и их автомодельные решения}

\section{Метод факторизации для уравнения Шредингера}

\subsection{Факторизационная цепочка}

Рассмотрим одномерное уравнение Шредингера
\begin{equation}\label{schr-eq}
L\psi(x)=-\psi_{xx}(x)+u(x)\psi(x)=\lambda\psi(x),
\end{equation}
описывающее движение нерелятивистской частицы на линии
$x\in \R$ в потенциале $u(x)$, задаваемом некоторой
функцией ограниченной снизу. Для удобства, перенормировкой
координаты $x$ и энергии $\lambda$ мы устранили из
оператора Шредингера массу частицы и постоянную Планка
$\hbar$. В рамках квантовой механики разрешенные значения энергии
определяются как собственные значения оператора $L$,
соответствующие условию квадратичной интегрируемости
модуля волновых функций $\psi(x)$ на каком-либо конечном или
бесконечном интервале $x\in \R $. С математической точки зрения
интересны все решения уравнения \re{schr-eq} при $x, \lambda\in \C $.

Дифференциальный оператор второго порядка $L=-d^2/dx^2+u(x)$
может быть факторизован в виде произведения двух операторов первого порядка
с точностью до некоторой константы $\lambda_0$:
\begin{equation}\label{fact}
L=A^+A^-+\lambda_0, \qquad A^\pm=\mp d/dx +f(x),
\end{equation}
где функция $f(x)$ есть решение уравнения Риккати
$f^2(x)-f_x(x)+\lambda_0=u(x)$. Подстановка
$f(x)=-\phi_{0,x}(x)/\phi_0(x)$ показывает, что $-\phi_{0,xx}(x)+u(x)\phi_0(x)
=\lambda_0\phi_0(x)$, то есть $\phi_0(x)$ является решением начального
уравнения (\ref{schr-eq}) при $\lambda=\lambda_0$.

Предположим, что $f(x)$ является гладкой функцией при конечных $x\in{\R}$.
Тогда $A^+$ является эрмитовым сопряжением оператора $A^-$ в
${\rm L}^2({\R})$, а $L$ --- самосопряжен. При таких условиях $\lambda_0$
не может быть больше минимального собственного значения $L$, так как
для любой нормированной функции $\psi(x)$ из области определения $L$ имеем
$\la \psi|L|\psi\ra=||A^-\psi||^2+\lambda_0\geq\lambda_0$.
Предположим, что $\lambda_0$ есть минимальное собственное значение $L$.
Обозначим $\psi_0(x)$ соответствующую нормированную собственную функцию
(не имеющую, как известно, нулей), $||\psi_0||^2=
\int_{-\infty}^\infty|\psi_0(x)|^2dx=1.$ Тогда, функция
$$
\phi_0(x)=a\psi_0(x)+b\psi_0(x)\int^x\frac{dy}{\psi_0^2(y)},
$$
где $a$ и $b$ произвольные параметры, является общим решением уравнения
$L\phi_0=\lambda_0\phi_0$. При $b\neq 0$, функция $f(x)=-d\log\phi_0(x)/dx$
сингулярна в некоторой точке и операторы $A^\pm$ (\ref{fact})
плохо определены. Поэтому мы исключим этот случай и положим $a=1, b=0$.

Используя собственную функцию, соответствующую наименьшему собственному
значению $L$, всегда можно факторизовать $L$ как указано в (\ref{fact}) с
хорошо определенными операторами $A^\pm$. И наоборот, если найти
факторизацию (\ref{fact}), такую что нулевая мода оператора
$A^-$, $A^-\psi_0=0$, принадлежит $L^2({\R})$ и $f(x)$ несингулярна,
то $\lambda_0$ является наименьшим собственным значением  и
\begin{equation}\label{zmode}
\psi_0(x)=\frac{e^{-\int_{0}^xf(y)dy}}
{\left(\int_{-\infty}^\infty e^{-2\int_{0}^x f(y)dy}dx\right)^{1/2}}
\end{equation}
является соответствующей нормированной собственной функцией.

Определим новый оператор Шредингера $\tilde L$ перестановкой операторных
факторов в выражении (\ref{fact})
\begin{equation}\label{new-op}
\tilde L=A^-A^+ + \lambda_0.
\end{equation}
Соответствующий потенциал имеет форму $\tilde u(x)=f^2(x)+f_x(x)+\lambda_0$.
Очевидно, что справедливо сплетающее соотношение
\begin{equation}\label{inter}
LA^+=A^+\tilde L, \qquad A^-L=\tilde LA^-,
\end{equation}
играющее ключевую роль в этом формализме. Из равенств (\ref{inter}) можно увидеть,
что если $\psi(x)$ удовлетворяет уравнению (\ref{schr-eq}),
то функция $\tilde\psi=A^-\psi$ определяет формальную собственную функцию
$\tilde L$. Действительно,
\begin{equation}\label{interpsi}
\tilde L\tilde\psi=\tilde L (A^-\psi)=A^-(L\psi)=\lambda(A^-\psi).
\end{equation}
На самом деле это дает общее решение дифференциального уравнения
$\tilde L\tilde\psi=\lambda\tilde\psi$ для всех $\lambda$, кроме
случая $\lambda=\lambda_0$, соответствующего нулевой моде $A^-$.
Недостающее решение уравнения $\tilde L\tilde\phi_0=\lambda_0\tilde\phi_0$
легко восстанавливается:
\begin{equation}\label{gensol2}
\tilde\phi_0(x)= \frac{g}{\psi_0(x)}+\frac{e}{\psi_0(x)}\int^x\psi_0^2(y)dy,
\end{equation}
где $g$ и $e$ произвольные параметры и $\psi_0$ обозначает основное состояние
$L$ (заметим, что справедливо равенство $A^+\psi_0^{-1}=0$).

Обозначим $\lambda_n, \tilde \lambda_n$ и $\psi_n,$ $\tilde\psi_n=A^-\psi_n$
дискретные собственные значения и соответствующие собственные функции операторов
$L$ и $\tilde L$. Из (\ref{interpsi}) следует, что спектры $\lambda_n$ и
$\tilde \lambda_n$ почти совпадают $\tilde\lambda_n=\lambda_n, n=1,2,\dots$.
Единственное отличие может касаться нулевой моды оператора $A^-$.
Действительно, точка $\lambda=\lambda_0$ не принадлежит спектру $\tilde L$:
 поскольку $\psi_0\in \text{L}^2({\R})$,
мы имеем $\psi_0^{-1}\notin \text{L}^2({\R})$, и, как следствие, для
любых $g$ и $e$ функция (\ref{gensol2}) не может быть нормируемой.

Следовательно, $\lambda_1$ является наименьшим собственным значением оператора
$\tilde L$ и точка $\lambda_0$ была ``стерта" из спектра $L$.
Повторяя эту процедуру еще раз, т.е. найдя подходящую факторизацию
$\tilde L=\tilde A^+\tilde A^-+\lambda_1$ и переставив операторные
множители $\tilde A^+$ и $\tilde A^-$, мы можем стереть точку
$\lambda=\lambda_1$ в спектре. Таким образом можно удалить произвольное число
собственных значений $L$ снизу.
Часто гораздо легче найти наименьшее собственное значение заданного оператора
$\tilde L$, чем остальные. Если каким-либо способом оно найдено,
то оно совпадает со вторым собственным значением $L$, и т.д. Это наблюдение
является центральным в методе факторизации \cite{Sc1} поскольку оно сводит
проблему поиска дискретного спектра данного оператора к проблеме поиска
наименьших собственных значений последовательности операторов, построенных из
$L$ алгоритмом ``факторизации и перестановки".

Можно обернуть процедуру устранения наименьшего собственного значения заданного
оператора Шредингера. Пусть оператор $\tilde L$ имеет наименьшее собственное значение
$\lambda_1$. Факторизуем его как $\tilde L=A^-A^++\lambda_0$, где
$\lambda_0<\lambda_1$. Если нулевая мода оператора $A^-$ нормируема,
то тогда оператор $L=A^+A^-+\lambda_0$ имеет тот же спектр, что и
$\tilde L$ с дополнительным вставленным собственным значением в произвольной
точке $\lambda= \lambda_0$. Если же факторизовать $L$ (или $\tilde L$) так
что  $f(x)$ несингулярно, а нулевые моды операторов $A^\pm$ ненормируемы,
то тогда дискретные спектры операторов $L$ и $\tilde L$ совпадают
(изоспектральная ситуация).

Имеются более сложные возможности для изменения спектральных данных оператора $L$.
Например, если при первой факторизации разрешить
сингулярную $f(x)$, то $\tilde L$ будет иметь сингулярный потенциал и его спектр,
вообще говоря, будет сильно отличаться от спектра $L$. Однако, можно потребовать,
чтобы после ряда таких ``неудачных" факторизаций возник оператор Шредингера
с несингулярным потенциалом. Оказывается, что таким образом можно устранять
или вставлять собственные значения в произвольном месте спектра \cite{Kr}.
В частности, двухшаговая процедура перефакторизации позволяет вставлять дискретные
уровни в непрерывный спектр \cite{D}.

Перепишем теперь описанную конструкцию в формализме  дискретного времени.
Обозначим $L\equiv L_j$ и $\tilde L\equiv L_{j+1}$, где $j\in {\Z }$
(по соглашению, мы можем рассматривать индекс $j$ как непрерывную переменную,
равно как мы могли бы обозначить $\tilde L \equiv L_{j-1}$ и т.д.).
Таким образом мы получаем бесконечную последовательность операторов Шредингера
$L_j=-d^2/dx^2+u_j(x)$ с формальной факторизацией
\begin{equation}
L_j=A_j^+A_j^- + \lambda_j,\qquad A_j^\pm=\mp\;d/dx + f_j(x).
\label{fact2}\end{equation}
Соседние $L_j$ связаны друг с другом посредством абстрактной
факторизационной цепочки
\begin{equation}\label{fchain}
L_{j+1}=A_{j+1}^+A_{j+1}^-+\lambda_{j+1}=A_j^-A_j^++\lambda_j,
\end{equation}
а сплетающие соотношения принимают вид
$$
A_j^-L_j=L_{j+1}A_j^-, \qquad L_jA^+_j=A_j^+L_{j+1}.
$$
Подставляя явный вид $A^\pm_j$ в (\ref{fchain}) мы получаем
дифференциально-разностное уравнение для $f_j(x)$:
\begin{equation}\label{chain}
(f_j(x)+f_{j+1}(x))_x+f_j^2(x)-f_{j+1}^2(x)=
\mu_j\equiv \lambda_{j+1}-\lambda_j.
\end{equation}
Эта цепочка была выведена Инфельдом в работе \cite{Inf}, где проблема поиска
точно решаемых спектральных задач была сформулирована как задача поиска
решений (\ref{chain}) таких что величины $\lambda_j$ определяют
дискретный спектр оператора $L$, скажем, $L\equiv L_0$. В качестве
одной из возможностей, можно попытаться найти  решение $f_j(x)$ в виде
Лорановского ряда по $j$, считающейся непрерывной переменной (при этом $\lambda_j$
рассматриваются как неизвестные функции $j$). Подробный анализ показывает
\cite{Inf,IH}, что в
таком разложении имеется конечное число членов только если
$f_j(x)=a(x)j+b(x)+c(x)/j$, где $a, b, c$ являются некоторыми элементарными
функциями $x$. Это приводит к гипергеометрической функции $_2F_1$
и хорошо известным ``старым" решаемым потенциалам квантовой механики.

По построению, константы $\lambda_j, j\geq 0,$ определяют наименьшие
собственные значения операторов $L_j$ при условии, что нулевые моды $A_j^-$
нормируемы и $f_j(x)$ несингулярны. В общей ситуации, сдвиги
$j\to j+1$ могут описывать все три возможности --- устранение или вставку
собственного значения или изоспектральное преобразование.

Для положительного целого числа $n$, определим операторы
$$
M_j^-=A_{j+n-1}^-\cdots A_{j+1}^-A_j^-, \qquad
M_j^+=A_j^+A_{j+1}^+\cdots A_{j+n-1}^+.
$$
Сплетающие соотношения
$$
L_{j+n}M_j^-=M_j^-L_j, \qquad M_j^+L_{j+n}=L_jM^+_j
$$
гарантируют, что почти для всех $\lambda$ решения уравнений
$$
L_j\psi^{(j)}=\lambda\psi^{(j)} \qquad \mbox{и} \qquad
L_{j+n}\psi^{(i+n)}=\lambda\psi^{(j+n)}
$$
связаны друг с другом как  $\psi^{(j+n)}\propto M_j^-\psi^{(j)}$ и
$\psi^{(j)}\propto M^+_j\psi^{(j+n)}$. Как следствие, произведение
$M^+_jM^-_j$ коммутирует с $L_j$, а $M^-_jM^+_j$
коммутирует с $L_{j+n}$. Простое вычисление приводит к соотношениям
\be
M_j^+M_j^-=\prod_{k=0}^{n-1}(L_j-\lambda_{j+k}), \qquad
M_j^-M_j^+=\prod_{k=0}^{n-1}(L_{j+n}-\lambda_{j+k}).
\lab{polrel}\end{equation}
Пусть $\psi^{(j)}(x)\in$ L$^2({\R})$ будет собственной функцией
$L_j$ с собственным значением $\lambda$ и конечной нормой $||\psi^{(j)}||$.
Если мы положим $\psi^{(j+n)}=M_j^-\psi^{(j)}$, тогда
$$
||\psi^{(j+n)}||^2=(\lambda-\lambda_j)\cdots(\lambda-\lambda_{j+n-1})\,
||\psi^{(j)}||^2.
$$
Нули выражения в правой части этого равенства при
$\lambda=\lambda_k$ для некоторых $k$ показывает, что соответствующие
собственные значения были стерты из спектра. При условии, что нулевые моды
операторов $A^-_j, j=0, \dots, n-1,$ нормируемы и не имеют нулей,
$A^-_j\psi_0^{(j)}=0$, $||\psi_0^{(j)}||=1$,
мы находим, что функции
$$
\psi^{(0)}_n(x)=\frac{A^+_0\cdots A^+_{n-1} \psi_0^{(n)}(x)}
{\sqrt{(\lambda_n-\lambda_{n-1})\cdots(\lambda_n-\lambda_0)}}
$$
определяют нормированные собственные функции оператора $L_0$ с собственными
значениями $\lambda_n$:
$L_0\psi^{(0)}_n=\lambda_n\psi^{(0)}_n$ и $||\psi_n^{(0)}||=1$.

Основное преимущество метода факторизации состоит в его чисто операторной
формулировке. Можно заменить операторы Шредингера любым другим (дифференциальным
оператором более высокого порядка, ко\-неч\-но-раз\-ност\-ным,
дифференциально-разностным, интегральным
и т.д.) оператором $L$, допускающим  факторизацию на хорошо определенные
операторы более простой структуры  $A^\pm$ (в некоторых ситуациях эти операторы не
обязаны быть эрмитово-сопряженными друг к другу). Во всех этих случаях мы имеем
дело с абстрактной операторной факторизационной цепочкой (\ref{fchain})
с соответствующей интерпретацией констант $\lambda_j$.

\subsection{Автомодельные потенциалы и квантовые алгебры}

Опишем ряд конкретных систем, удовлетворяющих операторной
факторизационной цепочке (\ref{fchain}). В качестве ведущего примера
можно взять стандартную проблему гармонического осциллятора---базовую модель
квантовой механики и квантовой теории поля, и попытаться образовать из
операторов $A^\pm_j$ некоторую нетривиальную алгебру симметрии.
Действительно, соотношения \re{polrel} уже выглядят как составляющие
некоторой алгебры, но они не образуют замкнутой системы---связь между
проблемами на собственные значения для операторов $L_j$ и $L_{j+n}$
слишком слаба (эти соотношения справедливы практически для произвольного
начального потенциала $u_j(x)$). Для того чтобы замкнуть систему необходимо
наложить дополнительную связь между $L_j$ и $L_{j+n}$,
скажем $L_{j+n}\propto L_j$, что заставит операторы $M_j^\pm$ преобразовывать
пространство собственных функций взятого оператора $L_j$ на себя.

В простейшем случае можно потребовать периодичности последовательности Гамильтонианов
$L_j$: $L_{j+N}=L_j$ для некоторого целого $N>0$. В результате, операторы $M^\pm_j$
при $n=N$ коммутируют с $L_j$, $[M^\pm_j,L_j]=0$, и
этот замечательный факт упрощает поиск собственных функций для $L_j$.
При реализации  $A^\pm_j$ дифференциальными операторами это приводит к
дифференциальным операторам $N$-го порядка, коммутирующим с изначальным
оператором Шредингера \cite{BC1}. Существует обобщение этого условия
чистой периодичности на периодичность с точностью до подкрутки
$L_{j+N}=UL_jU^{-1}$, где $U$ некоторый обращаемый оператор.
Это опять приводит к коммутирующим операторам $[B^\pm_j,L_j]=0$, где
$B^+_j=M^+_jU,\; B^-_j= U^{-1}M_j^-$, но теперь появляется существенная
дополнительная свобода в выборе $U$.

Другой возможный вариант ``замыкания" или редукции последовательности
операторов $L_j$ состоит в требовании их периодичности с точностью до
добавки константы и подкрутки, $L_{j+N}=UL_jU^{-1}+\mu$, где $\mu$ некоторый
произвольный параметр. Это приводит к лестничным соотношениям
$[L_j,B^\pm_j]=\pm \mu B^\pm_j$ и операторным тождествам
$$
B^+_jB^-_j=\prod_{k=0}^{N-1}(L_j-\lambda_{j+k}),\qquad
B^-_jB^+_j=\prod_{k=0}^{N-1}(L_j+\mu-\lambda_{j+k}),
$$
где операторы $B_j^\pm$ определены так же как и в предыдущем случае.
Обозначив $B^\pm\equiv B_0^\pm$ и $L\equiv L_0$, мы можем образовать
полиномиальную алгебру
\begin{equation}\label{pol}
[L,B^\pm]=\pm\mu B^\pm,\qquad [B^+,B^-]=P_{N-1}(L),
\end{equation}
где $P_{N-1}(x)$ некоторый полином степени $N-1$ от переменной $x$.
Теория представлений такой алгебры проанализирована, например, в статье \cite{Sm}.
Отметим, что при $N=1$ мы получаем алгебру Гейзенберга-Вейля, описывающую
гармонический осциллятор, а при $N=2$ возникает $su(1,1)$ алгебра.

Квантовые алгебры, или $q$-аналоги алгебр (\ref{pol}) возникают при следующем
операторном автомодельном ограничении, наложенном на цепочку (\ref{fchain}):
\begin{equation}\label{q-closure}
L_{j+N}=q^2UL_jU^{-1}+\mu.
\end{equation}
Если числовой множитель $q^2\neq 1$, то тогда можно устранить $\mu$ простым
однородным сдвигом $L_j\to L_j+\mu/(1-q^2)$. Поэтому ниже мы будем подразумевать,
что $\mu=0$. Подставим условие замыкания (\ref{q-closure}) с $\mu=0$ в (\ref{polrel}).
Тогда операторы $L=L_0$ и $B^+=M^+_0U,\; B^-= U^{-1}M_0^-$
удовлетворяют следующим соотношениям
\begin{equation}\label{q-alg}
LB^\pm=q^{\pm 2}B^\pm L,\qquad
B^+B^-=\prod_{j=0}^{N-1}(L-\lambda_j),\qquad
B^-B^+=\prod_{j=0}^{N-1}(q^2L-\lambda_j).
\end{equation}
При $N=1$ эти равенства определяют некоторую реализацию алгебры
$q$-гармоничес\-ко\-го осциллятора
\begin{equation}\label{qosc}
B^-B^+-q^2B^+B^-=\rho, \qquad [B^\pm,\rho]=0,
\end{equation}
с $\rho=\lambda_0(q^2-1)$. Можно положить $\rho=1$, выбрав условие нормировки
$\lambda_0=1/(q^2-1)$. Такой $q$-аналог алгебры Гейзенберга-Вейля
появился в физике очень давно \cite{Coo,FB}. В современное время эта алгебра стала
весьма популярной благодаря связи с точно решаемыми моделями статистической
механики и квантовыми группами (см., например, \cite{Mac}). Классификация
неприводимых представлений этой алгебры дана в работе \cite{kul:q-osc}.
При $N=2$, соотношения (\ref{q-alg}) определяют специфический $q$-аналог
конформной алгебры $su(1,1)$, допускающий структуру алгебры Хопфа, и
так далее \cite{S0,spi:symmetries,spi:universal}.

Рассмотрим теперь реализацию указанных соотношений на уравнении Шредингера.
Начнем с рассмотрения замыкания $\tilde L=ULU^{-1}+\mu$,
где $L=-\partial_x^2+u(x)$ и $U$ есть оператор сдвига $Uf(x)=f(x+a)$.
В результате, операторы $B^\pm$ принимают вид $B^+=(-\partial_x+f(x))U$,
$B^-=U^{-1}(\partial_x+f(x))$ и мы получаем алгебру Гейзенберга Вейля
$[L,B^\pm]=\pm\mu B^\pm, [B^-,B^+]=\mu$. Потенциал, входящий в $L$,
определяется как $u(x)=f^2(x)-f_x(x)+\lambda_0$,
где $f(x)$ удовлетворяет следующему нелинейному дифференциальному уравнению
с запаздывающим аргументом
\begin{equation}\label{nonlocal}
\left(f(x)+f(x+a)\right)_x+f^2(x)-f^2(x+a)=\mu.
\end{equation}
Для $a=0$ мы получаем из этого уравнения стандартную модель гармонического
осциллятора $f(x)=\mu (x-x_0)/2$ и $u(x)=\mu^2 (x-x_0)^2/4-\mu/2+\lambda_0$,
связанную с полиномами Эрмита. В работе \cite{spi:symmetries}
автор нашел мероморфное решение (\ref{nonlocal}) при $a\neq 0$ и $\mu=0$
в терминах $\cal P$-функции Вейерштрасса
$$
f(x)=-\frac{1}{2}\,\frac{{\cal P}'(x-x_0)-{\cal P}'(a)}
{{\cal P}(x-x_0)-{\cal P}(a)}.
$$
При $\mu\neq 0$ крайне трудно описать решения этого уравнения
аналитические в какой-либо конечной области. Некоторый анализ мероморфных
решений проведен в работе \cite{tov}. Заметим, что уравнение (\ref{nonlocal})
так же может быть получено наложением ограничения $f_{j+1}(x)=f_j(x+a)$
на цепочку (\ref{chain}).

Перейдем теперь прямиком к рассмотрению соотношений (\ref{q-alg}).
В этом случае мы можем взять в качестве $U$ оператор растяжений
(или оператор $q$-разности), $Uf(x)=|q|^{1/2}f(qx)$. Для действительных
 $q\neq 0$ этот оператор унитарен $U^\dagger=U^{-1}$. Определив $L\equiv L_0,$
$B^-\equiv U^{-1}(\partial_x+f_{N-1})\cdots(\partial_x+f_0)$ и
$B^+\equiv (-\partial_x+f_0)\cdots(-\partial_x+f_{N-1})U$,
мы сможем реализовать соотношения (\ref{q-alg}) при условии, что $f_j(x)$
удовлетворяют следующей системе нелинейных смешанных дифференциальных
и $q$-разностных уравнений:
$$
(f_0(x)+f_1(x))_x+f_0^2(x)-f_1^2(x)=\mu_0, \quad \dots\dots
$$
\begin{equation}\label{q-pereq}
(f_{N-1}(x)+qf_0(qx))_x+f_{N-1}^2(x)-q^2f_0^2(qx)=\mu_{N-1}.
\end{equation}
Эквивалентным образом, эти уравнения появляются из цепочки (\ref{chain})
после наложения ограничений
\begin{equation}
f_{j+N}(x)=qf_j(qx),\qquad \mu_{j+N}=q^2\mu_j,
\label{q-per}\end{equation}
имеющих простую теоретико-групповую интерпретацию.

Автомодельное замыкание (\ref{q-per}) легко вывести с помощью Лиевской
техники симметрийных редукций. В контексте дифференциальных уравнений
эта теория обрисована, например, в монографии \cite{Mi3}, а её расширение на
конечно-разностные уравнения рассмотрена, например, в работах \cite{LW,Mae}.
Во первых, необходимо заметить что если функции $f_j(x)$ и параметры $\mu_j$
представляют собой некоторые решения цепочки (\ref{chain}), то их дискретное
масштабное преобразование $qf_j(qx)$ и $q^2\mu_j$, при некотором $q\in\C$
также задают решения (\ref{chain}). Аналогично, сдвиг в нумерации этих
переменных $f_j(x)\to f_{j+N}(x),\; \mu_j\to\mu_{j+N}$ отображает
решения (\ref{chain}) на решения. Рассмотрим семейство  автомодельных
решений (\ref{chain}), которое инвариантно относительно комбинации этих
преобразований симметрии. Так, требуя, чтобы эти два преобразования
были эквивалентны друг другу, мы приходим к (\ref{q-per}). Для $N=1$, эта
редукция может быть переписана в виде $f_j(x)=q^jf_0(q^jx),\;
\lambda_j=q^{2j}\lambda_0$ и в этой форме она была открыта в работе
\cite{Sh}. Общая $q$-периодическая редукция (\ref{q-per})
была найдена автором \cite{S0} исходя из идеи $q$-деформации
алгебры обобщенной суперсимметрии, связанной с парастатистикой \cite{RS}.

Решения уравнения (\ref{q-pereq}) имеют очень сложную структуру и общие методы
решения дифференциально-разностных уравнений соответствующего типа дают
весьма слабые результаты. Например, для $N=1$ в работе \cite{sko-spi:self}
были доказаны существование и единственность решений, аналитичных вблизи
точки $x=0$, при некоторых дополнительных ограничениях. В диссертации \cite{Liu}
доказано существование несингулярных при $x\in\R $ решений (асимптотика
этих функций при $|x|\to \infty$ всё ещё полностью не определена, в отличие от
ряда других похожих потенциалов \cite{nov}).
Некоторое представление о степени сложности общего решения возникает
при рассмотрении специальных значений параметров  $q, \mu_j$,
когда $f_j(x)$ могут быть выражены через известные специальные функции.

Пусть $q$ будет произвольным и $f_j^2(x)=
\frac{1}{1-q^2}\sum_{m=0}^{N-1}\mu_m -\sum_{m=0}^{j-1}\mu_m.$
Это приводит к $L_0=-d^2/dx^2$ или к свободной квантовомеханической частице,
которая приобретает таким образом $q$-алгебраическую интерпретацию
\cite{spi:universal}.

Предположим, что $f_0(x)$ несингулярна при $x=0$. Тогда, в так называемом
``кристаллическом" пределе $q\to 0$, потенциал $u_0(x)$ сворачивается
к общему $N$-солитонному потенциалу уравнения Кортвега-де Вриза.

Подставим теперь $q=1$ в (\ref{q-pereq}) и предположим, что
$\sum_{m=0}^{N-1}\mu_m\neq 0$. Тогда для $N=1,2$ легко вывести вид потенциалов
$u_0(x)\propto x^2, ax^2+b/x^2$. Для $N=3$ соответствующая система уравнений
на $f_j(x)$ приводит к ``циклическому" представлению уравнения Пенлеве-IV
\cite{Bu,VS}. При $N=4$ возникает функция Пенлеве-V \cite{Ad}, и т.д.
Последние два случая приводят также к коммутаторному представлению уравнений Пенлеве,
впервые замеченному в работе \cite{Fl}. При $q=-1$ возникают аналогичные ситуации
с дополнительным условием, что функции $f_j(x)$ oбладают определенной
четностью \cite{spi:universal}. Если $\sum_{m=0}^{N-1}\mu_m= 0$, то тогда
при нечетных $N$ и некоторых случаях четного $N$ потенциал $u_0(x)$
выражается через гиперэллиптические функции \cite{VS,Weiss}.

Случай, когда $q$ является примитивным корнем единицы, $q^n=1,
q\neq \pm1,$ весьма интересен \cite{sko-spi:self,sko-spi:spectra}.
Для любых нечетных $n$ и четных $n$ с определенными ограничениями
потенциалы выражаются через гиперэллиптические функции с дополнительными
кристаллографическими симметриями. Например, при $N=1$ и $q^3=1$ или $q^4=1$
возникает уравнение Ламе с эквиангармонической или лемнискатической
$\cal P$-функцией Вейерштрасса. Таким образом, это классическое
дифференциальное уравнение второго порядка оказывается связанным с теорией
представления алгебры $q$-гармонического осциллятора при $q$ равном корням единицы.

Общее семейство автомодельных потенциалов объединяет некоторые из функций Пенлеве
с гиперэллиптическими функциями. Благодаря связи с квантовыми алгебрами,
оно может быть проинтерпретировано как новый класс нелинейных $q$-специальных
функций (``непрерывные $q$-трансценденты Пенлеве"). Стандартные
\cite{sla:generalized} и базисные
гипергеометрические функции \cite{gas-rah:basic}
возникают в этом формализме через когерентные состояния.

Что касается дискретного спектра операторов Шредингера с автомодельными
потенциалами, из метода факторизации и условия (\ref{q-per}) следует, что в простейшем
случае он состоит из $N$ независимых геометрических прогрессий:
$\lambda_{pN+k}=\lambda_kq^{2p}, \; k=0,\dots, N-1,\;
p\in {\N }$. Тот же вывод следует из теории унитарных представлений
алгебры (\ref{q-alg}). Единственным условием справедливости этого
формального вывода является требование, чтобы функции $f_j(x)$ не
были сингулярными для $x\in \R $ и имели положительную при
$x\to +\infty$ и отрицательную при $x\to-\infty$ асимптотики.

\subsection{Когерентные состояния}

Существует несколько определений квантовомеханических когерентных состояний
\cite{Per}. В контексте спектр-генерирующих алгебр, таких как (\ref{pol}) и
(\ref{q-alg}), они определяются как собственные функции соответствующих
понижающих операторов. Такие функции играют роль генерирующих функций для
векторов неприводимых представлений данной алгебры. Опишем кратко когерентные
состояния квантовой алгебры \re{q-alg}, сконструированные в работах
\cite{spi:coherent,spi:universal}.

Обозначим  $|\lambda\rangle$ собственные состояния абстрактного оператора
$L$ входящего в соотношения (\ref{q-alg}),
$L|\lambda\rangle=\lambda|\lambda\rangle.$
Предположим, что операторы $B^\pm$ сопряжены друг к другу по отношению к некоторому
скалярному произведению $\langle \sigma|\lambda\rangle$.
Для простоты положим, что спектр $L$ невырожденный. Тогда операторы
$B^\pm$ действуют на $|\lambda\rangle$ следующим образом:
$$
B^-|\lambda\rangle=
\prod_{k=0}^{N-1}\sqrt{\lambda-\lambda_k}|\lambda q^{-2}\rangle, \qquad
B^+|\lambda\rangle=
\prod_{k=0}^{N-1}\sqrt{\lambda q^2-\lambda_k}|\lambda q^{2}\rangle.
$$
Пусть $N$ нечетно, $0<q^2<1$ и $\lambda_0<\ldots <\lambda_{N-1}<q^2\lambda_0<0$.
Тогда для $\lambda<0$ оператор $L$ может иметь дискретный спектр образованный
только из геометрических прогрессий (не более $N$ штук), соответствующих
унитарным неприводимым представлениям алгебры \re{q-alg}. Это следует из того
факта, что $B^-$ является понижающим оператором для $\lambda<0$ состояний и
произведение $\prod_{k=0}^{N-1}(\lambda-\lambda_k)$ становится отрицательным
при $\lambda< \lambda_0$. Поскольку $B^-|\lambda_k\rangle=0$, эта проблема
не возникает для специальных значений $\lambda$, а именно, при
$\lambda=\lambda_{pN+k}\equiv \lambda_kq^{2p}, \; p\in {\N }$.
Для нормируемых $|\lambda_k\rangle$, мы получаем последовательности
нормируемых состояний вида $|\lambda_{pN+k}\rangle \propto
\left(B^+\right)^p|\lambda_k\rangle$.

Когерентные состояния первого типа определяются как собственные
состояния $B^-$:
\begin{equation}
B^-|\alpha\rangle_-^{(k)}=\alpha|\alpha\rangle_-^{(k)},
\quad k=0, \dots, N-1,
\label{coh1}\end{equation}
где $\alpha$ некоторое комплексное число. Представляя $|\alpha\rangle_-^{(k)}$
как суперпозицию состояний $|\lambda_{pN+k}\rangle$ с различными коэффициентами,
находим
$$
|\alpha\rangle_-^{(k)}= \sum_{p=0}^\infty C_p^{(k)}\alpha^p
|\lambda_{pN+k}\rangle =C^{(k)}(\alpha)\; {_N}\varphi_{N-1}\left(
{0, \dots, 0 \atop \lambda_0/\lambda_k, \dots,
\lambda_{N-1}/\lambda_k}; q^2, z\right)|\lambda_k\rangle,
$$
где ${_N}\varphi_{N-1}$ есть стандартный односторонний $q$-гипергеометрический ряд
с операторным аргументом $z=(-1)^N\alpha B^+/\lambda_0\cdots\lambda_{N-1}$
и $C^{(k)}(\alpha)$ есть нормировочная константа (среди аргументов ряда
$_N\varphi_{N-1}$ член $\lambda_k/\lambda_k=1$ отсутствует). Верхний
индекс $k$ просто нумерует неприводимые представления
алгебры (\ref{q-alg}) с низшим весом, каждое из которых имеет свое собственное
когерентное состояние (или генерирующую функцию). Состояния
$|\alpha\rangle_-^{(k)}$ нормируемы при
$|\alpha|^2< |\lambda_0\cdots\lambda_{N-1}|$. Когерентные состояния такого типа
были впервые рассмотрены в работе \cite{AC} для специальной реализации алгебры
$q$-гармонического осциллятора (\ref{qosc}) (см. также \cite{ASu}).

Нетрудно увидеть, что нулевые моды оператора  $L$ определяют специальное
вырожденное представление алгебры (\ref{q-alg}). Поскольку операторы $B^\pm$
коммутируют с $L$ в подпространстве этих нулевых мод, то они сами по себе
могут рассматриваться как когерентные состояния.

Необычные когерентные состояния возникают из представлений без старшего веса
алгебры (\ref{q-alg}), соответствующих собственным функциям $L$ с $\lambda>0$.
Пусть $\lambda_+>0$ будет точкой дискретного спектра $L$. Тогда операторы
$B^\pm$ генерируют из $|\lambda_+\rangle$ часть дискретного спектра $L$
в виде двусторонней геометрической прогрессии $\lambda_+ q^{2n}, n\in {\Z }$.
Поскольку $\prod_{k=0}^{N-1}(\lambda_+-\lambda_k)>0 $ для произвольного $N$ и
$\lambda_+>0$, число таких неприводимых представлений в спектре $L$ не ограничено.
Более того, непрерывная прямая сумма таких представлений может образовывать
непрерывный спектр оператора $L$.

Важным фактом является то, что для состояний с $\lambda>0$ роль понижающего
оператора играет $B^+$. Поэтому мы так же можем определить новые когерентные
состояния как собственные вектора этого оператора \cite{spi:universal}:
\begin{equation}
B^+|\alpha\rangle_+=\alpha|\alpha\rangle_+.
\label{coh2}\end{equation}
Предположим, что двусторонняя последовательность $\lambda_+q^{2n}>0$ принадлежит
дискретному спектру $L$. Тогда состояния $|\alpha\rangle_+$ представляются
в виде суперпозиции состояний $|\lambda_+q^{2n}\rangle$:
$$
|\alpha\rangle_+= \sum_{n=-\infty}^\infty C_n\alpha^n|\lambda_+ q^{2n}
\rangle =C(\alpha)\; {_0\psi_N}\left({0, \dots, 0\atop \lambda_0/\lambda_+, \dots,
\lambda_{N-1}/\lambda_+}; q^2, z \right)|\lambda_+\rangle,
$$
где ${_0\psi_N}$ есть двусторонний $q$-гипергеометрический ряд с операторным
аргументом $z=\alpha B^-/(-\lambda_+)^N$, а $C(\alpha)$ есть нормировочная
константа. В этом случае состояния $|\alpha\rangle_+$ нормируемы при
$|\alpha|^2>|\lambda_0\cdots\lambda_{N-1}|$.

Предположим, что область  $\lambda>0$ занята непрерывным спектром.
Соответствующие состояния $|\lambda\rangle$ могут быть нормированы
так: $\langle \sigma|\lambda\rangle=\lambda\delta(\lambda-\sigma)$.
Возьмем $N=1$ и предположим, что состояния $|\lambda\rangle$ невырожденные.
Положив $\rho=1$ в (\ref{qosc}), мы получаем разложение
$$
|\alpha\rangle_+^{(s)}= C(\alpha)
\int_0^\infty \frac{\lambda^{\gamma_s}|\lambda\rangle d\lambda}
{\sqrt{(-\lambda q^2(1-q^2);q^2)_\infty}},
$$
где использовано обозначение $(a;q)_\infty=\prod_{j=0}^\infty(1-aq^j)$ и
$$
\gamma_s=\frac{2\pi i s-\ln(\alpha q^2\sqrt{1-q^2})}
{\ln q^2}, \qquad   s\in {\Z }.
$$
Существует счетное семейство таких когерентных состояний имеющих
единичную норму для $|\alpha|^2 > 1/(1-q^2)$ при следующем выборе
нормировочной константы $C(\alpha)$:
$$
|C(\alpha)|^{-2}=\int_0^\infty\frac{\lambda^{-\nu}d\lambda}
{(-\lambda q^2(1-q^2);q^2)_\infty}
=\frac{\pi}{\sin\pi\nu}\frac{(q^{2\nu};q^2)_\infty
(q^2(1-q^2))^{\nu-1}}{(q^2;q^2)_\infty},
$$
где $\nu=\ln|\alpha q\sqrt{1-q^2}|/\ln q.$
Этот интеграл является частным случаем $q$-бета интеграла Рамануджана
\cite{gas-rah:basic}.

Выше мы рассматривали только абстрактные когерентные состояния.
В уравнении Шредингера с автомодельными потенциалами состояния непрерывного
спектра $\lambda>0$ дважды вырождены. Структура состояний $|\alpha\rangle_\pm$
в этом случае является исключительно сложной. Даже для нулевого потенциала
(свободной частицы) возникает весьма нетривиальная ситуация.
Положим $L=-d^2/dx^2$ и
$$
B^-=U^{-1}(d/dx+1/\sqrt{1-q^2}),\qquad B^+=(-d/dx+1/\sqrt{1-q^2})U,
$$
где $U$ есть оператор растяжений, $Uf(x)=q^{1/2}f(qx), \; 0<q<1.$
Можно проверить прямым расчетом, что $B^-B^+-q^2B^+B^-=1$ и
$L=B^+B^- - 1/(1-q^2).$ Уравнение $B^-\psi_\alpha^-(x)=
\alpha\psi_\alpha^-(x)$ теперь совпадает с запаздывающим уравнением пантографа,
которое было детально изучено в работе \cite{KaM}. Как следует из ее результатов,
в запаздывающем случае уравнение пантографа не допускает решения
принадлежащих $\mbox{L}^2({\R })$. Однако, оператор $B^+$ имеет бесконечно
много нормируемых собственных состояний. Уравнение
$B^+\psi_\alpha^+(x)=\alpha\psi_\alpha^+(x)$
совпадает теперь с опережающим уравнением пантографа
\begin{equation}\label{pant}
\frac{d}{dx}\psi_\alpha^+(x)=-\alpha q^{-3/2}\psi_\alpha^+(q^{-1}x)
+\frac{q^{-1}}{\sqrt{1-q^2}}\psi_\alpha^+(x),
\end{equation}
которое имеет бесконечно много решений из  $\mbox{L}^2({\R })$ \cite{KaM}.
Для модели свободной частицы, $N>1$ алгебра симметрий приводит к
обобщенным уравнениям пантографа, рассмотренным в статье \cite{Is}.
Одной из открытых проблем в этой области остается задача описания
минимального набора решений уравнения  (\ref{pant}) составляющих полный
базис векторов Гильбертова пространства  $\mbox{L}^2({\R })$.

\subsection{Солитоны и модели статистической механики}

Рассмотрим $N$-солитонное решение уравнения Кортвега-де Вриза
$$u_N(x, t)=- \frac{d^2}{dx^2}\log \tau_N(x, t),$$
где $\tau_N$ есть Вронскиан $N$ различных решений
$\phi_1,\ldots,\phi_N$ свободного
уравнения Шредингера: $\tau_N=\det(\partial_x^{i-1}\phi_k)$.
Существует несколько детерминантных представлений этой тау-функции, например,
\cite{AS}
\be
\tau_N=\det \Phi, \qquad \Phi_{ij}=\delta_{ij}+ {2\sqrt{k_i k_j}
\over k_i+k_j} e^{(\theta_i+\theta_j)/2},
\lab{phase}\end{equation}
$$
\theta_i=k_ix-k_i^3 t +\theta_i^{(0)}, \qquad i, j=1, 2, \dots, N.
$$
Известно, что переменные $k_i$ описывают амплитуды солитонов (они связаны с
собственными значениями оператора $L_N=-d^2/dx^2+u_N(x)$ как $\lambda_i=-k_i^2/4$),
$\theta_i^{(0)}/k_i$ являются $t=0$ фазами солитонов, а $k_i^2$ -- их скоростями.

Некоторые автомодельные потенциалы, появляющиеся из $q$-периодических
замыканий (\ref{q-per}), могут рассматриваться как бесконечно-солитонные
системы. Действительно, рассмотрим $pM$-солитонные решения Кортевега-де Вриза
с параметрами $k_j$, удовлетворяющими ограничениям $k_{j+M}=qk_j, 0<q<1$, и
рассмотрим предел $p\to\infty$. Предельный потенциал имеет автомодельный дискретный
спектр $k_{pM+m}=q^pk_m, m=0, \dots, M-1,$ и бесконечное число свободных
параметров $\theta_j^{(0)}$.

Масштабно-преобразованный потенциал $\tilde u(x,t)=q^2u(qx,q^3t)$
имеет ту же самую солитонную интерпретацию с фазами
$\tilde \theta_j(x,t)=$ $\theta_j(qx,q^3t)=$
$k_{j+M}x-k_{j+M}^3t+\theta_j^{(0)}.$ Потребуем, чтобы имело место равенство
$\theta_{j+M}^{(0)}=\theta_j^{(0)}$, т.е. наложим дополнительное ограничение
$\theta_j(qx, q^3 t)=\theta_{j+M}(x,t)$, согласованное с условием (\ref{q-per}).
В этой картине дискретное преобразование $x\to qx, t\to q^3t$ просто стирает
$M$ солитонов, соответствующих наименьшим собственным значениям  оператора $L$.
В некотором смысле, эта интерпретация предоставляет наиболее простую характеризацию
автомодельных потенциалов.

Рассмотрим теперь приложения в статистической механике.
$N$-солитонная тау-фун\-кция уравнения Кортевега-де Вриза \re{phase}
может быть представлена в следующем виде, найденным Хиротой
\cite{AS,Hir}:
\begin{equation}
\tau_N = \sum_{\sigma_i=0,1} \exp \left( \sum_{0\leq i<j \le N-1 } A_{ij}
\sigma_i \sigma_j + \sum_{i=0}^{N-1} \theta_i\sigma_i \right).
\label{N_soliton}\end{equation}
Коэффициенты $A_{ij}$ описывают фазовые сдвиги при рассеянии $i$-го
и $j$-го солитонов и имеют вид
\begin{equation}
e^{A_{ij}}={ (k_i-k_j)^2 \over  (k_i+k_j)^2 }.
\label{KDV_phase}\end{equation}
Как замечено в работе \cite{lou-spi:self}, для $\theta_i=\theta^{(0)}$
эта тау-функция совпадает со статистической суммой большого канонического ансамбля
решеточного газа на линии (для двумерного Кулоновского газа, см. статью
\cite{lou-spi:soliton} и ниже). В этой интерпретации дискретные переменные
$\sigma_i$ описывают числа
заполнения молекулами ячеек решетки, $\theta^{(0)}$ есть химический потенциал,
а $A_{ij}$ есть энергия взаимодействия молекул \cite{bax:exactly}.

Известно, что простая замена переменных переводит модель решеточного
газа в модель Изинга:
\begin{equation}\label{Zising}
Z_N=\sum_{s_i=\pm1} e^{- \beta E}, \qquad
E =\sum_{0\leq i<j \leq N-1} J_{ij}s_is_j -\sum_{i=0}^{N-1} H_i s_i,
\end{equation}
где $N$ есть число спинов $s_i=\pm 1$, $J_{ij}$ являются обменными константами,
$H_i$ обозначает внешнее магнитное поле, а $\beta=1/kT$
обратную температуру. Действительно, замена чисел заполнения в \re{N_soliton}
спиновыми переменными $\sigma_i=(s_i+1)/2$ приводит к
\be
\tau_N=e^{\varphi} Z_N, \qquad \varphi=\frac{1}{4}\sum_{i<j} A_{ij}
+ \frac{1}{2}\sum_{j=0}^{N-1} \theta_j,
\lab{tauising}\end{equation}
где
\be
\beta J_{ij} =-\frac{1}{4} \; A_{ij},
\qquad
\beta H_i=\frac{1}{2}\theta_i +\frac{1}{4}\sum_{j=0, i\neq j}^{N-1} A_{ij}.
\lab{identif}\end{equation}

Аналогичная связь с нелокальными цепочками Изинга справедлива для всей
бесконечной иерархии уравнений Кадомцева-Петвиашвили и для некоторых
других дифференциальных и конечно-разностных нелинейных интегрируемых
эволюционных уравнений. Соответствующие тау-функции имеют ту же самую
форму \re{N_soliton} с более сложной структурой фазовых сдвигов $A_{ij}$
и фаз $\theta_i$ (некоторый список таких уравнений приведен в
монографии \cite{AS}).

Общая связь между $N$-солитонными решениями интегрируемых иерархий и
статистическими суммами решеточных моделей статистической механики,
установленная в работах \cite{lou-spi:self,lou-spi:spectral,lou-spi:soliton},
приводит и к некоторым новым точкам зрения на автомодельные потенциалы.
Конкретно, как показано в заметке \cite{lou-spi:self} автомодельные спектры
могут быть выведены из условия трансляционной инвариантности бесконечных
спиновых цепочек. Например, потребуем от спинового взаимодействия,
чтобы оно было инвариантно по отношению к сдвигу решетки на одну ячейку
$j\to j+1$, то есть $J_{i+1,j+1}=J_{ij}$. В результате, обменные
константы $J_{ij}$ (или солитонные сдвиги фаз $A_{ij}$) зависят только от
расстояния между ячейками $|i-j|$. Это очень простое и естественное физическое
ограничение однозначно приводит к тому, что параметры $k_i$ образуют
геометрическую прогрессию
\begin{equation}
k_i=k_0q^i, \qquad q=e^{-2\alpha}, \qquad A_{ij}=2\ln|\tanh\alpha (i-j)|,
\label{interactions} \end{equation}
где $k_0$ and $0<q<1$ являются некоторыми свободными параметрами. В более сложном
случае можно потребовать трансляционной инвариантности по отношению к сдвигам
на $M$ ячеек, или $J_{i+M,j+M}=J_{ij}$. Это приводит к автомодельным спектрам
общего вида $k_{j+M}=qk_j$. Основным недостатком описанной взаимосвязи между
солитонами и цепочками Изинга состоит в том, что при фиксированном $q$ температура
$T$ (или $\beta$) также фиксирована, что явно видно из сравнения (\ref{identif}) с
(\ref{interactions}) (в интерпретации решеточного Кулоновского газа фактически
оказывается, что $\beta=2$).

Для цепочек конечной длины $0\leq j\leq N-1$, трансляционная инвариантность
не является точной. Предел бесконечного числа солитонов $N\to\infty$ в статистической
механике соответствует термодинамическому пределу. В этой картине, все координаты
интегрируемых иерархий ($x,t$ и другие высшие ``времена") интерпретируются как
специфические параметры, описывающие неоднородное внешнее магнитное поле $H_i$.
Поскольку $0<q<1$, то $x$ и $t$-зависимые части $H_i$ вымирают экспоненциально быстро
при $i\to\infty$. Поэтому, в пределе $N\to\infty$ только постоянные $\theta_i^{(0)}$
существенны для функции распределения (точнее, они определяют лидирующую асимптотику
функции распределения при $N\to\infty$). Несмотря на то, что в термодинамическом пределе
зависимость от $x$ и $t$ вымывается и обменное взаимодействие трансляционно инвариантно,
наличие бесконечного числа параметров $\theta_i^{(0)}$ не позволяет найти замкнутое
выражение для лидирующего члена $Z_N$. Однако, для автомодельных бесконечно-солитонных
систем, характеризующихся периодичностью фаз $\theta_i^{(0)}$ или, эквивалентно,
периодичностью внешнего магнитного поля $H_{i+M}=H_i$, такие выражения могут
быть найдены.

Рассмотрим случай $M=1$ с однородным магнитным полем $H_i=H$.
Тогда с уравнением Кортевега-де Вриза ассоциируется антиферромагнитная цепочка
Изинга: $0<|\tanh \alpha(i-j)|<1$ и $J_{ij} = -A_{ij}/4\beta > 0$
(аналогичная картина справедлива и для произвольного периода $M$).
Обменное взаимодействие нелокально, но оно спадает экспоненциально быстро и, поэтому,
фазовые переходы отсутствуют для ненулевых температур.

Для вычисления функции распределения можно использовать детерминантные
представления для тау-функций. Как показано в работах
\cite{lou-spi:self,lou-spi:spectral}, для трансляционно-инвариантных
цепочек Изинга в однородном магнитном поле тау-функции становятся детерминантами
некоторых матриц Тёплица и их естественных $(M\times M)$-блочных обобщений.
Например, для $M=1$ был получен следующий результат:
$Z_N\to \exp(-N\beta f_I)$ при $N\to\infty$,
где удельная свободная энергия $f_I$ имеет вид
\begin{eqnarray}
&& -\beta f_I(q, H)=\ln \frac{2(q^4;q^4)_\infty \cosh \beta H}
{(q^2;q^2)_\infty^{1/2}}
+ \frac{1}{4\pi}\int_0^{2\pi}d\nu \ln (|\rho(\nu)|^2 - q\tanh^2 \beta H),
\lab{free}\\ && \makebox[2em]{}
|\rho(\nu)|^2=
\frac{(q^2e^{i\nu};q^4)_\infty^2(q^2e^{-i\nu};q^4)_\infty^2}
{(q^4e^{i\nu};q^4)_\infty^2(q^4e^{-i\nu};q^4)_\infty^2}\;
\frac{1}{4\sin^2(\nu/2)}=
q\frac{\theta_4^2(\nu/2, q^2)}{\theta_1^2(\nu/2, q^2)}
\nonumber\end{eqnarray}
(здесь $\theta_{1,4}$ обозначают стандартные тета-функции Якоби).
При выводе этой формулы использовалась $_1\psi_1$ сумма Рамануджана
на границе сходимости соответствующего двустороннего базисного
гипергеометрического ряда. Маленький вспомогательный  параметр
$\epsilon$, введенный в выражение для плотности $\rho(\nu)$,
\begin{equation}
\rho(\nu)=\frac{(q^2;q^4)_\infty^2}{(q^4;q^4)_\infty^2}
\sum_{k=-\infty}^\infty
\frac{e^{i\nu k-\epsilon k}}{1-q^{4k+2}}=
\frac{(q^2e^{i\nu-\epsilon};q^4)_\infty (q^2e^{-i\nu+\epsilon};q^4)_\infty}
{(e^{i\nu-\epsilon};q^4)_\infty (q^4e^{-i\nu+\epsilon};q^4)_\infty},
\label{nu}
\end{equation}
гарантирует абсолютную сходимость ряда при $\epsilon> 0$. Несмотря на то, что
в пределе $\epsilon\to 0$ возникает особенность по $\nu$
в правой части выражения (\ref{nu}), это не приводит к сингулярностям после
взятия интеграла в (\ref{free}). Полная намагниченность цепочки
$m(H)=-\partial_H f_I = \stackreb{\lim}{N\to\infty} N^{-1}
\sum_{i=0}^{N-1}\langle s_i \rangle$ принимает следующий привлекательный вид
\be
m(H)=\left(1-\frac{1}{\pi}
\int_0^{\pi}\frac{\theta_1^2(\nu, q^2)d\nu}{\theta_4^2(\nu, q^2)\cosh^2\beta H
-\theta_1^2(\nu, q^2) \sinh^2\beta H} \right)\tanh \beta H.
\lab{magkdv}\end{equation}

Для изменения фиксированного значения температуры можно попробовать заменить
(\ref{interactions}) выражением $A_{ij}=2 n\ln|(k_i-k_j)/(k_i+k_j)|,$
для некоторой последовательности целых чисел $n$, и поискать интегрируемые уравнения
допускающие $N$-солитонные решения с такими фазовыми сдвигами. Это крайне
нетривиальная задача и в работах \cite{lou-spi:self,lou-spi:spectral}
было найдено только еще одно допустимое значение температуры.
Оно соответствует $n=2$ и возникает из специальной редукции многосолитонного решения
уравнения Кадомцева-Петвиашвили B-типа \cite{DJKM}. Общая
тау-функция для указанного уравнения генерирует значительно более сложные
цепочки Изинга чем для уравнения Кортевега-де Вриза.

Известно, что функция распределения вероятностей в теории случайных матриц
описывает статистическую сумму $n$-частичного газа с логарифмическим взаимодействием,
которая определяется инте\-гра\-лом по
кон\-фи\-гу\-ра\-ци\-он\-но\-му про\-стран\-ст\-ву.
Поэтому описанные выше модели Изинга (или решеточного газа) связаны с моделями
случайных матриц с дискретным набором собственным значений. Опишем кратко это
соответствие сравнивая многосолитонные решения уравнения Кадомцева-Петвиашвили
и Дайсоновский ансамбль на круге.

Модель случайных $n\times n$ унитарных матриц была предложена Дайсоном
в качестве альтернативы Вигнеровскому ансамблю эрмитовых матриц.
Соответствующие собственные значения имеют вид $\epsilon_j=e^{i\phi_j}$,
$j=1, \dots, n$, и, после интегрирования по вспомогательным ``угловым" переменным,
распределение вероятностей по фазам $0\leq \phi_j<2\pi$
принимает вид $Pd\phi_1 \dots d\phi_n\propto
\prod_{i<j}|\epsilon_i-\epsilon_j|^2 d\phi_1 \dots d\phi_n.$
Можно ослабить некоторые из требований, использованных Дайсоном, и работать с
функциями распределения более общего вида. В работе \cite{Gaud1}
Годен предложил круговой ансамбль с вероятностным распределением
\begin{equation}
Pd\phi_1 \dots d\phi_n\propto
\prod_{i<j}\left|\frac{\epsilon_i-\epsilon_j}
{\epsilon_i-\omega\epsilon_j}\right|^2 d\phi_1 \dots d\phi_n,
\label{distr1}
\end{equation}
содержащем свободный непрерывный параметр $\omega$. Он интерполирует между
моделью Дайсона  ($\omega=0$) и однородным распределением ($\omega=1$).
Модель (\ref{distr1}) допускает также интерпретацию решеточного газа на круге
со статистической суммой
$$
Z_n\propto \int_0^{2\pi}\dots\int_0^{2\pi}d\phi_1\dots d\phi_n
\exp\left(-\beta \sum_{i<j}V(\phi_i-\phi_j)\right),
$$
где $\beta=1/kT=2 $ фиксирована и потенциальная энергия имеет вид
\begin{equation}
\beta V(\phi_i-\phi_j)=
\ln\left(1+\frac{\sinh^2\gamma}{\sin^2((\phi_i-\phi_j)/2)}\right),
\qquad \omega=e^{-2\gamma}.
\label{gas1}
\end{equation}

Как это было обнаружено в \cite{lou-spi:spectral}, статистическая сумма большого
канонического ансамбля для этой модели может быть получена из тау-функции
специаль\-ного бесконечно-соли\-тон\-но\-го решения иерархии уравнений Кадомцева-Петвиашвили.
При этом конечно-соли\-тон\-ные ре\-ше\-ния опре\-де\-ля\-ют дис\-кре\-ти\-за\-цию этой модели,
то есть решеточную модель случайных матриц с дискретным набором собственных значений.
Например, можно взять унитарную $n\times n$ матрицу с собственными значениями
равными $N$-м корням единицы, $\epsilon_j=\exp(2\pi im_j/N),\; m_j=0, \dots, N-1 $.
Вероятностная мера непрерывна по вспомогательным ``угловым" переменным унитарных матриц
и дискретна по фазам собственных значений $\phi_j$. При этом интегралы по $\phi_j$
заменяются на конечные суммы по $m_j$ и непрерывная модель воспроизводится в пределе
$m_j, N\to\infty$ при конечных $m_j/N$:
\begin{equation}
\left(\frac{2\pi}{N}\right)^n\sum_{m_1=0}^{N-1} \dots \sum_{m_n=0}^{N-1}
\stackreb{\to}{N\to\infty} \int_{0}^{2\pi}d\phi_1 \dots
\int_{0}^{2\pi}d\phi_n.
\label{sum}\end{equation}
$n$-Частичная статистическая сумма равна
$$
Z_n(N,\omega)=\left(\frac{2\pi}{N}\right)^n\sum_{m_1=0}^{N-1} \dots
\sum_{m_n=0}^{N-1} \prod_{1 \le i < j \le n}
\left|\frac{\epsilon_i-\epsilon_j}
{\epsilon_i/\sqrt{\omega}-\sqrt{\omega}\epsilon_j}\right|^2,
$$
а для большого канонического ансамбля статистическая сумма принимает вид
\begin{equation}
Z(\omega,\theta)=\sum_{n=0}^N \frac{Z_n(N,\omega)e^{\theta n}}{n!} =
\sum_{\sigma_m=0,1} \exp \left( \sum_{0 \le m<k \le N-1}A_{mk}
\sigma_m\sigma_k+(\theta+\eta) \sum_{m=0}^{N-1} \sigma_m \right),
\label{grand}\end{equation}
где $\eta=\ln({2\pi}/{N})$ есть химический потенциал и параметры
$$
A_{mk} = \ln \frac{\sin^2(\pi(m-k)/N)}{\sin^2(\pi(m-k)/N) +\sinh^2\gamma}
= \ln \frac{(a_m-a_k)(b_m-b_k)}{(a_m+b_k)(b_m+a_k)},
$$
равны сдвигам фаз солитонов уравнения Кадомцева-Петвиашвили с
\begin{equation}
a_m=e^{2\pi i m/N}, \qquad b_m=-\omega a_m, \quad m=0,1, \dots, N-1.
\label{momentums}\end{equation}
Дискретизацию случайных матриц на круге с помощью корней единицы рассматривал
ещё сам Годен \cite{Gaud2}, где также была наблюдена связь с моделями Изинга.
Связь с солитонными интегрируемыми уравнениями была установлена значительно
позднее \cite{lou-spi:spectral}.

Опишем теперь кратко связь вышеперечисленных конструкций с двумерным Кулоновским газом.
Решение уравнения Пуассона на плоскости,
$\Delta V({\bf r},{\bf r}^\prime)= - 2\pi\delta({\bf r}-{\bf r}^\prime),$
определяет электростатический потенциал $V(z,z^\prime)=-\ln |z-z^\prime|$,
созданный заряженной частицей, находящейся в точке $z'=x'+iy'$.
В областях с диэлектрической или проводящей границами потенциал имеет более
сложный вид поскольку нормальная компонента электрического поля
${\cal \bf E} = -{\bf \nabla} V$ должна зануляться на поверхности
диэлектриков, ${\cal E}_n = 0,$ в то время как касательная компонента
исчезает на поверхности металлов, ${\cal E}_t =0.$ Введение искусственных
изображений зарядов может упростить решение уравнения Пуассона при
простых геометрических конфигурациях границ.

Энергия электростатической системы из $N$ частиц, находящихся
в ограниченной области на плоскости равна
\begin{equation}
E_N=\sum_{1\le i<j\le N} q_iq_jV(z_i,z_j)+\sum_{1\le i\le N} q_i^2v(z_i)
+\sum_{1\le i\le N} q_i\phi(z_i),
\label{energy}
\end{equation}
где $z_j=x_j+i y_j$ и $q_j$ обозначают координаты и заряды частиц.
Первый член представляет собой стандартную Кулоновскую энергию,
второй описывает взаимодействие с границами (или взаимодействие зарядов
с их образами) и последний член соответствует вкладу внешних полей.

Предположим, что наша система состоит из частиц с равными зарядами,
$q_j=+1$ (в однокомпонентном случае) или $q_j=\pm 1$ (в двухкомпонентном случае), находящимися
на дискретной решетке $\Gamma$ на плоскости. Статистическая сумма большого
канонического  ансамбля такого двумерного решеточного Кулоновского газа выглядит
следующим образом
$$
G_N=\sum_{n=0}^{N}\frac{e^{\mu n}}{n!}\sum_{z_1\in \Gamma}
\dots \sum_{z_n\in \Gamma}  e^{-\beta E_n},
$$
где $\mu$ обозначает химический потенциал. В общем случае
$G_N$ можно представить в виде \cite{bax:exactly}
\begin{equation}\label{G}
G_N=\sum_{\{\sigma(z)\}}\exp\Bigl(\frac{1}{2}
\sum_{z\neq z^\prime}W(z,z^\prime)\sigma(z)\sigma(z^\prime)
+\sum_{z\in \Gamma}w(z)\sigma(z)\Bigr)
\end{equation}
с эффективным потенциалом
$$
W\left(z,z^\prime\right)=-\beta q(z)q(z^\prime)V(z,z^\prime), \quad
w(z)=\mu(z)-\beta\left(q^2(z)v(z)+q(z)\phi(z)\right),
$$
где $\sigma(z)=0$ или 1 в зависимости от того свободна или занята ячейка
с координатой $z$, а функции $q(z)$ и $\mu(z)$ характеризуют распределение
частиц различных типов. Так, в двухкомпонентном случае, когда
$q({z_\pm})={\pm 1}$ заряды занимают подрешетки $\{z_\pm\}$,
имеем $\mu({z_\pm})=\mu_{\pm}$.

Рассмотрим теперь $N$-солитонную тау-функцию (\ref{N_soliton})
и заменим в ней номер солитона $j$ переменной $z$, принимающей
$N$ дискретных значений. После этого, она может быть переписана в виде
\begin{equation}
\tau_N=\sum_{\sigma(z)=0,1}\exp\Bigl(\frac{1}{2}\sum_{z\neq z^\prime}A_{z
z^\prime}\sigma(z)\sigma(z^\prime)+\sum_{\{z\}}
\theta(z)\sigma(z)\Bigr).
\label{tau}\end{equation}
Это выражение совпадает с (\ref{G}) при отождествлении фазовых сдвигов
$A_{zz^\prime}$ с потенциалом Кулоновского взаимодействия
$W(z,z^\prime)$ и фаз $\theta(z)$ с функцией $w(z)$. Благодаря этому
наблюдению, сделанному в работе \cite{lou-spi:soliton}, можно построить
много новых точно решаемых моделей Кулоновского газа
в дополнение к уже известным случаям
(см., например, статьи \cite{FJT,Gaud1,Gaud2} и указанные в них ссылки).

Для солитонных решений иерархии уравнений Кадомцева-Петвиашвили
фазовые сдви\-ги имеют вид \cite{DJKM}
\begin{equation}
A_{zz^\prime}=\ln\frac{(a_z-a_{z^\prime})(b_z-b_{z^\prime})}
{(a_z+b_{z^\prime})(b_z+a_{z^\prime})}, \qquad
\theta(z)=\theta^{(0)}(z)+\sum_{p=1}^\infty(a_z^p-(-b_z)^p)t_p,
\label{kp}\end{equation}
где $t_p$ обозначает $p$-е ``время" иерархии,
а $a_z$ и $b_z$ произвольные функции $z$.
Если положить $a_z=z=x+i y, \; b_z=-z^*=-x +i y,$ то тогда
\begin{equation}\label{kpint}
A_{zz'}=W(z,z^\prime)=-2 V(z,z')=
2\ln |z-z^\prime| -2\ln |z^*-z^\prime|,
\end{equation}
где $V(z,z')$ обозначает потенциал, созданный единичным положительным
зарядом, помещенным в точку $z'$ над проводником, занимающим нижнюю
полуплоскость с границей вдоль $y= 0$ линии. В этом случае
$V(z,z')$ решает уравнение Пуассона с граничным условием ${\cal E}_t(y=0)=0$.
Тот же самый потенциал создается положительным зарядом
в точке $z'$ и его образом противоположного заряда, расположенном в
точке $\left(z'\right)^*$.

Аналогично ситуации с моделями случайных матриц, соответствие между
моделями Кулоновского газа и солитонными решениями справедливо только
при фиксированной температуре $\beta$, которая находится из сравнения
выражений (\ref{kpint}) с (\ref{G}) и равна $\beta=2.$
Поскольку $w(z)=\theta(z)$, мы имеем для фаз при нулевом времени
в (\ref{kp}):
\begin{equation}
\theta^{(0)}(z)=\mu-\beta(\ln |z^*-z|+ \phi(z)),\qquad \beta=2,
\label{fkp}
\end{equation}
где второй член соответствует взаимодействию ``заряд-образ". Функция
$\phi(z)$ описывает потенциал, созданный нейтрализующим фоном с некоторой
плотностью $\rho({\bf r})$, то есть $\Delta\phi({\bf r})=-2\pi\rho({\bf r})$
с граничным условием $\phi_x(y=0)=0.$
Гармонический член $\sum_{p=0}^\infty (z^p-(z^*)^p)t_p=-\beta \phi_{ext}(z)$
в (\ref{kp}) соответствует внешнему электрическому полю. Таким образом,
при мнимых временах эволюции солитонные решения иерархии Кадомцева-Пет\-виа\-шви\-ли
описывают поведение двумерного решеточного Кулоновского газа в переменном внешнем
электрическом поле. При этом редукция к уравнению Кортевега-де Вриза,
связанному с уравнением Шредингера, соответствует расположению заряженных частиц
на линии с фиксированным значением $x$.

Рассматривая другие солитонные уравнения можно получить модели
Кулоновского газа, помещенного в одной четверти плоскости с диэлектрическими
и проводящими границами, модели дипольного газа и т.д. С помощью конформных
преобразований $z\to f(z)$ можно поместить частицы в области достаточно
сложной конфигурации (углы с наклоном $\pi/n,\; n>2$, в трубки с параллельными
границами, прямоугольники, и тому подобное). Более детальное рассмотрение этих
вопросов проведено в работе \cite{lou-spi:soliton}.

\section{Конечно-разностное уравнение Шредингера}

\subsection{Дискретная факторизационная цепочка}

Линейное конечно-разностное уравнение второго порядка
\be
L\psi(x)\equiv a(x+1)\psi(x+1)+a(x)\psi(x-1)+b(x)\psi(x)=\lambda\psi(x),
\lab{jac} \end{equation}
может интерпретироваться как уравнение, определяющее собственные
частоты гармонических колебаний неоднородной дискретной струны
или как уравнения на собственные значения, определяющее допустимые энергии
частицы движущейся вдоль некоторой неоднородной решеткм (модель сильной связи).
Помимо этого, уравнение \re{jac} может рассматриваться как вспомогательная
спектральная задача, необходимая для интегрирования уравнений движения
цепочки Тоды.

Заменяя $x\pm1$ в \re{jac} на $x\pm h$ и взяв предел нулевого расстояния между
ячейками решетки $h\to 0$, мы получаем непрерывное уравнение Штурма-Лиувилля
$$h^2(a_0(x)\psi'(x))'+(a(x+h)+a(x)+b(x)-\lambda)\psi(x)+O(h^3)=0,$$
где штрихи обозначают производную $d/dx$ и $a_0(x)$ обозначает лидирующий член
асимптотики $a(x),\, a(x)=a_0(x)+O(h)$. Если $a_0=const$ и асимптотические
разложения $a(x)$ и $b(x)$ подобраны соответствующим образом, то
возникает стандартное уравнение Шредингера:
$-\psi''(x)+u(x)\psi(x)=\tilde\lambda\psi(x),$
рассмотренное в предыдущих параграфах с точки зрения метода факторизации.
В этом параграфе мы опишем конечно-разностные операторы  Шредингера с
автомодельными потенциалами, сконструированные в работе \cite{svz:difference}
с помощью того же самого метода.

Уравнение \re{jac} необходимо дополнить граничными условиями. Пусть
$\Gamma_c$ описывает некоторую координатную решетку, то есть набор
дискретных и непрерывных значений $x$, на которых определены $a(x)$ и $b(x)$, а
$\Gamma_s$ описывает решетку спектрального  параметра, то есть
набор значений $\lambda$ в задаче на собственные значения
$L\psi_n=\lambda_n\psi_n,\; n\in\Gamma_s,$ с выбранными граничными условиями.
Чаще всего встречаются $\Gamma_s$ состоящие из набора дискретных точек
и непрерывной части идущей от какой-либо точки до бесконечности,
а $\Gamma_c$ состоящие из бесконечного набора точек, отрезка, половины
или всей прямой линии. Для конечных $\Gamma_c$ стандартны граничные условия
$a(0)\psi(-1)=0,$  $\psi(0)\neq 0,$ $a(x_{max})\psi(x_{max})=0.$
Если $\Gamma_c$ идет от $0$ до бесконечности, то граничные условия в нуле
те же самые и дополнительно требуется ограниченность $\psi(x)$.
Если края $\Gamma_c$ определяются двумя ближайшими нулями $a(x)$,
тогда достаточно потребовать конечности $\psi(x)$ в этих точках.
Если $\Gamma_c=\R$, то рассматриваются только ограниченные $\psi(x)$.

Если $x\in \N$, то есть $\Gamma_c$ дискретна и имеет границу, то
уравнение \re{jac} определяет трехчленное рекуррентное соотношение
для ортогональных полиномов степени $x$ с аргументом равным $\lambda$,
$\psi(x)=p_x(\lambda)$. Возможна и обратная ситуация, когда для
некоторой последовательности $\lambda=\lambda_n,\, n\in \N$,
решения \re{jac} генерируют полиномы степени $n$ от некоторой функции
$z(x)$, $\psi_n(x)=p_n(z(x))\psi_0(x)$. Для задач с чисто дискретным
спектром, когда $\psi_n(x)\in \ell^2(\Gamma_c)$ и $x\in \Z$,
соотношения ортогональности и полноты имеют вид:
\be
\sum_{x\in \Gamma_c} \psi_n^*(x)\psi_m(x)=\delta_{nm}, \qquad
\sum_{n\in \Gamma_s}\psi^*_n(x)\psi_n(y)=\delta_{xy},
\lab{comp}\end{equation}
где $\delta_{xy}=0$ если $x\neq y$ или $1$ при $x=y$. Если $\Gamma_c=\Gamma_s$
и $x$ и $n$ входят симметричным образом (их перестановка эквивалентна
перестановке параметров), то такие системы называются самодуальными.

В теории ортогональных полиномов дискретные аналоги преобразований Дарбу рассматривались
Кристоффелем \cite{Sz} и Геронимусом \cite{Ger1,Ger2}.
Факторизация конеч\-но-раз\-ност\-ных операторов рассматривалась в статьях
\cite{Mi1,Mi2,ata-sus:realization} и некоторых других исследованиях. Подход работы
\cite{svz:difference} (см. также \cite{nov-dyn}) отличается от предыдущих тем,
что в нем не находятся симметрии заданных систем, а описываются целые семейства
спектральных задач с заданными свойствами симметрии.
Простейшие системы, возникающие таким образом, описываются ортогональными
полиномами Шарлье, Кравчука, Мейкснера или их $q$-аналогами; полиномы
Стильтьеса-Вигерта и $q$-Эрмита также легко воспроизводятся в этой схеме
(см. \cite{koe-swa:askey} для определения этих полиномов).
Более сложные системы связаны с дискретными трансцендентами Пенлеве.

Рассмотрим некоторое множество уравнений, имеющих вид \re{jac}:
\be
L_j\psi^{(j)}(x)=\lambda\psi^{(j)}(x), \qquad j\in {\Z },
\lab{set} \end{equation}
где
\be
L_j=a_j(x+1)T^++a_j(x)T^- +b_j(x), \qquad T^\pm\psi(x)=\psi(x\pm 1).
\lab{ham}\end{equation}
Мы подразумеваем, что переменная $x$ действительная, так что операторы
$L_j$ формально являются эрмитовыми. Гамильтонианы $L_j$
можно факторизовать следующим образом
\be
L_j=A_j^+A_j^-+\lambda_j,
\lab{fac} \end{equation}
где
\be
A_j^+=p_j(x)T^-+f_j(x),\qquad A_j^-=p_j(x+1)T^+ + f_j(x).
\lab{oper} \end{equation}
Предположим, что эти операторы эрмитово сопряжены друг с другом по отношению
к некоторому скалярному произведению, $(A_j^\pm)^\dagger=A_j^\mp$.  Из
определений \re{ham}-\re{oper} находится связь между ``потенциалами"
$a_j(x)$, $b_j(x)$ и ``суперпотенциалами" $p_j(x), f_j(x)$:
\be
a_j(x)=p_j(x)f_j(x-1),\qquad b_j(x)=p^2_j(x)+f^2_j(x)+\lambda_j.
\lab{rel} \end{equation}

Факторизационная цепочка
\be
L_{j+1}=A^+_{j+1}A^-_{j+1}+\lambda_{j+1}=
(-1)^{\s_j}(A^-_jA^+_j+\lambda_j),
\lab{f-chain} \end{equation}
связывает соседние операторы $L_j$ и $L_{j+1}$. Она отличается от
\re{fchain} наличием знаков $(-1)^{\s_j}$. В принципе они могут быть
устранены изменением знаков $a_{j+1}(x),$ $b_{j+1}(x)$ и $\lambda_{j+1}$,
но иногда удобно считать, что $a_j(x)>0$ (например, при рассмотрении
операторов с ограниченным снизу спектром). В непрерывном случае такие
изменения знаков были запрещены, так как вид кинетического члена в
$L_j$ был жестко зафиксирован.

Подставив \re{oper} в \re{f-chain}, мы получаем (1+1)-мерную полностью
дискретную цепочку спектральных преобразований:
\ba
p_j(x)f_j(x)&=&(-1)^{\s_j} p_{j+1}(x)f_{j+1}(x-1),
\lab{d1} \\
p_j^2(x+1)+f^2_j(x)&=&(-1)^{\s_j} (p^2_{j+1}(x)+f^2_{j+1}(x))+\mu_j,
\lab{d2} \\
\mu_j&\equiv&(-1)^{\s_j}\lambda_{j+1}-\lambda_j.
\nonumber
\ea

Как указывалось ранее, при реализации $A_j^\pm$ дифференциальными
операторами  первого порядка дискретные спектры операторов $A_j^+A_j^-$ и
$A_j^-A_j^+$ могут отличаться только наименьшим собственным  значением.
В дискретном случае пространство нулевых мод операторов $A^\pm$
одномерно, поэтому это утверждение по прежнему остается справедливым.
Операторы
\be
M^+_j=A_j^+ A^+_{j+1}\cdots A^+_{j+N-1}, \qquad  M_j^-=(M^+_j)^\dagger,
\lab{prod} \end{equation}
для произвольного целого $N>0$ удовлетворяют соотношениям
\ba
L_j M^+_j&=&(-1)^{s_j}M^+_jL_{j+N},
\qquad s_j=\sum_{k=0}^{N-1}\s_{j+k}
\nonumber \\
M^-_j L_j&=&(-1)^{s_j}L_{j+N}M^-_j,
\lab{lad} \ea
которые показывают, что $M^\pm_j$ формально отображают собственные функции
операторов $L_j$ и $L_{j+N}$ друг на друга. Легко увидеть так же, что
\ba
M_j^+ M_j^-&=&\prod_{k=0}^{N-1}((-1)^{s_j-s_{jk}} L_j-\lambda_{j+k}), \qquad
s_{jk}=\sum_{l=k}^{N-1} \s_{j+l}, \nonumber \\
M_j^-M_j^+&=&\prod_{k=0}^{N-1}((-1)^{s_{jk}} L_{j+N}-\lambda_{j+k}).
\lab{alg} \ea
После наложения условий замыкания,
\be
L_{j+N}=q U L_jU^{-1}+\lambda_{j+N}-q\lambda_j,
\lab{closure} \end{equation}
где $q$ обозначает произвольный параметр и $U$ есть унитарный
оператор, мы получим некоторую полиномиальную алгебру симметрии.
Комбинации $B^+_j\equiv M_j^+ U$ и $B^-_j\equiv U^{-1}M_j^-$
становятся операторами симметрии Гамильтониана  $L_j$.
Выбор оператора $U$ достаточно произволен, но, по аналогии с
непрерывным случаем, ограничимся генераторами аффинных преобразований.

Полагая $B^\pm\equiv B^\pm_1,$  $\omega\equiv \lambda_{N+1}-
q\lambda_1,$ $H\equiv L_1, s_1=\sum_{l=1}^N \s_l=0$ и
подставляя \re{closure} в \re{lad} и \re{alg}, получаем
\begin{eqnarray}
&& HB^+ -q B^+H=\omega B^+,  \qquad
B^- H-q HB^-=\omega B^-,
\lab{pol1} \\
&& B^+B^-=(-1)^S\prod_{k=1}^{N}(H- \tau_k),
\lab{pol2}
\qquad
B^-B^+=(-1)^S\prod_{k=1}^{N}(q H+\omega -\tau_k),
\end{eqnarray}
где
$$
\tau_k=(-1)^{\sum_{l=k}^N \s_l}\lambda_k, \qquad
S=\sum_{k=1}^N\sum_{l=k}^N \s_l.
$$

Предположим, что $\tau_k<\tau_{k+1}$, $k=1,\dots, N,$ тогда уравнение
$B^-\psi_0^{(k)}=0$ определяет $N$ ``вакуумных" состояний с энергиями
$\tau_k$ (при условии, что они удовлетворяют необходимым граничным условиям).
Действуя повышающим оператором $B^+$ на эти состояния,
$\psi_m^{(k)}=(B^+)^m\psi_0^{(k)}$, мы получаем физические связанные состояния.
Для $q=1$, спектр $H$ состоит из $N$ независимых арифметических прогрессий
с шагом $\omega$. Для $0<q<1$, спектр может иметь достаточно сложный вид,
в частности, его дискретная часть может состоять из $N$ геометрических
прогрессий с точкой накопления $\lambda_1-\omega/(1-q)$.
Для $q>1$, спектр чисто дискретный и растет экспоненциально, что невозможно
в случае обычного уравнения Шредингера.

Выберем $U\equiv T_{\delta}^+$,
$T_{\delta}^\pm\psi(x)=\psi(x\pm\delta)$, где $\delta$ произвольный
действительный параметр. Тогда замыкание \re{closure} эквивалентно
следующим условиям:
\begin{equation}
p_{j+N}(x)=\sqrt{q} p_j(x+\delta), \quad
f_{j+N}(x)=\sqrt{q} f_j(x+\delta), \quad
\mu_{j+N}=q\mu_j,  \quad \s_{j+N}=\s_j.
\lab{clos} \end{equation}
Потребовав, чтобы $\delta$ было дробным числом, мы получаем систему обычных
разностных уравнений, методы интегрирования которых хорошо разработаны.

\subsection{Автомодельные редукции, дискретные уравнения Пенлеве и ортогональные полиномы}

Рассмотрим дискретную цепочку \re{d1}, \re{d2} с замыканием \re{clos}
для некоторых простых значений $N$ и $\delta$. Для $N=2$, в работе \cite{svz:difference}
было найдено только одно значение $\delta$, а именно, $\delta=1$, для которого
автомодельные потенциалы описываются в терминах элементарных функций.
Соответствующая система уравнений имеет вид:
\ba
p_1(x)f_1(x)&=&\pm \; p_2(x)f_2(x-1),\nonumber \\
p_2(x)f_2(x)&=&\pm \; q\; p_1(x+1)f_1(x), \lab{syst1} \\
p_1^2(x+1)+f^2_1(x)&=&\pm\; (p^2_2(x)+f_2^2(x))+\mu_1, \nonumber \\
p_2^2(x+1)+f^2_2(x)&=&\pm\; q\; (p^2_1(x+1)+f_1^2(x+1))+\mu_2.
\lab{syst2} \ea
Обозначив $F_j\equiv f_j^2, P_j\equiv p_j^2$, можно получить решение:
\be
P_1(x)={\mu_1\pm\mu_2+\gamma^2(\mu_2\pm\mu_1q)q^{2x-1}+cq^x\over
(1-q)(1\mp\gamma^2q^{2x-1})(1\mp\gamma^2q^{2x})}, \lab{mei1} \end{equation}
\be
F_1(x)=\gamma^2q^{2x}\;{\mu_2\pm\mu_1q+\gamma^2(\mu_1\pm\mu_2)q^{2x+1}\pm
cq^x\over (1-q)(1\mp\gamma^2q^{2x+1})(1\mp\gamma^2q^{2x})},
\lab{mei2} \end{equation}
\smallskip
\be
F_2(x)=\gamma^2 q^{2(x+1)}P_1(x+1),\qquad P_2(x)=\gamma^{-2}q^{-2x}F_1(x),
\lab{mei3} \end{equation}
где $\gamma^2$ и $c$ являются двумя константами интегрирования.
В действительности, $\gamma^2$ и $c$ могут быть произвольными периодическими
функциями с периодом $1$, но мы пренебрегаем такой свободой.
Для общего значения параметров $\mu_{1,2}, \gamma^2$ в \re{mei1}-\re{mei3}
постоянная $c$ определяется из условия $a(0)=0$, подразумевающего, что
$\Gamma_c$ занимает половину линии или конечный интервал.

При выборе верхнего (положительного) знака, волновые функции
$\psi(x)$ связаны с $q$-аналогами полиномов Мейкснера и спектр генери\-рую\-щей
квантовой алгеброй $su_q(1,1)$ (эта реализация была найдена еще в работе \ci{alexei}).
Для нижнего знака возникает $su_q(2)$ алгебра и формулы \re{mei1}-\re{mei3} приводят
к $q$-аналогам дискретных полиномов Кравчука, которые могут быть выражены через
$_3\varphi_2$ базисные гипергеометрические ряды \ci{stanton}. В пределе $q\to 1$,
возникают стандартные полиномы Мейкснера и Кравчука.

Частный случай \re{mei1}-\re{mei3} связан с моделью $q$-гармонического осциллятора.
Действительно, подстановка $q\to q^2,\;
\mu_2\to q\mu_1,$ приводит к системе с $N=1,\; \delta=1/2$:
\be
P_1(x)={\mu_1(1+\gamma^2q^{4x-1})+cq^{2x}(1+q)^{-1}\over
(1-q)(1-\gamma^2q^{4x-2})(1-\gamma^2q^{4x})},
\qquad F_1(x)=\gamma^2 q^{4x+1}P_1(x+\fac12).
\lab{qhermite}\end{equation}
Существует специальный выбор параметра $c$,
а именно, $c=-\mu_1\gamma(1+q)^2/q,$  при котором сингулярности знаменателей
$P_1(x)$ и $F_1(x)$ сокращаются с нулями числителей.
Потенциалы $a(x)$ и $b(x)$ становятся определенными на всей линии ($\Gamma_c=\R$)
и приводят к так называемым непрерывным полиномам $q$-Эрмита.

Наиболее простая система возникает для $N=1,\; \delta=0$:
\be
F_1(x)=\gamma^2 q^{2x},\qquad
P_1(x)=(1-q)^{-1}(\mu_1+\gamma^2(1-q)q^{2x-1}+cq^x).
\lab{char}\end{equation}
Имеется две возможности. Если $P_1(x)$ не имеет нулей непрерывного аргумента
$x$, тогда выбор $\Gamma_c=\Z$ является естественным. В противном случае параметр
$c$ зафиксирован требованием $P_1(0)=0$ и тогда $\Gamma_c= \N$ и эта система
связана с $q$-полиномами Шарлье.

Система уравнений, возникающая из $N=2,\; \delta=0$, $q\neq 1$, допускает только один
интеграл. При $q=1$, появляется второй интеграл и мы получаем:
$$
F_2(x)={\gamma^2\over F_1(x)},\qquad
P_2(x)=\mu_\pm x+c \mp P_1(x),\qquad
P_1(x)={\gamma^2 (\mu_\pm x+c)\over F_1(x)F_1(x-1)\pm\gamma^2},
$$
и уравнение:
\be
{\gamma^2(\mu_\pm x+c)\over F_1(x)F_1(x-1)\pm\gamma^2} +
{\gamma^2(\mu_\pm (x+1)+c)\over F_1(x)F_1(x+1)\pm\gamma^2} =
\mu_1 \pm(\mu_\pm x +c+ {\gamma^2\over F_1(x)})-F_1(x),
\lab{unknown} \end{equation}
где $\mu_\pm=\mu_2\pm \mu_1$, а $\gamma^2$ и $c$ обозначают
константы интегрирования. Неизвестно, можно ли проинтегрировать
уравнение \re{unknown} еще раз в общем виде. Однако, если положить
$\mu_\pm=0$, то тогда легко понизить его порядок:
$$
I={\gamma^2 c\over (F_1(x)F_1(x+1)\pm\gamma^2)^2}+
{F_1(x)+F_1(x+1)-\mu_1 \mp c\over F_1(x)F_1(x+1)\pm \gamma^2}=const.
$$
Общее решение этого уравнения выражается через эллиптические функции.
Благодаря условию $\mu_\pm=0$, операторы $B^\pm$
коммутируют с Гамильтонианом, то есть они являются интегралами движения.

Интересная система возникает при $N=1$, $\delta=1/3$:
$$
P_1(y)=\gamma^2q^{-2y}F_1(y-1)F_1(y-2),
$$
\be
\gamma^2 q^{-2(y+3)}(F_1(y+2)F_1(y+1)-q^2F_1(y)F_1(y-1))
+F_1(y)-qF_1(y+1)=\mu_1,
\lab{qpain} \end{equation}
где $y=3x$. Это уравнение дальше не упрощается, но при $q=1$
оно допускает дополнительный интеграл и приводит к
\be
\gamma^2 F(y)(F(y-1)+F(y+1)-\gamma^{-2})=\mu y+c,
\lab{pain} \end{equation}
где $F\equiv F_1,\; \mu\equiv \mu_1$.
Это уравнение очень похоже на дискретный аналог первого уравнения
Пенлеве (PI), рассмотренного в работе \ci{fokas}, и в непрерывном пределе
\be
F(y)={1-3h^2u(\xi)\over4\gamma^2},\quad
\xi={h\over\mu}(\mu y+c+{1\over 8\gamma^2}),\qquad
h\to 0,\qquad \mu\to -\; {3h^5\over 16\gamma^2},
\lab{cont}\end{equation}
оно сводится к стандартному уравнению Пенлеве-I:
$d^2 u(\xi)/d\xi^2=6u^2(\xi)+\xi.$

При $N=1$, $\delta=2,$ имеем
\be
F_1(x)=\gamma^2\beta(x)\beta(x+1),\qquad  P_1(x)=q^{-x}\beta^{-1}(x),
\lab{n1d2}\end{equation}
где $\beta(x)$ определяется уравнением $(\mu\equiv\mu_1>0$):
\be
q^{-x}({1\over\beta(x)}-{1\over\beta(x+1)})+
\gamma^2 q(\beta(x-1)\beta(x)-q\beta(x+1)\beta(x+2)) =\mu.
\lab{ric} \end{equation}
В работе \cite{svz:difference} было найдено одно элементарное решение этого
уравнения: $\beta^2=\mu/\gamma^2 q(1-q)$,
приводящее к полиномам Стильтьеса-Вигерта и алгебре $q$-гармонического осциллятора
\cite{ata-sus:realization}. Это решение вырождается при $q\to 1$ если держать
размер решетки конечным. Действительно, при $q=1$, порядок уравнения \re{ric}
легко понижается:
\be
\beta(x-1)+\beta(x+1)=-\; {1+(\mu x+c)\beta(x) \over \gamma^2\beta^2(x)},
\lab{pain2}
\end{equation}
где $c$ обозначает константу интегрирования. Это уравнение представляет собой
специальный случай дискретного Пенлеве-II уравнения, рассмотренного в
статье \ci{papa} и оно не допускает решения в виде постоянной. Поскольку
$\beta(x)>0$, решетка $\Gamma_c$ должна быть полубесконечной,
распространяющейся от $-\infty$ до некоторой точки.
Если $\mu=0$, то \re{pain2} интегрируется в терминах эллиптических функций.

Имеются также более сложные примеры, связанные с дискретными уравнениями Пенлеве
\cite{svz:difference}. Таким образом, при автомодельных редукциях цепочки спектральных
преобразований для конечно-разностного оператора Шредингера в качестве
``нелинейных" специальных функций возникают дискретные трансценденты Пенлеве и
их $q$-аналоги, причем $q$-деформированные уравнения, как правило, имеют более высокий
порядок чем уравнения на линейной решетке.
Некоторые другие примеры дискретных автомодельных систем, появляющихся в рамках
разработанного формализма, рассмотрены в работе \cite{gram}.

\subsection{Цепочки Тоды и Вольтерра с дискретным временем и их симметрии}

Рассмотрим последовательность спектральных задач \re{set}, заменив
$x$ на $n\in\Z$ и перейдя к операторам $H_j=g_jL_jg_j^{-1}$ с
некоторыми подходящими функциями $g_j$, такими что
\begin{equation}
H_j\psi_n^{(j)}\equiv \psi_{n+1}^{(j)}+u_n^j \psi_{n-1}^{(j)}+b_n^j \psi_n^{(j)}
=\lambda\psi_n^{(j)}, \qquad n, j\in \Z.
\label{dse}\end{equation}
Наложим граничные условия $\psi_{-1}^{(j)}(\lambda)=0,$
$\psi_0^{(j)}(\lambda)=1$, и будем считать, что $u_n^j\geq 0 $ и $b_n^j$ действительны
и несингулярны. Тогда по теореме Фавара $\psi_n^{(j)} (\lambda)$,
$n=0, 1, \dots$, образуют систему ортогональных полиномов $\lambda$ с положительной
мерой. Положим
\begin{equation}
\psi_n^{(j+1)}=\frac{\psi_{n+1}^{(j)}+C_n^{j+1} \psi_n^{(j)}}{\lambda-\lambda_{j+1}}
\equiv \frac{L_{j+1}\psi_n^{(j)}}{\lambda-\lambda_{j+1}},
\label{forwdse}\end{equation}
что определяет спектральное преобразование Кристоффеля, так как оно
связано с ядерными полиномами \cite{chi}.
Обратное преобразование, изученное Геронимусом \cite{Ger1,Ger2}, имеет вид:
\begin{equation}
\psi_n^{(j-1)}=\psi_n^{(j)}+A_n^j \psi_{n-1}^{(j)}\equiv R_j\psi_n^{(j)}.
\label{backdse}\end{equation}
В преобразованиях \re{forwdse}, \re{backdse} коэффициенты $A_n^j $ и $C_n^j $
аналогичны дискретным суперпотенциалам рассмотренным в предыдущих параграфах.
Условия совместности этих двух соотношений приводят к факторизации
$H_j=L_jR_j+\lambda_j$ и равенствам
\begin{equation}
u_n^j =A_n^j C_n^j , \qquad b_n^j =A_{n+1}^j +C_n^j +\lambda_j.
\label{factrec}\end{equation}
Абстрактная факторизационная цепочка $L_{j+1}R_{j+1}+\lambda_{j+1}=
R_jL_j+\lambda_j$ в этом случае эквивалентна связанной паре
дискретных нелинейных уравнений:
\begin{equation}
A_n^{j}C_{n-1}^{j}=A_n^{j-1}C_n^{j-1}, \qquad
A_n^{j}+C_n^{j}+\lambda_{j} = A_{n+1}^{j-1} +C_n^{j-1}+\lambda_{j-1}.
\label{factdse}\end{equation}
Эти уравнения определяют цепочку Тоды с дискретным временем \cite{SZ0}.
Очевидно, что после подходящих переобозначений они совпадают с
\re{d1}, \re{d2}. Заметим, что первые интегралы $\lambda_j$ определяют
нарушение изоспектральности операторов $H_j$.

$q$-Периодическое замыкание  \re{clos} в новых обозначениях принимает вид:
\begin{equation}
A_n^{j+N}=qA_{n+k}^j , \qquad C_n^{j+N}= qC_{n+k}^j , \qquad
\lambda_{j+N}=q\lambda_j,
\label{qperdse}\end{equation}
где $k\in\Z$ в соответствии с принятым ограничением на $n$. Если $n$ рассматривается
как непрерывная переменная, тогда $k$ может принимать произвольные значения.

Цепочка Вольтерра с дискретным временем вида
\begin{equation}
D_n^j  \left( D_{n-1}^j -\beta_j \right)=D_n^{j-1}\left( D_{n+1}^{j-1}
- \beta_{j-1} \right),
\label{dtvl}\end{equation}
была введена в работе \cite{SZ0,SZ1} как цепочка преобразований
Кристоффеля-Геронимуса для двухдиагональных матриц Якоби. Ранее
такое уравнение появлялось при анализе алгоритмов для
численного расчета собственных значений матриц в виде
$g$-алгоритма Бауэра \cite{Bau}.

Решения уравнения \re{dtvl} определяют решения цепочки Тода с дискретным временем
\re{factdse} с помощью соотношений $\lambda_j=const -\beta_j^2$ и
\begin{equation}
A_n^j =D_{2n}^j D_{2n+1}^j , \qquad
C_n^j =(D_{2n+1}^j -\beta_j)(D_{2n+2}^j -\beta_j),
\label{map1}\end{equation}
или
\begin{equation}
A_n^j =D_{2n-1}^j D_{2n}^j , \qquad
C_n^j =(D_{2n}^j -\beta_j)(D_{2n+1}^j -\beta_j).
\lab{map2}\end{equation}
В частности, в работе \cite{SZ0} был представлен следующий результат.
\begin{theorem}\label{sep-aw}
Обобщенное разделение переменных в цепочке \re{dtvl} вида
\begin{equation}
D_{2n}^j=\frac{\phi_jr_n}{g_{2n+j}}, \qquad
D_{2n+1}^j=\frac{\sigma_jp_n}{g_{2n+j+1}},
\label{sepvar}\end{equation}
интегрируемо в терминах элементарных функций. Оно приводит к рекуррентным
коэффициентам для двух широких классов ортогональных полиномов:
(ассоциированных) полиномов Аски-Вильсона и симметричных полиномов
Поллачека \cite{askeyismail}.
\end{theorem}
Детали доказательства этой теоремы приведены в статье \cite{SZ1}.

Интересная дискретная симметрия, связанная с аффинной группой Вейля, для факторизационной цепочки
дифференциального уравнения Шредингера  была найдена в работе \cite{Ad}. Она возникает благодаря
тому, что последовательность $K$ преобразований Дарбу с параметрами
$\lambda_j, \lambda_{j+1}, \dots, \lambda_{j+K-1}$ дает один и тот же результат
независимо от порядка проведения преобразований. Этот факт очевиден из Вронскианного
представления $K$-шагового преобразования \ci{Crum}. В результате, появляется
нетривиальная симметрия цепочки \re{chain}. В простейшем случае $K=2$ она имеет
вид:
\begin{equation}
\tilde f_{j-1}=f_{j-1}-\frac{\lambda_{j}-\lambda_{j-1}}{f_{j}+f_{j-1}}, \qquad
\tilde f_{j}=f_j+\frac{\lambda_{j}-\lambda_{j-1}}{f_{j}+f_{j-1}},
\label{adcont1}\end{equation}
\begin{equation}
\tilde\lambda_{j-1}=\lambda_{j}, \qquad\tilde\lambda_{j}=\lambda_{j-1},
\label{adcont2}\end{equation}
а все остальные переменные остаются неизменными
$\tilde f_k(x)=f_k(x),\; \tilde\lambda_k=\lambda_k, \; k\neq j, j-1.$
Действительно, требование
\begin{equation}
\frac{L_{j}L_{j-1}}{(\lambda-\lambda_{j})(\lambda-\lambda_{j-1})}=
\frac{\tilde L_{j}\tilde L_{j-1}}{(\lambda-\tilde\lambda_{j})
(\lambda-\tilde\lambda_{j-1})},
\label{adrule}\end{equation}
допускает только две возможности для преобразований $\lambda_j$:
$\tilde\lambda_j=\lambda_j$ и $\tilde \lambda_{j-1}=\lambda_{j-1}$,
что не приводит к нетривиальным преобразованиям, или (\ref{adcont2}).
Соотношения \re{adcont1} выводятся из (\ref{adrule}) или прямым разрешением
ограничений, вытекающих из цепочки \re{chain}. Если обозначить генераторы указанной
симметрии как $B_j$, т.е. $B_jf_k(x)=\tilde f_k(x)$, и т.д., то они образуют
аффинную группу Вейля \cite{Ad}:
\begin{equation}
B^2_j=1, \qquad B_iB_j=B_jB_i, \quad i\neq j\pm 1, \qquad (B_{j-1}B_j)^3= 1.
\label{weyl}\end{equation}

Используя ту же аргументацию, как и в непрерывном случае \re{adrule},
в работе \cite{spi:jphysa} автор нашел дискретный аналог симметрии \re{adcont1},
\re{adcont2}.

\begin{theorem}
Дискретная факторизационная цепочка (цепочка Тоды с дискретным временем) (\ref{factdse})
инвариантна относительно преобразований
$$
\tilde A_n^{j-1}=A_n^{j-1}+\frac{(\lambda_{j}-\lambda_{j-1})
(A_n^{j}+A_n^{j-1})}{C_{n-1}^{j}+C_n^{j-1}-A_n^{j}-A_n^{j-1}},
$$
\begin{equation}
\tilde A_n^{j}=A_n^{j}-\frac{(\lambda_{j}-\lambda_{j-1})(A_n^{j}+A_n^{j-1})}
{C_{n-1}^{j}+C_n^{j-1}-A_n^{j}-A_n^{j-1}},\label{lawAj}\end{equation}
$$
\tilde C_n^{j-1}=C_n^{j-1}-\frac{(\lambda_{j}-\lambda_{j-1})
(C_{n-1}^{j}+C_n^{j-1})}{C_{n-1}^{j}+C_n^{j-1}-A_n^{j}-A_n^{j-1}},
$$
\begin{equation}
\tilde C_n^{j}=C_n^{j}+\frac{(\lambda_{j}-\lambda_{j-1})
(C_n^{j}+C_{n+1}^{j-1})}{C_n^{j}+C_{n+1}^{j-1}-A_{n+1}^{j}-A_{n+1}^{j-1}},
\label{lawCj}\end{equation}
\begin{equation}
\tilde\lambda_{j-1}=\lambda_{j}, \quad \tilde \lambda_{j}=\lambda_{j-1},
\lab{lambdalaw}\end{equation}
при условии что все остальные $A_n^k, C_n^k, \lambda_k$
при $k\neq j, j-1$ остаются неизменными.
\end{theorem}

В непрерывном пределе $x=nh, h\to 0$ (если он существует) полученные преобразования
симметрии переходят в соотношения \re{adcont1}, \re{adcont2}.
Как и в случае стандартного уравнения Шредингера, преобразования $j\to j+k$ в
\re{factdse} и \re{lawAj}-\re{lambdalaw} описывают преобразования Бэклунда
для дискретных трансцендентов Пенлеве и их $q$-аналогов, появляющихся
из автомодельной редукции \re{qperdse}.

Исходя из связи функций $A_n^j , C_n^j $ с $D_n^j $, можно найти
аналогичную симметрию и для цепочки Вольтерра с дискретным временем \re{dtvl}.

\begin{theorem}
Для фиксированного $j$, дискретная симметрия, связанная со свободой промежуточного
шага при сдвигах на два шага дискретного времени в цепочке Вольтерра
(\ref{dtvl}) с $\beta_j,\beta_{j-1}\neq 0$, имеет вид:
\begin{equation}
\tilde D_{n}^j =\frac{1}{\beta_{j-1}}\left( \beta_jD_{n}^j  +
{(\beta_j^2-\beta_{j-1}^2)(D_n^j D_{n+1}^j +D_n^{j-1}D_{n+1}^{j-1}) \over
\beta_j(\beta_j-D_{n-1}^j -D_{n+1}^j ) +
\beta_{j-1}(\beta_{j-1} - D_n^{j-1}-D_{n+2}^{j-1})}
\right),
\lab{dnj}\end{equation}
\begin{equation}
\tilde D_{n}^{j-1}=\frac{1}{\beta_j} \left( \beta_{j-1}D_n^{j-1} -
{(\beta_j^2-\beta_{j-1}^2)(D_{n-1}^j D_n^j +D_{n-1}^{j-1}D_n^{j-1})
\over \beta_j(\beta_j-D_{n-2}^j -D_n^j )+
\beta_{j-1}(\beta_{j-1}-D_{n-1}^{j-1}-D_{n+1}^{j-1}) }\right).
\lab{dnj-1}\end{equation}
Преобразования параметров $\beta_j$ имеют вид:
\begin{equation}
\tilde \beta_j=\beta_{j-1}, \qquad \tilde\beta_{j-1}=\beta_j.
\lab{betalaw}\end{equation}
Все другие величины $D_n^k$ и $\beta_k$ с $k\neq j, j-1$ остаются неизменными.
\end{theorem}

Как видно из структуры цепочки \re{dtvl}, одновременное изменение знаков
$\tilde D^{j}, \tilde D^{j-1}_n$ и $\tilde \beta_{j}, \tilde\beta_{j-1}$
оставляет эту цепочку инвариантной. Можно проверить, что преобразования
суперпотенциалов $A_n^j , C_n^j $ для цепочки Тоды с дискретным временем
\re{lawAj}-\re{lambdalaw} следуют из связей \re{map1} и \re{map2} и
преобразований \re{dnj}-\re{betalaw}.

Если $\beta_{j-1}=0$ (при фиксированном $j$), тогда появляется любопытная
свобода:
\begin{equation}
\tilde D_n^j =(a_j\pm (-1)^n\sqrt{a_j^2-1})\;\frac{D_n^{j-1}(D_n^j -\beta_j)}
{D_{n-1}^j -\beta_j}, \qquad \tilde D_n^{j-1}=D_{n-1}^j ,
\lab{bj-1=0}\end{equation}
где $a_j$ обозначают произвольные параметры. Предел $\beta_{j-1}\to 0$ в \re{dnj},
\re{dnj-1} соответствует выбору $a_j=1$ в \re{bj-1=0}.
Если $\beta_j=0$, то имеем
\begin{equation}
\tilde D_n^j =D_{n+1}^{j-1}, \qquad
\tilde D_n^{j-1}= (b_j\pm (-1)^n\sqrt{b_j^2-1})\; \frac{D_n^j (D_n^{j-1}-
\beta_{j-1})}{D_{n+1}^{j-1}-\beta_{j-1}}
\lab{bj=0}\end{equation}
с другими свободными параметрами $b_j$. Предел $\beta_j\to 0$ в
\re{dnj}, \re{dnj-1} соответствует выбору $b_j=1$.
Если $\beta_j=0$ для всех $j$, тогда цепочка \re{dtvl} допускает общий интеграл вида
\begin{equation}
D_n^j  = (\gamma_j\pm (-1)^n \sqrt{\gamma_j^2-1}) D_{n+1}^{j-1},
\lab{triv}\end{equation}
где $\gamma_j$ --- произвольные параметры. Аналог интересующей нас симметрии имеет
теперь двухпараметрическую свободу
$$
\tilde D_n^j  = (a_j \pm (-1)^n\sqrt{a_j^2-1})D_{n+1}^{j-1}, \quad
\tilde D_n^{j-1} = (b_j \pm (-1)^n\sqrt{b_j^2-1})D_{n-1}^j .
$$
Пределы $\beta_j, \beta_{j-1}\to 0$ в \re{dnj}, \re{dnj-1}
дают преобразования, определенные на подпространстве решений
с $\gamma_j=\pm 1$ в \re{triv}.

\section{Обобщенная задача на собственные значения для двух матриц Якоби}

\subsection{Обобщенная факторизация и $R_{II}$-цепочка}

Изложение этого параграфа основано на результатах работы \cite{spi-zhe:spectral}.
Обозначим $P_n^{(j)}(z)$, $n,\; j\in {\Z }$, бесконечный набор функций
независимой переменной $z\in {\C }$. Предположим, что эти функции удовлетворяют
соотношениям
\begin{equation}
\label{r1}
P_n^{(j+1)}(z)=\frac{D_n^{j+1} P_{n+1}^{(j)}(z) +C_n^{j+1} \,
(z-\alpha_n^{j+1} )P_n^{(j)}(z)}{z-\lambda_{j+1}},
\end{equation}
\begin{equation}
\label{r2}
P_n^{(j-1)}(z)=B_n^j P_n^{(j)}(z)+A_n^j \,(z-\beta_n^j )\, P_{n-1}^{(j)}(z),
\end{equation}
где коэффициенты $A_n^j , B_n^j , C_n^j , D_n^j $
и спектральные переменные $\alpha_n^j ,
\beta_n^j , \lambda_j$ не зависят от $z$. На языке теории интегрируемых
систем уравнения \re{r1} и \re{r2} образуют некоторую дискретную пару Лакса со
спектральным параметром $z$.

Сдвинув дискретную переменную $j\to j-1$ в \re{r1} и избавившись от $P_n^{(j-1)}(z)$
с помощью \re{r2}, мы получаем трехчленное рекуррентное соотношение
\begin{equation}
\label{recrel}
 P_{n+1}^{(j)}(z)+r_n^j \,(v_n^j -z)P_n^{(j)}(z)
+u_n^j \, (z-\alpha_n^j )(z-\beta_n^j )\, P_{n-1}^{(j)}(z)=0,
\end{equation}
в котором потенциалы $u_n^j , r_n^j , v_n^j $ имеют вид
$$
u_n^j =\frac{A_n^j \,C_n^j }{D_n^j \,B_{n+1}^j }, \qquad
r_n^j =\frac{1-D_n^j \,A_{n+1}^j -C_n^j \,B_n^j }{D_n^j \,B_{n+1}^j },
$$
\begin{equation}
r_n^j \,v_n^j = \frac{\lambda_j-\beta_{n+1}^j \,D_n^j \,A_{n+1}^j
-\alpha_n^j \,C_n^j \,B_n^j }{D_n^j \,B_{n+1}^j }.
\label{super1} \end{equation}
Аналогично, сдвинув $j\to j+1$ в \re{r2} и устранив $P_n^{(j+1)}(z)$
с помощью \re{r1}, мы опять получаем соотношение \re{recrel}, но
уже с другими рекуррентными коэффициентами. Условия совместности
этих рекуррентных соотношений приводят к ограничениям
$$
\beta_n^j =\beta_n, \qquad \alpha_n^j =\alpha_{n+j}
$$
и системе трех нелинейных конечно-разностных уравнений
\begin{eqnarray}
&&\frac{A_n^j C_n^j }{B_{n+1}^j D_n^j } = \frac{A_n^{j+1} \, C_{n-1}^{j+1} }
{B_n^{j+1} \, D_n^{j+1} }, \label{rel1} \\
&&\frac{C_n^j \,B_n^j  + A_{n+1}^j D_n^j  - 1}{B_{n+1}^j D_n^j } =
\frac{C_n^{j+1} B_n^{j+1} +A_n^{j+1} D_{n-1}^{j+1} -1}
{B_n^{j+1} \,D_n^{j+1} },\label{rel2}\\
&&\frac{\alpha_{n+j}\, C_n^j \,B_n^j  +\beta_{n+1}A_{n+1}^j D_n^j
-\lambda_j}{B_{n+1}^j D_n^j }
= \frac{\alpha_{n+j+1} C_n^{j+1}  B_n^{j+1} +\beta_n \,A_n^{j+1}
D_{n-1}^{j+1} -\lambda_{j+1}}{B_n^{j+1} \,D_n^{j+1} }.
\label{rel3} \end{eqnarray}
Эта система определяет (1+1)-мерную дискретную интегрируемую цепочку,
поскольку она возникла из условия совместности дискретной пары Лакса.
Она играет ключевую роль в последующих рассмотрениях. Переменная $j$
может интерпретироваться как дискретное время, поскольку полученные
уравнения обобщают
цепочку Тоды с дискретным временем, описанную в предыдущих параграфах.

Если на рекуррентное соотношение \re{recrel} наложить граничные условия
$$
P_0^{(j)} (z)=p_j, \qquad P_1^{(j)} (z)=r_0^j (z-v_0^j ),
$$
где $p_j, r_0^j$ ненулевые числа, то $P_n^{(j)}(z)$, $n\geq 0$, представляют собой
полиномы $n$-й степени по $z$. Для того, чтобы оборвать рекуррентное соотношение
\re{r2} слева при $n=0$, мы налагаем ограничения $A_0^j =A_0^j \beta_0=0.$
Цепные дроби, связанные с трехчленным рекуррентным соотношением типа
\re{recrel}, названы в работе \cite{ism-mas:general} $R_{II}$-дробями.
Поэтому мы будем называть $P_n^{(j)}(z)$-функции $R_{II}$-полиномами,
а уравнения \re{rel1}-\re{rel3} $R_{II}$-цепочкой.

Заметим, что $R_{II}$-полиномы могут быть редуцированы к так называемым
$R_I$ и Лорановским биортогональным полиномам или стандартным ортогональным
полиномам устранением из \re{recrel} квадратичной зависимости от $z$ и
последующими упрощениями зависимости этого рекуррентного соотношения от $z$ различными
предельными переходами. Преобразования \re{r1} и \re{r2} аналогичны преобразованиям
Кристоффеля и Геронимуса в теории ортогональных полиномов, описанным
в предыдущих параграфах. Взятые вместе, они могут рассматриваться как
дискретные спектральные преобразования для $R_{II}$-рекуррентного соотношения \re{recrel}.

Как показано в работе \cite{ism-mas:general}, для заданного множества $R_{II}$-полиномов,
удовлетворяющих ограничениям $P_n^{(j)}(\alpha_{n+j})\neq 0,$
$P_n^{(j)}(\beta_n)\neq 0$
и $u_n^j \neq 0$, всегда существует линейный функционал ${\cal L}_j$
(дискретная переменная $j$ может рассматриваться как вспомогательная
величина в равенстве \re{recrel}) такой, что
\be
{\cal L}_j\left[\frac{z^mP_n^{(j)}(z)}{\prod_{k=1}^n(z-\alpha_{j+k})
(z-\beta_k)}\right]=0, \qquad 0\leq m <n.
\lab{bio_1}  \end{equation}
Это соотношение появилось впервые в теории многоточечных Паде
аппроксимаций \cite{gon:speed,GL}.
На основе очень простых аргументов из линейной алгебры, в работе \cite{zhe:bio}
это соотношение было переписано как условие биортогональности для двух
рациональных функций, построенных из $P_n^{(j)}(z)$:
\be
{\cal L}_j\left[H_m^{(j)}(z)\, R_n^{(j)}(z) \right]= 0
\quad \mbox{при} \quad n\neq m,
\lab{biort}\end{equation}
где рациональные функции $R_n^{(j)}(z)$ и $H_m^{(j)}(z)$ определяются ниже.
Положим $R_0^{(j)}(z)=S_0^{(j)}(z)=1$ и обозначим
\be
R_n^{(j)}(z)=\frac{P_n^{(j)}(z)}{\prod_{k=1}^n(z-\alpha_{j+k})}, \qquad
S_n^{(j)}(z)=\frac{P_n^{(j)}(z)}{\prod_{k=1}^nu_k^j (z-\beta_k)},
\lab{RS}\end{equation}
для $n>0$. Эти функции удовлетворяют трехчленному рекуррентному соотношению
с линеаризованной зависимостью от $z$:
\be
(z-\alpha_{n+j+1})R_{n+1}^{(j)}(z) + r_n^j (v_n^j -z)R_n^{(j)}(z)
+ u_n^j (z-\beta_n)R_{n-1}^{(j)}(z)=0,
\lab{recrelR} \end{equation}
\be
u_{n+1}^j (z-\beta_{n+1})S_{n+1}^{(j)} (z)+r_n^j (v_n^j -z)S_n^{(j)}(z)+
(z-\alpha_{n+j})S_{n-1}^{(j)}(z)=0.
\lab{recrelS}\end{equation}
Уравнения \re{recrelR} и \re{recrelS} могут рассматриваться как
обобщенные спектральные задачи вида \cite{Wilk}
$$L\psi(z)=zM\psi(z),$$
где операторы $L$ и $M$ задаются двумя общими трехдиагональными матрицами Якоби.

Верхний индекс $j$ не играет существенной роли в определении соотношений
биортогональности. Поэтому для простоты положим временно $j=0$ и устраним верхние индексы.
Тогда соотношение \re{recrelR} может быть переписано в виде $L R_n(z)=zM R_n(z)$,
где
\begin{eqnarray*}
&& LR_n(z)\equiv \alpha_{n+1}R_{n+1}(z)-r_nv_nR_n(z)+\beta_nu_nR_{n-1}(z),
\\
&& MR_n(z)\equiv R_{n+1}(z)-r_nR_n(z)+u_nR_{n-1}(z).
\end{eqnarray*}
В этих обозначениях уравнение \re{recrelS} принимает вид
$L^T S_n(x)=xM^T S_n(x)$, где $L^T$ и $M^T$ обозначают матрицы $L$ и $M$
транспонированные по отношению к формальному скалярному произведению
$$
(S(x),R(z))\equiv\sum_{k=0}^\infty S_n(x)R_n(z)
$$
определенному на пространстве рациональных функций. В явном виде мы имеем
\begin{eqnarray*}
&& L^TS_n(z)=u_{n+1}\beta_{n+1}S_{n+1}(z)-r_nv_nS_n(z)+\alpha_nS_{n-1}(z),
\\ &&
M^TS_n(z)=u_{n+1}S_{n+1}(z)-r_nS_n(z)+S_{n-1}(z).
\end{eqnarray*}
Из последовательности очевидных равенств
\begin{eqnarray}\nonumber
&& 0=(S(x),LR(z))-z(S(x), MR(z))
\\ && \makebox[2em]{}
=(L^T S(x), R(z))-z(M^T S(x), R(z))=(x-z)(M^T S(x),R(z)),
\lab{biortMSR}\end{eqnarray}
легко заключить, что функции $H_n(x)\equiv M^T S_n(x)$
ортогональны $R_n(z)$ при различных значениях собственных значений $x\neq z$.
Восстановив верхний индекс $j$ можно увидеть, что функции $H_n^{(j)}(z)$
имеют следующий явный вид
\be
H_n^{(j)}(z)\equiv u_{n+1}^j S_{n+1}^{(j)}(z)-r_n^j S_n^{(j)}(z)+S_{n-1}^{(j)}(z)
\lab{H_n}\end{equation}
для $n=1, 2, \dots,$ а при $n=0$ имеем
$$
H_0^{(j)}(z)=u_1^j S_1^{(j)}(z)-r_0^j =\frac{r_0^j (\beta_1-v_0^j )}{z-\beta_1}.
$$

Поскольку мы имеем дело с матрицами и их собственными значениями,
ортогональность для различных собственных значений \re{biortMSR}
предполагает, что имеет место еще и дуальное соотношение ортогональности
для функций  $H_m^{(j)}(z)$ и $R_n^{(j)}(z)$ при равных значениях собственных
значений $z$. Оно определяется с помощью функционала ${\cal L}_j$,
отображающего рациональные функции $z$ на комплексные числа $\C $ \re{bio_1}.
В результате, биортогональность $H_m^{(j)}(z)$ и $R_n^{(j)}(z)$ для $m\neq n$
может быть проверена прямой подстановкой соответствующих выражений в
равенство \re{biort} и применением условий \re{bio_1}.

Любое нетривиальное решение $R_{II}$-цепочки с подходящими граничными условия\-ми
при $n=0$ соответствует определенной системе биортогональ\-ных рациональ\-ных функций (БРФ).
Опишем кратко процедуру построения полиномов $P_n^{(j)}(z)$ при заданных
коэффициентах $A_n^j , \dots, \lambda_j.$ Определим два вспомогательных
полинома $n$-й степени:
\be
Y_n^j =\prod_{k=1}^n(z-\lambda_{j+k}), \qquad
Z_n^j =\prod_{k=1}^n(z-\alpha_{j+k}), \quad n>0,
\lab{YZ}\end{equation}
и положим $Y_0^j =Z_0^j =1$.
Тогда $P_n^{(j)}(z)$ могут быть представлены в виде
\be
P_n^{(j)}(z)=Z_n^j (z)\sum_{k=0}^n\zeta_n^j (k)\frac{Y_k^j (z)}{Z_k^j (z)}
\lab{PYZ}\end{equation}
с некоторыми неизвестными коэффициентами $\zeta_n^j (k), k\leq n$. Подставляя
это выражение в \re{r1}, мы приходим к системе уравнений
\begin{eqnarray*}
&& \zeta_{n+1}^{j}(0)\, D_n^{j+1} + \zeta_n^{j}(0)\, C_n^{j+1} =0, \qquad
\zeta_n^{j+1} (n)=D_n^{j+1} \,\zeta_{n+1}^{j}(n+1),
\\ &&
D_n^{j+1} \, \zeta_{n+1}^{j}(k)+C_n^{j+1} \,\zeta_n^{j}(k)
=\zeta_n^{j+1} (k-1), \quad k=1, 2, \dots, n.
\end{eqnarray*}
Из первых двух уравнений находим $\zeta_n^{j}(0)$ и $\zeta_n^{j}(n)$
(эти коэффициенты однозначно фиксируются граничными условиями $\zeta_0^{j}(0)=1$):
\be
\zeta_n^{j}(0)=(-1)^n\prod_{m=0}^{n-1}\frac{C_m^{j+1} }{D_m^{j+1} }, \qquad
\zeta_n^{j}(n)=\prod_{m=0}^{n-1}\frac{1}{D_m^{j+n-m}}.
\lab{zeta0} \end{equation}
Определив нормированные коэффициенты $\eta_n^{j}(k)=\zeta^j _n(k)/\zeta_n^j (0)$,
можно переписать оставшуюся часть уравнений в следующем виде
\be
\eta^j _{n+1}(k)=\eta_n^j (k)-\frac{\zeta_n^{j+1} (0)}
{\zeta_n^{j}(0) C_n^{j+1} }\eta_n^{j+1} (k-1),\quad
k=1, 2, \dots, n.\lab{rec_eta} \end{equation}
Поскольку $\eta_n^j (0)=1$ и $\eta_k^j (k)$ уже были определены,
это рекуррентное соотношение позволяет однозначно определить все
коэффициенты $\eta_n^{j}(k)$ итерационным образом.

Покажем, что $R_{II}$-цепочка ассоциирует с каждым заданным трехчленным
рекур\-рент\-ным соотношением \re{recrel} другое рекуррентное
соотношение той же самой природы. Дей\-ст\-ви\-тель\-но, из уравнения \re{r1} мы находим
\be
P_{n+1}^{(j)}(z)=\frac{z-\lambda_{j+1}}{D_n^{j+1} }P_n^{(j+1)} (z)-
\frac{C_n^{j+1} (z-\alpha_{n+j+1})}{D_n^{j+1} }P_n^{(j)}(z).
\lab{Pn+1} \end{equation}
Аналогичным образом, пользуясь соотношением \re{r2} мы можем выразить
$P_{n-1}^{(j)}(z)$ через $P_n^{(j-1)}(z)$ и $P_n^{(j)}(z)$:
\be
P_{n-1}^{(j)}(z)=\frac{P_n^{(j-1)}(z)-B_n^j P_n^{(j)}(z)}{A_n^j (z-\beta_n)},
\quad n> 0.
\lab{cont2}\end{equation}
Подставляя выражения \re{Pn+1} и \re{cont2} в равенство \re{recrel},
мы получаем трехчленное рекуррентное соотношение по дискретному времени $j$:
\begin{eqnarray}\nonumber
&& \frac{z-\lambda_{j+1}}{D_n^{j+1} }P_n^{(j+1)}(z)+
\left( r_n^j (v_n^j -z)+\frac{C_n^{j+1} (\alpha_{n+j+1}-z)}{D_n^{j+1} }
+\frac{u_n^j B_n^j (\alpha_{n+j}-z)}{A_n^j }
\right)P_n^{(j)}(z)
\\ && \makebox[6em]{}
+\frac{u_n^j (z-\alpha_{n+j})}{A_n^j }P_n^{(j-1)}(z)=0.
\lab{contig}\end{eqnarray}
Это уравнение опять принадлежит тому же классу обобщенных задач на собственные
значения. Заменив $P_n^{(j)}(z)$ на $S_n^{(j)}(z)$ в равенстве \re{contig} и сравнив
результат с \re{recrelS}, можно увидеть, что $S_n^{(j)}(z)$ удовлетворяет
рекуррентным соотношениям $R_{II}$-типа по обоим дискретным индексам $n$ и $j$.
Заметим, однако, что в контексте $R_{II}$-полиномов мы имеем
ограничение $n\ge 0$, в то время как значения $j$ не ограничены.
Предположим, что зависимость от $j$ входит в $P_n^{(j)}(z)$ через некоторые непрерывные
параметры. Тогда уравнение \re{contig} определяет соотношение сопряжения для соответствующей
системы полиномов, то есть связывает эти полиномы при различных значениях параметров.

\subsection{Сопряженные полиномы}

Рассмотрим подробнее функции $H_n^{(j)}(z)$. Удобно представить их в виде
\be
H_n^{(j)}(z)=\frac{Q_n^{(j)}(z)}{(z-\beta_{n+1})\prod_{k=1}^n u_k^j (z-\beta_k)},
\lab{HQ} \end{equation}
где $Q_n^{(j)}(z)$ полиномы $n$-й степени, которые будут называться
сопряженными с $P_n^{(j)}(z)$. Их явный вид находится из определения \re{H_n}:
\be
Q_n^{(j)}(z)=P_{n+1}^{(j)}(z)-r_n^j (z-\beta_{n+1})P_n^{(j)}(z) +
u_n^j (z-\beta_n)(z-\beta_{n+1})P_{n-1}^{(j)}(z),
\lab{Q0} \end{equation}
для $n>0$ и при $n=0$ имеем $Q_0^{(j)}(z)=r_0^j (\beta_1-v_0^j )$.
Используя рекуррентное соотношение \re{recrel}, мы можем представить $Q_n^{(j)}(z)$
в виде
\be
Q_n^{(j)}(z)= r_n^j (\beta_{n+1}-v_n^j )P_n^{(j)}(z) +
u_n^j (z-\beta_n)(\alpha_{n+j}-\beta_{n+1})P_{n-1}^{(j)}(z)
\lab{Q1} \end{equation}
или, альтернативно,
\be Q_n^{(j)}(z)= \frac{(\beta_{n+1}-\alpha_{n+j})P_{n+1}^{(j)}(z) +
r_n^j (\alpha_{n+j}-v_n^j )(z-\beta_{n+1})P_n^{(j)}(z)}{z-\alpha_{n+j}}.
\lab{Q2} \end{equation}
Из представления \re{Q1} очевидно, что $Q_n^{(j)}(z)$ действительно являются полиномами
$n$-й степени. С помощью формул \re{Q1} и \re{Q2} можно выразить $P_n^{(j)}(z)$
через $Q_n^{(j)}(z)$:
\be
 P_n^{(j)}(z)=\gamma_n^j Q_{n}^{(j)}(z)+\delta_n^j (z-\alpha_{n+j-1})Q_{n-1}^{(j)}(z)
\lab{PQ1} \end{equation}
или
\be
P_n^{(j)}(z)=\frac{\sigma_n^j Q_{n+1}^{(j)}(z) + \tau_n^j (z-\alpha_{n+j})Q_n^{(j)}}
{z-\beta_{n+1}}, \lab{PQ2} \end{equation}
где
\ba
&&\gamma_n^j =\frac{r_{n-1}^j (\alpha_{n+j-1}-v_{n-1}^j )}{\epsilon_n^j },
\qquad \delta_n^j =\frac{u_n^j (\beta_{n+1}-\alpha_{n+j})}{\epsilon_n^j },
\nonumber \\ &&\sigma_n^j =\frac{\alpha_{n+j}
-\beta_{n+1}}{\epsilon_{n+1}^j },  \qquad
\tau_n^j =\frac{r_{n+1}^j (\beta_{n+2}-v_{n+1}^j )}{\epsilon_{n+1}^j },
\nonumber \ea
\be
\epsilon_n^j =
r_n^j r_{n-1}^j (\alpha_{n+j-1}-v_{n-1}^j )(\beta_{n+1}-v_n^j )-
u_n^j (\alpha_{n+j}-\beta_{n+1})(\beta_n-\alpha_{n+j-1}).
\nonumber \end{equation}
Здесь подразумевается, что $\epsilon_n^j \neq 0$.
Подставляя \re{PQ1} и \re{PQ2} в \re{r1} и \re{r2}, мы находим, что
полиномы $Q_n^{(j)}(z)$ удовлетворяют соотношениям
\begin{eqnarray*}
&& Q_n^{(j+1)}=\frac{\tilde D_n^{j+1} Q_{n+1}^{(j)} +\tilde C_n^{j+1} (z-\tilde
\alpha_{n+j+1})Q_n^{(j)}} {z-\tilde\lambda_{j+1}}, \quad
\\ &&
Q_n^{(j-1)}=\tilde B_n^j Q_n^{(j)}+\tilde A_n^j (z-\tilde\beta_n)Q_{n-1}^{(j)}
\end{eqnarray*}
со спектральными коэффициентами и суперпотенциалами вида
\begin{eqnarray}
&& \tilde \alpha_{n+j}=\alpha_{n+j-1}, \quad \tilde \beta_n = \beta_{n+1},
\quad  \tilde \lambda_j= \lambda_j,
\lab{compar}
\\ &&
\tilde A_n^j =\frac{D_{n-1}^j }{D_n^j }A_n^j  , \qquad \tilde B_n^j  =B_{n+1}^j .
\lab{com_AB}\end{eqnarray}
Два других оставшихся суперпотенциала $\tilde D_n^j $ и $\tilde C_n^j $ имеют
намного более сложный вид:
\be
\tilde D_n^j =\frac{\epsilon_n^j }{\epsilon_{n+1}^{j-1}}D_{n-1}^j , \qquad
\tilde C_n^j =\frac{\epsilon_n^j }{\alpha_{n+j-1}-\beta_{n+1}}
\left(D_{n-1}^j \tau_n^{j-1} + \frac{\lambda_j-\beta_{n+1}}
{\epsilon_n^j B_{n+1}^j }\right).
\lab{com_D}\end{equation}
Очевидно, что условие совместности $j\to j\pm 1$
преобразований для $Q_n^{(j)}(z)$ полиномов порождает
$R_{II}$-рекуррентное соотношение и $R_{II}$-цепочку с новыми
неизвестными функциями $\tilde A_n^j , \dots, \tilde \lambda_j$.
Таким образом, мы нашли нетривиальный автоморфизм $R_{II}$-цепочки.

\begin{theorem} Преобразования \re{compar}-\re{com_D}
определяют симметрию $R_{II}$-цепочки \re{rel1}-\re{rel3}, то есть
эта система уравнений не изменяется при соответствующих подстановках,
индуцированных переходом от $R_{II}$-полиномов $P_n^{(j)}(z)$ к их
партнерам $Q_n^{(j)}(z)$.
\end{theorem}

\section{Симметрии $R_{II}$-цепочки спектральных преобразований и
автомодельные редукции}

\subsection{Дискретные симметрии $R_{II}$-цепочки}

Опишем теперь некоторые симметрии $R_{II}$-цепочки. Для начала кратко
рассмотрим нормировочную (или калибровочную) свободу.
Можно преобразовать рекуррентные коэффициенты в соотношении \re{recrel}
умножением полиномов на произвольный весовой множитель
$\xi_n^j $ не зависящий от $z$, $P_n^{(j)}(z) = \xi_n^j \, \tilde P_n^{(j)}(z)$.
Это приводит к рекурсиям \re{r1}, \re{r2} с перенормированными ингредиентами
\be
\tilde A_n^j =\frac{A_n^j }{t_{n-1}^j }, \quad \tilde B_n^j = B_n^j  w_n^j ,
\quad \tilde C_n^j = \frac{C_n^j }{w_n^j }, \quad \tilde D_n^j =D_n^j  \,
t_n^j ,
\lab{tr_sup} \end{equation}
где $t_n^j =\xi_{n+1}^{j-1}/\xi_n^j$ и $w_n^j = \xi_n^j /\xi_n^{j-1}.$
Коэффициенты $t_n^j$ и $w_n^j $ подчиняются соотношению
$t_n^j \, w_{n+1}^j =t_n^{j+1} \, w_n^{j+1}.$ Преобразованные рекуррентные
коэффициенты имеют вид
\be
\tilde r_n^j  = \frac{\xi_n^j }{\xi_{n+1}^j }r_n^j ,\qquad \tilde u_n^j  =
\frac{\xi_{n-1}^j }{\xi_{n+1}^j }u_n^j ,
\lab{tr_rec} \end{equation}
а все остальные объекты, входящие в \re{recrel}, не изменяются.  Таким образом,
имеется широкая свобода в представлении рекуррентных коэффициентов для полиномов
$P_n^{(j)}(z)$.

В статье \cite{ism-mas:general} использовалась калибровка $r_n^j =1$.
Можно также выбрать другую калибровку $r_n^j -u_n^j =1, r_0^j =1,$
приводящую к моническим полиномам $P_n^{(j)}(z)=z^n+O(z^{n-1})$,
которая может быть удобна для упрощения выражений. Из уравнений
\re{r1} и \re{r2} видно, что условие моничности подразумевает следующее
ограничение на суперпотенциалы: $A_n^j +B_n^j =C_n^j +D_n^j =1.$
В этой нормировке в действительности имеется всего два независимых уравнения
\re{rel1} и \re{rel3}, поскольку уравнение \re{rel2} выполняется автоматически.

У монической калибровки имеется существенный технический недостаток ---
в ней трудно конструировать явные решения $R_{II}$-цепочки
$A_n^j , \dots, D_n^j $. Для этой цели необходимо редуцировать число
суперпотенциалов и, поэтому, удобно выбрать калибровку:  $B_n^j =1$.
Как видно из \re{tr_sup}, этот выбор оставляет свободу в преобразованиях
суперпотенциалов:
\be
A_n^j  \to A_n^j /t_{n-1}, \qquad D_n^j  \to D_n^j  t_n,
\lab{ga} \end{equation}
где множитель $t_n$ не зависит от $j$.

Опишем теперь более сложные свойства $R_{II}$-цепочки.
Пусть $A_n^j , B_n^j , C_n^j , D_n^j $ удовлетворяют уравнениям \re{rel1}, \re{rel2}.
Эти функции определяли бы решение полной $R_{II}$-цепочки, если бы
выполнялись равенства $\lambda_j=\alpha_{n+j}=\beta_n=const$,
поскольку в этом случае уравнение \re{rel3} совпадает с \re{rel2}.
Сдвигом аргумента полиномов $z\to z-const$ такую ситуацию можно редуцировать
к случаю $\lambda_j=\alpha_{n+j}=\beta_n=0.$ Однако решения, полученные при
таких ограничениях слишком тривиальны $P_n^{(j)}(z)= \gamma_n^j  z^n,$
где  $\gamma_n^j $ не зависят от $z$. Действительно, при $n=0$ в \re{r1}
из начального условия $P_0^{(j)}(z)=1$ мы находим $P_1^{(j)}(z)= \gamma_1^j z$
и приведенное утверждение вытекает по рекурсии. Поэтому мы предполагаем, что
такая простая ситуация не имеет места в дальнейших рассмотрениях.

Нетрудно заметить, что аффинные преобразования аргумента $z$, $z \to \xi z+\eta,$
могут быть компенсированы подходящими преобразованиями параметров
$\alpha_{n+j}, \beta_n, \lambda_j$ и рекуррентных коэффициентов,
аналогично случаю ортогональных полиномов. Однако, БРФ ассоциируются с
обобщенной спектральной задачей $L\psi(z)=zM\psi(z)$,
которая допускает также преобразование инверсии $z\to 1/z$,
приводящего к перестановке операторов $L$ и $M$. В результате,
общие рациональные преобразования аргументов $R_{II}$-полиномов,
комбинированные с компенсирующим калибровочным преобразованием,
\be
\tilde P_n^{(j)}(z) = (\zeta \, z +\sigma)^n \, P_n^j \left(\frac{\xi\, z +
\eta}{\zeta\, z + \sigma}\right), \lab{rat_P}\end{equation}
где $\xi, \eta, \zeta, \sigma$ обозначают произвольные параметры не зависящие
от $j$, не выводят из класса $R_{II}$-полиномов.

\begin{theorem} Полиномы \re{rat_P} удовлетворяют трехчленному рекуррентному
соотношению \re{recrel} с коэффициентами
\begin{eqnarray}\nonumber
&& \tilde r_n^j =r_n^j (\xi-\zeta v_n^j ), \quad
\tilde v_n^j =\frac{\sigma v_n^j -\eta}{\xi - \zeta v_n^j },\quad
\tilde u_n^j =u_n^j (\xi-\zeta \alpha_{n+j})(\xi-\zeta \beta_n),
\\ &&
\tilde\alpha_{n+j}=\frac{\sigma \alpha_{n+j}-\eta}{\xi-\zeta\alpha_{n+j}},
\qquad \tilde\beta_n=\frac{\sigma\beta_n-\eta}{\xi-\zeta\beta_n}.
\lab{partr}\end{eqnarray}
Как следствие, $R_{II}$-цепочка инвариантна относительно
преобразований \re{partr} и
\begin{eqnarray}\nonumber
&& \tilde \lambda_j = \frac{\sigma\lambda_j-\eta}{\xi-\zeta \lambda_j},\quad
\tilde A_n^j  = A_n^j (\xi-\zeta \beta_n), \quad
\tilde B_n^j  =B_n^j , \quad
\\ &&
\tilde C_n^j  = \frac{C_n^j (\xi- \zeta\alpha_{n+j})}{\xi-\zeta\lambda_j},
\qquad \tilde D_n^j = \frac{D_n^j }{\xi-\zeta \lambda_j}.
\lab{rat_pot} \end{eqnarray}
\end{theorem}

\smallskip

Это утверждение доказывается простой подстановкой указанных выражений в
требуемые уравнения.

Другой тип симметрий индуцируется дискретными преобразованиями
двумерной дискретной решетки, образованной переменными $n$ и $j$.
Рассмотрим отражения
1) $j\to -j, n\to -n;$ 2) $n\to j, j\to n;$ 3) $j\to -j-n;$
4) $n\to -n-j$ и построим симметрии $R_{II}$-цепочки, связанные с ними.

\begin{theorem}
Следующие инволютивные преобразования определяют дискретные симметрии
$R_{II}$-цепочки:
\ba
1.&& \tilde A_n^j =D_{-n}^{-j}, \; \tilde D_n^j = A_{-n}^{-j},\;
\tilde B_n^j = C_{-n}^{-j},\; \tilde C_{n}^j =B_{-n}^{-j}, \;
\nonumber\\
&& \tilde \beta_n= \beta_{1-n}, \; \tilde \alpha_{n+j}=\alpha_{-n-j},
\; \tilde \lambda_j= \lambda_{-j};
\nonumber \\
2.&& \tilde A_n^j =\frac{1}{A_j^n}, \; \tilde B_n^j =\frac{B_j^n}{A_j^n},\;
\tilde C_n^j =\frac{C_{j-1}^{n+1}}{D_{j-1}^{n+1}},\;
\tilde D_n^j =\frac{1}{D_{j-1}^{n+1}}, \;
\nonumber \\
&& \tilde \alpha_{n+j}=\alpha_{n+j}, \;
\tilde \beta_n = \lambda_n, \; \tilde\lambda_j= \beta_j;
\lab{invol} \\
3.&& \tilde A_n^j =\frac{A_n^{1-j-n}}{B_n^{1-j-n}}, \;
\tilde B_n^j = \frac{1}{B_n^{1-j-n}}, \;
\tilde C_n^j =\frac{1}{C_n^{1-j-n}}, \;
\tilde D_n^j =\frac{D_n^{1-j-n}}{C_n^{1-j-n}},\;
\nonumber \\
&& \tilde \lambda_j= \alpha_{1-j}, \; \tilde \alpha_{j+n}=\lambda_{1-j-n}, \;
\tilde\beta_n =\beta_n;
\nonumber \\
4.&& \tilde A_n^j =B_{1-n-j}^j , \; \tilde B_n^j =A_{1-n-j}^j , \;
\tilde C_n^j =D_{-n-j}^j , \; \tilde D_n^j = C_{-n-j}^j , \;
\nonumber \\
&& \tilde \alpha_{n+j}=\beta_{1-n-j}, \; \tilde \beta_n =
\alpha_{1-n}, \; \tilde\lambda_j= \lambda_j.
\nonumber\ea
\end{theorem}

\smallskip

Доказательство этого утверждения состоит в проверке того,
 что после подстановки тильдованных выражений в уравнения \re{rel1}-\re{rel3}
возникает $R_{II}$-цепочка с отраженными точками решетки, так как
это указано выше. В некотором смысле эта теорема показывает, что
спектральные параметры $\lambda_j, \alpha_{n+j}, \beta_n$
эквивалентны друг другу несмотря на то, что они входят несимметричным
образом в начальные формулы \re{r1}, \re{r2}.

Эти четыре преобразования не исчерпывают все возможные инволюции
$R_{II}$-цепочки. Например, должны существовать инволюции, порождаемые
свободой в промежуточных шагах двойных сдвигов по дискретному времени,
которые обобщают аналогичные симметрии для цепочки Тоды с дискретным временем
\re{lawAj}-\re{lambdalaw}.

Предположим, что суперпотенциалы $A_n^j , B_n^j , C_n^j , D_n^j $ и
спектральные параметры $\alpha_n,$  $\beta_n,$  $\lambda_n$ описываются
мероморфными функциями непрерывных переменных $n$ и $j$. Существование
таких решений, вообще говоря, не гарантировано, так как разностные уравнения
определяются на числовых последовательностях, которые не обязаны допускать
аналитическое продолжение на макроскопические области. Обычно мероморфные
решения интегрируемых цепочек появляются из автомодельных решений
соответствующих уравнений. В общем случае,  преобразования
\re{invol} существенно изменяют вид заданного решения. Однако существует
специальный класс решений, для которых происходит только изменение параметров
решений. Попробуем построить такие решения в явном виде.

Прежде всего отметим, что существуют специфические комбинации дискретных
переменных $n$ и $j$, а именно,
$$
u_1=n,\; u_2= j,\; u_3=n+j,\; u_4=n-j,\; u_5=2n+j,\; u_6=2j+n,
$$
которые выражаются друг через друга при указанных отражениях
решетки с точностью до изменения знаков. Поэтому, симметричные
произведения некоторых функций этих комбинаций не будет изменять
своей формы при отражениях базовой решетки. Это наводящее соображение
позволяет наложить конструктивное ограничение на суперпотенциалы,
заключающегося в требовании что $A_n^j , B_n^j, C_n^j , D_n^j $
расщепляются на произведения функций, каждая из которых зависит только
от одной из этих шести комбинаций:
$$
A_n^j =\prod_{k=1}^6 A^{(k)}(u_k),
\quad B_n^j =\prod_{k=1}^6 B^{(k)}(u_k), \quad
C_n^j =\prod_{k=1}^6C^{(k)}(u_k), \quad
D_n^j =\prod_{k=1}^6 D^{(k)}(u_k).
$$
Заранее не известно будет ли этот анзац совместимым с уравнениями
 \re{rel1}-\re{rel3}. Упростим максимально вид суперпотенциалов
пользуясь калибровочной свободой. Наложим условие
$B_n^j =1,$ которое позволяет нормировать полиномы
$P_0^{(j)}=1$. Будем подразумевать так же, что $D_n^j $ не зависит от
переменной $u_1=n$, то есть $D^{(1)}(u)=1$, чего всегда можно добиться
 преобразованием \re{ga}.

После этого, первое уравнение \re{rel1} может быть полностью
разрешено. Оно  приводит к следующим связям между функциями
$A^{(k)}, C^{(k)}, D^{(k)}$:
\begin{eqnarray*}
&& C^{(1)}(u)=1, \quad C^{(6)}(u)=\frac{D^{(6)}(u)D^{(6)}(u+1)}
{A^{(6)}(u)A^{(6)}(u+1)}, \quad D^{(2)}(u)=A^{(2)}(u)C^{(2)}(u),
\\ &&
D^{(3)}(u)=A^{(3)}(u), \quad D^{(4)}(u)=A^{(4)}(u)C^{(4)}(u)C^{(4)}(u-1),
\quad D^{(5)}(u)=\frac{A^{(5)}(u)}{C^{(5)} (u-1)}.
\end{eqnarray*}
Тем не менее, сохраняется одиннадцать неизвестных функций, оставляющих
слишком большую свободу в решениях. Некоторые наводящие соображения для
поиска дальнейших возможных ограничений дает обобщенное разделение
переменных для цепочки Вольтерра с дискретным временем, приводящее
к полиномам Аски-Вильсона описанное в работах \ci{SZ0,SZ1} и теореме \ref{sep-aw}.
Используя эти и другие эвристические аргументы можно прийти к следующему
анзацу обобщенного разделения переменных:
\begin{eqnarray}\nonumber
&& A_n^j =\frac{d(n)\rho(2j+n)}{g(2n+j)g(2n+j-1)\phi(n-j)\phi(n-j-1)},
\qquad B_n^j =1,
\\ &&
C_n^j =\frac{c(n+j)\phi(n-j)\phi(n-j+1)}{\sigma(j)g(2n+j)g(2n+j+1)},
\quad D_n^j =\frac{\rho(2j+n)\phi(n-j)\phi(n-j+1)}{\sigma(j)},
\label{sep}\end{eqnarray}
где $d(0)=0$. Уравнение \re{rel1} удовлетворяется автоматически для
произвольных функций $d(x),\dots,$ $\sigma(x)$. Отметим, что первое,
вторая и третья инволюции нарушают условие $B_n^j =1$ и необходимо
произвести некоторое калибровочное преобразование \re{tr_sup} для
его восстановления. Тогда можно убедиться, что применение инволюций
\re{sep} просто переставляет функции $d(x), c(x), \sigma(x)$ с точностью
до простого изменения их аргументов. Аналогичная ситуация имеет место и
для $g(x), \rho(x), \phi(x)$. Поэтому естественно ожидать, что соответствующие
функции имеют одинаковую форму.

Нам осталось решить уравнения \re{rel2} и \re{rel3}. В некоторых специальных
случаях можно редуцировать \re{rel3} к \re{rel2}.

\begin{proposition}\label{trick}
Предположим, что суперпотенциалы \re{sep}
определяют некоторое решение уравнений \re{rel1} и \re{rel2}, такое
что функции $d(x), \sigma(x), c(x)$ содержат некоторый набор свободных
параметров, не входящих в функции $g(x), \rho(x), \phi(x)$. Тогда
уравнение \re{rel3} удовлетворяется следующим выбором спектральных
параметров
\begin{equation}
\lambda_j=\tilde \sigma(j)/\sigma(j),
\quad \beta_n=\tilde d(n)/d(n), \quad
\alpha_{n+j}=\tilde c(n+j)/c(n+j),
\label{specpar}\end{equation}
где тильдованные функции отличаются от нетильдованных только
выбором свободных параметров.
\end{proposition}

Подставив анзац \re{sep} в уравнение \re{rel2}, мы перепишем его в виде
\ba
&& \frac{c(n +j)}{g(2n +j + 1)g(2n +j)\rho (2j + n)\rho (2j + n + 1)}
\nonumber \\
&& +\frac{d(n + 1)}{g(2n+2+j)g(2n +j + 1)\phi(n + 1 - j)\phi(n - j)}
\nonumber \\
&& - \frac {\sigma(j)}{\rho (2\,j + n)
\,\phi(n + 1 - j)\,\phi(n - j)\,\rho (2\,j + n + 1)} \nonumber \\
&&  = \frac {c(n +j + 1)}{g(2n + 2 + j)g(2n + j + 1)\rho (2j + 2 + n)
\rho(2j + n + 1)}
\nonumber \\ &&
+  \frac {d(n)}{g(2n + j+1)g(2n+ j)\phi(n - j - 1)\phi(n - j)}
\nonumber \\ &&
- \frac {\sigma (j+ 1)}{\rho (2j+ 2 + n)\phi(n - j)\phi(n - j - 1)
\rho(2j + n + 1)}.
\lab{Ans} \ea
Общее решение этого уравнения неизвестно, но его богатый подкласс
может быть найден после наложения дополнительных естественных
ограничений. Потребуем, чтобы функции $g(x), \rho(x), \phi(x)$
имели простые нули в точках $x=x_2, x_1, x_0$ соответственно, где
$x_2, x_1, x_0$ обозначают произвольные постоянные. Потребуем также, чтобы
$g(x)\neq 0$ при $x=x_2-1, x_2-2$,
$\rho(x)\neq 0$ при $x=x_1-1, x_1-2$ и $\phi(x)\neq 0$
при $x=x_0\pm 1$. Условие сокращения полюсов в уравнении \re{Ans} приводит
к соотношениям
\begin{eqnarray*}
&& \frac{c(x_1-x)}{\sigma(x)}=\frac{g(2x_1-3x)\,g(2x_1-3x+1)}
{\phi(x_1-3x)\, \phi(x_1-3x+1)},
\\  &&
\frac{c(x_2-x)}{d(x)}=\frac{\rho(2x_2-3x)\,\rho(2\,x_2-3x+1)}
{\phi(3x-x_2-1)\, \phi(3x-x_2)},
\\  &&
\frac{\sigma(x-x_0)}{d(x)}=\frac{\rho(3x-2x_0-1)\,
\rho(3x-2x_0)} {g(3x-x_0-1)\, g(3x-x_0)}.
\end{eqnarray*}
Эти ограничения удовлетворяются, если положить
\ba\nonumber
&& \phi(x)=\psi(x-x_0), \quad g(x)=\psi(x-x_2), \quad \rho(x)=\psi(x-x_1)
\\ && \sigma(x)=d(x+x_0), \quad c(x)= d(x_2-x),
\lab{csa} \ea
где $\psi(x)$ произвольная нечетная функция $\psi(x) = -\psi(-x)$
с небольшим ограничением на положение своих нулей, упомянутым выше.
При этом необходимо, чтобы параметры $x_0,x_1,x_2$ удовлетворяли
равенству $x_2=x_0+x_1.$
В дальнейшем ограничимся этим частным выбором функций, входящих в анзац \re{sep}.

Очевидно, что теперь остались только две неизвестные функции
$d(x)$ и $\psi(x)$, а уравнение \re{Ans} принимает вид:
\begin{eqnarray}\nonumber
&& \frac{d(x_2-n-j)} {\psi(2n+j-x_2) \psi(2n+j+1-x_2)
\psi(2j+n-x_1)\psi(2j+n+1-x_1)}
\\ \nonumber  && \makebox[1em]{}
+\frac{d(n+1)}{\psi(2n+j+1-x_2) \psi(2n+j+2-x_2) \psi(n-j-x_0)
\psi(n-j+1-x_0)}
\\ \nonumber && \makebox[1em]{}
 - \frac{d(j+x_0)} {\psi(2j+n-x_1)
\psi(2j+n+1-x_1)\psi(n-j-x_0) \psi(n-j+1-x_0)}
\\ \nonumber && \makebox[1em]{}
=\frac{d(x_2-n-j-1)}{ \psi(2n+j+1-x_2) \psi(2n+j+2-x_2)
\psi(2j+n+1-x_1) \psi(2j+n+2-x_1)}
\\ \nonumber && \makebox[1em]{}
+\frac{d(n)}
{\psi(2n+j-x_2)\psi(2n+j+1-x_2)\psi(n-j-1-x_0)\psi(n-j-x_0)}
\\  && \makebox[1em]{}
- \frac{d(j+1+x_0)} {\psi(2j+n+1-x_1) \psi(2j+n+2-x_1) \psi(n-j-1-x_1)
\psi(n-j-x_1)}.
\lab{main}\end{eqnarray}
В дальнейшем будем называть \re{main} основным уравнением.
Предположим, что функции $\psi(x)$ и $d(x)$ являются целыми, то есть
они не имеют сингулярностей при конечных значениях аргумента $x$.
Тогда, из нашего рассмотрения очевидно, что уравнение \re{main}
не содержит полюсов при конечных значениях $n$ и $j$ для произвольных
$\psi(x)$ и $d(x)$, при условии, что $\psi(x)$ имеет только простые нули.

\subsection{Рациональные и тригонометрические решения}

Рассмотрим решения основного уравнения  \re{main}  в классе
рациональных функций, а затем в классе элементарных функций.
Взяв $\psi(x)$ в каком-либо простом виде, можно
анализировать получающееся уравнение на $d(x)$.
Самый простой допустимый выбор соответствует $\psi(x)=x$ и полиномиальному
$d(x)$. Используя компьютерную систему аналитических расчетов легко
проверить, что в этом случае $d(x)$ может быть полиномом шестой степени
\be\lab{polynom}
d(x)=x\prod_{k=1}^5(x-d_k)
\end{equation}
с любопытным ограничением на положение его нулей:
\begin{equation} \label{poly}
\sum_{k=1}^5 d_k=1+2(x_0+x_2).
\end{equation}
Имеется очевидная свобода в умножении $d(x)$ на произвольный множитель,
который мы не указали. Так же один из нулей $d(x)$ был
помещен в точку $x=0$ для того, чтобы иметь $d(0)=0$. Поэтому, осталось
только четыре свободных параметра в $d(x)$. Взяв в предложении \ref{trick}
в качестве $\tilde d(x)$ полином с той же структурой что и $d(x)$:
$$
\tilde d(x)=x\prod_{k=1}^k(x-e_k), \qquad \sum_{k=1}^5e_k=1+2(x_0+x_2),
$$
содержащий четыре других свободных параметра, мы находим соответствующие
спектральные параметры
\begin{equation}
\lambda_j=\prod_{k=1}^5\frac{j+x_0-e_k}{j+x_0-d_k},
\quad \beta_n=\prod_{k=1}^5\frac{n-e_k}{n-d_k}, \quad
\alpha_n=\prod_{k=1}^5\frac{n-x_2 + e_k}{n-x_2+d_k}.
\label{sppar}\end{equation}

Удобно обозначить $s\equiv j+2-x_2$ и $a\equiv 2j+1+x_0-x_2$.
Тогда получаем результат, анонсированный в работе \cite{zhe-sip:umn}.

\begin{theorem}
Рекуррентное соотношение (\ref{recrel}) для полученного рационального
решения $R_{II}$-цепочки \re{polynom}-\re{sppar}
приводит к $R_{II}$-полиномам $P_n^{(j)}(z)$, которые выражаются
через совершенно уравновешенный  2-сбалансированный обобщенный
гипергеометрический ряд $_9F_8$:
$$
P_n^{(j)}(z)= f_n^j (z)
{_9}F_8\left({a,a/2+1,-n,s+n-1,a+2-s-y_1, \dots,a+2-s-y_5\atop a/2,
a+n+1, a+2-s-n, s-1+y_1, \dots, s-1+y_5};1\right),
$$
\begin{equation}\label{9F8}
f_n^j (z)=\frac{(1-z)^n\prod_{k=1}^5(s-1+y_k)_n}{(n+s-1)_n(a+1)_n},
\end{equation}
где $y_1(z), \dots, y_5(z)$ являются корнями следующего алгебраического
уравнения пятой степени:
$$z\prod_{k=1}^5(y-d_k)=\prod_{k=1}^5(y-e_k).$$
\end{theorem}

\smallskip

Согласно общепринятой классификации \cite{gas-rah:basic}, обобщенный
гипергеометрический ряд
$$
{_{r+1}F_r}\left( {a_1, \dots, a_{r+1} \atop b_1, \dots, b_r}; z\right)
=\sum_{n=0}^\infty\frac{(a_1)_n\dots (a_{r+1})_n}{n!(b_1)_n\dots(b_r)_n}z^n
$$
называется вполне уравновешенным если $1+a_1=a_{k+1}+b_k, \; k=1,\dots, r$.
Он становится совершенно уравновешенным если, в дополнение к этим условиям,
выполняется равенство $a_2=a_1/2+1$. Кроме того, он называется
$k$-сбалансированным, если выполняются ограничения
$k+a_1+\dots+ a_{r+1}=b_1+\dots+b_r$ и $z=1$. Такие ряды
обладают некоторыми специальными свойствами.

Известное семейство рациональных функций Вильсона \cite{wil:hypergeometric}
соответствует случаю, ког\-да $\beta_n,$ $\alpha_n,$ $\lambda_n$
сводятся к полиномам второй степени. Такой ситуации можно добиться, потребовав
чтобы функция $d(x)$ была полиномом четвертой степени делящим полином
$\tilde d(x)$. Ключевое новые свойство полиномов (\ref{9F8})
состоит в том, что они содержат восемь независимых свободных параметров
(в \ci{wil:hypergeometric} их было только пять) и что для их
представления в виде гипергеометрического ряда необходимо решать
алгебраическое уравнение выше второго порядка. Фактически, имеется десять
свободных параметров в \re{9F8} в дополнение к степени полиномов
$n$ и их аргументу $z$. Однако, два из них могут быть включены в определение
аргумента $z$ с помощью дробно линейного преобразования \re{rat_P},
которое сохраняет фиксированные ведущие $j,n\to\infty$ асимптотики
$\lambda_j, \alpha_n, \beta_n\to 1$.

Рассмотрим сразу же $q$-обобщение этого решения, которое и будет
строго доказано. Естественно заменить $\psi(x)=x$ в \re{main}
следующей нечетной функцией, определяющей $q$-числа
$$
\psi(x)=\frac{q^{x/2}-q^{-x/2}}{q^{1/2}-q^{-1/2}},
$$
где $q$ обозначает произвольный комплексный параметр.
С помощью компьютерной программы символьных вычислений MAPLE
легко проверить, что $q$-аналог  полинома $d(x)$ имеет вид
(с точностью до общего множителя)
\be
d(x)=\psi(x)\prod_{k=1}^5\psi(x-d_k),
\lab{d(x)}\end{equation}
где необходимо наложить то же самое ограничение \re{poly}
на нули функции $d(x)$. Заметим, что различными предельными переходами
по параметрам $d_k$ можно редуцировать общее число множителей в
произведении \re{d(x)} с 6 до 4,3,2 и 1.

Зафиксировав $\tilde d(x)=\psi(x)\prod_{k=1}^5\psi(x-e_k)$, где $e_k$
удовлетворяют тем же самым ограничениям, что и в рациональном случае,
и подставив это выражение в \re{specpar}, находим
\be
\lambda_j=
\prod_{k=1}^5\frac{1-q^{j+x_0-e_k}}{1-q^{j+x_0-d_k}},
\quad
\beta_n=\prod_{k=1}^5\frac{1-q^{n-e_k}}{1-q^{n-d_k}},
\quad
\alpha_n=\prod_{k=1}^5\frac{1-q^{x_2-n-e_k}}{1-q^{x_2-n-d_k}}.
\lab{lba}\end{equation}
Для полноты, мы приведем также явный вид суперпотенциалов
\be
A_n^j =-\;\frac{(q^{1/2}-q^{-1/2})^{-3}
(1-q^n)(1-aq^{n-1})\prod_{k=1}^5(1-q^{n-d_k})}
{a^{1/2}q^{n/2}(1-sq^{2n-2})(1-sq^{2n-3})(1-sq^{n-1}/a)(1-sq^{n-2}/a)},
\lab{A}\end{equation}
\be
C_n^j =\frac{a^2(1-sq^{n-2})(1-sq^n/a)
(1-sq^{n-1}/a)\prod_{k=1}^5(1-sq^{n+d_k-2})}
{s^2q^{2n-1}(1-sq^{2n-2})(1-sq^{2n-1})(1-aq/s)\prod_{k=1}^5
(1-aq^{1-d_k}/s)},
\lab{C}\end{equation}
\be
D_n^j =-\;(q^{1/2}-q^{-1/2})^3\frac{a^{5/2}(1-aq^{n-1})(1-sq^n/a)(1-sq^{n-1}/a)}
{s^2q^{(3n-1)/2}(1-aq/s)\prod_{k=1}^5(1-aq^{1-d_k}/s)},
\lab{D}\end{equation}
где использованы обозначения
$$
a\equiv q^{2j+1-x_1}, \qquad s\equiv q^{j+2-x_2}.
$$

Согласно стандартному определению \cite{gas-rah:basic}, $q$-гипергеометрический
ряд $_{r+1}\varphi_r$ имеет вид
$$
{_{r+1}\varphi_r}\left( {a_1, \dots, a_{r+1}\atop b_1, \dots , b_r}; q, z
\right)=\sum_{k=0}^{\infty}\frac{(a_1,\dots, a_{r+1};q)_k}
{(q,b_1,\dots,b_r;q)_k}z^k,
$$
где $(a;q)_n$ есть общепринятое обозначение для $q$-факториалов
$$
(a;q)_0=1, \quad (a;q)_n=\prod_{k=1}^n(1-aq^{k-1}), \quad
(a_1, \dots, a_r;q)_n=(a_1;q)_n\dots(a_r;q)_n.
$$
Эти ряды называются вполне уравновешенными, если выполняются ограничения
 $qa_1=a_2b_1=\dots=a_{r+1}b_r$
и совершенно уравновешенными если в дополнение к этим условиям имеют место
равенства $a_2=qa_1^{1/2}$ и $a_3=-qa_1^{1/2}$. Аналогично случаю
$_{r+1}F_r$ рядов, $_{r+1}\varphi_r$ ряд называется сбалансированным
если $qa_1\dots a_{r+1}=b_1\dots b_r$ и $z=q$.

\begin{theorem}
Трехчленное рекуррентное соотношение для $R_{II}$-полиномов \re{recrel}
с рекуррентными коэффициентами, соответствующими решению $R_{II}$-цепочки
в элементарных функциях \re{lba}-\re{D}, имеет следующее явное решение:
\be
P_n^{(j)}(z)= f_n^j (z)
{_{10}}\varphi_9\left({a,qa^{1/2},-qa^{1/2},q^{-n},sq^{n-1},
aq^2/sy_1, \dots,aq^2/sy_5\atop a^{1/2},-a^{1/2},aq^{n+1},
aq^{2-n}/s, sy_1/q, \dots, sy_5/q}; q, q\right),
\lab{10phi9}\end{equation}
где
$$
f_n^j (z)=\frac{(q^{1/2}-q^{-1/2})^{-3n}(z-1)^n\prod_{k=1}^5(sy_k/q; q)_n}
{a^{n/2}q^{n(n+1)/4}(sq^{n-1},aq;q)_n},
$$
и $y_1(z), \dots, y_5(z)$ обозначают решения алгебраического уравнения
пятой степени
\be
(z-1)\prod_{k=1}^5(y_k(z)-y)=
z\prod_{k=1}^5(q^{d_k}-y)-\prod_{k=1}^5(q^{e_k}-y).
\lab{algeq}\end{equation}
\end{theorem}

\noindent
{\bf Доказательство.}
Для того, чтобы найти явный вид $P_n^{(j)}(z)$ мы используем
представление \re{PYZ}. Прежде всего мы найдем коэффициенты
$\zeta_n^j(0)$ в новых обозначениях:
$$
\zeta_n^j (0)=\frac{(q^{1/2}-q^{-1/2})^{-3n}\prod_{k=1}^5(sq^{d_k-1};q)_n}
{a^{n/2}q^{n(n+1)/4}(sq^{n-1},aq;q)_n}.
$$
После этого необходимо вычислить $q$-факториальный вид $Y_n^j $ и
$Z_n^j$, для чего необходимы решения алгебраического уравнения \re{algeq}:
$$
Y_n^j (z)=(z-1)^n\prod_{k=1}^5\frac{(aq^2/sy_k;q)_n}{(aq^{2-d_k}/s;q)_n},
\qquad
Z_n^j (z)=(z-1)^n\prod_{k=1}^5\frac{(sy_k/q;q)_n}{(sq^{d_k-1};q)_n}.
$$
Наконец, решая рекуррентное соотношение \re{rec_eta}, что представляет собой
наиболее сло\-жную часть доказательства, мы находим
$$
\eta_n^j (k)=\frac{({a,qa^{1/2},-qa^{1/2},q^{-n},sq^{n-1},
aq^{2-d_1}/s, \dots,aq^{2-d_5}/s};q)_kq^k}
{(q,a^{1/2},-a^{1/2},aq^{n+1},aq^{2-n}/s, sq^{d_1-1},
\dots, sq^{d_5-1};q)_k}.
$$
Теперь достаточно просто подставить эти выражения в формулу \re{PYZ},
которая и приводит к представлению $P_n^{(j)}(z)$ в указанном виде
совершенно уравновешенного сбалансированного $_{10}\varphi_9$
базисного гипергеометрического ряда.
\hfill{Q.E.D.}

Специальный подкласс полученных $R_{II}$-полиномов соответствует
биортогональным рациональным функциям Рах\-ма\-на-Виль\-со\-на, которые изучались в работах
\cite{gup-mas:contiguous,rah:integral,rah:biorthogonality,
rah-sus:classical,wil:orthogonal}. Для того чтобы получить его,
необходимо выродить $d(x)$ к полиному четвертой степени с нулями
при $d_3=e_3, d_4=e_4, d_5=e_5$. Например, взяв $d_1\to\infty,
d_2\to -\infty$ таким образом, что  $d_1+d_2=e_1+e_2$ остается конечной
величиной. Расходимость в $d(x)$ проявляется только в виде общего
множителя, который может быть устранен преобразованием растяжения.
Тогда $d(x)$ делит $\tilde d(x)$ и получается
$$
\lambda_j=q^{n+x_0-t}+q^{-n-x_0+t} - v, \quad
\beta_n=q^{n-t}+q^{-n+t}-v, \quad
\alpha_n=q^{n-x_2+t}+q^{-n+x_2-t}-v,
$$
где $t=(e_1+e_2)/2,\; v=q^{(e_2-e_1)/2}+q^{(e_1-e_2)/2}.$
Можно показать, что в этой ситуации сопряженные полиномы
$Q_n^{(j)}(z)$ отличаются от $P_n^{(j)}(z)$ только заменой
параметров $x_0, e_1, e_2$ на $x_0-1, e_1-1, e_2-1$, соответственно.
В результате, мы имеем соотношение биортогональности между двумя
$_{10}\varphi_9$ функциями, отличающимися друг от друга только
выбором параметров.

Из соотношения \re{biort} следует, что
общая $_{10}\varphi_9$ функция \re{10phi9} биортогональна
линейной комбинации двух аналогичных $_{10}\varphi_9$ функций.
Для вывода этой комбинации необходимо использовать
трехчленное рекуррентное соотношение и пока неизвестно можно ли
ее свести к одному базисному гипергеометрическому ряду.
В этом случае, суперпотенциалы для сопряженных полиномов
$Q_n^{(j)}(z)$ зависят от параметров $e_k$ (чего нет для
$P_n^{(j)}(z)$) и имеют значительно более сложный
вид, чем \re{A}--\re{D}. Отметим, что мы можем построить сопряженные
полиномы для $Q_n^{(j)}(z)$ точно так же как это было сделано для
$P_n^{(j)}(z)$ и они не будут совпадать с $P_n^{(j)}(z)$ или
$Q_n^{(j)}(z)$. Это следует из того факта, что в общем случае
изменения спектральных параметров $\beta_n\to \beta_{n+1}$ и
$\alpha_n\to\alpha_{n-1}$, вызванные переходом к сопряженным полиномам
\re{compar}, не могут быть скомпенсированы переопределением параметров
системы. Очевидно, что эти переходы к сопряженным полиномам могут
итерироваться до бесконечности. При этом на каждом шаге мы имеем дело с
новым элементарным решением  $R_{II}$-цепочки и специфическим условием
биортогональности между линейными комбинациями  $_{10}\varphi_9$ рядов.

\subsection{Эллиптические решения, простейший случай}

Описанные в предыдущем параграфе решения основного уравнения
\re{main} допускают обобщение, описываемое эллиптическими тета-функциями.
Тета-функция Якоби $\theta_1(u)$ имеет вид \cite{whi-wat:course}
\ba\nonumber
&& \theta_1(u)= 2\sum_{n=0}^\infty (-1)^n p^{(n+1/2)^2/2}\sin \pi (2n+1)u
\\ && \makebox[2em]{}
=2p^{1/8} \sin \pi u\: \prod_{n=1}^{\infty}\left(1-2p^{n} \cos 2\pi u + p^{2n}
\right) (1-p^{n}),
\lab{ell_H} \ea
где $p$ обозначает произвольный комплексный параметр, $|p|<1$.
Модулярный параметр $\tau$ вводится обычным образом $p=\exp (2\pi i\tau)$.
Перечислим наиболее важные свойства этой функции:

(i) $\theta_1(u)$ нечетна, $\theta_1(-u)=-\theta(u)$;

(ii) $\theta_1(u)$ квазипериодична по отношению к сдвигам аргумента
на $1$ и $\tau$
\be
\theta_1(u+1)=-\theta(u), \quad \theta(u+\tau)= -p^{-1}\:
\exp(-2i\pi u)\: \theta_1(u); \lab{qp} \end{equation}

(iii) она удовлетворяет алгебраическому соотношению
\ba
&&\theta_1(x+z)\theta_1(x-z)\theta_1(y+w)\theta_1(y-w)-
\theta_1(x+w)\theta_1(x-w)\theta_1(y+z)\theta_1(y-z) \nonumber \\
&&\makebox[4em]{}
=\theta_1(x+y)\theta_1(x-y)\theta_1(z+w)\theta_1(z-w)
\lab{Riem} \ea
при произвольных $x,y,z,w\in\C $.

Следуя работам \cite{djkmo:exactly1,djkmo:exactly2,fre-tur:elliptic},
определим ``эллиптические числа"
\be
[x;h,\tau]= \frac{\theta_1(hx)}{\theta_1(\pi h)},
\lab{def_e} \end{equation}
где $h$ обозначает произвольный параметр. Ясно, что эти числа зависят от
трех переменных $x, h$ и $\tau$, но для удобства в обозначениях мы будем
опускать зависимость от $h$ и $\tau$: $[x]\equiv [x; h,\tau]$.
Они обладают свойствами
\ba
(i)&& \quad [-x]=-[x]; \nonumber \\
(ii)&&\quad [x+1/h]= -[x], \quad [x+ \tau/h]=-\exp(-i\pi\tau -2\pi i hx)\: [x];
\lab{quasi_e} \\
(iii)&& \quad
[x+z][x-z][y+w][y-w]-[x+w][x-w][y+z][y-z] \nonumber \\
&& \makebox[4em]{} = [x+y][x-y][z+w][z-w];
\lab{Riem_e} \\
(iv)&& \quad
\lim_{Im(\tau) \to +\infty}[x; h,\tau] = \frac{\sin(\pi h x)}{\sin(\pi h)};
\lab{lim} \\
(v) && \quad
\lim_{h \to 0} [x;h,\tau]=x.
\lab{lim2} \ea

Свойство (iv) означает, что в пределе $Im(\tau) \to +\infty$ эллиптические
числа становятся $q$-числами с $q=e^{2\pi i h}$, использовавшимися в
предыдущем параграфе. Свойство (v) устанавливает связь с обычными числами.

Определим эллиптические аналоги символа Похгаммера и факториала
$[x]_n=[x][x+1]\cdots [x+n-1], \; [n]!= [1]_n$
и приступим теперь к построению эллиптических решений основного уравнения
\re{main}. Функции $[x]$ и $\psi(x)$ обладают одинаковыми свойствами:
обе нечетны и в пределах  (iv) и (v) $[x]$ совпадает с $\psi(x)$
для элементарных и рациональных решений. Поэтому естественно положить
$\psi(x) = [x]$ и для функции $d(x)$ выбрать следующий анзац
\be
d(x)=[x]\prod_{k=1}^5[x-d_k],
\lab{An_d} \end{equation}
в котором параметры $d_k$ ограничены условием \re{poly}.
Тогда предельные переходы к предыдущим решениям очевидны.
Для того чтобы доказать, что \re{An_d} определяет решение $R_{II}$-цепочки
\re{main}, необходимо привести получающееся уравнение
к виду двояко-периодической функции без сингулярностей в фундаментальном
параллелограмме квазипериодов функции $\theta_1(x)$. По теореме Лиувилля
такая функция должна быть постоянной, значение которой легко устанавливается.

Рассмотрим комбинацию тета-функций
$$
R(n)=\Biggl( \frac{d(x_2-n-j)} {[2n+j-x_2][2n+j-x_2+1][2j+n+x_0-x_2]
[2j+n+x_0-x_2+1]}
$$
$$
+\frac{d(n+1)}{[2n+j-x_2+1][2n+j-x_2+2][n-j-x_0][n-j-x_0+1]}
$$
$$
-\frac{d(x_2-n-j-1)}{[2n+j-x_2+1][2n+j-x_2+2][2j+n+x_0-x_2+1]
[2j+n+x_0-x_2+2]}
$$
$$
-\frac{d(n)}{[2n+j-x_2][2n+j-x_2+1][n-j-x_0-1][n-j-x_0]}
$$
$$
+ \frac{d(j+x_0+1)}
{[2j+n+x_0-x_2+1][2j+n+x_0-x_2+2][n-j-x_0-1][n-j-x_0]}\Biggr)
$$
$$
\times [2j+n+x_0-x_2][2j+n+x_0-x_2+1][n-j-x_0]
[n-j-x_0+1],
$$
с $d(x)$ заданной выражением \re{An_d}. Будем рассматривать $n$ как
непрерывную комплексную переменную. Тогда нетрудно убедиться, что
$R(n+1/h)=R(n)$ и $R(n+\tau/h)=R(n)$ благодаря специальному ограничению
\re{poly}. По построению, все полюсы $R(n)$ как функции $n$
были устранены заранее специальным выбором функций $c(x), \sigma(x), \phi(x)$ и
$\rho(x)$, то есть $R(n)$ оказывается целой двояко-периодической функцией.
По теореме Лиувилля $R(n)=C_1$ есть величина, не зависящая от $n$, которая
однако может зависеть от других переменных $j, x_2, x_0$. Для того, чтобы
доказать, что $C_1=d(j+x_0)$, что будет подразумевать справедливость
равенства \re{main}, необходимо рассмотреть комбинацию
$$
S(j)=\Biggl( \frac{d(x_2-j)}
{[j-x_2][j-x_2+1][2j+x_0-x_2][2j+x_0-x_2+1]} $$ $$
-\frac{d(j+x_0)}{[2j+x_0-x_2][2j+x_0-x_2+1][-j-x_0][-j-x_0+1]}
$$
$$
-\frac{d(x_2-j-1)}{[j-x_2+1][j-x_2+2][2j+x_0-x_2+1][2j+x_0-x_2+2]}
$$
$$
-\frac{d(0)}{[j-x_2][j-x_2+1][-j-x_0-1][-j-x_0]}
$$
$$
+ \frac{d(j+x_0+1)}
{[2j+x_0-x_2+1][2j+x_0-x_2+2][-j-x_0-1][-j-x_0]}\Biggr)
$$
$$
\times [j-x_2+1][j-x_2+2][-j-x_0][-j-x_0+1].
$$
Опять, рассматривая $j$ как непрерывную переменную, можно проверить, что
$S(j+1/h)=S(j)$ и $S(j+\tau/h)=S(j)$. По построению $S(j)$ не имеет полюсов
по $j$. Поэтому, $S(j)=C_2$ есть величина, не зависящая от $j$.
Взяв предел $j\to -x_0$, можно заключить, что $C_2=-d(1)$, откуда вытекает
$C_1=d(j+x_0)$.

Таким образом, функция $d(x)$ \re{An_d} удовлетворяет
основному уравнению \re{main} при выполнении ограничения \re{poly}.
Отметим важное отличие полученного эллиптического решения от предыдущих:
в этом случае нельзя удалять параметры
в бесконечную точку и таким образом понижать число сомножителей
$[x]$ в $d(x)$ из-за квазипериодичности тета-функций.
Можно взять другое аналогичное решение уравнения \re{main}
$$
\tilde d(x)=[x]\: \prod_{k=1}^5[x-e_k], \quad \sum_{k=1}^5 e_k=1+2(x_0+x_2),
$$
и с его помощью построить спектральные параметры $\alpha_n, \beta_n,
\lambda_n$ так же как и в предыдущих случаях.

Рассмотрим подробнее решения, соответствующие ограничениям
\be
e_3 = d_3, \; e_4 = d_4, \; e_5 = d_5, \quad e_1 + e_2 =
d_1 + d_2. \lab{t_par} \end{equation}
Используя соотношения \re{specpar}, мы получаем выражения
\ba
&&\beta_k = \frac{[k-e_1][k- e_2]}{[k- d_1][k- d_2]}, \lab{e_bt} \\
&&\alpha_k = \frac{[k-x_2+ e_1][k-x_2+e_2]}{[k-x_2+ d_1][k-x_2+ d_2]},
\lab{e_al} \\
&&\lambda_k = \frac{[k+x_0-e_1][k+x_0- e_2]}{[k+x_0- d_1][k+x_0- d_2]}.
\lab{l_bt} \ea

Рассмотрим теперь выражение \re{PYZ} для $R_{II}$-полиномов.
Прежде всего выберем следующую параметризацию аргумента $z$:
\be
z(\xi) = \frac{[\xi][\xi+e_2 - e_1]}{[\xi +d_2 - e_1][\xi + d_1 - e_1]}.
\lab{z_xi} \end{equation}
Пользуясь тождеством \re{Riem_e}, можно получить
\ba
&&z(\xi)- \lambda_{k}=
\frac{[k+\xi+x_0-e_1][k-\xi+x_0-e_2][d_2-e_1][e_1 -d_1]}
{[\xi+d_2-e_1][\xi+d_1-e_1][k+x_0-d_1][k+x_0-d_2]}, \lab{zl} \\
&&z(\xi) -\alpha_{k} =\frac{[k+\xi-x_2+e_2]
[k-\xi-x_2+e_1][d_2-e_1][e_1 -d_1]} {[\xi+d_2-e_1]
[\xi+d_1-e_1][k-x_2+d_1][k-x_2+d_2]}. \lab{za} \ea
Следовательно,
\ba
&&Z_n^j (z)= \prod_{k=1}^n(z-\alpha_{j+k}) =
\left(\frac{[d_2-e_1][e_1 - d_1]}{[\xi+d_2-e_1][\xi + d_1 -e_1]}\right)^n
\nonumber \\
&&\makebox[4em]{}
\times \frac{[1+j+\xi+e_2 - x_2]_n [1+ j - \xi + e_1 - x_2]_n}
{[1+j-x_2+d_1]_n[1+j-x_2+d_2]_n}, \lab{Z_e} \\
&&Y_n^j (z)= \prod_{k=1}^n(z-\lambda_{j+k}) =
\left(\frac{[d_2-e_1][e_1 - d_1]}{[\xi+d_2-e_1][\xi + d_1 - e_1]}\right)^n
\nonumber \\
&&\makebox[4em]{}\times \frac{[1+j+\xi-e_1 + x_0]_n [1+ j - \xi - e_2 + x_0]_n}
{[1+j+x_0-d_1]_n[1+j+x_0-d_2]_n}. \lab{Y_e} \ea

Определим теперь коэффициенты $\eta_n^j (k)$ из разностного уравнения
\re{rec_eta}. Коэффициенты $\zeta_n^j (0)$ имеют вид
\be
\zeta_n^j (0)=(-1)^n \: \prod_{m=0}^{n-1} \frac{C_m^{j+1} }{D_m^{j+1} } =
(-1)^n \frac{[j+1-x_2]_n \: \prod_{k=1}^5 [j+1-x_2+d_k]_n}
{[j+1-x_2]_{2n} [2j+2-x_1]_n}. \lab{e_zeta} \end{equation}
Следовательно,
\ba
&&\frac{\zeta_n^{j+1} (0)}{\zeta_n^j (0)C_n^{j+1} } =
\frac{[j+x_0+1][2j+2-x_1][2j+3-x_1][2n+j+2-x_2]}
{[n-j-x_0][n-j-x_0-1][2j+n+2-x_1][2j+n+3-x_1]}
\nonumber \\
&&\makebox[6em]{}\times \prod_{k=1}^5{\frac{[j+1+x_0-d_k]}{[j+1-x_2+d_k]}}.
\lab{L_e} \ea
Рассмотрим теперь анзац
\be
\eta_n^j (k)=G(k;j) \: \frac{[-n]_k [1-x_2+j+n]_k}{[x_0+1-n+j]_k \:
[2-x_1+n+2j]_k},
\lab{eta_G} \end{equation}
где $G(k,j)$ неизвестные коэффициенты. Подставив выражения \re{L_e} и \re{eta_G}
в соотношение \re{rec_eta} и используя тождество \re{Riem_e}, можно увидеть, что
выражения, содержащие зависимость от $n$, сокращаются и остается
уравнение на $G(k;j)$:
\be
\frac{G(k;j)}{G(k-1;j+1)} = \frac{[2j+2-x_1][2j+3-x_1]}{[k][k+2j-x_1+1]}
\prod_{m=1}^5\frac{[j+1+x_0-d_m]}{[j+1-x_2+d_m]} \lab{eq_G} \end{equation}
с начальным условием $G(0;j)=1$. Легко проверить, что единственным решением
уравнения \re{eq_G} с заданными начальными условиями является функция
\be
G(k;j)= \frac{[-x_1+1+2j]_k[1-x_1+2j+2k]}{[k]![1-x_1+2j]}
\prod_{m=1}^5 \frac{[j+1+x_0-d_m]_k}{[j+1-x_2+d_m]_k} .
\lab{sol_G} \end{equation}
Таким образом, коэффициенты $\eta_n^j (k)$ равны \re{eta_G},
где $G(k;j)$ определены в \re{sol_G}.

Подставляя  \re{Y_e}, \re{Z_e}, \re{eta_G} в определение \re{PYZ},
мы получаем выражение
\ba\nonumber
\lefteqn{ P_n^{(j)}(z)= Z_n^j (z) \zeta_n^j (0) \: {_{12}}v_{11}
(2j+1-x_1; -n, 1+j-x_2+n, j+1+x_0-d_3, j+1  }&& \\
&&  +x_0-d_4,j+1+x_0-d_5, j+1+ \xi + x_0 -
e_1, j+1 - \xi + x_0 - e_2;h, \tau; 1 ), \lab{PW} \ea
где символ $_{12}v_{11}$ обозначает совершенно уравновешенный эллиптический
гипергеометрический ряд со специальным значением аргумента.
Общий ряд такого типа  имеет вид
\be
{_{r+1}}v_r(a_0;\: a_1, a_2, \dots, a_{r-4}; h, \tau; z) =
\sum_{k=0}^{\infty} \frac{[a_0+2k]}{[a_0]} \:
\prod_{m=0}^{r-4} \frac{[a_{m}]_k} {[1+a_0-a_{m}]_k}z^k,
\lab{def_W} \end{equation}
где параметры $a_0, \dots, a_{r-4}$ удовлетворяют условию балансировки
\be
 \sum_{m=1}^{r-4} a_{m} =\frac{r-7}{2}+ \frac{r-5}{2}a_0.
 \lab{cond_a} \end{equation}
Для того, чтобы избежать проблем со сходимостью ряда \re{def_W},
необходимо положить один из параметров $a_k$ равным отрицательному
целому числу. В нашем случае условие \re{cond_a} выполняется
благодаря ограничениям $x_2=x_0+x_1$ и \re{poly}. Для
$Im\:(\tau)\to+\infty$ этот $_{r+1}v_r$ ряд переходит в
совершенно уравновешенный сбалансированный $q$-гипер\-гео\-ме\-три\-чес\-кий ряд
$_{r-1}\varphi_{r-2}$.

\begin{remark}
В оригинальной работе Френкеля и Тураева \cite{fre-tur:elliptic}
ряд $_{12}v_{11}(\ldots;h,\tau;1)$
обозначался как $_{10}\omega_9(\ldots)$. В статье \cite{spi-zhe:spectral}
использовалось обозначение $_{10}E_9(\ldots)$. Однако, в работах
\cite{spi:theta,spi:bailey} было прояснено происхождение условия совершенной
уравновешенности и символ $_{10}E_9(\ldots)$ был заменен на логически более
оправданные $_{12}E_{11}(\ldots)$ для рядов общего типа и $_{12}V_{11}(\ldots)$
для совершенно уравновешенных рядов. В новом издании книги \cite{gas-rah2}
большие буквы зарезервированы за рядами в мультипликативной системе обозначений,
а в аддитивной системе, такой как в выражении \re{def_W},
предложено пользоваться маленькими буквами. Подробнее об этих соглашениях
рассказано в следующей главе.
\end{remark}

Так же как и в рациональном и тригонометрическом случаях, из
\re{Z_e}, \re{Y_e} и \re{PW} видно, что зависимость от $j$ входит
только в комбинациях $x_0 +j, x_1 -2j, x_2-j$, то есть сдвиги
$j\to j\pm 1$ эквивалентны простым переопределениям параметров.
Таким образом, функции $P_n^{(j)}(z)$ \re{PW} определяют $R_{II}$-полиномы от $z$,
удовлетворяющие рекуррентным соотношениям \re{r1}, \re{r2} и
условиям биортогональности \re{bio_1} по отношению к некоторому
функционалу ${\cal L}_j$. Предполагается, что для этой системы существует
обобщение теоремы Леонарда о самодуальных (или биспектральных \cite{DG})
ортогональных полиномах \cite{leo:orthogonal}, полностью характеризующей
полиномы Аски-Вильсона \cite{ask-wil:some}, на случай биортогональных рациональные
функций. С точки зрения обобщения теоремы Магнуса \cite{Mag}, полученного
в работах \cite{spi-zhe:gevp,spi-zhe:grids}, эллиптические сетки
\re{z_xi} фиксируются однозначно из условия существования понижающих
операторов в пространстве рациональных функций.
Интересно было бы прояснить связь этих функций с эллиптическими
рациональными функциями Золотарева \cite{zol:sur}.

\subsection{Общее эллиптическое решение}

Мы рассматривали специальное эллиптическое решение $R_{II}$-цепочки
\re{t_par}. Обсудим теперь общую ситуацию когда все параметры $e_i$
полинома $\tilde  d(x)=[x]\prod_{i=1}^5[x-e_i]$
отличаются от $d_i$ и удовлетворяют единственному ограничению
\be
\sum_{i=1}^5 e_i= \sum_{i=1}^5 d_i=1+2(x_0+x_2). \lab{res_ed} \ee
Тогда,
\be
\alpha_k= \prod_{i=1}^5 {\frac{[k-x_2+e_i]}{[k-x_2+d_i]}}, \quad
\beta_k= \prod_{i=1}^5 {\frac{[k-e_i]}{[k-d_i]}}, \quad
\lambda_k= \prod_{i=1}^5 {\frac{[k+x_0-e_i]}{[k+x_0-d_i]}}
\lab{gen_abm}\ee
и удобная параметризация выражений $z-\lambda_k$ и $z-\alpha_k$ определяется
равенством
\be z-\prod_{i=1}^5 \frac{[x-e_i]}{[x-d_i]} = \kappa(z)\:  \prod_{i=1}^5
\frac{[x-\nu_i(z)]} {[x-d_i]}, \lab{fund_1} \ee
где $\kappa(z), \nu_i(z)$ не зависят от $x$.

Используя соотношение \re{fund_1} и свойство $[-x]=-[x]$, находим
\ba
&&Z_n^j(z)=\prod_{k=1}^n (z-\alpha_{j+k}) = \kappa^n(z) \: \prod_{i=1}^5
\frac{[j+1-x_2+\nu_i(z)]_n}{[j+1-x_2+d_i]_n}, \nonumber \\
&&Y_n^j(z)=\prod_{k=1}^n (z-\lambda_{j+k}) = \kappa^n(z) \: \prod_{i=1}^5
\frac{[j+1+x_0-\nu_i(z)]_n}{[j+1+x_0-d_i]_n}. \lab{ZY_gen} \ea
Биортогональные рациональные функции $R_n^j(z)$
имеют форму (в несколько другой нормировке, чем в \re{biort}):
\be
R_n^j(z)= \frac{P_n^j(z)}{Z_n^j(z)\zeta_n^j(0)}=
\sum_{k=0}^n \eta_n^j(k) \frac{Y_k^j(z)}{Z_k^j(z)}, \lab{R_gen} \ee
где коэффициенты $\eta_n^j(k)$ определены в \re{eta_G}, \re{sol_G}.
Подставляя выражения \re{ZY_gen} и \re{eta_G} в определение \re{R_gen}, получаем
\ba
&&R_n^j(z)= {_{12}}v_{11}(2j+1-x_1;-n, j+1-x_2+n, j+1+x_0-\nu_1(z),
\nonumber \\ && \makebox[6em]{}
j+1+x_0-\nu_2(z), \dots, j+1+x_0-\nu_5(z); h, \tau ;1).
\lab{R_e_g} \ea
Таким образом, $R_n^j(z)$ опять выражаются через совершенно уравновешенный
эллиптический  гипергеометрический ряд, но с другой параметризацией и очень
сложной формой нулей $\nu_i(z)$. Решение $R_{II}$-цепочки, приводящее к
функции \re{R_e_g} содержит 12 параметров: $x_0, x_2, d_1, \dots,
d_4, e_1, \dots, e_4, h, \tau$. Еще один свободный параметр появляется
как отношение коэффициентов пропорциональности полиномов $d(x)$ и $\tilde  d(x)$.
Дробно-линейные преобразования $z$ позволяют зафиксировать три параметра, т.е.
остается всего 10 независимых параметров. К сожалению, для этой общей системы
еще не удалось найти замкнутого выражения для функционала ${\cal L}_j$,
определяющего условие биортогональности. Более того, сопряженные рациональные
функции имеют теперь намного более сложный вид чем в предыдущем
случае \re{t_par}.

Укажем кратко на другую автомодельную редукцию $R_{II}$-цепочки.
Она происходит из инвариантности относительно преобразований Мебиуса и
дискретных сдвигов по решеткам переменных $n$ и $j$.
Определим это автомодельное решение, потребовав чтобы сдвиг
 $j\to j+M$, $M\in\N$, был эквивалентен сдвигу по другой дискретной
переменной $n\to n+k, k\in \Z,$ в комбинации с преобразованием
\re{partr}, \re{rat_pot} для фиксированных $\xi,
\zeta, \sigma, \eta$:
$$
\lambda_{j+M}=\frac{\sigma \lambda_j -\eta}{\xi -\zeta \lambda_j}, \quad
\alpha_{j+M}=\frac{\sigma \alpha_{j+k} -\eta}{\xi -\zeta \alpha_{j+k}}, \quad
\beta_n=\frac{\sigma \beta_{n+k} -\eta}{\xi -\zeta \beta_{n+k}},
$$
\be
A_n^{j+M}=A_{n+k}^j(\xi-\zeta b_{n+k}), \quad
B_n^{j+M}=B_{n+k}^j,
\lab{ssr}\ee
$$
C_n^{j+M}=C_{n+k}^j\frac{\xi-\zeta\alpha_{n+j+k}}{\xi-\zeta\lambda_j}, \quad
D_n^{j+M}=\frac{D_{n+k}^j}{\xi-\zeta\lambda_j}.
$$
Тогда спектральные коэффициенты $\lambda_j$ состоят из суперпозиции $M$
независимых последовательностей вида:
\be
\lambda_{Mi+m}=\frac{a_mq^i+b_m}{c_mq^i+d_m}, \quad m=1, 2, \dots M,
\lab{newseries}\ee
с некоторыми постоянными $a_m,b_m,c_m,d_m,q$.
Коэффициенты  $\beta_n$ и $\alpha_{n+j}$ состоят из $|k|$ и $|M-k|$
последовательностей чисел аналогичного вида.
Это решение $R_{II}$-цепочки является естественным обобщением автомодельных
потенциалов уравнения Шредингера, описанных в начале этой главы.

\section{Нелинейная цепочка для симметричных $R_{II}$-по\-ли\-но\-мов}

Известно, что матрицы Якоби с нулевыми диагональными элементами, называемые
в дальнейшем двухдиагональными матрицами, приводят к нелинейным интегрируемым
системам с дополнительными симметриями. Наложим на основное трехчленное
рекуррентное соотношение (\ref{recrel}) ограничения $b_n=0$,
$\alpha_n=-\beta_n\equiv\gamma_n$. В новых обозначениях оно принимает вид
\begin{equation}
V_{n+1}(z)+v_n(z^2-\gamma_n^2)V_{n-1}(z)=z\rho_n V_n(z).
\label{sympol}\end{equation}
Для начальных условий $V_0=1, V_1=\rho_0 z$ это рекуррентное соотношение
генерирует симметричные полиномы, обладающие фиксированной четностью
$V_n(-z)=(-1)^nV_n(z)$. Калибровочное преобразование
$V_n=\prod_{k=0}^{n-1} \rho_k\: U_n$ устраняет $\rho_n$
$$
U_{n+1}(z)+\tilde v_n(z^2-\gamma_n^2)U_{n-1}(z)=zU_n(z),
\quad \tilde v_n=v_n/\rho_n\rho_{n-1}.
$$
Нетрудно видеть, что множества полиномов нечетной и четной степеней
$U_{2n+1}$ и $U_{2n}$ раздельно удовлетворяют трехчленным рекуррентным
соотношениям. Так,
\begin{eqnarray}
&U_{2n+2}(z)+r_n(b_n-z^2)U_{2n}(z)+
u_n(z^2-\alpha_n)(z^2-\beta_n)U_{2n-2}(z)=0, &
\nonumber \\
&r_n=1-\tilde v_{2n}-\tilde v_{2n+1}, \quad
r_nb_n=-\gamma_{2n}^2\tilde v_{2n}-\gamma_{2n+1}^2\tilde v_{2n+1}, &
\nonumber \\
&u_n=\tilde v_{2n-1}\tilde v_{2n}, \quad \alpha_n=\gamma_{2n-1}^2,
\quad \beta_n=\gamma_{2n}^2. &
\nonumber \end{eqnarray}
Видно, что это есть общее рекуррентное соотношение (\ref{recrel})
для полиномов $P_n(z^2)\equiv U_{2n}(z)$. Аналогично, в нечетном случае имеем
\begin{eqnarray}
&U_{2n+3}(z)+r_n(b_n-z^2)U_{2n+1}(z)+
u_n(z^2-\alpha_n)(z^2-\beta_n)U_{2n-1}(z)=0, &
\nonumber \\
&r_n=1-\tilde v_{2n+1}-\tilde v_{2n+2}, \quad
r_nb_n=-\gamma_{2n+1}^2\tilde v_{2n+1}-\gamma_{2n+2}^2\tilde v_{2n+2}, &
\nonumber \\
&u_n=\tilde v_{2n}\tilde v_{2n+1}, \quad \alpha_n=\gamma_{2n}^2,
\quad \beta_n=\gamma_{2n+1}^2. &
\nonumber\end{eqnarray}
Это есть рекуррентное соотношение  (\ref{recrel}) для полиномов
$P_n(z^2)\equiv U_{2n+1}(z)/z$. Эта конструкция сформулирована в работе
автора \cite{spi:solitons}, в которой предложено называть цепные дроби,
порожденные соотношением (\ref{sympol}), $V$-дробями. Они связаны с
$R_{II}$-дробями работы \cite{ism-mas:general} точно так же как $S$-дроби
связаны с $J$-дробями в теории ортогональным полиномов.

Определим симметричные БРФ
$R_n(z,\gamma_k)=V_n(z)/\prod_{l=1}^n(z-\gamma_l)$. Они удовлетворяют
равенству $LR_n(z)=zMR_n(z)$, где $LR_n\equiv \gamma_{n+1}R_{n+1}-
v_n\gamma_nR_{n-1}$ и $MR_n\equiv R_{n+1}-\rho_nR_n+v_nR_{n-1}$.
Преобразования четности принимают вид $R_n(-z,\gamma_k)=R_n(z,-\gamma_k)$.

Применение аналога преобразований Кристоффеля (\ref{r1}) к симметричным БРФ
нарушает симметрию четности. Это очевидно поскольку оно преобразует
$\alpha_n\to\alpha_{n+1}$ и не меняет $\beta_n$. Однако, $\alpha_n$
и $\beta_n$ входят в первоначальное рекуррентное соотношение (\ref{recrel})
симметричным образом. Поэтому можно произвести второе преобразование
Кристоффеля с переставленными $\alpha_n$ и $\beta_n$. Это приводит
к $\beta_n\to\beta_{n+1}$ и не меняет $\alpha_n$. В результате
восстанавливается необходимая симметрия между $\alpha_n$ и $\beta_n$
так что опять можно редуцировать рекуррентное соотношение для БРФ к симметричной
форме. После ряда технических вычислений автором был получен
следующий аналог преобразований Кристоффеля, отображающий одни симметричные БРФ
на другие симметричные БРФ. Обозначим $R_n^{(j)}(z)$ симметричные БРФ, удовлетворяющие
соотношению $L_jR_n^{(j)}(z)=zM_jR_n^{(j)}(z)$, или
\begin{equation}\label{ttr}
(z-\gamma_{n+1}^j )R_{n+1}^{(j)}(z)+v_n^j (z+\gamma_n^j )R_{n-1}^{(j)}(z)=z\rho_n^j R_n^{(j)}(z).
\end{equation}
Сдвиг дискретного времени $j\to j+1$ (что не должно быть перепутано со сдвигом
$j$ в соотношении (\ref{r1})) определяется следующим образом:
\begin{eqnarray}\nonumber
&& R_n^{(j+1)}(z)=\frac{z-\gamma_1^j }{z^2-\kappa_j^2}
\Bigl((z-\gamma_{n+2}^j )R_{n+2}^{(j)}(z)
- \Bigl(1+\frac{\gamma_{n+2}^j }{\gamma_{n+1}^j }\Bigr)\frac{R_{n+2}^{(j)}(\kappa_j)}
{R_{n+1}^{(j)}(\kappa_j)}\: zR_{n+1}^{(j)}(z)
\\ && \makebox[6em]{}
+\frac{\gamma_{n+2}^j R_{n+2}^{(j)}(\kappa_j)}{\gamma_{n+1}^j R_n^{(j)}(\kappa_j)}
(z+\gamma_{n+1}^j )R_n^{(j)}(z) \Bigr).
\label{dtv} \end{eqnarray}
Дадим явное выражение для коэффициентов перед $R_{n}^{(j)}(z)$ в терминах
решения уравнения $L_jR_n^{(j)}(\kappa_j)=\kappa_jM_jR_n^{(j)}(\kappa_j)$,
где $\kappa_j$ произвольный параметр. Из условия совместности определения (\ref{dtv}) с
рекуррентным соотношением (\ref{ttr}) имеем $\gamma_n^{j+1}=\gamma_{n+1}^j $,
т.е. $\gamma_n^j  \equiv\gamma_{n+j}$, и
\begin{eqnarray}
&v_n^{j+1}=v_n^j \frac{\gamma_{n+j}\gamma_{n+j+2} R_{n-1}^{(j)}(\kappa_j)
M_jR_{n+2}^{(j)}(\kappa_j)}{\gamma_{n+j+1}^2R_{n+1}^{(j)}(\kappa_j)
M_jR_n^{(j)}(\kappa_j)},&
\nonumber \\
&\rho_n^{j+1}=\frac{R_{n+1}^{(j)}(\kappa_j)M_jR_{n+2}^{(j)}(\kappa_j)}
{R_{n+2}^{(j)}(\kappa_j)M_jR_{n+1}^{(j)}(\kappa_j)}
\Bigl(\rho_n^j +\frac{(\gamma_{n+j+2}-\gamma_{n+j+1})L_jR_{n+1}^{(j)}(\kappa_j)}
{\gamma_{n+j+1}\gamma_{n+j+2}R_{n+1}^{(j)}(\kappa_j)}\Bigr).&
\nonumber\end{eqnarray}

Определим $G_n^j \equiv v_n^j R_{n-1}^{(j)}(\kappa_j)/\rho_n^j R_n^j (\kappa_j)$.
Тогда из соотношения (\ref{ttr}) можно вывести факторизацию
$$
\frac{v_n^j }{\rho_n^j \rho_{n-1}^j }=
\frac{G_n^j (\kappa_j-(\kappa_j+\gamma_{n+j-1})G_{n-1}^j )}
{\kappa_j-\gamma_{n+j}}.
$$
Аналогично находим
$$
\frac{R_{n+1}^{(j)}}{\rho_n^j R_n^{(j)}}=
\frac{\kappa_j-(\kappa_j+\gamma_{n+j})G_n^j }{\kappa_j-\gamma_{n+j+1}},
\quad \frac{M_jR_n^{(j)}}{\rho_n^j R_n^{(j)}}=
\frac{\gamma_{n+j+1}-(\gamma_{n+j+1}
+\gamma_{n+j})G_n^j }{\kappa_j-\gamma_{n+j+1}},
$$
где $R_n^{(j)}=R_n^{(j)}(\kappa_j)$.
Применяя эти формулы к калибровочно инвариантной комбинации
$v_n^{j+1}/\rho_n^{j+1}\rho_{n-1}^{j+1}$, мы получаем нелинейную
цепочку с дискретным временем
\begin{eqnarray}
&\gamma_{n+j}^2\gamma_{n+j+2}^2\frac{(\kappa_j-\gamma_{n+j+1})G_n^j
(\kappa_j-(\kappa_j+\gamma_{n+j+1})G_{n+1}^j  )}
{G_n^{j+1}(\kappa_{j+1}-(\kappa_{j+1}+\gamma_{n+j})G_{n-1}^{j+1})}
=\frac{\chi_n^j \chi_{n+1}^j }{\kappa_{j+1}-\gamma_{n+j+1}},&
\nonumber \\
&\chi_n^j \equiv
{(\gamma_{n+j}^2-\gamma_{n+j+1}^2)\kappa_jG_n^j +\gamma_{n+j+1}^2
(\kappa_j-\gamma_{n+j})}.&
\nonumber \end{eqnarray}
Это уравнение обобщает цепочку Вольтерра с дискретным временем \re{dtvl},
выведенную в работе \cite{SZ1}, и $g$-алгоритм Бауэра \cite{Bau}.

%% file: CHAPTER2.TEX
\chapter[Общая теория тета-гипергеометрических рядов]
{Общая теория тета-гипер\-гео\-ме\-три\-чес\-ких рядов}

В этой главе мы опишем общую схему построения рядов гипергеометрического
типа, основанных на тета-функциях Якоби. Основное внимание
будет уделяться эллиптическим гипергеометрическим рядам.
Возможности обобщения на Римановы поверхности произвольного
рода обсуждаются в последнем параграфе.

\section{Тета-гипергеометрические ряды}

\subsection{Мультипликативная система обозначений}

Опишем мультипликативную систему обозначений для тета-функций,
которая будет использоваться в дальнейшем. Возьмем  комплексный параметр
$p$, лежащий внутри единичного круга, $|p|<1$, и $q\in\C$, $|q|\leq 1$.
Введем также модулярные параметры $\tau,\, \mbox{Im} (\tau)>0,$ и
$\sigma,\, \mbox{Im} (\sigma ) \geq 0,$ связанные с $p$ и $q$,
\begin{equation}\label{mod-par}
p=e^{2\pi i\tau }, \qquad q=e^{2\pi i\sigma}.
\end{equation}
Определим $p$-сдвинутые факториалы
$$
(a;p)_\infty=\prod_{n=0}^\infty(1-ap^n),\qquad
(a;p)_s=\frac{(a;p)_\infty}{(ap^s;p)_\infty}, \quad s\in\C.
$$
Для положительных целых $n$ имеем
$$ (a;p)_n=(1-a)(1-ap)\cdots(1-ap^{n-1}), \qquad
(a;p)_{-n}=\frac{1}{(ap^{-n};p)_n}. $$
Удобно использовать также компактные обозначения
$$
(a_1, \dots, a_k;p)_\infty \equiv (a_1;p)_\infty\cdots(a_k;p)_\infty.
$$

Определим укороченную тета-функцию Якоби
\begin{equation}\label{theta}
\theta(z;p)=(z,pz^{-1};p)_\infty.
\end{equation}
Она обладает следующими  свойствами симметрии
\begin{equation}\label{fun-rel}
\theta(pz;p)=\theta(z^{-1};p)=-z^{-1}\theta(z;p).
\end{equation}
Как следствие, выполняется также соотношение $\theta(p^{-1}z;p)=-p^{-1}z\theta(z;p)$.
Очевидно, что $\theta(z;p)=0$ при $z=p^{-M},\, M\in\mathbb{Z}$,
и $\theta(z;0)=1-z$.

Стандартная $\theta_1$-функция Якоби \cite{whi-wat:course}
выражается через $\theta(z;p)$ так
\begin{eqnarray} \nonumber
\theta_1(u;\sigma,\tau) \equiv \theta_1(\sigma u|\tau)&=&
-i\sum_{n=-\infty}^\infty (-1)^np^{(2n+1)^2/8}q^{(n+1/2)u}
\\ \nonumber
&=& 2\sum_{n=0}^\infty(-1)^n e^{\pi i\tau(n+1/2)^2}\sin\pi(2n+1)\sigma u
\\ \nonumber
&=& 2p^{1/8}\sin\pi \sigma u\:(p,pe^{2\pi i\sigma u},
pe^{-2\pi i\sigma u};p)_\infty  \\
&=& p^{1/8} iq^{-u/2}\: (p;p)_\infty\: \theta(q^{u};p), \quad
u\in\mathbb{C}.
\label{theta1}\end{eqnarray}
Мы искусственно ввели второй модулярный параметр $\sigma$
в определение $\theta_1$-функции. Дело в том, что часто будет возникать
произведение тета-функций с эквидистантными значениями аргументов и
удобно с самого начала произвести соответствующее растяжение аргумента
и придавать $u$ целые значения. Тета-функции Якоби с другими индексами
$\theta_{2,3,4}(u)$ могут быть получены из $\theta_1(u)$ простым сдвигом
переменной $u$ \cite{whi-wat:course}, то есть общая структура этих
функций по существу не сильно отличается от $\theta_1(u)$.

В дальнейшем удобно так же заменять $\theta_1$-символ
обозначениями эллиптических чисел $[u]$, отличающимися от
использовавшихся в предыдущей главе простой перенормировкой:
\be
[u] \equiv [u;\sigma,\tau]= \theta_1(u;\sigma,\tau), \qquad
[u_0,\ldots,u_k]= \prod_{m=0}^k[u_m].
\nonumber\ee
Зависимость от $\sigma$ и $\tau$ будет явно указываться только
при необходимости. Функция $[u]$ целая, нечетная,
$[-u]=-[u]$, и двояко-квазипериодическая
\be
[u+\sigma^{-1}]     = -[u],   \qquad
[u+\tau\sigma^{-1}] = -e^{-\pi i\tau-2\pi i\sigma u}[u].
\label{quasi}\ee
Известно, что тета-функция $[u]$ однозначно восстанавливается
(с точностью до постоянного множителя) из свойств
(\ref{quasi}) и условия целости функции.

Модулярные преобразования, образующие $PSL(2,\mathbb{Z})$
группу, меняют параметры $\sigma$ и $\tau$ дробно линейным образом
\begin{equation}\label{sl2}
\tau\to \frac{a\tau+b}{c\tau+d},\qquad \sigma\to \frac{\sigma}{c\tau+d},
\end{equation}
где $a,b,c,d \in\mathbb{Z}$ и $ad-bc=1$. Эта группа генерируется двумя
простыми преобразованиями: $\tau\to\tau+1$,
$\sigma\to\sigma,$ и $\tau\to -\tau^{-1},$ $\sigma\to\sigma\tau^{-1}$.
В этих двух случаях имеем
\begin{eqnarray}\label{modular-tr1}
&& [u;\sigma,\tau+1] = e^{\pi i/4} [u;\sigma,\tau],\qquad
\\  \label{modular-tr2}
&& [u;\sigma/{\tau},-1/\tau]
= -i(-i\tau)^{1/2}e^{\pi i\sigma^2u^2/\tau} [u;\sigma,\tau],
\end{eqnarray}
где знак квадратного корня в выражении $(-i\tau)^{1/2}$ фиксируется
из условия положительности его реальной части.

Обозначим эллиптические символы Похгаммера в аддитивной форме
$$
[u_1,\ldots,u_k]_{\pm n}=\prod_{m=1}^k[u_m]_{\pm n}, \qquad
[u]_n=\prod_{k=0}^{n-1}[u+k],\qquad [u]_{-n}=\frac{1}{[u-n]_n},
\; n\in\N,
$$
а в мультипликативной форме
\begin{eqnarray} &&
\theta(t;p;q)_n=\prod_{m=0}^{n-1}\theta(tq^m;p),  \qquad
\theta(t;p;q)_{-n}=\prod_{m=1}^{n}\frac{1}{\theta(tq^{-m};p)}
\nonumber \\ &&
\theta(t_0,\ldots,t_k;p;q)_n=\prod_{m=0}^k\theta(t_m;p;q)_n.
\label{ell-shift}\end{eqnarray}

\subsection{Односторонние $_rE_s$ и двусторонние $_rG_s$ ряды}

В этом параграфе вводятся формальные степенные ряды $_rE_s$ и $_rG_s$,
построенные из тета-функций Якоби. Они обобщают обычные и базисные
односторонние гипергеометрические ряды $_rF_s$ и $_r\varphi_s$ и их
двусторонние аналоги $_rH_s$ и $_r\psi_s$. Определение, приведенное
ниже, следует по духу качественному, но конструктивному определению
обычных и базисных гипергеометрических рядов, уходящему идейно еще к
Похгаммеру и Хорну \cite{aar:special,gas-rah:basic,ggr:general}.
В случае тета-гипергеометрических рядов мы используем
в качестве ключевого свойства квазипериодичность тета-функций
Якоби (\ref{quasi}).

\begin{definition}
Формальные ряды $\sum_{n\in\mathbb{N}} c_n$ и $\sum_{n\in\mathbb{Z}} c_n$
называются {\em тета-гипер\-гео\-ме\-три\-чес\-ки\-ми рядами эллиптического типа}
если отношения $h(n)=c_{n+1}/{c_n}$ и $1/h(n)$ являются мероморфными
двояко-квазипериодическими функциями $n$, рассматриваемой как комплексная
переменная. Точнее, функция $h(x)$, $x\in \mathbb{C}$, должна обладать
свойствами:
\begin{equation}\label{h(n)-period}
h(x+\sigma^{-1})=ah(x),\qquad
h(x+\tau\sigma^{-1})=be^{2\pi i\sigma \gamma x}h(x),
\end{equation}
где $\sigma^{-1}$ и $\tau\sigma^{-1}$ квазипериоды тета-функции
$[u]$ (\ref{quasi}) и $a, b, \gamma$ некоторые комплексные числа.
\end{definition}

\begin{theorem}
Пусть мероморфная функция $h(n)$ удовлетворяет условиям (\ref{h(n)-period}).
Тогда она должна иметь вид
\begin{equation}
h(n)=\frac{[n+u_1,\ldots, n+u_r]}{[n+v_1,\ldots, n+v_s]}\, q^{\beta n} y,
\label{h(n)-quasi}\end{equation}
где $r, s$ произвольные неотрицательные целые числа, а
$u_1,\ldots, u_r, v_1,\ldots,v_s, \beta, y$ произвольные комплексные
параметры, ограниченные условием несингулярности $h(n)$ и
связанные с параметрами квазипериодичности $a, b, \gamma$ следующим образом:
\begin{eqnarray}\nonumber
&& a = (-1)^{r-s} e^{2\pi i\beta}, \qquad \gamma= s-r, \\
&& b=(-1)^{r-s}e^{\pi i\tau (s-r+2\beta)}
e^{2\pi i\sigma(\sum_{m=1}^sv_m-\sum_{m=1}^ru_m)}.
\label{multipliers}\end{eqnarray}
\end{theorem}
{\bf Доказательство.}
Рассмотрим решетку на комплексной плоскости $n$, образованную периодами
$\sigma^{-1}$ and $\tau\sigma^{-1}$. Поскольку множители квазипериодичности в
(\ref{h(n)-period}) являются целыми функциями переменной $n$ без нулей,
мероморфная
функция $h(n)$ имеет то же самое конечное число нулей и полюсов в каждом
параллелограмме периодов. Обозначим через $-u_1,\ldots, -u_r$ положение
нулей $h(n)$ и $-v_1,\ldots, -v_s$ положение полюсов в одном из таких
параллелограммов. Для простоты будем подразумевать, что эти нули и полюсы
простые (для создания кратных полюсов и нулей достаточно положить
некоторые $u_m$ или $v_m$ равными друг другу). Конечность числа нулей $h(n)$
следует из меромофности $1/h(n)$.

Представим отношение соседних членов тета-гипергеометрического ряда
в виде $c_{n+1}/c_n$ $=h(n)g(n)$, где $h(n)$ имеет форму
(\ref{h(n)-quasi}) с некоторым свободным параметром $\beta$.
Поскольку нули и полюсы $c_{n+1}/c_n$ уже содержатся в
$h(n)$, функция $g(n)$ должна быть целой функцией без нулей,
удовлетворяющей ограничениям $g(n+\sigma^{-1})=a'g(n),$
$g(n+\tau\sigma^{-1})=b'e^{2\pi i\sigma \gamma' n}g(n)$ для некоторых
комплексных чисел $a', b', \gamma'$. Однако, единственной функцией,
удовлетворяющей этим требованиям, является экспоненциальная функция
$q^{\beta' n}$, где $\beta'$ обозначает произвольный параметр.
Поскольку такой множитель уже присутствует в (\ref{h(n)-quasi}),
то мы можем положить $g(n)=1$ и это доказывает, что наиболее общая
форма $c_{n+1}/c_n$ имеет вид (\ref{h(n)-quasi}) (заметим, что $y$
представляет собой произвольный коэффициент пропорциональности). Прямая
подстановка $h(n)$ в (\ref{quasi}) определяет связь между параметрами
$a, b, \gamma$ и $u_m,\, m=1,\ldots,r,$ $v_k,\, k=1,\ldots,s,$
$\beta$, как это указано в (\ref{multipliers}).
\hfill{Q.E.D.}

\smallskip
Решая рекуррентное соотношение первого порядка для
$c_n$ с нормировкой $c_0=1$, мы получаем следующий явный вид
тета-гипергеометрического ряда в ``аддитивной" форме
\begin{equation}
\sum_{n\in\mathbb{N}\, \mbox{или}\, \mathbb{Z}}
\frac{[u_1,\ldots,u_{r}]_n}
{[v_1,\ldots,v_s]_n}\,q^{\beta n(n-1)/2} y^n.
\label{_rE_s-1}\end{equation}

Для того, чтобы упростить предельный переход $\mbox{Im}(\tau)\to+\infty$
(или $p\to 0$) в ряде (\ref{_rE_s-1}) (тригонометрическое вырождение),
перенормируем $y$ и введем в качестве аргумента ряда переменную $z$:
$$
y \equiv (ip^{1/8}(p;p)_\infty)^{s-r}q^{(u_1+\ldots+u_r-v_1-\ldots-v_s)/2}z.
$$
Заменим параметр $\beta$ другим параметром $\alpha$:
$\beta\equiv \alpha + (r-s)/2$. Тогда, мы можем переписать функцию
(\ref{h(n)-quasi}) в ``мультипликативном" виде, основанном на
укороченной тета-функции $\theta(t;p)$:
\begin{equation}
h(n)=\frac{\theta(t_1q^n,\ldots,t_{r}q^n;p)}
{\theta(w_1q^n,\ldots, w_{s}q^n;p)}\, q^{\alpha n}z,
\label{h(n)-2}\end{equation}
где $t_m=q^{u_m},\, m=1,\ldots, r,$  $w_k=q^{v_k},\, k=1,\ldots, s,$ и
$\theta(t_1,\ldots,t_k;p)=\prod_{m=1}^k\theta(t_m;p).$

Теперь мы в состоянии привести явный вид одностороннего
тета-гипергеометрического ряда $_rE_s$. При этом мы принимаем соглашения,
аналогичные общепринятым для  $_r\varphi_s$ и $_rF_s$ рядов.
Точнее, в выражении (\ref{h(n)-2}) мы заменяем $s$ на $s+1$ и полагаем
$u_r\equiv u_0$ и $v_{s+1}=1$. Это не ограничивает общности ситуации,
так как подобное условие может быть устранено фиксацией одного из
параметров в числителе $u_m$ равным единице. После этого мы полагаем
\be
_rE_s\left({t_0,\ldots, t_{r-1}\atop w_1,\ldots,w_s};q,p;\alpha,
z\right) %} && \\ && \makebox[4em]{}
= \sum_{n=0}^\infty \frac{\theta(t_0,t_1,\ldots,t_{r-1};p;q)_n}
{\theta(q,w_1,\ldots,w_s;p;q)_n}\, q^{\alpha n(n-1)/2} z^n.
\label{_rE_s-2}\ee
Порядок $q$ и $p$, фигурирующих в обозначениях $_rE_s$ ряда, важен в виду
отсутствия перестановочной симметрии между этими параметрами.
Ряд \re{_rE_s-2} носит формальный характер, так как мы не анализируем
условий его сходимости (в приложениях будут использоваться только обрывающиеся
ряды).

Важно отметить, что тета-гипергеометрические ряды не допускают
конфлюэнтных предельных переходов. Действительно, из-за квазипериодичности
тета-функций пределы параметров
$t_m, w_m\to 0$ или $t_m, w_m \to\infty$ не определены и поэтому нельзя
таким образом перейти от $_rE_s$ рядов к аналогичным рядам с
меньшими значениями индексов $r$ и $s$ (что было возможно в
случае $_rF_s$ и $_r\varphi_s$ рядов.

Для двусторонних тета-гипергеометрических рядов мы введем другие обозначения:
\be
 _rG_s\left({t_1,\ldots, t_{r}\atop w_1,\ldots,w_{s}};q,p;\alpha, z
\right) % } && \\ && \makebox[4em]{}
=\sum_{n=-\infty}^\infty\frac{\theta(t_1,\ldots,t_{r};p;q)_n}
{\theta(w_1,\ldots,w_{s};p;q)_n}\, q^{\alpha n(n-1)/2} z^n.
\label{_rG_s-2}\ee
Это выражение было получено с помощью выражения (\ref{h(n)-2}), использованного
без каких-либо изменений. Эллиптические символы Похгаммера при отрицательных
значениях индекса имеют вид
$$
\theta(t;p;q)_{-n}=\frac{1}{\theta(tq^{-n};p;q)_n},\quad n\in\mathbb{N}.
$$
Благодаря свойству $\theta(q;p;q)_{-n}=0$ (или $[1]_{-n}=0$) при $n>0$,
выбор $t_{s+1}=q$ (или $v_{s+1}=1$) в $_rG_{s+1}$ ряде приводит к его обрыву
с одной стороны. После переобозначения $t_r\equiv t_0$ (или
$u_{r}\equiv u_0$), можно прийти таким образом к общему $_rE_s$ ряду.
Поскольку двусторонние ряды более общие чем односторонние, достаточно изучить
основные свойства тета-гипергеометрических рядов в этом случае без спецификации
к односторонним рядам.

Рассмотрим предел $\mbox{Im}(\tau)\to+\infty$, при котором $p\to 0$.
Прямым почленным предельным переходом получаем
\be
\lim_{p\to 0} {_rE_s} ={_r\varphi_s}\left({t_0,t_1,\ldots, t_{r-1}
\atop w_1,\ldots,w_s }; q;\alpha, z\right) %} && \\ && \makebox[3em]{}
=\sum_{n=0}^\infty \frac{(t_0,t_1,\ldots,t_{r-1};q)_n}
{(q,w_1,\ldots,w_s;q)_n}\, q^{\alpha n(n-1)/2} z^n.
\label{Phi-2}\ee
Этот базисный гипергеометрический ряд отличается от стандартного наличием
дополнительного параметра $\alpha$. Определение $_r\varphi_s$ ряда,
предложенное в монографии \cite{sla:generalized}, использует $\alpha=0$.
Определение, предложенное в книге \cite{gas-rah:basic},
\begin{equation}
{_r\varphi_s}=\sum_{n=0}^\infty \frac{(t_0,t_1,\ldots,
t_{r-1};q)_n}{(q,w_1,\ldots,w_s;q)_n}
\left((-1)^nq^{n(n-1)/2}\right)^{s+1-r}\, z^n,
\label{Phi-3}\end{equation}
согласуется с (\ref{Phi-2}) при $\alpha=s+1-r$ после замены $z$ на
$(-1)^{s+1-r}z$. На самом деле, случаи $\alpha=0$
и $\alpha=s+1-r$ связаны друг с другом инверсией $q\to q^{-1}$ с последующим
переопределением параметров $t_m, w_m, z$. Структурно, базисные гипергеометрические
ряды могут быть определены как суммы  $\sum_nc_n$, для которых
отношение $c_{n+1}/c_n$
равно общей рациональной функции $q^n$ и (\ref{Phi-2}) удовлетворяет
этому требованию только для целых $\alpha$. Для того, чтобы иметь в пределе $p\to 0$
стандартные $q$-гипергеометрические ряды, мы полагаем $\alpha=0$.
Вообще говоря, неясно насколько существенно данное ограничение,  это
прояснится только после обнаружения каких-либо
реальных приложений рядов $_rE_s$.

В случае двусторонних рядов мы также полагаем $\alpha=0$, так что в $p\to 0$
пределе $_rG_s$ ряд сводится к общему $_r\psi_s$ ряду:
\be
{_r\psi_s}\left({t_1,\ldots, t_{r}
\atop w_1,\ldots,w_s }; q;z\right) % } && \\ && \makebox[3em]{}
=\sum_{n=-\infty}^\infty \frac{(t_1,\ldots,t_{r};q)_n}
{(w_1,\ldots,w_s;q)_n}\, z^n.
\label{Psi}\ee

\begin{definition}
Ряды $_{r+1}E_r$ и $_{r}G_{r}$ (с $\alpha=0$) называются сбалансированными,
если их параметры удовлетворяют ограничениям (в аддитивной форме)
\begin{equation}
u_0+\ldots+u_{r}=1+v_1+\ldots+v_r +k/\sigma, \quad k\in\Z,
\label{balance1}\end{equation}
и
\begin{equation}
u_1+\ldots+u_{r}=v_1+\ldots+v_r +k/\sigma, \quad k\in\Z,
\label{balance-g}\end{equation}
соответственно.
В мультипликативной форме эти ограничения принимают следующий вид:
$\prod_{m=0}^r t_m=q\prod_{k=1}^rw_k$ и $\prod_{m=1}^r t_m=$
$\prod_{k=1}^rw_k$, соответственно.
\end{definition}

При $r=s+1$ число тета-функций в числителе и знаменателе $h(n)$ совпадает.
Поэтому все параметры $u_m$ и $v_m$ в $_{r+1}E_r$ рядах определены
с точностью до добавления  $k/\sigma, \quad k\in\Z$. Именно эта свобода и
приводит к добавлению аналогичного члена в условиях \re{balance1} и
\re{balance-g}. В дальнейшем мы не будем указывать эту неоднозначность,
отсутствующую в мультипликативных обозначениях, но необходимо помнить, что
выражения типа $u_m/2$ содержат уже неоднозначность в добавлении
$\Z/2\sigma$.

\begin{remark}
В пределе $p\to 0$ ряд $_{r+1}E_r$ переходит в $_{r+1}\varphi_r$
при условии, что параметры $u_m$ (или $t_m$), $m=0,\ldots,r,$ и
$v_k$ (или $w_k$), $k=1,\ldots,r,$  {\em остаются фиксированными}.
При этом наше условие балансировки не совпадает с приведенным в
\cite{gas-rah:basic}, где $_{r+1}\varphi_r$ называется сбалансированным при
$q\prod_{m=0}^r t_m=\prod_{k=1}^rw_k$ (одновременно полагается $z=q$,
но мы опускаем это условие). Противоречие в этих определениях
снимается после наложения условия совершенной уравновешенности, описанного ниже.
\end{remark}

\section{Эллиптические гипергеометрические ряды одной переменной}

Концепция эллиптического гипергеометрического ряда, предложенная ниже,
играет важную роль для всей теории рядов гипергеометрического типа,
так как она дает объяснение происхождения особенностей некоторых
обычных и базисных гипергеометрических рядов.

\begin{definition}
Формальные ряды $\sum_{n\in\mathbb{N}} c_n$ и $\sum_{n\in\mathbb{Z}} c_n$
называются  эллиптическими гипергеометрическими рядами
если отношения $h(n)=c_{n+1}/{c_n}$ являются дискретной последовательностью
значений при $x\in \N$ или $x\in\Z$ некоторой эллиптической (т.е. мероморфной
двояко-периодической) функции $h(x),\, x\in\mathbb{C}$.
\end{definition}

\begin{theorem}
Пусть $\sigma^{-1}$ и $\tau\sigma^{-1}$ обозначают два периода эллиптической
функции $h(x)$, т.е. $h(x+\sigma^{-1})=h(x)$ и $h(x+\tau\sigma^{-1})=h(x)$.
Пусть порядок $h(x)$ (т.е. число нулей или полюсов в параллелограмме периодов
с учетом их кратностей) равен $r+1, r>0$. Тогда односторонний (или двусторонний)
эллиптический гипергеометрический ряд совпадает со сбалансированным
тета-гипергеометрическим рядом $_{r+1}E_r$ (или $_{r+1}G_{r+1}$).
\end{theorem}
{\bf Доказательство.}
Хорошо известно, что любую эллиптическую функцию $h(x),\, x\in\mathbb{C}$,
порядка $r+1$ с периодами $\sigma^{-1}$ и $\tau\sigma^{-1}$
можно представить в виде отношения произведений $\theta_1$-функций
\cite{whi-wat:course}:
\begin{equation}
h(x)=z\prod_{m=0}^r\frac{[x+\alpha_m;\sigma,\tau]}{[x+\beta_m;\sigma,\tau]},
\label{factorization}\end{equation}
где координаты нулей $\alpha_0,\ldots,\alpha_r$ и полюсов
$\beta_0,\ldots,\beta_r$ удовлетворяют ограничению:
\begin{equation}\label{balance2}
\sum_{m=0}^r\alpha_m=\sum_{m=0}^r\beta_m.
\end{equation}
Это представление позволяет легко идентифицировать односторонний эллиптический
гипергеометрический ряд со сбалансированным $_{r+1}E_r$ рядом. Для этого достаточно
сдвинуть $x\to x-\beta_0+1$, положить $x\in \mathbb{N}$, обозначить
$u_m=\alpha_m-\beta_0+1,\, v_m=\beta_m-\beta_0+1$ и решить рекуррентное соотношение
$c_{n+1}=h(n)c_n$. После этого, условие (\ref{balance2}) становится условием
балансировки для $_{r+1}E_r$ рядов. Очевидно, что аналогичная ситуация имеет
место и для двустороннего ряда $_{r+1}G_{r+1}$, для которого $\alpha_m$ и
$\beta_m$ просто совпадают с $u_m$ и $v_m$, соответственно.
\hfill{Q.E.D.}

Отметим, что благодаря условию (\ref{balance2}),
функция (\ref{factorization}) может быть переписана в виде отношения
укороченных $\theta(t;p)$-функций:
$$
h(x)=z\prod_{m=0}^r\frac{\theta(t_mq^x;p)}{\theta(w_mq^x;p)},
$$
где $t_m=q^{\alpha_m}$ и $w_m=q^{\beta_m}$.

\begin{definition}
Тета-гипергеометрические ряды называются модулярными
гипергеометрическими рядами если они инвариантны относительно
действия $PSL(2,\mathbb{Z})$ группы (\ref{sl2}).
\end{definition}

Рассмотрим какого типа ограничения необходимо наложить на параметры $_rE_s$
и $_rG_s$ рядов, для того чтобы получить модулярный гипергеометрический
ряд. Очевидно, что достаточно установить модулярность
$h(n)=c_{n+1}/c_n$. Из явного вида этой функции (\ref{h(n)-quasi})
и преобразований (\ref{modular-tr1}) и (\ref{modular-tr2})
легко увидеть, что в случае одностороннего ряда необходимо иметь
$$
\sum_{m=0}^{r-1}(x+u_m)^2=(x+1)^2+\sum_{m=0}^s(x+v_m)^2,
$$
что возможно только если а) $s=r-1$, б) параметры удовлетворяют
условию балансировки (\ref{balance1}), и в) имеет место равенство:
\begin{equation}
u_0^2+\ldots+u_{r-1}^2=1+v_1^2+\ldots+v_{r-1}^2.
\label{modular-E}\end{equation}
При этих условиях $_{r}E_{r-1}$ ряд модулярно инвариантен.
Отметим, что модулярность подразумевает эллиптичность тета-гипергеометрического
ряда, но не наоборот. Однако, более сильное требование эллиптичности,
сформулированное ниже, автоматически приводит к модулярной инвариантности.
Модулярные гипергеометрические ряды представляют частные примеры
(мероморфных) модулярных форм Якоби в смысле Эйхлера и Загира \cite{eic-zag:theory}.

В случае двусторонних рядов, необходимо иметь $r=s$, условие балансировки
(\ref{balance-g}), и ограничение
\begin{equation}
u_1^2+\ldots+u_{r}^2=v_1^2+\ldots+v_r^2
\label{modular-G}\end{equation}
для модулярности $_{r}G_{r}$ ряда.

\begin{definition}
Тета-гипергеометрический ряд $_{r+1}E_r$ называется  вполне уравновешенным
если его параметры удовлетворяют ограничениям
\begin{equation}\label{well-poised-1}
u_0+1=u_1+v_1=\ldots=u_{r}+v_r
\end{equation}
в аддитивной форме или
\begin{equation}\label{well-poised-2}
qt_0=t_1w_1=\ldots=t_{r}w_r
\end{equation}
в мультипликативной форме. Аналогично, ряды $_{r}G_{r}$ называются вполне
уравновешенными если $u_1+v_1=\ldots=u_r+v_r$ или $t_1w_1=\ldots=t_rw_r$.
\end{definition}

Это определение вполне уравновешенных рядов полностью соответствует
аналогичному понятию для обычных и $q$--гипергеометрических рядов
\cite{gas-rah:basic}. Заметим, что оно не подразумевает выполнения
условия балансировки.

\begin{definition}
Ряд $_{r+1}E_r$ называется совершенно уравновешенным
если, в дополнение к условиям (\ref{well-poised-1})
или (\ref{well-poised-2}), выполняются ограничения
\begin{eqnarray}\nonumber
&& u_{r-3}=\frac{1}{2}u_0+1,\quad u_{r-2}=\frac{1}{2}u_0+1-\frac{1}{2\sigma}, \\
&& u_{r-1}=\frac{1}{2}u_0+1-\frac{\tau}{2\sigma},
\quad u_{r} = \frac{1}{2}u_0+1+\frac{1+\tau}{2\sigma},
\label{very-well-poised-1}\end{eqnarray}
или, в мультипликативной форме,
\be
t_{r-3}=t_0^{1/2}q,\quad t_{r-2}=-t_0^{1/2}q, \quad
 t_{r-1}=t_0^{1/2}qp^{-1/2}, \quad t_{r} =- t_0^{1/2}qp^{1/2}.
\label{very-well-poised-2}\ee
\end{definition}

Совершенно уравновешенные ряды имеют специфический вид, получающийся после
упрощения явного вида его коэффициентов. Пользуясь формулами
$$
\theta(zp^{-1/2};p)=-zp^{-1/2}\theta(zp^{1/2};p)
$$
и
$$
\theta(z,-z,zp^{1/2},-zp^{1/2};p)=\theta(z^2;p),
$$
находим
$$
\frac{\theta(t_{r-3},\ldots,t_{r};p;q)_n}
{\theta(qt_0/t_{r-3},\ldots,qt_0/t_{r};p;q)_n}=
\frac{\theta(t_0q^{2n};p)}{\theta(t_0;p)}\, (-q)^n.
$$
В результате,
\begin{eqnarray}\nonumber
\lefteqn{
_{r+1}E_r\left({t_0,t_1,\ldots, t_{r-4},qt_0^{1/2},-qt_0^{1/2},
qp^{-1/2}t_0^{1/2},-qp^{1/2}t_0^{1/2} \atop
qt_0/t_1,\ldots,qt_0/t_{r-4},t_0^{1/2},-t_0^{1/2},
p^{1/2}t_0^{1/2},-p^{-1/2}t_0^{1/2} };q,p;z\right) } &&
\\ && \makebox[5em]{}
= \sum_{n=0}^\infty \frac{\theta(t_0q^{2n};p)}{\theta(t_0;p)}
\prod_{m=0}^{r-4}\frac{\theta(t_m;p;q)_n}{\theta(qt_0/t_m;p;q)_n}\, (-qz)^n.
\label{vwp-1}\end{eqnarray}

Для удобства, в работе \cite{spi:bailey} было введено отдельное обозначение для
совершенно уравновешенного тета-гипергеометрического ряда, поскольку
он содержит значительно меньше свободных параметров чем общий $_{r+1}E_r$ ряд.
Заменим $z$ на $-z$ и определим
\be
_{r+1}V_r(t_0;t_1,\ldots,t_{r-4};q,p;z)
\equiv \sum_{n=0}^\infty \frac{\theta(t_0q^{2n};p)}{\theta(t_0;p)}
\prod_{m=0}^{r-4}\frac{\theta(t_m;p;q)_n}{\theta(qt_0t_m^{-1};p;q)_n}\,
(qz)^n.  \label{vwp-2}\ee
В целях согласования с монографией \cite{gas-rah2}, для обозначения
этого ряда в терминах эллиптических чисел используем другой символ:
\begin{eqnarray} \nonumber
\lefteqn{
_{r+1}v_r(u_0;u_1,\ldots,u_{r-4};\sigma,\tau;z) }&&
\\ && \makebox[3em]{}
\equiv \sum_{n=0}^\infty \frac{[u_0+2n]}{[u_0]}\prod_{m=0}^{r-4}
\frac{[u_m]_n}{[u_0+1-u_m]_n}\, z^n
q^{n(\sum_{m=1}^{r-4}u_m-(r-7)/2-(r-5)u_0/2)}.
\label{vwp-3}\end{eqnarray}

При $p=0$, тета-гипергеометрический ряд (\ref{vwp-2}) сводится к совершенно
уравновешенному базисному гипергеометрическому ряду
\be
_{r-1}W_{r-2}(t_0;t_1,\ldots,t_{r-4};q;qz)
= \sum_{n=0}^\infty \frac{1-t_0q^{2n}}{1-t_0}
\prod_{m=0}^{r-4}\frac{(t_m;q)_n}{(qt_0t_m^{-1};q)_n}\, (qz)^n,
\label{vwp-phi}\ee
который отличается от $_{r-1}W_{r-2}$ ряда, описанного в \cite{gas-rah:basic},
заменой $z$ на $qz$. Накладывая условие балансировки
$\sum_{m=0}^r u_m=1+\sum_{m=1}^r(u_0+1-u_m)$
на ряд (\ref{vwp-3}), получаем
$\sum_{m=1}^{r-4} u_m= \frac{r-7}{2}+\frac{r-5}{2}u_0+\frac{\Z}{2\sigma}.$
В мультипликативной форме это условие принимает вид
$$
\prod_{m=1}^{r-4}t_m=\pm q^{(r-7)/2}t_0^{(r-5)/2}.
$$
Оказывается, что при выборе положительного знака мы получаем в точности
условие балансировки для совершенно
уравновешенных базисных гипергеометрических рядов в стандартном
виде \cite{gas-rah:basic}. Таким образом, для совершенно уравновешенных рядов
можно согласовать условие балансировки приведенные выше с определением данным в
\cite{gas-rah:basic}. Это происходит благодаря тому, что ограничения
(\ref{very-well-poised-1}) по отдельности не имеют предела
$\mbox{Im} (\tau)\to+\infty$ и только их комбинация в тета-функциях
приводит к осмысленному результату. Заметим, что для сбалансированных
рядов дополнительный множитель, стоящий справа от $z^n$ в
(\ref{vwp-3}), сводится к $(\pm 1)^n$.
Для нечетных $r$ выбор положительного знака в условии балансировки
однозначен --- известно, что только в этом случае возникают нетривиальные тождества
для рядов. В случае четного $r$ неоднозначность остается (необходимо
указывать ветви корней $q^{1/2}$ и $t_0^{1/2}$),
но для четных $r$ неизвестны какие-либо формулы суммирования или
преобразования для совершенно уравновешенных рядов гипергеометрического
типа. Для нечетных  $r$ необходимый знак фиксируется с помощью
дополнительной симметрии (см. ниже).

Суммируя приведенный анализ можно прийти к заключению,
что условие эллиптичности функции $h(n)=c_{n+1}/c_n$
в тета-гипергеометрических рядах придает глубокий смысл (неестественному)
условию балансировки для обычных и $q$-ги\-пер\-гео\-ме\-три\-чес\-ких рядов.

Аналогично, двусторонние ряды $_{r}G_{r}$ называются совершенно
уравновешенными если удовлетворены условия (\ref{very-well-poised-1}) или
(\ref{very-well-poised-2}) с произвольным параметром $u_0$ или $t_0=q^{u_0}$.
Чтобы сделать эти ряды похожими на односторонние, положим
$t_kw_k=qt_0$ и заменим $z$ на $-z$. Тогда,
\be
_{r}G_{r}(t_0;t_1,\ldots,t_{r-4};q,p;z)
= \sum_{n=-\infty}^\infty \frac{\theta(t_0q^{2n};p)}{\theta(t_0;p)}
\prod_{m=1}^{r-4}\frac{\theta(t_m;p;q)_n}{\theta(qt_0t_m^{-1};p;q)_n}\,
(qz)^n  \label{vwp-g2}\ee
или в аддитивной форме
\begin{eqnarray} \nonumber
\lefteqn{
_{r}g_{r}(u_0;u_1,\ldots,u_{r-4};\sigma,\tau ;z) }&&
\\ && \makebox[3em]{}
= \sum_{n=-\infty}^\infty \frac{[u_0+2n]}{[u_0]}\prod_{m=1}^{r-4}
\frac{[u_m]_n}{[u_0+1-u_m]_n}\, z^n
q^{n(\sum_{m=1}^{r-4}u_m-(r-8)/2-(r-4)u_0/2)}.
\label{vwp-g3}\end{eqnarray}
Условие балансировки для этого ряда принимает вид
$$
\sum_{m=1}^{r-4}u_m=\frac{r-8}{2} +\frac{r-4}{2}u_0+
\frac{\Z}{2\sigma}, \quad
\mbox{или} \quad \prod_{m=1}^{r-4}t_m=\pm q^{(r-8)/2}t_0^{(r-4)/2}.
$$
При условии $u_{r-3}=u_0$ (или $t_{r-3}=t_0$)
$_{r+1}G_{r+1}$ ряд переходит в $_{r+1}E_r$ ряд.

\begin{remark}
В рамках описанной классификации, предложенной в работе \cite{spi:theta},
эллиптическое обобщение $_{r+1}W_r$ ряда Френкеля и Тураева \cite{fre-tur:elliptic}
совпадает с совершенно уравновешенным сбалансированным тета-гипергеометрическим
рядом $_{r+3}V_{r+2}$ с единичным аргументом $z=1$. Такие ряды впервые появились
при описании эллиптических решений уравнения Янга-Бакстера \cite{bax:exactly,
tf,abf:eight,djkmo:exactly1,djkmo:exactly2,djmo:fusion}.
В совершенно другом контексте -- при изучении автомодельных биортогональных рациональных
функций, эти ряды были выведены в работах \cite{spi-zhe:spectral,spi-zhe:classical,
spi-zhe:gevp}, как это уже было описано в предыдущей главе.
\end{remark}

\begin{definition}
Ряды $\sum_{n\in\mathbb{N}}c_n$ и $\sum_{n\in\mathbb{Z}}c_n$
называются полностью эллиптическими гипергеометрическими рядами если
$h(n)=c_{n+1}/c_n$ является эллиптической функцией  всех свободных параметров
входящих в него (за исключением параметра $z$, на который всегда можно умножить
$h(n)$), с равными периодами.
\end{definition}

\begin{theorem}
Наиболее общий (в смысле максимального числа независимых свободных
параметров среди $u_k$ и $v_k$)
полностью эллиптический гипергеометрический ряд совпадает с совершенно уравновешенным
сбалансированным тета-гипергеометрическим рядом $_{r}v_{r-1}$
(в одностороннем случае) и $_rg_r$ (в двустороннем случае) при $r>3$.
Полная эллиптичность гарантирует модулярную инвариантность.
\end{theorem}
{\bf Доказательство.}
Достаточно доказать это утверждение для двустороннего ряда, легко сводимого к
одностороннему. Эллиптичность по $n$ приводит к $h(n)$ вида
$$
h(n)=\frac{[n+u_1,\ldots, n+u_r]}{[n+v_1,\ldots, n+v_r]}\, z,
$$
с произвольными $u_1,\ldots,u_r$ и $v_1,\ldots,v_{r}$, удовлетворяющими
условию $u_1+\ldots+u_r=v_1+\ldots+v_{r}$. Из этого представления очевидно,
что имеется свобода в сдвиге параметров на произвольную константу:
$u_m\to u_m+u_0$, $v_m\to v_m+u_0$,
$m=1,\ldots, r$, без нарушения условия балансировки.

Определим теперь максимально возможное число независимых переменных
в полностью эллиптических гипергеометрических рядах.
Предположим, что параметры $u_1,\ldots,$ $u_{r-1}$ линейно независимы
и функция $h(n)$ симметрична и двояко периодична по ним.
Поскольку минимальный
порядок эллиптической функции равен двум, то $h(n)$ должна иметь как
минимум два нуля или полюса (с учетом их кратностей) по $u_1$.
Двукратные нули или полюсы редуцируют общее число свободных параметров
и поэтому мы предполагаем их отсутствие.
Поэтому, $u_r$ зависит линейно от $u_1$ ($u_2,\ldots,u_{r-1}$
по предположению не зависят от этой переменной):
$u_r=\alpha\sum_{k=1}^{r-1}u_k+\beta$,
где $\alpha, \beta$ некоторые неизвестные коэффициенты
(очевидно, что $\alpha$ должно быть целым).

При этом параметры $v_m,\, m=1,\dots,r,$ могут зависеть от $u_k$
двумя способами. Если они зависят только от суммы $\sum_{k=1}^{r-1}u_k$,
то тогда нельзя обеспечить инвариантность относительно сдвигов
$u_k\to u_k+\tau\sigma^{-1}$ для любого $k$.
Поэтому остается единственный выбор, обеспечивающий перестановочную
симметрию произведения $\prod_{m=1}^r[x+v_m]$ по $u_1,\ldots,u_{r-1}$,
$$
v_m=\gamma\sum_{k=1}^{r-1}u_k+\delta u_m+\rho,
\quad m=1,\ldots,r-1,
$$
и $v_r=\mu\sum_{k=1}^{r-1}u_k+\nu$,
где $\gamma, \delta, \rho, \mu, \nu$ обозначают некоторые числовые коэффициенты
($\gamma, \delta, \mu$ должны быть целыми). Другие симметричные комбинации
$u_m$ требуют перемножения более чем $r$  тета-функций. Подстановка этих выражений
в условие балансировки дает $1+\alpha=(r-1)\gamma+\delta+\mu$ и
$\beta=(r-1)\rho+\nu$, что гарантирует инвариантность $h(x)$ относительно
сдвигов $u_k\to u_k+\sigma^{-1}$ и сокращает частично факторы, появляющиеся
из сдвигов $u_k\to u_k+\tau\sigma^{-1}$. Несколько громоздкий, но технически
прямой анализ условия сокращения множителей вида $e^{-2\pi i\sigma u}$
приводит к уравнениям
$$
\delta^2=1,\qquad \alpha^2=(r-1)\gamma^2+2\gamma\delta+\mu^2,\qquad
\alpha\beta=(r-1)\gamma\rho+\rho\delta+\mu\nu.
$$
При этом необходимо учитывать, что уравнения на $\nu$ и $\beta$
содержат произвол в добавлении члена вида $\Z\sigma^{-1}$.
Ограничение, появляющееся из условия сокращения множителей вида
$e^{-\pi i\tau}$, оказывается несущественным и оно опускается.

Пусть $\delta=1$. Тогда уравнения на $\alpha,\gamma,
\mu$ фиксируют $\gamma=0$ или $\gamma(r-1)(r-2)+2\mu(r-1)=2.$ Поскольку
$\gamma$ и $\mu$ целые числа, второй случай исключается
(целые числа в левой части пропорциональны $r-1$, в то время как
 это невозможно обеспечить в правой части при $r>3$).
Выбор $\gamma=0$ приводит к $\alpha=\mu$ и, как следует из двух других уравнений,
с необходимостью $\rho=0$ и $\beta=\nu$. При этом $h(n)=1$,
т.е. мы приходим к тривиальной ситуации, которую отбрасываем.

Пусть $\delta=-1$. Решение взятых уравнений приводит к
$$
\alpha=\frac{\gamma r}{2}-1, \quad
\mu=1-\frac{\gamma (r-2)}{2}, \quad
\beta=\frac{\rho r}{2}, %+\frac{m\sigma^{-1}}{2-\gamma(r-1)},
\quad \nu=-\frac{\rho (r-2)}{2}, %+\frac{m\sigma^{-1}}{2-\gamma(r-1)},
$$
где $\gamma$ есть целое число (четное для нечетных $r$) и
$\rho$ --- произвольный параметр. Легко увидеть, что формальным сдвигом
$\rho\to \rho -\gamma\sum_{k=1}^{r-1}u_k$ параметр $\gamma$
устраняется из $h(n)$. Этот сдвиг разрешен при условии, что функция $h(n)$
эллиптична и по параметру $\rho$, но это условие зависит от четности $r$.
При четных $r$, оно выполняется автоматически для произвольных $r$
только при указанном выборе параметров.
Таким образом, в этом случае
\begin{equation}\label{vwp-even}
h(n)=\prod_{m=1}^{r}\frac{[n+u_m]}{[n+\rho-u_m]}\, z,\qquad
\sum_{k=1}^ru_k=\frac{\rho r}{2},
\end{equation}
что соответствует правильному выбору знака в условии балансировки
для вполне уравновешенных эллиптических
гипергеометрических рядов при четных $r$ упомянутому выше.

При нечетных $r$ возможна инвариантность при сдвигах $\rho$ только
на $2\tau\sigma^{-1}$. Обозначим теперь $\rho\equiv 2g_0$
и сдвинем $u_m\to u_m+g_0$, $m=1,\ldots,r-1$. После этого получаем
\begin{equation}\label{vwp-theorem}
h(n)=\prod_{m=1}^{r-1}\frac{[n+g_0+u_m]}{[n+g_0-u_m]}\,
\frac{[n+g_0-\sum_{k=1}^{r-1}u_k-l/2\sigma]}
{[n+g_0+\sum_{k=1}^{r-1}u_k+l/2\sigma]}\, z, \qquad l\in\Z,
\end{equation}
т.е. параметр $g_0$ играет ту же самую роль, что и $n$, и эллиптичность по
нему очевидна. Нетрудно узнать в (\ref{vwp-theorem})
наиболее общее выражение для отношения $c_{n+1}/c_n$ во вполне уравновешенных
и сбалансированных тета-гипергеометрических рядах и выбор $l$
соответствует неоднозначности в знаке условия балансировки.
Таким образом мы доказали,
что эллиптичность по свободным параметрам в $h(n)$ приводит к условиям
балансировки и вполне уравновешенности. Причем для четных $r$ удается
зафиксировать даже неоднозначность в знаке условия балансировки
благодаря дополнительной симметрии---эллиптичности по $\rho$ с теми же
периодами что и для других параметров.

Докажем, что полностью эллиптические ряды модулярно инвариантны. Для этого
достаточно проверить, что суммы квадратов параметров $u$, входящих
в эллиптические числа $[n+u]$ в числителе $h(n)$ (\ref{vwp-theorem})
и знаменателе, совпадают. Параметры числителя генерируют сумму
$$
\sum_{k=1}^{r-1}u_k^2 +\Bigl(-\sum_{k=1}^{r-1} u_k-l/2\sigma\Bigr)^2
$$
которая элементарным образом равна сумме для знаменателя:
$$
\sum_{k=1}^{r-1}(-u_k)^2+\Bigl(\sum_{k=1}^{r-1} u_k+l/2\sigma\Bigr)^2,
$$
т.е. модулярная инвариантность выполняется автоматически.
Отметим, что вполне уравновешенные ряды без условия балансировки
не обладают этой симметрией. Аналогичная ситуация имеет место и для
функции (\ref{vwp-even}).

Все рассмотрения, приведенные выше, проводились в расчете на двусторонние
ряды $_{r}G_r$, но переход к вполне уравновешенным и сбалансированным
$_rE_{r-1}$ рядам осуществляется простым ограничением на один из параметров.
\hfill{Q.E.D.}

\smallskip

Мы установили интересное происхождение условий балансировки и вполне
уравновешенности для обычных и $q$-гипергеометрических рядов, требующее
переосмысления этих концепций. Другая любопытная интерпретация
условия вполне уравновешенности описана в статье \cite{and:thread}.
Однако, происхождение условия совершенной уравновешенности остается не
до конца понятым. Возможно, что эллиптическая функция $h(n)$ с таким
ограничением обладает дополнительными арифметическими свойствами.

Бесконечные ряды $_{r}E_s$ и $_rG_s$ носят формальный характер.
Из-за квазипериодичности тета-функций, определение условий сходимости
этих рядов представляет собой серьезную проблему и она не будет рассматриваться
здесь. При определенном выборе параметров в сбалансированном ряде
$_{r+1}E_r$, его радиус сходимости равен 1. Строгий смысл этим рядам придается
наложением условий их обрыва. Тета-гипергеометрический ряд обрывается, если
для некоторого $m$,
\begin{equation}\label{e-trunc1}
u_m=-N-K\sigma^{-1}-M\tau\sigma^{-1}, \qquad N\in\mathbb{N},\quad
K, M\in\mathbb{Z},
\end{equation}
или в мультипликативной форме
\begin{equation}\label{e-trunc2}
t_m=q^{-N}p^{-M}, \qquad N\in\mathbb{N},\quad M\in\mathbb{Z}.
\end{equation}
Совершенно уравновешенные эллиптические гипергеометрические ряды
дважды периодичны по параметрам с периодами $\sigma^{-1}$ и $\tau\sigma^{-1}$.
Поэтому обрывающиеся ряды такого типа не зависят от целых чисел $K$ и $M$.

Самым общим тождеством в теории $q$-гипергеометрических рядов является
четырехчленное тождество Бэйли для необрывающихся  $_{10}\varphi_9$
совершенно уравновешенных сбалансированных рядов со специальным значением
аргумента \cite{gas-rah:basic}. В случае обрыва рядов остается только два
члена. В работе \cite{fre-tur:elliptic}, Френкель и Тураев доказали
эллиптическое обобщение тождества Бэйли для обрывающихся рядов.
В наших обозначениях, это тождество имеет вид
\begin{eqnarray}\nonumber
{_{12}V_{11}}(t_0;t_1,\dots,t_7;q,p;1) &=&
\frac{\theta(qt_0,qs_0/s_4,qs_0/s_5,qt_0/t_4t_5;p;q)_N}
{\theta(qs_0,qt_0/t_4,qt_0/t_5,qs_0/s_4s_5;p;q)_N}   \\
&& \times {_{12}V_{11}}(s_0;s_1,\dots,s_7;q,p;1),
\label{ft-bailey}\end{eqnarray}
где подразумевается, что $\prod_{m=1}^7t_m=q^2t_0^3$, $t_6=q^{-N}$,
$N\in\mathbb{N}$, и
\ba\nonumber
&& s_0=\frac{qt_0^2}{t_1t_2t_3},\quad s_1=\frac{s_0t_1}{t_0},\quad
s_2=\frac{s_0t_2}{t_0},\quad s_3=\frac{s_0t_3}{t_0},
\\ &&
s_4=t_4,\quad s_5=t_5,\quad s_6=t_6,\quad s_7=t_7.
\nonumber\ea

Если положить $t_2t_3=qt_0$, то левая часть (\ref{ft-bailey}) перейдет в
обрывающийся $_{10}V_9$ ряд, а в ряде в правой части получим $s_1=1$,
т.е. только его первый член будет отличен от нуля. Это приводит
к формуле суммирования Френкеля-Тураева, описывающей эллиптическое
обобщение формулы суммирования Джексона для обрывающегося совершенно уравновешенного
$_8\varphi_7$ ряда \cite{gas-rah:basic}. После уменьшения номеров параметров
$t_{4,5,6,7}$ на  два, получаем
\be
{_{10}V_9}(t_0;t_1,\dots, t_5;q,p;1) % } && \\ && \makebox[4em]{}
= \frac{\theta (qt_0;p;q)_N\prod_{1\leq r<s\leq 3}
      \theta (qt_0/t_rt_s;p;q)_N}{\theta (qt_0/t_1t_2t_3;p;q)_N
\prod_{r=1}^3\theta (qt_0/t_r;p;q)_N},
\label{ft-sum}\ee
где $t_r$ удовлетворяет условию балансировки $\prod_{r=1}^5 t_r =qt_0^2$
и условию обрыва $t_4=q^{-N}$, $N\in \mathbb{N}$.

\begin{remark}
Несмотря на наличие двойной периодичности по параметрам, бесконечные
полностью эллиптические гипергеометрические ряды (если они сходятся)
не являются эллиптическими функциями $u_m$ поскольку в параллелограмме
периодов имеется бесконечно много полюсов. Действительно, некоторые
из полюсов параметров $t_s,\; s=1,\ldots,r-4$,
находятся в точках $t_s=t_0q^{n+1}p^m,$ где $n\in\mathbb{N}$ и
$m\in\mathbb{Z}$. Если $q^k\neq p^l$ для любых $k, l\in\mathbb{N},$
то тогда существует бесконечно много целых $n$ и $m$, таких что $t_s$
остается внутри кольца $|p|< |t_s|< 1.$ Это означает, что в этом кольце
имеется бесконечно много полюсов для $t_s$.
\end{remark}

\section{Эллиптическая цепочка Бэйли}

Совершенно уравновешенные эллиптические гипергеометрические ряды
\begin{equation}
_{r+1}V_r(t_0;t_1,\ldots,t_{r-4};q,p;z)
=\sum_{n=0}^\infty \frac{\theta(t_0q^{2n})}{\theta(t_0)}
\prod_{m=0}^{r-4}\frac{\theta(t_m)_n}{\theta(qt_0t_m^{-1})_n}
\, (qz)^n
\label{vwp-2'}\end{equation}
для нечетных $r$ имеют условие балансировки вида
$$
\prod_{m=1}^{r-4}t_m= q^{(r-7)/2}t_0^{(r-5)/2}.
$$
При $z=1$, обрывающийся сбалансированный $_{10}V_9$ ряд суммируем
\cite{fre-tur:elliptic}.
Если мы возьмем предел $t_3\to t_0/t_1t_2$ в \re{ft-sum}, то условие
обрыва ряда сокращается и мы получаем замкнутое выражение для сбалансированного
$_8V_7$ ряда с неопределенным пределом суммирования
\begin{equation}
_{8}V_7^{(N)}(t_0;t_1,t_2,t_0/t_1t_2;q,p;1)
=\frac{\theta (qt_0,qt_1,qt_2,qt_0/t_1t_2)_N}
{\theta(q,qt_0/t_1,qt_0/t_2,qt_1t_2)_N},
\label{8V7}\end{equation}
где $_8V_7^{(N)}$ означает, что рассматривается
сумма первых $N+1$ членов $_8V_7$ ряда. Если положить $t_2=q^{-N}$,
то правая часть \re{8V7} равна нулю за исключением случая $N=0$.

Определим матрицу, которая будет играть ключевую роль при определении
цепочки преобразований Бэйли,
\begin{equation}
M_{nr}(a,k)=\frac{\theta(k/a)_{n-r}\theta(k)_{n+r}}
{\theta(q)_{n-r}\theta(aq)_{n+r}}\frac{\theta(aq^{2r})}
{\theta(a)}a^{n-r}
\end{equation}
и диагональную матрицу
\begin{equation}
D_{nr}(a;b,c)=\frac{\theta(b,c)_r}{\theta(aq/b,aq/c)_r}
\left(\frac{aq}{bc}\right)^r\delta_{nr}.
\end{equation}

\begin{definition}
Две последовательности $\alpha_n$ и $\beta_n,\, n\in\N,$ образуют
пару Бэйли по отношению к параметрам $a$ и $k$, если они связаны
матричным соотношением
\begin{equation}
\beta(a,k)=M(a,k)\alpha(a,k),
\label{ell-pair}\end{equation}
где $\alpha(a,k)$ и $\beta(a,k)$ обозначают столбцы, образованные из
элементов $\alpha_n$ и $\beta_n.$ Равенство \re{ell-pair}
имеет покомпонентный вид $\beta_n=\sum_{m=0}^n M_{nm}(a,k)\alpha_m$.
\end{definition}

С помощью формулы \re{8V7} легко увидеть, что
обратная матрица $M^{-1}(a,k)=M(k,a)$ (см., например,
\cite{war:extensions}). Другие полезные свойства матриц $M$ и $D$:
$M_{nr}(a,a)=\delta_{nr}$ и $D_{nr}(bc/q;b,c)=\delta_{nr}$.

\begin{theorem}
{\em (Эллиптическая лемма Бэйли 1 \cite{spi:bailey}).}
Если $\alpha(a,k)$ и $\beta(a,k)$ образуют пару Бэйли, то
\begin{eqnarray}\nonumber
&& \alpha'(a,k)=D(a;b,c)\alpha(a,m),\quad  \\
&& \beta'(a,k)=D(a;b,c)D^{-1}(m;b,c)M(m,k)D(m;b,c)\beta(a,m),
\label{a1}\end{eqnarray}
где $m=kbc/aq$, так же образуют пару Бэйли.
\end{theorem}

\begin{theorem}
{\em (Эллиптическая лемма Бэйли 2).}
Если $\alpha(a,k)$ и $\beta(a,k)$ образуют па\-ру Бэй\-ли, то
\begin{eqnarray}\nonumber
&& \alpha'(a,k)=D(k;d,e)D^{-1}(f;d,e)M(f,a)D(f;d,e)\alpha(f,k),\quad
\\
&& \beta'(a,k)=D(k;d,e)\beta(f,k),
\label{a2-ser}\end{eqnarray}
где $f=ade/kq$, так же образуют пару Бэйли.
Параметры $d$  и  $e$ полностью независимы от переменных
 $b $ и $ c$ в первой эллиптической лемме Бэйли.
\end{theorem}

\noindent
{\bf Доказательство.} С одной стороны, мы имеем
$$
\beta'(a,k)= D(k;d,e)\beta(f,k)=D(k;d,e)M(f,k)\alpha(f,k).
$$
С другой стороны,
\begin{eqnarray*}
&& \beta'(a,k)=M(a,k)\alpha'(a,k)
\\ && \makebox[4em]{}
= M(a,k)D(k;d,e)D^{-1}(f;d,e)M(f,a)D(f;d,e)\alpha(f,k).
\end{eqnarray*}
Эти два выражения равны благодаря ключевому равенству
\begin{equation}
D^{-1}(a;b,c)M(a,k)D(a;b,c)=D^{-1}(m;b,c)M(m,k)D(m;b,c)M(a,m),
\label{key}\end{equation}
эквивалентному $_{10}V_9$ формуле суммирования Френкеля-Тураева, что легко проверяется
прямой подстановкой. Аналогичным образом доказывается и первая лемма.
\hfill{Q.E.D.}

\smallskip
Последовательное применение первого и второго преобразований дает
\begin{eqnarray}\label{tr-1}
&&
\alpha^{(12)}(a,k)=D(a;b,c)D(m;d,e)D^{-1}(f;d,e)M(f,a)D(f;d,e)\alpha(f,m),
\\ &&
\beta^{(12)}(a,k)=D(a;b,c)D^{-1}(m;b,c)M(m,k)D(m;b,c)D(m;d,e)\beta(f,m),
\label{tr-2}\end{eqnarray}
где $m=kbc/aq$ и $f=ade/mq$. Перестановка преобразований дает другой
результат
\begin{eqnarray}\label{tr-3}
&&
\alpha^{(21)}(a,k)=D(k;d,e)D^{-1}(f;d,e)M(f,a)D(f;d,e)D(f;b,c)\alpha(f,m),
\\ &&
\beta^{(21)}(a,k)=D(k;d,e)D(f;b,c)D^{-1}(m;b,c)M(m,k)D(m;b,c)\beta(f,m),
\label{tr-4}\end{eqnarray}
где $f=ade/kq$ и $m=kbc/fq$. При $f=a$ или $m=k$ эти двойные
преобразования сводятся к первичным двум. Простейшие пары Бэйли
$\alpha_n=M_{n0}(a,k)$ и $\beta_n=\delta_{n0}$ порождают двоичное дерево
тождеств для эллиптических гипергеометрических рядов, включающее \eqref{ft-bailey}.
Соответствующие ряды содержат,
вообще говоря, произвольное число суммирований, то есть мы неизбежно
приходим к многократным рядам. Решетка (или дерево) Бэйли
таких тождеств на уровне $q$-гипергеометрических рядов обсуждалась в
работах \cite{aab:bailey} и \cite{and:bailey}.

В сравнительно старой работе \cite{bre:some}, Брессуд построил
некоторые примеры того, что Эндрюс назвал в \cite{and:bailey} вполне
уравновешенными парами Бэйли. Оказалось, что один из примеров Брессуда,
\begin{eqnarray*}
&& \alpha_n^{Br}(a,k)=\frac{1-aq^{2n}}{1-a}
\frac{(a,qa^2/k^2;q)_n}{(q,k^2/a;q)_n}
\left(\frac{k^2}{qa^2}\right)^n,
\\
&& \beta_n^{Br}(a,k)=\frac{(k^2,q^2a^2/k^2;q^2)_n
(-k^2/a;q)_{2n}}{(q^2,k^4/a^2;q^2)_n(-aq;q)_{2n}}
\left(\frac{k^2}{qa^2}\right)^n,
\end{eqnarray*}
допускает эллиптическое обобщение. Рассмотрим последовательности
\begin{eqnarray}\label{alpha4}
&& \alpha_n(a,k)=\frac{\theta(aq^{2n};p)}{\theta(a;p)}
\frac{\theta(a,qa^2/k^2;p;q)_n}{\theta(q,k^2/a;p;q)_n}
\left(\frac{k^2}{qa^2}\right)^n,
\\
&& \beta_n(a,k)=\frac{\theta(k^2,q^2a^2/k^2;p^2;q^2)_n
\theta(-k^2/a;p;q)_{2n}}{\theta(q^2,k^4/a^2;p^2;q^2)_n\theta(-aq;p;q)_{2n}}
\left(\frac{k^2}{qa^2}\right)^n.
\label{beta4}\end{eqnarray}
Покажем, что они удовлетворяют соотношению (\ref{ell-pair})
после возведения в квадрат всех параметров, входящих в него
(для $p=0$ это было показано в работе \cite{and-ber:bailey}).

Подставляя (\ref{alpha4}) в (\ref{ell-pair})
с $q, p, a, k$ замененными на $q^2, p^2, a^2, k^2$, получаем
\begin{eqnarray*}
&& \beta_n(a,k)=\sum_{j=0}^n\frac{\theta(aq^{2j};p)\theta(k^2;p^2;q^2)_{n+j}
\theta(k^2/a^2;p^2;q^2)_{n-j}}{\theta(a;p)\theta(a^2q^2;p^2;q^2)_{n+j}
\theta(q^2;p^2;q^2)_{n-j}}
\\ &&\makebox[6em]{}
\times \frac{\theta(a,qa^2/k^2;p;q)_j}
{\theta(q,k^2/a;p;q)_j}\left(\frac{k^2}{qa^2}\right)^j
\\
&& \makebox[3em]{}
=\frac{\theta(k^2/a^2,k^2;p^2;q^2)_n}{\theta(q^2,a^2q^2;p^2;q^2)_n}
{_{10}V_9}(a;kq^n,-kq^n,q^{-n},-q^{-n},qa^2/k^2;q,p;1),
\end{eqnarray*}
где мы использовали факторизацию
$$
\theta(z^2;p^2;q^2)_n=\theta(z;p;q)_n\theta(-z;p;q)_n.
$$
Суммируя этот $_{10}V_9$ ряд, получаем
$$
\beta_n(a,k)=
\frac{\theta(k^2,q^2a^2/k^2,k^2/a^2;p^2;q^2)_n\theta(-aq^{1-2n}/k^2;p;q)_n}
{\theta(q^2,a^2q^{2-2n}/k^2;p^2;q^2)_n\theta(-aq;p;q)_{2n}
\theta(aq^{1-n}/k^2;p;q)_n}.
$$
Последующее упрощение этого выражения с помощью соотношения
$\theta(z;p)=-z\theta(z^{-1};p)$ позволяет получить
(\ref{beta4}), что и требовалось. Таким образом мы доказали, что
(\ref{alpha4}) и (\ref{beta4}) действительно образуют эллиптическую
пару Бэйли.

Подставляя  (\ref{alpha4}) и (\ref{beta4}) в соотношения
(\ref{a'}), (\ref{b'}), где $q, p, a, k,\rho_1,\rho_2$ необходимо
возвести во вторую степень, получаем новую пару Бэйли:
\begin{eqnarray}\label{a5}
&& \alpha_n'(a,k)=\frac{\theta(aq^{2n};p)\theta(\rho_1^2,\rho_2^2;p^2;q^2)_n
\theta(a,qa^2/k^2;p;q)_n}
{\theta(a;p)\theta(a^2q^2/\rho_1^2,a^2q^2/\rho_2^2;p^2;q^2)_n
\theta(q,k^2/a;p;q)_n}\left(\frac{k}{\rho_1\rho_2}\right)^{2n},
\\ \nonumber
&& \beta_n'(a,k)= \frac{\theta(k^2\rho_1^2/a^2,k^2\rho_2^2/a^2,k^2,
k^2/c^2;p^2;q^2)_n}{\theta(a^2q^2/\rho_1^2,a^2q^2/\rho_2^2,q^2c^2,
q^2;p^2;q^2)_n}\sum_{j=0}^n\frac{\theta(c^2q^{4j};p^2)}
{\theta(c^2;p^2)}
\\ \nonumber && \makebox[6em]{}
\times\frac{\theta(\rho_1^2,\rho_2^2,k^2q^{2n},q^{-2n},c^2,a^2q^2/c^2;
p^2;q^2)_j}{\theta(k^2\rho_2^2/a^2,k^2\rho_1^2/a^2,c^2q^{2-2n}/k^2,
c^2q^{2+2n},q^2,c^4/a^2;p^2;q^2)_j}
\\ \nonumber && \makebox[6em]{}
\times \frac{\theta(-c^2/a;p;q)_{2j}}
{\theta(-aq;p;q)_{2j}}\left(\frac{qc^2}{a^2}\right)^j
\\ \nonumber && \makebox[3em]{}
=\frac{\theta(k^2\rho_1^2/a^2,
k^2\rho_2^2/a^2,k^2,k^2/c^2;p^2;q^2)_n}
{\theta(a^2q^2/\rho_1^2,a^2q^2/\rho_2^2,q^2c^2,q^2;p^2;q^2)_n}\;
{_{14}V_{13}}(c^2;\rho_1^2,\rho_2^2,q^{-2n},
\\ && \makebox[6em]{}
k^2q^{2n},a^2q^2/c^2,-c^2/a,-c^2q/a,-c^2/ap,-c^2qp/a;q^2,p^2;1),
\label{b5}\end{eqnarray}
где в последнем равенстве использовалась факторизация
\begin{equation}
\frac{\theta(-c^2/a;p;q)_{2j}}{\theta(-aq;p;q)_{2j}}
=\left(\frac{a^2q}{c^2}\right)^j\frac{\theta(-c^2/a,-c^2q/a,-c^2/ap,-c^2qp/a;
p^2;q^2)_j}{\theta(-aq^2,-aq,-aq^2p,-aq/p;p^2;q^2)_j}.
\label{fact-prod}\end{equation}
Подставляя теперь последовательность (\ref{a5}) в соотношение
(\ref{ell-pair}) (с параметрами, возведенными во вторую степень),
нетрудно увидеть что
\begin{eqnarray}\nonumber
&& \beta_n'(a,k)=\frac{\theta(k^2/a^2,k^2;p^2;q^2)_n}
{\theta(q^2,a^2q^2;p^2;q^2)_n}\;{_{14}V_{13}}(a;qa^2/c^2,
\\ && \makebox[6em]{}
\rho_1,-\rho_1,\rho_2,-\rho_2, kq^n,-kq^n,q^{-n},-q^{-n};q,p;1).
\label{b6}\end{eqnarray}
Сравнивая выражения (\ref{b5}) и (\ref{b6}), можно получить
новое тождество для эллиптических гипергеометрических рядов.
Таким образом, мы доказали следующую теорему.

\begin{theorem}\label{imp-th}
Для произвольных комплексных параметров $a,c,k,\rho_1,\rho_2$
справедливо следующее соотношение между двумя $_{14}V_{13}$
рядами с единичным аргументом:
\begin{eqnarray}\nonumber
\lefteqn{ _{14}V_{13}(c^2;\rho_1^2,\rho_2^2,k^2q^{2n},a^2q^2/c^2,q^{-2n}, } &&
\\ \nonumber && \makebox[4em]{}
-c^2/a,-c^2q/a,-c^2/ap,-c^2qp/a;q^2,p^2;1)
\\ \nonumber
&& =\frac{\theta(k^2/a^2,q^2c^2,a^2q^2/\rho_1^2,a^2q^2/\rho_2^2;p^2;q^2)_n}
{\theta(a^2q^2,k^2/c^2,k^2\rho_1^2/a^2,k^2\rho_2^2/a^2;p^2;q^2)_n}
\\ && \makebox[1em]{}
\times {_{14}V_{13}}(a;\rho_1,-\rho_1,\rho_2,-\rho_2,kq^n,-kq^n,q^{-n},-q^{-n},
qa^2/c^2;q,p;1),
\label{new-id}\end{eqnarray}
где $c=k\rho_1\rho_2/aq$.
\end{theorem}

\begin{corollary}
Рассмотрим предел $p\to 0$. Факторизация  (\ref{fact-prod})
принимает вид
$$
\frac{(-c^2/a;q)_{2j}}{(-aq;q)_{2j}}=\frac{(-c^2/a,-c^2q/a;q^2)_j}
{(-aq^2,-aq;q^2)_j}.
$$
В результате, соотношение (\ref{new-id}) сводится к
\begin{eqnarray}\nonumber
\lefteqn{ _{10}W_{9}(c^2;\rho_1^2,\rho_2^2,k^2q^{2n},a^2q^2/c^2,q^{-2n},
-c^2/a,-c^2q/a;q^2;qc^2/a^2)  }
\\ \nonumber
&& =\frac{(k^2/a^2,q^2c^2,a^2q^2/\rho_1^2,a^2q^2/\rho_2^2;q^2)_n}
{(a^2q^2,k^2/c^2,k^2\rho_1^2/a^2,k^2\rho_2^2/a^2;q^2)_n}
\\ && \makebox[1em]{}
\times {_{12}W_{11}}(a;\rho_1,-\rho_1,\rho_2,-\rho_2,kq^n,-kq^n,q^{-n},-q^{-n},
qa^2/c^2;q;q),
\label{and-ber}\end{eqnarray}
которое было получено в работах \cite{nas-rah:series,and-ber:bailey}.
\end{corollary}

Сделаем несколько замечаний. Во-первых, оба $_{14}V_{13}$ ряда в
(\ref{new-id}) сбалансированы, т.е. мы имеем дело с эллиптическими
гипергеометрическими рядами. Это резко отличается от формулы
(\ref{and-ber}), где $_{10}W_9$ ряд не сбалансирован, согласно определению
из \cite{gas-rah:basic}, в то время как ряд $_{12}W_{11}$ сбалансирован.
Это противоречие показывает, что наше определение балансировки,
которое напрямую связано с условием эллиптичности, более гибкое, чем
$q$-гипергеометрическое определение. Во-вторых, оба $_{14}V_{13}$
ряда имеют единичный аргумент, т.е. ограничение $z=1$ все еще сопровождает
условия совершенной уравновешенности и сбалансированности для появления
нетривиальных тождеств в рядах гипергеометрического типа с большими
значениями индексов (в рамках $q$-определений \cite{gas-rah:basic}
это условие соответствует выбору $z=q$). В третьих, соотношение
(\ref{new-id}) можно интерпретировать как правило квадратичных
преобразований базисных переменных $q\to q^2,\; p\to p^2$
для специфических $_{r+1}V_r$ рядов (ср., например, с
\cite{bis:change,war:summation}). Для последнего замечания, обозначим
$a=q^{\alpha}, k=q^\kappa, \rho_1=q^{u_1},
\rho_2=q^{u_2}$. Поскольку оба $_{14}V_{13}$ ряда полностью эллиптические,
т.е. они эллиптичны по $\alpha, \kappa, u_1, u_2$ переменным, то соотношение
(\ref{new-id}) описывает не что иное как бесконечную последовательность
тождеств для эллиптических функций, пронумерованных целым числом
 $n\in \mathbb{N}$. Другие примеры аналогичных преобразований и
связанных формул суммирования для эллиптических гипергеометрических
рядов обсуждались в работах \cite{war:extensions,war:summation2}.

Приведенная выше конструкция, обобщающая  результаты Эндрюса и Берковича
 \cite{and:bailey,and-ber:bailey}, должна допускать расширение
на многократные ряды. Некоторые преобразования типа Бэйли
для многократных рядов обсуждались в статьях \cite{den-gus:beta,mil-lil:alcl,
mil-lil:consequences,BS,war:summation,ros:elliptic,ros:bailey,ros-sch1,
ros-sch2}. Было бы интересно прояснить связаны ли
эллиптические пары Бэйли с концепцией суперсимметрии, по аналогии
со связью описанной в работах \cite{bms:susy,and-ber:trinomial}.

Подчеркнем еще раз, что найденное дерево тождеств
(\ref{ell-pair}), (\ref{a1}), (\ref{a2-ser}) не допускает конфлюэнтных
пределов, т.е. нельзя устремить параметры $t_m$ к нулю или бесконечности при
$p\neq 0$, и этот факт сильно ограничивает общее число тождеств, которые
могут быть получены на эллиптическом уровне. В частности, не ясно
существуют ли эллиптические аналоги триномиальных пар Бэйли
\cite{and-ber:trinomial, war:note}. Более того, ``альтернативная лемма
Бэйли'' работ \cite{and:bailey,and-ber:bailey} не допускает прямого
эллиптического обобщения, так как нет прямого эллиптического аналога необходимой для
этого $_3\varphi_2$ формулы суммирования $q$-Пфаффа-Саальшутца.
Вводя зависимость параметров от $p$, можно прямым образом спуститься
с $_{10}V_9$ формулы суммирования на этот уровень, но использование этого предела не
приводит к новым цепочкам Бэйли.

\section{Многократные ряды}

\subsection{Общее определение}

Следуя определению тета-гипергеометрических рядов одной переменной
можно определить формальные многократные суммы с коэффициентами,
которые связаны определенным образом с квазипериодическими тета-функциями
Якоби. Однако, мы ограничимся рассмотрением  только многократных
эллиптических гипергеометрических рядов.

\begin{definition}
Формальные ряды
$ \sum_{\lambda_1,\ldots,\lambda_n\in \Z} c(\lambda_1,\ldots,\lambda_n) $
называются  многократными эллиптическими гипергеометрическими
рядами если для всех $k=1,\ldots,n$ функции
$$
h_k(\mathbf{\lambda})=\frac{
c(\lambda_1,\ldots,\lambda_k+1,\ldots,\lambda_n)
}{c(\lambda_1,\ldots,\lambda_k,\ldots,\lambda_n)},
$$
являются эллиптическими функциями переменных $\lambda_k,\, k=1,\ldots,n,$
рассматриваемыми как комплексные переменные. Эти ряды называются полностью
эллиптическими если функции $h_k(\mathbf{\lambda})$ дополнительно
эллиптичны по всем свободным параметрам за исключением мультипликативных
факторов.
\end{definition}

В этом определении рассмотрен только общий случай двусторонних рядов,
который включает в себя случай суммирования по многомерному кубу
$\sum_{\lambda_1,\ldots,\lambda_n\in \N}$ и разбиениям
$ \sum_{0\leq\lambda_1\leq \ldots\leq \lambda_n}^\infty$, которые появляются
при подходящих ограничениях на коэффициенты $c(\mathbf{\lambda})$.

Предположим, что $h_k(\mathbf{\lambda})$ симметричны по $\lambda_1,\ldots,
\lambda_{k-1},\lambda_{k+1},\ldots,\lambda_n$. Тогда, коэффициенты
$c(\mathbf{\lambda})$ могут иметь вид
\begin{equation}\label{c-general}
c(\mathbf{\lambda})=
\prod_{k=1}^n\Biggl(\prod_{1\leq i_1<\ldots<i_k\leq n} \prod_{m=1}^{r_k}
\frac{[u_{km}]_{\lambda_{i_1}+\ldots+\lambda_{i_k} } }
{[v_{km}]_{\lambda_{i_1}+\ldots+\lambda_{i_k} } }\Biggr)
z_1^{\lambda_1}\cdots z_n^{\lambda_n},
\end{equation}
где
$ \sum_{k=1}^n C_{n-1}^{k-1} \sum_{m=1}^{r_k} (u_{km}-v_{km})=0. $
Если действие симметрической группы ${\cal S}_n$ переставляет
$\lambda_1,\ldots, \lambda_n$ одновременно со свободными параметрами,
входящими в $c(\mathbf{\lambda})$ (отличными от $z_1,\ldots,z_n$),
то ситуация богаче и допускает более общие комбинации $\lambda_k$
чем указано в выражении (\ref{c-general}). Дальнейший перебор возможностей
здесь не столь принципиален и мы переходим к рассмотрению конкретных
интересных примеров.

\subsection{Формулы суммирования для эллиптических гипергеометрических
рядов на корневой системе $C_n$}

В настоящее время известно три типа многократных эллиптических
гипергеометрических рядов, связанных с системой корней  $C_n$ и приводящих
к конструктивным тождествам (многократным обобщениям $_{10}V_9$-формулы суммирования
Френкеля-Тураева). Первый пример соответствует эллиптическому обобщению
$q$-гипергеометрических примеров, рассмотренных Аомото, Ито и Макдональдом
\cite{aom:elliptic,ito:theta,mac:constant}.

Структура таких рядов  видна из следующего многомерного
обобщения формулы суммирования Френ\-ке\-ля-Ту\-рае\-ва, предложенного Варнааром
в статье \cite{war:summation}, частично доказанного ван Диехеном и автором
в заметке \cite{die-spi:elliptic} и полностью доказанного Розенгреном в
работе \cite{ros:proof}. Пусть $N\in\mathbb{N}$ и семь параметров
$t, t_r\in\mathbb{C},\, r=0,\ldots,5,$ ограничены условием балансировки
$t^{2n-2}\prod_{r=0}^5t_r=q$ и условием обрыва $t^{n-1}t_0t_4=q^{-N}.$
Тогда выполняется следующее тождество для тета-функций
\begin{eqnarray}\nonumber
&&
\sum_{0\leq \lambda_1\leq \lambda_2\leq \cdots \leq \lambda_n\leq N}
q^{\sum_{j=1}^n\lambda_j} t^{2 \sum_{j=1}^n (n-j)\lambda_j}
\\  \nonumber && \makebox[8em]{}\times
\prod_{1\leq j<k\leq n} \Biggl(
\frac{\theta(\tau_k\tau_jq^{\lambda_k+\lambda_j},
\tau_k\tau_j^{-1}q^{\lambda_k-\lambda_j};p)}
{\theta(\tau_k\tau_j,\tau_k\tau_j^{-1};p)} \\  \nonumber
&& \makebox[8em]{}\times
\frac{\theta(t\tau_k\tau_j;p;q)_{\lambda_k+\lambda_j}}
     {\theta(qt^{-1}\tau_k\tau_j;p;q)_{\lambda_k+\lambda_j}}
\frac{\theta(t\tau_k\tau_j^{-1};p;q)_{\lambda_k-\lambda_j}}
     {\theta(qt^{-1}\tau_k\tau_j^{-1};p;q)_{\lambda_k-\lambda_j}}
\Biggr) \\  \nonumber
&&\makebox[3em]{} \times \prod_{j=1}^n
\Biggl( \frac{\theta(\tau_j^2q^{2\lambda_j};p)}{\theta(\tau_j^2;p)}
\prod_{r=0}^5 \frac{\theta(t_r\tau_j;p;q)_{\lambda_j}}
     {\theta(qt_r^{-1}\tau_j;p;q)_{\lambda_j}}\Biggr)  \\
&&\makebox[3em]{}
=\prod_{j=1}^n\frac{\theta(qt^{n+j-2}t_0^2;p;q)_N
      \prod_{1\leq r <s \leq 3} \theta(qt^{1-j} t_r^{-1}t_s^{-1};p;q)_N}
      {\theta(qt^{2-n-j}\prod_{r=0}^3t_r^{-1};p;q)_N
      \prod_{r=1}^3 \theta(qt^{j-1}t_0t_r^{-1};p;q)_N}.
\label{multi-1}\end{eqnarray}
Здесь параметры $\tau_j$ связаны с $t_0$ и $t$ так:
$\tau_j=t_0t^{j-1}$, $j=1,\ldots ,n$. Заметим, что члены ряда
$c(\mathbf{\lambda})$ симметричны при одновременной перестановке
переменных $\lambda_j$ и $\lambda_k$ и параметров $\tau_j$ и $\tau_k$
для произвольных $j\neq k$ (при этом необходимо считать $\tau_j$
независимыми переменными).

\begin{theorem}
Ряд в левой части формулы (\ref{multi-1}) является многократным
полностью эллиптическим гипергеометрическим рядом.
\end{theorem}
{\bf Доказательство.}
Отношение соседних членов ряда имеет вид
\begin{eqnarray}\nonumber
h_l(\mathbf{\lambda})=\prod_{j=1}^{l-1}
\frac{\theta(\tau_j\tau_lq^{\lambda_j+\lambda_l+1},\tau_j^{-1}\tau_l
q^{\lambda_l+1-\lambda_j}, t\tau_j\tau_lq^{\lambda_j+\lambda_l},
t\tau_j^{-1}\tau_lq^{\lambda_l-\lambda_j};p)}
{\theta(\tau_j\tau_lq^{\lambda_j+\lambda_l},\tau_j^{-1}\tau_l
q^{\lambda_l-\lambda_j}, t^{-1}\tau_j\tau_lq^{\lambda_j+\lambda_l+1},
t^{-1}\tau_j^{-1}\tau_lq^{\lambda_l+1-\lambda_j};p)}
\\  \nonumber
\times \prod_{k=l+1}^n\frac{\theta(\tau_k\tau_lq^{\lambda_k+\lambda_l+1},
\tau_k\tau_l^{-1}q^{\lambda_k-\lambda_l-1},t\tau_k\tau_lq^{\lambda_k+
\lambda_l},t^{-1}\tau_k\tau_l^{-1}q^{\lambda_k-\lambda_l};p)}
{\theta(\tau_k\tau_lq^{\lambda_k+\lambda_l},\tau_k\tau_l^{-1}
q^{\lambda_k-\lambda_l}, t^{-1}\tau_k\tau_lq^{\lambda_k+\lambda_l+1},
t\tau_k\tau_l^{-1} q^{\lambda_k-\lambda_l-1};p)}
\\
\times qt^{2(n-l)}
\frac{\theta(\tau_l^2q^{2\lambda_l+2};p)}{\theta(\tau_l^2;p)}
\prod_{m=0}^5\frac{\theta(t_m\tau_lq^{\lambda_l};p)}
{\theta(t_m^{-1}\tau_lq^{\lambda_l+1};p)}.
\label{h-AIM}\end{eqnarray}

Используя равенства (\ref{fun-rel}), легко проверить эллиптичность
этих $h_l(\mathbf{\lambda})$ по $\lambda_i$ для $i< l$ и $i>l$
(для этого достаточно убедиться, что $h_l$ не меняется после
замены $q^{\lambda_i}$ на $pq^{\lambda_i}$). Для проверки эллиптичности
по $\lambda_l$ необходим существенно более сложный расчет.
Замена  $q^{\lambda_l}$ на $pq^{\lambda_l}$ в произведении
$\prod_{j=1}^{l-1}$ выдает множитель $t^{-4(l-1)}$.
Произведение $\prod_{k=l+1}^n$ дает множитель $t^{-4(n-l)}$.
Оставшаяся часть $h_l$ генерирует множитель $q^{-4}\prod_{m=0}^5
qt_m^{-2}$. Произведение всех этих трех множителей принимает вид
$q^2t^{-4(n-1)}\prod_{m=0}^5 t_m^{-2}$ и равно единице благодаря
условию балансировки.

Таким образом мы показали, что взятый ряд действительно является
многократным эллиптическим гипергеометрическим рядом. Докажем теперь
его полную эллиптичность, т.е. $p$-мультипликативную инвариантность
для параметров $t_m,\, m=0,\ldots,4,$ и $t$.
Инвариантность при $t_m\to pt_m,\, m=1,\ldots,4,$
следует из условия балансировки, так же как это было в случае однократного
ряда. При $t_0\to pt_0$ преобразовании, произведение
$\prod_{j=1}^{l-1}$ дает множитель $t^{-4(l-1)}$, а
произведение $\prod_{k=l+1}^n$ приводит к фактору $t^{-4(n-l)}$.
Оставшаяся часть $h_l$ генерирует $q^2\prod_{m=0}^5t_m^{-2}$.
Произведение всех множителей равно единице.

Наконец, замена $t\to pt$ приводит к наиболее сложному вычислению.
Произведение $\prod_{j=1}^{l-1}$ выдает сложный множитель
$(q^{2\lambda_l+1}t^{2(l-1)}\tau_l^2p^{2l-3})^{2(1-l)}$.
Произведение $\prod_{k=l+1}^n$ генерирует не менее сложное выражение
$(q^{2\lambda_l+1}t^{2(l-1)}\tau_l^2p^{2(l-1)})^{2(l-n)}$.
Оставшаяся часть $h_l$ приводит к множителю
$(q^{2\lambda_l+1}t^{2(l-1)}\tau_l^2)^{2(n-1)}p^{(4n-6)(l-1)}$
(после использования условия балансировки). Общий фактор равен
единице, что и доказывает теорему. \hfill{Q.E.D.}

\smallskip

Второй пример многократных рядов соответствует эллиптическому
обобщению базисных гипергеометрических рядов типа Милна-Густафсона
\cite{mil:multidimensional,gus:macdonald,den-gus:beta},
которые, в свою очередь, являются $q$-аналогами рядов Биденхарна-Лаука-Холмана
\cite{hbl:hypergeometric}. Их качественную структуру можно увидеть из
следующей формулы суммирования, предложенной ван Диехеном и автором
в заметке \cite{die-spi:modular}. Пусть $q^n\neq p^m$ при $n,m\in\mathbb{N}$.
Тогда для параметров $t_0,\ldots ,t_{2n+3}$, удовлетворяющих условиям
балансировки $q^{-1}\prod_{r=0}^{2n+3}t_r=1$ и обрыва
$q^{N_j}t_jt_{n+j}=1,\, j=1,\ldots ,n,$ где $N_j\in\mathbb{N},$
имеет место тождество
\begin{eqnarray}\nonumber
&& \sum_{\stackrel{0\leq \lambda_j\leq N_j}{j=1,\ldots ,n}}
q^{\sum_{j=1}^n j\lambda_j}
\prod_{1\leq j<k\leq n}
\frac{\theta (t_jt_kq^{\lambda_j+\lambda_k},
              t_jt_k^{-1}q^{\lambda_j-\lambda_k} ;p)}
     {\theta (t_jt_k,t_jt_k^{-1} ;p)}  \\
\label{multi-2}
&& \quad\qquad\qquad\qquad\times \prod_{1\leq j\leq n} \Biggl(
\frac{\theta (t_j^2q^{2\lambda_j};p)}{\theta (t_j^2;p)}
\prod_{0\leq r\leq 2n+3} \frac{\theta (t_jt_r;p;q)_{\lambda_j}}
     {\theta (qt_jt_r^{-1};p;q)_{\lambda_j}}\Biggr)  \\
&&= \theta (qa^{-1}b^{-1},qa^{-1}c^{-1},qb^{-1}c^{-1};p;q)_{N_1+\cdots +N_n}
\nonumber\\
&&\qquad\times\prod_{1\leq j<k\leq n}
\frac{\theta (qt_jt_k;p;q)_{N_j}\theta (qt_jt_k ;p;q)_{N_k}}
     {\theta (qt_jt_k;p;q)_{N_j+N_k}} \nonumber\\
&&\qquad \times \prod_{1\leq j\leq n}
\frac{\theta (qt_j^2;p;q)_{N_j}}
     {\theta (qt_ja^{-1},qt_jb^{-1},
 qt_jc^{-1},q^{1+N_1+\cdots+N_n-N_j}t_j^{-1}a^{-1}b^{-1}c^{-1};p;q)_{N_j}},
\nonumber\end{eqnarray}
где $a\equiv t_{2n+1}$, $b\equiv t_{2n+2}$, $c\equiv t_{2n+3}$.
Заметим, что для этого ряда коэффициенты $c(\mathbf{\lambda})$ симметричны
при одновременной перестановке $\lambda_j$ и $\lambda_k$ совместно с
параметрами $t_j$ и $t_k$ для произвольных $j,k=1,\ldots,n,\, j\neq k$.

\begin{theorem}
Ряд гипергеометрического типа, стоящий в левой части равенства (\ref{multi-2}),
полностью эллиптичен.
\end{theorem}
{\bf Доказательство.}
Отношения соседних членов ряда (\ref{multi-2}) имеют вид
\begin{eqnarray}\nonumber
\lefteqn{
h_l(\mathbf{\lambda})= \prod_{j=1}^{l-1}
\frac{\theta(t_jt_lq^{\lambda_j+\lambda_l+1},t_jt_l^{-1}
q^{\lambda_j-\lambda_l-1};p)}{\theta(t_jt_lq^{\lambda_j+\lambda_l},
t_jt_l^{-1}q^{\lambda_j-\lambda_l};p)}
}&& \\  \nonumber &&
\times \prod_{k=l+1}^{n}
\frac{\theta(t_lt_kq^{\lambda_l+\lambda_k+1},t_lt_k^{-1}
q^{\lambda_l+1-\lambda_k};p)}{\theta(t_lt_kq^{\lambda_l+\lambda_k},
t_lt_k^{-1}q^{\lambda_l-\lambda_k};p)}
%  \\  \times
q^l\frac{\theta(t_l^2q^{2\lambda_l+2};p)}
{\theta(t_l^2q^{2\lambda_l};p)} \prod_{m=0}^{2n+3}
\frac{\theta(t_lt_mq^{\lambda_l};p)}{\theta(t_lt_m^{-1}q^{\lambda_l+1};p)}.
\label{h-MG}\end{eqnarray}
Легко проверить эллиптичность $h_l(\mathbf{\lambda})$ по $\lambda_j$
при $j<l$ и $j>l$. Эллиптичность по $\lambda_l$ следует из
более сложного вычисления: замена $q^{\lambda_j}\to
pq^{\lambda_j}$ приводит к множителям $q^{2(1-l)}$ (из произведения
$\prod_{j=1}^{l-1}$), $q^{2(l-n)}$  (из произведения $\prod_{k=l+1}^n$),
и $q^{-4}\prod_{m=0}^{2n+3}qt_m^{-2}$ (из остатка $h_l(\mathbf{\lambda})$).
Перемножение этих выражений дает единицу благодаря условию балансировки.

Эллиптичность по параметрам $t_m,\, m=0,\ldots,2n+3,$ проверяется
отдельно для $m<l, m>l$ и $m=l$. Первые два случая достаточно легки
и не заслуживают отдельного рассмотрения. Замена
$t_l\to pt_l$ приводит к множителям:
$q^{2(1-l)}$ (из произведения $\prod_{j=1}^{l-1}$),
$q^{2(l-n)}$ (из произведения $\prod_{k=l+1}^n$), и
$q^{-2n}\prod_{m=0}^{2n+3}t_m^{-2}$ (из остатка $h_l(\mathbf{\lambda})$).
Условие балансировки опять гарантирует, что полный множитель равен единице,
что завершает доказательство полной эллиптичности взятого ряда. \hfill{Q.E.D.}

\smallskip

В работе \cite{war:summation}, Варнаар доказал следующий эллиптический
аналог многомерной $C_n$ формулы суммирования Шлоссера \cite{sch:summation}:
\begin{eqnarray}\nonumber
\lefteqn{
\sum_{\lambda_1,\ldots, \lambda_n=0}^N
q^{\sum_{j=1}^nj\lambda_j}
\prod_{1\leq i<j \leq n} \frac{\theta(t_it_j^{-1}q^{\lambda_i-\lambda_j},
t_it_jq^{\lambda_i+\lambda_j};p)}
{\theta(t_it_j^{-1},t_it_j;p)}  } &&
\\ \label{C_n}
&& \times
\prod_{j=1}^n\frac{\theta(t_j^2q^{2\lambda_j};p)
\theta(t_j^2,bt_j, ct_j, dt_j,et_j,q^{-N};p;q)_{\lambda_j}}
{\theta(t_j^2;p)\theta(q,t_jq/b, t_jq/c,t_jq/d,t_jq/e,
t_j^2q^{N+1};p;q)_{\lambda_j}}
\\  \nonumber &&
=\prod_{1\leq i<j\leq n}\frac{\theta(t_it_jq^N;p)}{\theta(t_it_j;p)}
\prod_{j=1}^n\frac{\theta(t_j^2q,q^{2-j}/bc,q^{2-j}/bd,q^{2-j}/cd;p;q)_N}
{\theta(t_j^{-1}q^{2-n}/bcd, t_jq/b, t_jq/c, t_jq/d;p;q)_N},
\end{eqnarray}
где параметры удовлетворяют условию балансировки $bcde=q^{N-n+2}$.

\begin{remark}
Формулы суммирования, приведенные в работах \cite{sch:summation,war:summation},
содержат дополнительный параметр $a$, но он может быть устранен
переопределением переменных $t_j,b,c,d,e$.
\end{remark}

\begin{theorem}
Тета-гипергеометрический ряд, стоящий в левой части равенства (\ref{C_n}), является
пол\-нос\-тью эллиптическим и модулярно инвариантным.
\end{theorem}
{\bf Доказательство.}
Рассмотрим отношения соседних членов ряда (\ref{C_n}):
\begin{eqnarray}\nonumber
\lefteqn{
h_\ell(\mathbb{\lambda})=q^\ell \frac{\theta(t_\ell^2q^{2\lambda_\ell+2};p)}
{\theta(t_\ell^2q^{2\lambda_\ell};p)} \prod_{i=1}^{\ell-1}
\frac{\theta(t_it_\ell^{-1}q^{\lambda_i-\lambda_\ell-1},t_it_\ell
q^{\lambda_i+\lambda_\ell+1};p)}
{\theta(t_it_\ell^{-1}q^{\lambda_i-\lambda_\ell},t_it_\ell
q^{\lambda_i+\lambda_\ell};p)}  } &&
\\  \label{hlcn} &&
\times \prod_{j=\ell+1}^n
\frac{\theta(t_\ell t_j^{-1}q^{\lambda_\ell-\lambda_j+1},t_\ell t_j
q^{\lambda_\ell+\lambda_j+1};p)}
{\theta(t_\ell t_j^{-1}q^{\lambda_\ell-\lambda_j},t_\ell t_j
q^{\lambda_\ell+\lambda_j};p)}
\\ \nonumber &&  \times
\frac{\theta(t_\ell^2q^{\lambda_\ell}, bt_\ell q^{\lambda_\ell},
ct_\ell q^{\lambda_\ell},dt_\ell q^{\lambda_\ell},et_\ell q^{\lambda_\ell},
q^{-N+\lambda_\ell};p)}
{\theta(q^{\lambda_\ell+1}, t_\ell q^{\lambda_\ell+1}/b,
t_\ell q^{\lambda_\ell+1}/c,t_\ell q^{\lambda_\ell+1}/d,
t_\ell q^{\lambda_\ell+1}/e,t_\ell^2q^{N+1+\lambda_\ell};p)}.
\label{hcn}\end{eqnarray}
Докажем сначала, что эта функция эллиптична по $\lambda_k$.
Для этого достаточно заменить $q^{\lambda_k}$ на $pq^{\lambda_k}$,
применить правило $\theta(pz;p)=-z^{-1}\theta(z;p)$ и проверить
сокращение всех дополнительных факторов, появляющихся в результате этого.
Для $k<\ell$ или $k>\ell$, ситуация достаточно проста и мы не рассматриваем
ее. Для $k=\ell$, первая дробь из тета-функций в (\ref{hlcn}) дает
множитель $q^{-4}$, следующие два произведения, $q^{2(1-n)}$.
Последняя дробь, содержащая параметры $b,\ldots, e$, генерирует
произведение $q^{2N+6}(bcde)^{-2}$. В итоге, полный результирующий множитель
равен 1 благодаря условию балансировки.

Для проверки эллиптичности по параметрам и модулярной инвариантности
$h_\ell(\mathbf{\lambda})$, произведем замену $t_j,b,c,d,e \to q^{g_j},
q^\beta, q^\gamma, q^\delta, q^\epsilon$. После этого мы можем переписать
(\ref{hlcn}) в терминах эллиптических чисел следующим образом
\begin{eqnarray}\nonumber
&& h_\ell(\mathbb{\lambda})=\frac{[2g_\ell+2\lambda_\ell+2]}
{[2g_\ell+2\lambda_\ell]}
\prod_{i=1,i\neq \ell}^{n}
\frac{[g_i-g_\ell+\lambda_i-\lambda_\ell-1,g_i+ g_\ell
+\lambda_i+\lambda_\ell+1]}
{[g_i-g_\ell+\lambda_i-\lambda_\ell,g_i+g_\ell
+\lambda_i+\lambda_\ell]}
\\  && \qquad\qquad
\times \frac{[2g_\ell+\lambda_\ell,\beta+g_\ell+\lambda_\ell,
\gamma + g_\ell+\lambda_\ell]}
{[\lambda_\ell+1, g_\ell+\lambda_\ell+1-\beta,
g_\ell+\lambda_\ell+1-\gamma]}
\nonumber \\  &&  \qquad\qquad
\times \frac{[\delta+g_\ell+\lambda_\ell,\epsilon+g_\ell
+\lambda_\ell,\lambda_\ell -N]}
{[g_\ell+\lambda_\ell+1-\delta,
g_\ell+\lambda_\ell+1-\epsilon,
2g_\ell+N+1+\lambda_\ell]}.
 \label{hlcn-ell}\end{eqnarray}
Из этого представления следует, что при переписывании (\ref{C_n})
в терминах эллиптических чисел первый множитель
$q^{\sum_{j=1}^nj\lambda_j}$ исчезает. В итоге одновременная
перестановка переменных $\lambda_j, \lambda_k$ и параметров
$g_j, g_k$ не меняет коэффициентов $c(\mathbf{\lambda})$
(имеет место одинаковое число изменений знаков в числителе и
знаменателе $c(\mathbf{\lambda})$, происходящих из-за нечетности
$[-u]=-[u]$). Таким образом, инвариантность относительно
симметрической группы естественным образом реализована в этих рядах.

Эллиптичность по $g_k$ в (\ref{hlcn-ell}) следует из инвариантности
(\ref{hlcn}) относительно сдвига $t_k\to pt_k$. Для $k<\ell$ или $k>\ell$
это легко проверяется. Для $k=\ell$, первая дробь тета-функций
порождает множитель $q^{-4}$, следующие два произведения дают
$q^{2(1-n)}$. Остаток в $h_\ell$ приводит к множителю
$q^{2N+6}(bcde)^{-2}$, и полный множитель равен 1 благодаря условию
балансировки. Для проверки эллиптичности по $\beta$ (или, эквивалентно, по
$\gamma$ или $\delta$), заменим $b\to pb, e\to p^{-1}e$ в (\ref{hlcn}).
Только выражение в последней строчке дает некоторые вклады и, поскольку
$\theta(p^{-1}z;p)=-p^{-1}z\theta(z;p)$, они сокращают друг друга.

Обратимся теперь к модулярной инвариантности. Поскольку числитель
и знаменатель $h_\ell$ содержит одинаковое число $\theta_1$-функций,
численные факторы, появляющиеся из-за преобразований
$\tau\to \tau+1, \tau\to -1/\tau$ в (\ref{modular-tr1}),
(\ref{modular-tr2}), сокращаются.  Таким образом достаточно проверить,
что разность квадратов аргументов эллиптических чисел в числителе и
знаменателе равна нулю. Для первой дроби тета-функций в (\ref{hlcn-ell}),
эта разность равна $4+8(g_\ell+\lambda_\ell)$. Следующее произведение
дает $2(1+2(g_\ell+\lambda_\ell))(n-1)$. Остаток $h_\ell$ генерирует
разность
$$
2(\beta+\gamma+\delta+\epsilon)(1+2(g_\ell+\lambda_\ell))
-4 - 4(g_\ell+\lambda_\ell)(N+3)-2(N+1).
$$
Полная сумма этих выражений равна нулю, что доказывает модулярную
инвариантность функции $h_\ell$. По терминологии Эйхлера-Загира
\cite{eic-zag:theory}, это означает, что все функции
$h_\ell(\mathbf{\lambda}),$ $\ell=1,\ldots, n,$
в (\ref{hlcn-ell}) являются мероморфными формами Якоби веса 0 и
индекса 0. Все вместе они образуют некоторую структуру, связанную
с системой корней $C_n$.

Хотя совпадение левой и правой частей (\ref{C_n}) было доказано
Варнааром \cite{war:summation}, для полноты мы проверили, что выражение справа в
(\ref{C_n}) является эллиптической функцией всех параметров
$g_j, \beta, \gamma, \delta$ и оно модулярно инвариантно.
\hfill{Q.E.D.}

\begin{remark}
Модулярная инвариантность рядов (\ref{multi-1}) и (\ref{multi-2})
была установлена в работах \cite{die-spi:elliptic} и \cite{die-spi:modular},
соответственно. Мы выдвигаем гипотезу, что все полностью эллиптические
гипергеометрические ряды автоматически модулярно инвариантны так же
как это было для однократного ряда. По аналогии с однократным случаем,
можно ввести многократные тета-гипергеометрические
ряды и понятия вполне уравновешенности и совершенной уравновешенности
для них, но это здесь не делается, так как общий вид таких рядов
еще не установлен.
\end{remark}

\subsection{Аналоги суммы Френкеля-Тураева для $A_n$ и $D_n$ рядов}

В статье Бхатнагара и Шлоссера \cite{BS} изучались различные преобразования
Бэйли для многократных $_{10}\varphi_9$ рядов. Цепочка Бэйли, построенная автором
в работе \cite{spi:bailey}, генерирует бесконечную последовательность таких
преобразований для одномерных эллиптических гипергеометрических рядов. Сравнение
исследований \cite{BS} и \cite{spi:bailey} естественным образом привело к
предположению, что многие результаты работы \cite{BS} допускают эллиптическое обобщение.
В частности, должны быть справедливы следующие
эллиптические аналоги многократных формул суммирования Джексона, найденных
Милном \cite{mil:multidimensional},  Шлоссером \cite{Sch1} и Бхатнагаром \cite{B}.

Эллиптическая $A_n$-формула суммирования Милна:
\begin{eqnarray} \nonumber
\lefteqn{
\sum_{0\leq \lambda_j \leq N_j \atop j=1,\ldots, n}
q^{\sum_{j=1}^nj\lambda_j}
\prod_{j=1}^n\frac{\theta(t_jq^{\lambda_j+|\lambda|};p)}{\theta(t_j;p)}
\prod_{1\leq i<j \leq n} \frac{\theta(t_it_j^{-1}q^{\lambda_i-\lambda_j};p)}
{\theta(t_it_j^{-1};p)}  } &&
\\ \nonumber
&& \times \prod_{i,j=1}^n\frac{\theta(t_it_j^{-1}q^{-N_j};p;q)_{\lambda_i}}
{\theta(qt_it_j^{-1};p;q)_{\lambda_i}}
\prod_{j=1}^n\frac{\theta(t_j;p;q)_{|\lambda|}}
{\theta(t_jq^{1+N_j};p;q)_{|\lambda|}}
\\ \nonumber
&& \times \frac{\theta(b,c;p;q)_{|\lambda|}}
{\theta(q/d, q/e;p;q)_{|\lambda|}}
\prod_{j=1}^n \frac{\theta(dt_j, et_j;p;q)_{\lambda_j}}
{\theta(t_jq/b, t_jq/c;p;q)_{\lambda_j}}
\\
&&
=\frac{\theta(q/bd,q/cd;p;q)_{|N|}}{\theta(q/d,q/bcd;p;q)_{|N|}}
\prod_{j=1}^n\frac{\theta(t_jq,t_jq/bc;p;q)_{N_j}}
{\theta(t_jq/b,t_jq/c;p;q)_{N_j}},
\label{A_n}\end{eqnarray}
где $|\lambda|=\lambda_1+\cdots+\lambda_n$,
$|N|=N_1+\cdots+N_n$ и $bcde=q^{1+|N|}$.

Эллиптическая $D_n$-формула суммирования Бхатнагара-Шлоссера (как заявлено в \cite{BS},
$D_n$-суммы работ \cite{B,Sch1} существенно эквивалентны друг другу):
\begin{eqnarray}\nonumber
\lefteqn{
\sum_{0\leq \lambda_j \leq N_j \atop j=1,\ldots, n}
q^{\sum_{j=1}^nj\lambda_j}
\prod_{j=1}^n\frac{\theta(t_jq^{\lambda_j+|\lambda|};p)}{\theta(t_j;p)}
\prod_{1\leq i<j \leq n} \frac{\theta(t_it_j^{-1}q^{\lambda_i-\lambda_j};p)}
{\theta(t_it_j^{-1};p)}  } &&
\\ \nonumber &&
\times \prod_{1\leq i<j \leq n}
\theta(t_it_jq/bcd;p;q)_{\lambda_i+\lambda_j}^{-1}
\prod_{i,j=1}^n
\theta(t_it_jq^{1+N_j}/bcd;p;q)_{\lambda_i}
\\ \nonumber &&
\times \prod_{i,j=1}^n\frac{\theta(t_it_j^{-1}q^{-N_j};p;q)_{\lambda_i}}
{\theta(qt_it_j^{-1};p;q)_{\lambda_i}}
\prod_{j=1}^n\frac{\theta(t_j;p;q)_{|\lambda|}
\theta(bcdt_j^{-1};p;q)_{|\lambda|-\lambda_j}}
{\theta(t_jq^{1+N_j}, bcdt_j^{-1}q^{-N_j};p;q)_{|\lambda|}}
\\ \nonumber &&
\times \frac{\theta(b,c,d;p;q)_{|\lambda|}}
{\prod_{j=1}^n\theta(t_jq/b, t_jq/c, t_jq/d;p;q)_{\lambda_j}}
\\  &&
=\prod_{j=1}^n
\frac{\theta(t_jq,t_jq/bc,t_jq/bd, t_jq/cd;p;q)_{N_j}}
{\theta(t_jq/bcd,t_jq/b,t_jq/c,t_jq/d;p;q)_{N_j}}.
\label{D_n}\end{eqnarray}

\begin{remark} В представлении $A_n$-суммы Милна, приведенном
в приложении статьи \cite{BS}, имеется дополнительный параметр $a$,
но он может быть устранен переопределением переменных $t_j, d$ и $e$.
Аналогично, в Шлоссеровской версии $D_n$-суммы, приведенной там же,
имеется дополнительный параметр $a$, который может быть устранен переопределением $t_j$
(то же самое справедливо для Бхатнагаровской версии этой формулы суммирования).
\end{remark}

\begin{theorem}
$A_n$ и $D_n$ тета-гипергеометрические ряды, стоящие слева в равенствах
(\ref{A_n}) и (\ref{D_n}), являются полностью эллиптическими и модулярно
инвариантными. Соответствующие выражения в правых частях обладают теми
же свойствами.
\end{theorem}
{\bf Доказательство.}
Рассмотрим функции $h_\ell(\mathbf{\lambda})$ в $A_n$-случае:
\begin{eqnarray}\nonumber
\lefteqn{
h_\ell(\mathbb{\lambda})=q^\ell
\prod_{j=1\atop j\neq\ell}^n\frac{\theta(t_jq^{\lambda_j+|\lambda|+1};p)}
{\theta(t_jq^{\lambda_j+|\lambda|};p)}\;
\frac{\theta(t_\ell q^{\lambda_\ell+|\lambda|+2};p)}
{\theta(t_\ell q^{\lambda_\ell+|\lambda|};p)} } &&
\\ \nonumber &&  \times
\prod_{i=1}^{\ell-1}
\frac{\theta(t_it_\ell^{-1}q^{\lambda_i-\lambda_\ell-1};p)}
{\theta(t_it_\ell^{-1}q^{\lambda_i-\lambda_\ell};p)}
\prod_{j=\ell+1}^n
\frac{\theta(t_\ell t_j^{-1}q^{\lambda_\ell-\lambda_j+1};p)}
{\theta(t_\ell t_j^{-1}q^{\lambda_\ell-\lambda_j};p)}
\\ \nonumber &&  \times
\prod_{j=1}^n\frac{\theta(t_\ell t_j^{-1}q^{\lambda_\ell-N_j},
t_jq^{|\lambda|};p)}{\theta(t_\ell t_j^{-1}q^{\lambda_\ell+1},
t_jq^{1+N_j+|\lambda|};p)}
\frac{\theta(bq^{|\lambda|}, cq^{|\lambda|},
dt_\ell q^{\lambda_\ell},et_\ell q^{\lambda_\ell};p)}
{\theta(q^{|\lambda|+1}/d, q^{|\lambda|+1}/e,
t_\ell q^{\lambda_\ell+1}/b,t_\ell q^{\lambda_\ell+1}/c;p)}.
\label{han}\end{eqnarray}

Докажем сначала эллиптичность $h_\ell(\mathbf{\lambda})$ по $\lambda_k$.
Для $k<\ell$ или $k>\ell$, сдвиг $q^{\lambda_k}\to pq^{\lambda_k}$
в $h_\ell$ приводит к следующим дополнительным множителям: $q^{-n}$ и
$q^{-2}$ от первых двух дробей тета-функций, соответственно;
следующие два произведения дают $q$; из произведения $\prod_{j=1}^n$
в третьей строке получается вклад $q^{n+|N|}$, а остаток $h_\ell$
генерирует фактор $q^{1-|N|}$. Полный множитель равен 1.

Аналогично, заменяя $q^{\lambda_\ell}$ на $pq^{\lambda_\ell}$,
мы получаем последовательно множители $q^{1-n}$, $q^{-4}$,
$q^{1-\ell}$, $q^{\ell-n}$, $q^{2n+2|N|}$ и $q^{2-2|N|}$,
и их  перемножение дает 1. Таким образом, функции
$h_\ell(\mathbf{\lambda})$, порождающие левую часть (\ref{A_n}),
эллиптичны по всем $\lambda_k$.

Для того, чтобы увидеть естественное действие симметрической группы,
заменим $t_j,$ $b,$ $c,$ $d,$ $e \to q^{g_j},$ $q^{\beta},$ $q^\gamma,$
$q^\delta,$  $q^\epsilon$,
и перепишем формулу (\ref{A_n}) в терминах эллиптических чисел:
\begin{eqnarray} \nonumber
\lefteqn{
\sum_{0\leq \lambda_j \leq N_j \atop j=1,\ldots, n}
\prod_{j=1}^n\frac{[g_j+\lambda_j+|\lambda|]}{[g_j]}
\prod_{1\leq i<j \leq n} \frac{[g_i-g_j+\lambda_i-\lambda_j]}
{[g_i-g_j]}  } &&
\\ \nonumber
&& \times \prod_{i,j=1}^n\frac{[g_i-g_j-N_j]_{\lambda_i}}
{[1+g_i-g_j]_{\lambda_i}}
\prod_{j=1}^n\frac{[g_j]_{|\lambda|}}
{[g_j+1+N_j]_{|\lambda|}}
\frac{[\beta,\gamma]_{|\lambda|}}
{[1-\delta, 1 -\epsilon]_{|\lambda|}}
\prod_{j=1}^n \frac{[\delta+g_j,
\epsilon +g_j]_{\lambda_j}}
{[g_j+1-\beta, g_j+1-\gamma]_{\lambda_j}}
\\ &&
=\frac{[1-\beta-\delta,1 -\gamma-\delta]_{|N|}}
{[1-\delta,1-\beta-\gamma-\delta]_{|N|}}
\prod_{j=1}^n\frac{[1+g_j,1+g_j-\beta-\gamma]_{N_j}}
{[g_j+1-\beta,g_j+1-\gamma]_{N_j}}.
\label{A_n-ell}\end{eqnarray}
Рассмотрим перестановку $(\lambda_k, g_k, N_k)$ с
$(\lambda_m, g_m, N_m)$ для некоторых фиксированных $k$ и $m$.
Имеется только два множителя, которые требуют некоторого рассмотрения.
Первый, заданный произведением $\prod_{1\leq i<j\leq n}$ оказывается
симметричным из-за равного числа изменений знака в числителе и знаменателе
(точно так же, как это было для формулы (\ref{C_n})). Произведение
$\prod_{i,j=1}^n$ оказывается симметричным, поскольку оно учитывает
все возможные значения обоих индексов $i$ и $j$. Таким образом,
согласно определению, данному в начале статьи, мы имеем дело с
эллиптическим гипергеометрическим рядом.

Докажем теперь полную эллиптичность функций $h_\ell$.
Эллиптичность по $\beta$ (или $\gamma$), или инвариантность (\ref{han})
при сдвигах $b\to pb, e\to p^{-1}e$, следует из сокращения множителей
$-b^{-1}q^{-|\lambda|}$, $-p^{-1}et_\ell q^{\lambda_\ell}$ и
$-q^{|\lambda|+1}e^{-1},$ $-pbt_\ell^{-1}q^{-\lambda_\ell-1}$,
порожденных числителем и знаменателем выражения в последней строке,
соответственно. Аналогичная ситуация имеет место для сдвигов
$d\to pd, e\to p^{-1}e$, что доказывает эллиптичность по $\delta$.

Существенно более сложное вычисление необходимо для проверки модулярной
инвариантности $c(\mathbf{\lambda})$ в (\ref{A_n-ell}). Рассмотрим
для этого разность квадратов аргументов эллиптических чисел в
числителе и знаменателе $h_\ell$. Первая дробь тета-функций дает
$n-1+2(|g|-g_\ell+n|\lambda|-\lambda_\ell)$,
где $|g|=g_1+\ldots+g_n.$ Вторая дробь дает $4(1+|\lambda|
+\lambda_\ell+g_\ell)$. Следующие два произведения генерируют
разность $n-1+2n(g_\ell+\lambda_\ell)-2|g|-2|\lambda|.$
Произведение $\prod_{j=1}^n$ отношений тета-функций дает вклад
$-2(n+|N|)(1+g_\ell+\lambda_\ell+|\lambda|).$ Остаток $h_\ell$ дает
\begin{eqnarray*}
&& (\beta+|\lambda|)^2+(\gamma+|\lambda|)^2+(\delta+g_\ell+\lambda_\ell)^2
+(\epsilon+g_\ell+\lambda_\ell)^2-(|\lambda|+1-\delta)^2
\\ &&
-(|\lambda|+1-\epsilon)^2
-(g_\ell+\lambda_\ell+1-\beta)^2
-(g_\ell+\lambda_\ell+1-\gamma)^2.
\end{eqnarray*}
Сумма всех этих громоздких выражений равна нулю, то есть
$h_\ell$ инвариантна относительно преобразований (\ref{modular-tr1}),
(\ref{modular-tr2}), порождающих полную модулярную груп\-пу.

Полная эллиптичность и модулярность выражения,
стоящего справа в равенстве (\ref{A_n-ell}), легко проверяются.
Действительно, инвариантность правой части (\ref{A_n}) при сдвигах
$t_j\to pt_j$ или $d\to pd$ (то есть, эллиптичность по $g_j$ или
$\delta$) очевидна. Сдвиг $b\to pb$ генерирует фактор $c^{|N|}$
из эллиптических сдвинутых факториалов $\theta(\ldots)_{|N|}$,
и он сокращается дополнительными множителями, происходящими из
эллиптических сдвинутых факториалов $\theta(\ldots)_{N_j}$.
По симметрии, это означает эллиптичность и по $\beta$ и по $\gamma$.

Что касается разности квадратов аргументов эллиптических чисел
в числителе и знаменателе правой части (\ref{A_n-ell}),
отношение эллиптических сдвинутых факториалов $[\ldots]_{|N|}$
дает (после упрощений) $-2|N|\beta\gamma$. Произведение отношений
$[\ldots]_{N_j}$ дает вклад
$$
\sum_{j=1}^nN_j((1+g_j)^2+(1+g_j-\beta-\gamma)^2
-(1+g_j-\beta)^2-(1+g_j-\gamma)^2),
$$
который сокращает предыдущее
выражение. Это доказывает заявленную модулярную инвариантность.

Обратимся теперь к $D_n$-формуле (\ref{D_n}). Отношения коэффициентов ряда
в левой части выглядят следующим образом
\begin{eqnarray}\nonumber
\lefteqn{
h_\ell(\mathbb{\lambda})=q^\ell
\prod_{j=1\atop j\neq\ell}^n\frac{\theta(t_jq^{\lambda_j+|\lambda|+1},
bcdt_j^{-1}q^{|\lambda|-\lambda_j};p)}
{\theta(t_jq^{\lambda_j+|\lambda|},
t_\ell t_jq^{\lambda_j+\lambda_\ell+1}/bcd;p)}\;
\frac{\theta(t_\ell q^{\lambda_\ell+|\lambda|+2};p)}
{\theta(t_\ell q^{\lambda_\ell+|\lambda|};p)} } &&
\\ \label{hdn} &&  \times
\prod_{i=1}^{\ell-1}
\frac{\theta(t_it_\ell^{-1}q^{\lambda_i-\lambda_\ell-1};p)}
{\theta(t_it_\ell^{-1}q^{\lambda_i-\lambda_\ell};p)}
\prod_{j=\ell+1}^n
\frac{\theta(t_\ell t_j^{-1}q^{\lambda_\ell-\lambda_j+1};p)}
{\theta(t_\ell t_j^{-1}q^{\lambda_\ell-\lambda_j};p)}
\\ \nonumber &&  \times
\prod_{j=1}^n
\frac{\theta(t_\ell t_jq^{1+N_j+\lambda_\ell}/bcd, t_\ell t_j^{-1}
q^{-N_j+\lambda_\ell}, t_jq^{|\lambda|};p)}
{\theta(t_\ell t_j^{-1}q^{\lambda_\ell+1}, t_j q^{1+N_j+|\lambda|},
bcdt_j^{-1} q^{|\lambda|-N_j};p)}
\frac{\theta(bq^{|\lambda|},cq^{|\lambda|},dq^{|\lambda|};p)}
{\theta(t_\ell q^{\lambda_\ell+1}/b,t_\ell q^{\lambda_\ell+1}/c,
t_\ell q^{\lambda_\ell+1}/d;p)}.
\end{eqnarray}
Сдвиги $q^{\lambda_k}\to pq^{\lambda_k}$ не меняют эти функции, то есть
они эллиптичны по переменным $\lambda_k$. Для $k<\ell$ или $k>\ell$,
это следует из сокращения множителей квазипериодичности
$$
q^{-n}, \quad -\frac{t_kt_\ell}{bcd} q^{1+\lambda_k+\lambda_\ell},\quad
\prod_{j=1\atop j\neq k,\ell}^n
\left(-\frac{t_j}{bcd}q^{\lambda_j-|\lambda|}\right),
$$
появляющихся из произведения $\prod_{j=1,j\neq\ell}^n$,
фактора $q^{-2}$, происходящего из последующего отношения
тета-функций, $q$ (из двух следующих произведений), и
$
\prod_{j=1}^n (-bcdt_j^{-1}q^{|\lambda|+1}),$ $ \qquad
-q^{-3|\lambda|}/bcd
$
(из выражений двух последних строк, соответственно). Для $k=\ell$
мы имеем множители $q^{1-n}$ (из первой дроби),
$\prod_{j=1,j\neq\ell}^n (t_j^2t_\ell
q^{2\lambda_j+\lambda_\ell-|\lambda|+1}/b^2c^2d^2) $
(из второй дроби), $q^{-4}$ (третья дробь), $q^{1-n}$
(следующие два произведения) и
$$
\prod_{j=1}^n\frac{b^2c^2d^2}{t_\ell t_j^2}q^{1+|\lambda|-\lambda_\ell},
\qquad \frac{t_\ell^3}{b^2c^2d^2}q^{3(1+\lambda_\ell-|\lambda|)}
$$
(из последних двух строк, соответственно). Их произведение равно 1,
то есть, функции $h_\ell$ эллиптичны по всем $\lambda_k$.

После переписывания формулы (\ref{D_n}) в терминах эллиптических чисел
\begin{eqnarray}\nonumber
\lefteqn{
\sum_{0\leq \lambda_j \leq N_j \atop j=1,\ldots, N_j}
\prod_{j=1}^n\frac{[g_j+\lambda_j+|\lambda|]}{[g_j]}
\prod_{1\leq i<j \leq n} \frac{[g_i-g_j+\lambda_i-\lambda_j]}
{[g_i-g_j]}  } &&
\\ \nonumber &&
\times \frac{\prod_{i,j=1}^n
[g_i+g_j+1+N_j-\beta-\gamma-\delta]_{\lambda_i}}
{\prod_{1\leq i<j \leq n}
[g_i+g_j+1-\beta-\gamma-\delta]_{\lambda_i+\lambda_j}}
\\ \nonumber &&
\times \prod_{i,j=1}^n\frac{[g_i-g_j-N_j]_{\lambda_i}}
{[1+g_i-g_j]_{\lambda_i}}
\prod_{j=1}^n\frac{[g_j]_{|\lambda|}
[\beta+\gamma+\delta-g_j]_{|\lambda|-\lambda_j}}
{[g_j+1+N_j,\beta+\gamma+\delta-g_j-N_j]_{|\lambda|}}
\\ \nonumber &&
\times \frac{[\beta,\gamma,\delta]_{|\lambda|}}
{\prod_{j=1}^n[g_j+1-\beta,g_j+1-\gamma,g_j+1-\delta]_{\lambda_j}}
\\  &&
=\prod_{j=1}^n
\frac{[g_j+1,g_j+1-\beta-\gamma,g_j+1-\beta-\delta,
g_j+1-\gamma-\delta]_{N_j}}
{[g_j+1-\beta-\gamma-\delta,g_j+1-\beta,g_j+1-\gamma,g_j+1-\delta]_{N_j}},
\label{D_n-ell}\end{eqnarray}
перестановочная симметрия между $(\lambda_k, g_k, N_k)$ и
$(\lambda_m, g_m, N_m)$, для фиксированных $k$ and $m$, становится
очевидной.

Эллиптичность по переменным $g_k,\ k\neq \ell,$ следует из сокращения
множителей в выражении (\ref{hdn}), возникающих
после замены $t_k\to pt_k$, в порядке появления  $q^{-1}$,
$-\frac{t_kt_\ell q^{1+\lambda_k+\lambda_\ell}}{bcd}$,
$-p^{-1}bcdq^{|\lambda|-\lambda_k}t_k^{-1}$, $q$ и, наконец,
$pq^{-1-|\lambda|-\lambda_\ell}t_\ell^{-1}$. Аналогично, для $g_\ell$
имеем множители (в порядке появления)
$$
\prod_{j=1,j\neq \ell}^n (-t_jt_\ell q^{1+\lambda_j+\lambda_\ell}/bcd),\;
q^{-2}, \; q^{1-n},\;
\prod_{j=1,j\neq \ell}^n (-bcdq^{-\lambda_\ell}/t_jt_\ell), \;
-bcdq^{-1-|\lambda|-2\lambda_\ell}/t_\ell^3, \;
$$
и $-t_\ell^3q^{3+3\lambda_\ell}/bcd $, произведение которых равно 1.
Эллиптичность по $\beta$ (или $\gamma,\delta$) подтверждается
аналогичным вычислением, которое мы опускаем для краткости.

Доказательство модулярности $h_\ell$ требует существенно более длинного
расчета. Пропустим начальное выражение для разности квадратов аргументов
эллиптических чисел в $h_\ell$. Первая цепочка упрощений в нем
приводит к сумме
\begin{eqnarray*} &&
2(1+|\lambda|+g_\ell+\lambda_\ell)(\lambda_\ell-|\lambda|-N_\ell+1)
+(2g_\ell-\beta-\gamma-\delta+N_\ell+1+\lambda_\ell)^2
\\ &&
-(g_\ell-\beta-\gamma-\delta+N_\ell-|\lambda|)^2
+\sum_{\rho=\beta,\gamma,\delta}\left(
(\rho+|\lambda|)^2-(g_\ell+\lambda_\ell+1-\rho)^2\right),
\end{eqnarray*}
которая равна нулю.

Работа с правой частью (\ref{D_n}) существенно проще.
Ее эллиптичность по параметрам следует из того, что отношение
$$
\frac{[g_j,g_j-\beta-\gamma,g_j-\beta-\delta,g_j-\gamma-\delta]}
{[g_j-\beta-\gamma-\delta,g_j-\beta,g_j-\gamma,g_j-\delta]}
$$
эллиптично по всем переменным. Модулярная инвариантность
гарантирована следующим равенством
\begin{eqnarray*} &&
g_j^2+(g_j -\beta-\gamma)^2+(g_j-\beta-\delta)^2
+(g_j-\gamma-\delta)^2
\\ &&
=(g_j-\beta-\gamma-\delta)^2+(g_j-\beta)^2+(g_j-\gamma)^2+(g_j-\delta)^2.
\end{eqnarray*}

Таким образом, мы доказали полную эллиптичность и модулярность
обеих сторон обеих $A_n$ и $D_n$ формул суммирования. \hfill{Q.E.D.}

\begin{corollary}
Формулы суммирования (\ref{C_n}) ($C_n$-сумма Варнаара), (\ref{A_n})
(эллиптическая  $A_n$-сумма Милна) и  (\ref{D_n}) (эллиптическая
$D_n$-сумма Бхатнагара-Шлоссера) справедливы при разложении по малому
параметру $\sigma$ до степени $\sigma^{12}$.
\end{corollary}
{\bf Доказательство}
Для произвольно маленького $\epsilon>0$, можно подобрать достаточно
маленькое абсолютное значение параметра $\sigma$ такое, что
обе стороны равенств (\ref{C_n}), (\ref{A_n}), (\ref{D_n})
будут голоморфны по модулярному параметру $\tau$ в области
$\mbox{Im}(\tau) > \epsilon$. Являясь $PSL(2,\mathbb{Z})$ инвариантными
функциями, они имеют разложения в ряды по $\sigma$ вида
$\sum_{k=0}^\infty m_k(\tau)\sigma^{2k}$, где $m_k(\tau)$ суть
модулярные формы веса $2k$. Для $\mbox{Im}(\tau)\to\infty$ (т.е. $p\to 0$),
разности выражений левых и правых частей этих равенств зануляются.
Поэтому, в разложении этих разностей по $\sigma$ могут появиться
только параболические формы. Поскольку не существует параболических
модулярных форм с весом меньшим чем 12, рассматриваемые
$C_n, A_n, D_n$ формулы суммирования справедливы в $\sigma$-разложении
до членов степени $\sigma^{12}$.
\hfill{Q.E.D.}

\smallskip

Полные рекурсивые доказательства $C_n$-формулы суммирования, предложенной
ван Диехеном и автором \cite{die-spi:modular},
и $A_n$, $D_n$ формул суммирования (\ref{A_n}), (\ref{D_n}) получены Розенгреном в
подробной статье \cite{ros:elliptic}. Эти формулы суммирования
выводятся также  с помощью анализа вычетов
в точных формулах интегрирования, доказанных в следующей главе.
Эти результаты усиливают аргументы в пользу гипотезы, высказанной выше, что
полностью эллиптические гипергеометрические ряды автоматически модулярно
инвариантны. Другие способы доказательства тождеств для многократных
эллиптических гипергеометрических рядов приведены в работах
\cite{ros-sch2,kaj-nou,ros-sch1}. Таким образом, практически все
известные формулы суммирования для многократных $_8\varphi_7$
обрывающихся $q$-гипергеометрических рядов имеют эллиптические аналоги
(одна из недоказанных еще формул будет приведена в конце главы 3).
Что касается бесконечных сумм для односторонних \cite{sch:nonterminating}
и двусторонних \cite{mil:multidimensional,gus:macdonald,olshaski}
$q$-гипергеометрических рядов, то их эллиптические аналоги находятся
пока вне досягаемости ввиду проблем с их сходимостью.

\section{Некоторые результаты для тета-функций на Римановых поверхностях}

На основании успешного построения теории эллиптических гипергеометрических
рядов, изложенной выше, естественно попытаться обобщить её на случай
Римановых поверхностей произвольного рода. Очевидно, что роль тета-функций
Якоби должны играть многомерные тета-функции Римана, ассоциированные
с алгебраическими кривыми. В этом случае мы можем взять в качестве
 $h(n)$ мероморфную функцию на Римановой поверхности $S$ произвольного
рода $g$. Однако, все становится намного менее явным в этой ситуации.
В данном параграфе излагаются результаты работ автора в этом направлении
\cite{spi:modularity,spi:macd}, в которых многопараметрическая сумма
Макдональда, впервые описанная в статье \cite{BM}, обобщена на Римановы поверхности.

Тета-функция Римана $g$ переменных $u_1,\ldots,u_g\in \mathbb{C}$
с характеристиками $\mathbf{\alpha}=(\alpha_1,\ldots,\alpha_g)$ и
$\mathbf{\beta}=(\beta_1,\ldots,$ $\beta_g)$ определяется $g$-кратным
рядом вида
\begin{equation}
\Theta_{\mathbf{\alpha},\mathbf{\beta}}(\mathbf{u};\Omega)
=\sum_{\mathbf{n}\in \mathbb{Z}^g}
\exp \{ \pi i\sum_{j,k=1}^g (n_j+\alpha_j)\Omega_{jk}(n_k+\alpha_k)
+2\pi i\sum_{j=1}^g(u_j+\beta_j)(n_j+\alpha_j) \},
\label{theta-R}\end{equation}
где $\Omega_{jk}$ симметрическая матрица периодов.
Для общих абелевых многообразий матрица $\Omega_{jk}$ произвольна,
но для Якобианов, ассоциированных с Римановыми поверхностями,
она генерируется базисом голоморфных дифференциалов
${\bf \omega}= (\omega_1,\ldots,\omega_g)$
(см., например, \cite{Mu}). Очевидно, что для произвольного $j$ мы имеем
\begin{equation}
\Theta_{\mathbf{\alpha},\mathbf{\beta}}(u_1,\ldots,u_j+1,\ldots, u_g;\Omega)=
e^{2\pi i \alpha_j}\Theta_{\mathbf{\alpha},\mathbf{\beta}}(\mathbf{u};\Omega).
\label{per-1}\end{equation}
Аналогично,
\begin{equation}
\Theta_{\mathbf{\alpha},\mathbf{\beta}}(u_1+\Omega_{1k},\ldots, u_g
+\Omega_{gk};\Omega)=
e^{-\pi i\Omega_{kk}-2\pi i(\beta_k+u_k)}\Theta_{\mathbf{\alpha},\mathbf{\beta}}
(\mathbf{u};\Omega),
\label{per-2}\end{equation}
где $k=1,\ldots, g$.

Обозначим как $\Gamma_{1,2}$ подгруппу симплектической модулярной группы
$Sp(2g,\mathbb{Z})$ генерируемой матрицами вида
$$
\gamma=\left(\begin{array}{cc} a & b \\ c & d \end{array}\right)
\in Sp(2g,\mathbb{Z}),
$$
таких что $diag(a^tb)=diag(c^td)=0\,\mbox{mod}\, 2$.
Действие этой группы на матрицу периодов $\Omega$ и аргументы тета-функции
$\mathbf{u}$ выглядит так
\begin{equation}
\Omega'=(a\Omega+b)(c\Omega+d)^{-1},\qquad
\mathbf{u}'=\mathbf{u}^t(a\Omega+b)^{-1}.
\label{sp}\end{equation}
Аналогично определяются преобразованные характеристики
$$
\left(\begin{array}{c} \mathbf{\alpha}' \\ \mathbf{\beta}' \end{array}\right)
= \left(\begin{array}{cc} d & -c \\ -b & a \end{array}\right)
\left(\begin{array}{c} \mathbf{\alpha} \\ \mathbf{\beta} \end{array}\right)
+\frac{1}{2}\left(\begin{array}{c} diag(c^td) \\ diag(a^td) \end{array}\right).
$$
Тогда, закон преобразования $\Gamma_{1,2}$ группы на тета-функциях имеет вид
\begin{equation}
\Theta_{\mathbf{\alpha}',\mathbf{\beta}'}(\mathbf{u}';\Omega')
=\zeta\sqrt{\det(c\Omega+d)}e^{\pi i \mathbf{u}^t(c\Omega+d)^{-1}c\mathbf{u}}
\Theta_{\mathbf{\alpha},\mathbf{\beta}}(\mathbf{u};\Omega),
\label{sp-mod}\end{equation}
где $\zeta$ есть некоторый восьмой корень единицы \cite{Mu}.

Обозначим
\begin{equation}
v_j(a,b)\equiv \int_a^b\omega_j,\qquad a,b\in S,
\label{ab-int}\end{equation}
абелевы интегралы первого рода. Характеристики $\alpha_j,
\beta_j\in (0, 1/2)$, такие что $4\sum_{j=1}^g\alpha_j\beta_j$
$=1\, (\mbox{mod}\; 2)$, называются нечетными. Обозначим тета-функции с произвольными
(несингулярыми) нечетными характеристиками $[{\bf u}],\; {\bf u}\in \mathbb{C}^g,$
и примем соглашение
$ [\mathbf{u}_1,\ldots,\mathbf{u}_k]= \prod_{j=1}^k[\mathbf{u}_j]. $
Для таких функций имеем $[-{\bf u}]=-[{\bf u}]$. Известно, что тета-функции на
Римановых поверхностях удовлетворяют тождеству Фэя \cite{fay}:
\begin{eqnarray} \nonumber
&& [{\bf u}+\mathbf{v}(a,c),{\bf u}+\mathbf{v}(b,d),
\mathbf{v}(c,b),\mathbf{v}(a,d)]
\\ \nonumber && \makebox[2em]{}
+[{\bf u}+\mathbf{v}(b,c),{\bf u}+\mathbf{v}(a,d),
\mathbf{v}(a,c),\mathbf{v}(b,d)]
\\ && \makebox[4em]{}
= [{\bf u},{\bf u}+\mathbf{v}(a,c)+\mathbf{v}(b,d),
\mathbf{v}(c,d),\mathbf{v}(a,b)],
\label{Fay}\end{eqnarray}
справедливому для произвольных ${\bf u}\in \mathbb{C}^g$ и $a,b,c,d\in S$
(при этом $\mathbf{v}(a,c)+\mathbf{v}(b,d)=
\mathbf{v}(b,c)+\mathbf{v}(a,d)$). Здесь мы можем заменить тета-функции, не
зависящие от $\bf u$, на главную форму $E(x,y)$ поскольку отношения этих
функций равны отношениям соответствующих главных форм и, поэтому, не зависят от
выбора нечетных характеристик. Однако более удобно работать с одной тета-функцией
$[{\bf u}]$. Для функции $[{\bf u}+\mathbf{v}(a,b)]$ с определенной
нормировкой
базиса $\omega_j$ преобразования \re{per-1} и \re{per-2} соответствуют циклическим
перемещениям точки $b$ (или $a$) на Римановой поверхности вдоль $A_i$ и $B_i$
контуров, таких что $\int_{A_i}\omega_j=\delta_{ij}$ и
$\int_{B_i}\omega_j=\Omega_{ij}$, соответственно \cite{Mu}.

Опишем простую формулу суммирования для некоторого ряда, построенного из
Римановых тета-функций.

\begin{theorem}
Для неотрицательного целого $n$ рассмотрим $n+1$ произвольных переменных
$\mathbf{z}_k\in\mathbb{C}^g$ и $4n+4$ точки на Римановой поверхности
$a_k, b_k, c_k, d_k\in S$, $k=0,\ldots, n$. Имеет место следующая многопараметрическая
формула суммирования для Римановых тета-функций на алгебраических кривых
\begin{eqnarray}\nonumber
&& \sum_{k=0}^n
[\mathbf{z}_k+\mathbf{v}(b_k,c_k),\mathbf{z}_k+\mathbf{v}(a_k,d_k),
\mathbf{v}(a_k,c_k),\mathbf{v}(b_k,d_k)]
\\ \nonumber && \makebox[2em]{} \times
\prod_{j=0}^{k-1}
[\mathbf{z}_j,\mathbf{z}_j+\mathbf{v}(a_j,c_j)+\mathbf{v}(b_j,d_j),
\mathbf{v}(c_j,d_j),\mathbf{v}(a_j,b_j)]
\\ \nonumber && \makebox[2em]{} \times
\prod_{j=k+1}^{n}
[\mathbf{z}_j+\mathbf{v}(a_j,c_j),\mathbf{z}_j+\mathbf{v}(b_j,d_j),
\mathbf{v}(c_j,b_j),\mathbf{v}(a_j,d_j)]
\\ \nonumber &&
=\prod_{k=0}^n
[\mathbf{z}_k,\mathbf{z}_k+\mathbf{v}(a_k,c_k)+\mathbf{v}(b_k,d_k),
\mathbf{v}(c_k,d_k),\mathbf{v}(a_k,b_k)]
\\ && \makebox[2em]{}
-\prod_{k=0}^n [\mathbf{z}_k+\mathbf{v}(a_k,c_k),
\mathbf{z}_k+\mathbf{v}(b_k,d_k),
\mathbf{v}(c_k,b_k),\mathbf{v}(a_k,d_k)].
\label{sum-R}\end{eqnarray}
\end{theorem}
{\bf Доказательство.}
Для доказательства равенства \re{sum-R}, мы обозначим $f_l^{(n)}$ и $f_r^{(n)}$
его левую и правую части, соответственно. Определим также
\begin{eqnarray*}
&& g_k=[\mathbf{z}_k,\mathbf{z}_k+\mathbf{v}(a_k,c_k)+\mathbf{v}(b_k,d_k),
\mathbf{v}(c_k,d_k),\mathbf{v}(a_k,b_k)], \\
&& h_k= [\mathbf{z}_k+\mathbf{v}(a_k,c_k),
\mathbf{z}_k+\mathbf{v}(b_k,d_k),
\mathbf{v}(c_k,b_k),\mathbf{v}(a_k,d_k)],
\end{eqnarray*}
так что $f_r^{(n)}=\prod_{k=0}^n g_k -\prod_{k=0}^n h_k$.

Для $n=0$, равенство \re{sum-R} сводится к тождеству \re{Fay}.
Предположим, что \re{sum-R} справедливо для $n=0,\ldots, N-1$,
$N>1$. Тогда по индукции мы имеем
\begin{eqnarray*}
&& f_l^{(N)}=h_Nf_l^{(N-1)}
+[\mathbf{z}_N+\mathbf{v}(b_N,c_N),\mathbf{z}_N+\mathbf{v}(a_N,d_N),
\mathbf{v}(a_N,c_N),\mathbf{v}(b_N,d_N)]
\prod_{j=0}^{N-1}g_j
\\  && \makebox[4em]{}
=\xi_N\prod_{k=0}^{N-1}g_k -\prod_{k=0}^N h_k,
\end{eqnarray*}
где
$
\xi_N=h_N + [\mathbf{z}_N+\mathbf{v}(b_N,c_N),
\mathbf{z}_N+\mathbf{v}(a_N,d_N),
\mathbf{v}(a_N,c_N),\mathbf{v}(b_N,d_N)].
$
Используя тождество Фэя \re{Fay}, находим $\xi_N=g_N$.
Поэтому, $f_l^{(N)}=\prod_{k=0}^Ng_k-\prod_{k=0}^Nh_k=f_r^{(N)}$,
то есть формула \re{sum-R} справедлива для произвольного $n$.
\hfill{Q.E.D.}

\begin{remark}
Данная теорема является частным случаем общего утверждения для телескопирующихся
сумм:
\begin{equation}
\sum_{k=0}^{n} (x_k-y_k) \prod_{j=0}^{k-1} x_j
\prod_{j=k+1}^{n} y_j = \prod_{j=0}^{n} x_j - \prod_{j=0}^{n} y_j,
\end{equation}
которое доказывается аналогичным образом
$$
\sum_{k=0}^{n} (x_k-y_k) \prod_{j=0}^{k-1} x_j \prod_{j=k+1}^{n} y_j
 = \sum_{k=0}^{n} \prod_{j=0}^{k} x_j \prod_{j=k+1}^{n} y_j
 - \sum_{k=0}^{n} \prod_{j=0}^{k-1} x_j \prod_{j=k}^{n} y_j
$$
с очевидным сокращением членов в суммах правой части равенства.
\end{remark}

Для эллиптических кривых (т.е. для $g=1$) формула \re{sum-R} была
доказана Варнааром \cite{war:summation}. Ее дальнейшая редукция на
тригонометрический уровень приводит к тождеству Макдональда, которое
было впервые опубликовано Бхатнагаром и Милном \cite{BM} и которое обобщает
соотношения полученные Чу в статье \cite{C}. Как показано в работах \cite{BM,war:summation},
равенства подобного типа при специальном выборе параметров могут
быть использованы для получения более тонких формул суммирования
и преобразования для бибазисных и эллиптических гипергеометрических
рядов.

Предположим, что для некоторого целого $N>0$ точки $a_N,b_N,c_N,d_N\in S$
такие, что мы упираемся в ноль тета-функции: $[\mathbf{v}(a_N,b_N)]=0$
или $[\mathbf{v}(c_N,d_N)]=0$. Тогда, используя антисимметрию
$[\mathbf{v}(c_k,b_k)]=-[\mathbf{v}(b_k,c_k)]$,
мы можем переписать равенство \re{sum-R} в виде
\begin{eqnarray}\nonumber
&& \sum_{k=0}^N
\frac{[\mathbf{z}_k+\mathbf{v}(b_k,c_k),\mathbf{z}_k+\mathbf{v}(a_k,d_k),
\mathbf{v}(a_k,c_k),\mathbf{v}(b_k,d_k)]}
{[\mathbf{z}_k+\mathbf{v}(a_k,c_k),
\mathbf{z}_k+\mathbf{v}(b_k,d_k),
\mathbf{v}(b_k,c_k),\mathbf{v}(a_k,d_k)]}
\\ && \makebox[2em]{} \times
\prod_{j=0}^{k-1}
\frac{[\mathbf{z}_j,\mathbf{z}_j+\mathbf{v}(a_j,c_j)+\mathbf{v}(b_j,d_j),
\mathbf{v}(c_j,d_j),\mathbf{v}(a_j,b_j)]}
{[\mathbf{z}_j+\mathbf{v}(a_j,c_j),
\mathbf{z}_j+\mathbf{v}(b_j,d_j),
\mathbf{v}(c_j,b_j),\mathbf{v}(a_j,d_j)]}
= 1.
\label{sum2}\end{eqnarray}

Эквивалентно, для $[\mathbf{v}(c_0,b_0)]=0$ или $[\mathbf{v}(a_0,d_0)]=0$
мы получаем сумму
\begin{eqnarray}\nonumber
&& \sum_{k=0}^n
\frac{
[\mathbf{z}_k+\mathbf{v}(b_k,c_k),\mathbf{z}_k+\mathbf{v}(a_k,d_k),
\mathbf{v}(a_k,c_k),\mathbf{v}(b_k,d_k)]}
{[\mathbf{z}_0,\mathbf{z}_0+\mathbf{v}(a_0,c_0)+\mathbf{v}(b_0,d_0),
\mathbf{v}(c_0,d_0),\mathbf{v}(a_0,b_0)]}
\\ \nonumber &&  \times
\prod_{j=0}^{k-1}\frac{
[\mathbf{z}_j,\mathbf{z}_j+\mathbf{v}(a_j,c_j)+\mathbf{v}(b_j,d_j),
\mathbf{v}(c_j,d_j),\mathbf{v}(a_j,b_j)]}
{[\mathbf{z}_{j+1}+\mathbf{v}(a_{j+1},c_{j+1}),
\mathbf{z}_{j+1}+\mathbf{v}(b_{j+1},d_{j+1}),
\mathbf{v}(c_{j+1},b_{j+1}),\mathbf{v}(a_{j+1},d_{j+1})]}
\\ &&
=\prod_{k=1}^n\frac{
[\mathbf{z}_k,\mathbf{z}_k+\mathbf{v}(a_k,c_k)+\mathbf{v}(b_k,d_k),
\mathbf{v}(c_k,d_k),\mathbf{v}(a_k,b_k)]}
{[\mathbf{z}_{k}+\mathbf{v}(a_{k},c_{k}),
\mathbf{z}_{k}+\mathbf{v}(b_{k},d_{k}),
\mathbf{v}(c_{k},b_{k}),\mathbf{v}(a_{k},d_{k})]}.
\label{sum3}\end{eqnarray}

В эллиптическом случае и его дальнейших вырождениях, такие равенства
полезны для поиска пар матриц обратных друг к другу, которые задаются
простыми явными аналитическими выражениями. Естественно ожидать, что соотношения
полученные выше приведут к обобщениям некоторых результатов Варнаара
по многобазисным гипергеометрическим рядам и обращениям матриц  \cite{war:summation}.

Специальным выбором параметров можно попытаться придать этим суммам
вид сумм гипергеометрического типа. Мы опишем только два таких выбора.
Подставляя
$$
\mathbf{z}_k=\mathbf{u}_0+\mathbf{v}(x_k), \quad
b_k=c_0\equiv x_0,\quad c_k\equiv x_k,\quad d_k=d_0,
$$
где $\mathbf{u}_0\in\mathbb{C}^g$,
$\mathbf{v}(x_k)\equiv\mathbf{v}(x_0,x_k)$, $k=0,1,\ldots,$ и
$\mathbf{u}_2\equiv\mathbf{v}(d_0,x_0),$ $\mathbf{u}_1^{(k)}\equiv\mathbf{v}(x_k,a_k)$
в (\ref{sum3}), мы получаем следующее тождество.

\begin{corollary}
\begin{eqnarray}\nonumber
\lefteqn{
\sum_{k=0}^n\frac{[{\bf u}_0+2\mathbf{v}(x_k)]}{[\mathbf{u}_0]}
\frac{[\mathbf{u}_0-\mathbf{u}_1^{(k)}-\mathbf{u}_2,\mathbf{u}_1^{(k)}]}
{[\mathbf{u}_0-\mathbf{u}_1^{(0)}-\mathbf{u}_2,\mathbf{u}_1^{(0)}]}
\prod_{j=0}^{k-1}\Biggl(\frac{[{\bf u}_0+\mathbf{v}(x_j)]}{[\mathbf{v}(x_{j+1})]} }
&& \\ \nonumber && \times
\frac{[{\bf u}_0-{\bf u}_1^{(j)}-{\bf u}_2+\mathbf{v}(x_j),
{\bf u}_1^{(j)}+\mathbf{v}(x_j),{\bf u}_2+\mathbf{v}(x_j)]}
{[{\bf u}_1^{(j+1)}+{\bf u}_2+\mathbf{v}(x_{j+1}),
{\bf u}_0-{\bf u}_1^{(j+1)}+\mathbf{v}(x_{j+1}),
{\bf u}_0-{\bf u}_2+\mathbf{v}(x_{j+1})]}\Biggr)
\\ &&
= \prod_{k=1}^{n}\frac
{[{\bf u}_0+\mathbf{v}(x_k),
{\bf u}_0-{\bf u}_1^{(k)}-{\bf u}_2+\mathbf{v}(x_k),
{\bf u}_1^{(k)}+\mathbf{v}(x_k),{\bf u}_2+\mathbf{v}(x_k)]}
{[\mathbf{v}(x_k),{\bf u}_1^{(k)}+{\bf u}_2+\mathbf{v}(x_k),
{\bf u}_0-{\bf u}_1^{(k)}+\mathbf{v}(x_k),
{\bf u}_0-{\bf u}_2+\mathbf{v}(x_k)]}.
\label{theta-sum}\end{eqnarray}
\end{corollary}

Если бы существовала последовательность точек $a_k$ таких что
$\mathbf{u}_1^{(k)}=\mathbf{u}_1^{(0)}$, тогда четыре тета-функции
во втором множителе в левой части ({\ref{theta-sum}) сократили бы
друг друга и мы получили бы точный $g>1$ аналог суммы (\ref{8V7}).
Такое условие неявным образом подразумевалось в статье автора \cite{spi:modularity}
при выводе соответствующей формулы суммирования. Однако, в общем случае
это возможно только для $g=1$ (при выборе
$\omega= du, \, a_k=k+a_0, x_k=k+x_0$). Поэтому, правильная $g>1$ формула
суммирования \re{theta-sum}, опубликованная в работе \cite{spi:macd},
не обладает всеми структурными свойствами $g=1$ суммы.
Тем не менее, обратим внимание на тот факт, что правая часть
({\ref{theta-sum}) удовлетворяет аналогам условий балансировки и
вполне-уравновешенности: 1) суммы аргументов тета-функций
в числителе и знаменателе равны друг другу;
2) тета-функции в числителе и знаменателе, стоящие друг над/под другом,
имеют аргументы, суммы которых равны  ${\bf u}_0+2\mathbf{v}(x_k)$.

Другая любопытная формула суммирования получается после подстановки
$a_k=a,$ $c_k=b_0=c,$ $b_{k+1}=d_k\equiv x_k$ при $k=0,1,\ldots,$
и ${\bf z}_{k+1}={\bf z}_0+{\bf v}(c,x_k)$ при $k>0$ в (\ref{sum3}).
\begin{corollary}
\begin{eqnarray}\nonumber
\lefteqn{ \sum_{k=1}^{n}
\frac{[{\bf z_0},{\bf z}_0+{\bf v}(c,x_{k-1})+{\bf v}(a,x_k),
{\bf v}(a,c),{\bf v}(x_{k-1},x_k)]}
{[{\bf z}_0+{\bf v}(c,x_{k-1}),{\bf z}_0+{\bf v}(c,x_k),
{\bf v}(a,x_{k-1}),{\bf v}(a,x_k)]} } &&
\\ &&
 = \frac{[{\bf z}_0+{\bf v}(a,x_n),{\bf v}(c,x_n)]}
{[{\bf z}_0+{\bf v}(c,x_n),{\bf v}(a,x_n)]}
 - \frac{[{\bf z}_0+{\bf v}(a,x_0),{\bf v}(c,x_0)]}
{[{\bf z}_0+{\bf v}(c,x_0),{\bf v}(a,x_0)]}.
\label{sum4}\end{eqnarray}
\end{corollary}

Легко видеть, что все полученные суммы представляют собой ``полностью абелевы"
функции, то есть они инвариантны относительно произвольных перемещений
точек на Римановой поверхности вдоль циклов и соответствующих $2g$ сдвигов
переменных $\mathbf{z}_k$ (или $\mathbf{u}_0$). Аналогично, по построению
очевидна инвариантность относительно преобразований $Sp(2g,\mathbb{Z})$
модулярной группы (благодаря свойствам тождества Фэя).
Как видно из приведенного анализа, тета-гипергеометрические ряды
для Римановых поверхностей рода $g>1$ должны
обладать принципиально новыми свойствами по сравнению с $g=1$ случаем,
точная формулировка которых требует дополнительных исследований.

%% file: CHAPTER3.TEX
\chapter[Тета-гипергеометрические интегралы]
{Тета-гипергеометрические интегралы}

\section{Общее определение}

В соответствии с теорией общих тета-гипергеометрических рядов, представленной
в предыдущей главе, мы вводим аналогичную концепцию для интегралов.

\begin{definition}
Обозначим $\cal D$ некоторый $n$-мерный цикл в комплексном пространстве $\C^n$.
Пусть $\Delta(y_1,\ldots,$ $y_n)$ будет мероморфной функцией
аргументов $y_1,\ldots,y_n$. Рассмотрим многократные интегралы
\begin{equation}\label{definition}
I_n=\int_{\cal D} \; \Delta(y_1,\ldots,y_n)dy_1\cdots dy_n
\end{equation}
и отношения
\begin{equation}
h_\ell(\mathbf{y})=\frac{\Delta(y_1,\ldots,y_\ell+1,\ldots,y_n)
}{\Delta(y_1,\ldots,y_\ell,\ldots,y_n)}.
\label{hl}\end{equation}

Интегралы $I_n$ называются:

1)  чисто гипергеометрическими интегралами, если
\begin{equation}
h_\ell(\mathbf{y})=R_\ell(\mathbf{y})
\label{plain-int}\end{equation}
являются рациональными функциями $y_1,\ldots, y_n$ для всех
$\ell=1,\ldots, n;$

2)  $q$-гипергеометрическими интегралами, если
\begin{equation}
h_\ell(\mathbf{y})=R_\ell(q^\mathbf{y})
\label{q-int}\end{equation}
являются рациональными функциями $q^{y_1},\ldots, q^{y_n}$, $q\in\mathbb{C},$
для всех $\ell=1,\ldots, n;$

3)  эллиптическими гипергеометрическими интегралами, если
для всех $\ell=1,\ldots,n$ отношения $h_\ell(\mathbf{y})$
являются эллиптическими функциями переменных $y_1,\ldots,y_n$
с периодами $\sigma^{-1}$ и $\tau\sigma^{-1}, \, \mbox{Im}(\tau)>0$;

4) общими тета-гипергеометрическими интегралами, если $h_\ell(\mathbf{y})$
и $1/h_\ell(\mathbf{y})$ суть мероморфные функции, удовлетворяющие условиям
квазипериодичности
\begin{eqnarray}\nonumber
&& h_\ell(y_1,\ldots,y_k+\sigma^{-1},\ldots,y_n)
=e^{\sum_{j=1}^na_{\ell k}(j)y_j+b_{\ell k}} h_\ell(\mathbf{y}),
\\ && h_\ell(y_1,\ldots,y_k+\tau\sigma^{-1},\ldots,y_n)
=e^{\sum_{j=1}^n c_{\ell k}(j)y_j+d_{\ell k}} h_\ell(\mathbf{y}),
\label{quasi-per}\end{eqnarray}
с множителями квазипериодичности подобными множителям для сигма-функции
Вейерштрасса (связанной с $\theta_1(u)$ простым образом \cite{whi-wat:course}).
\end{definition}

Если считать переменные $y_1,\ldots,y_n$ дискретными, например,
$\mathbf{y}\in\mathbb{N}^n$, и заменить интегралы суммами
$\sum_{\mathbf{y}\in\mathbb{N}^n}$, то мы получим определения
обычных и $q$-гипергеометрических, эллиптических и тета-гипергеометрических рядов.
При $a_{\ell k}(j)=b_{\ell k}=c_{\ell k}(j)=d_{\ell k}=0$,
тета-гипергеометрические интегралы сводятся к эллиптическим.
Интегралы (или ряды), определенные так, не образуют алгебру поскольку
сумма двух гипергеометрических интегралов, вообще говоря,
не согласуется со взятыми определениями.

Сдвиги $y_\ell\to y_\ell+1$ в (\ref{hl}) могут быть заменены сдвигами на
произвольную константу $y_\ell\to y_\ell +\omega_1,$
$\omega_1=const.$ Однако, можно сделать $\omega_1$ равным единице простым
масштабированием $y_\ell$, что приводит к простой деформации цикла
интегрирования ${\cal D}$ в (\ref{definition}). Ниже будут использоваться
в основном циклы ${\cal D}=C^n$, где $C$ есть некоторый
контур на комплексной плоскости.

Рассмотрим детально случай $n=1$. Общая рациональная функция
$y$ может быть представлена в виде
$$
R(y)=\frac{\prod_{j=1}^n(1-a_j+y)\prod_{j=n+1}^r(a_j-1-y)}
{\prod_{j=1}^m(b_j-1-y)\prod_{j=m+1}^s(1-b_j+y)}\, x,
$$
где $n,r,m,s\in\N$, $x\in\C$ и $a_j, b_j$ описывает положения нулей и полюсов
функции $R(y)$. Общее решение уравнения $\Delta(y+1)=R(y)\Delta(y)$
имеет вид
\begin{equation}
\Delta(y)=\frac{\prod_{j=1}^m\Gamma(b_j-y)
\prod_{j=1}^n\Gamma(1-a_j+y)}{\prod_{j=m+1}^s\Gamma(1-b_j+y)
\prod_{j=n+1}^r\Gamma(a_j-y)}\, x^{y}\varphi(y),
\label{mej}\end{equation}
где $\Gamma(y)$ обозначает стандартную гамма-функцию, а
$\varphi(y)$ произвольную периодическую функцию, $\varphi(y+1)=\varphi(y)$.
Если положить $\varphi(y)=1$, то тогда для подходящего выбора контура $C$
интеграл  $I_1$ в (\ref{definition}) определяет не что иное как
функцию Мейера \cite{erd:higher}. В рассмотренной ситуации не имеется естественных
дополнительных ограничений, позволяющих зафиксировать бесконечномерную
(функциональную) свободу фигурирующую в решении $\Delta(y)$.

В $q$-случае, аналогичным образом можно записать
$$
R(q^y)=\frac{\prod_{j=1}^n(1-t_jq^y)\prod_{j=n+1}^r(1-t_j^{-1}q^{-y})}
{\prod_{j=1}^m(1-w_j^{-1}q^{-y})\prod_{j=m+1}^s(1-w_jq^y)}\, x.
$$
Для $0<|q|<1$, общее мероморфное решение уравнения
$\Delta(y+1)=R(q^y)\Delta(y)$ есть
\begin{equation}
\Delta(y)= \frac{\prod_{j=n+1}^r (t_j^{-1}q^{1-y};q)_\infty
\prod_{j=m+1}^s(w_jq^y;q)_\infty}
{\prod_{j=1}^n(t_jq^y;q)_\infty
\prod_{j=1}^m(w_j^{-1}q^{1-y};q)_\infty}\, x^{y}\varphi(y),
\label{q-mej}\end{equation}
где $\varphi(y)$ опять обозначает произвольную периодическую функцию,
$\varphi(y+1)=\varphi(y)$. В этом случае, для $\varphi(y)=1$
и подходящего контура $C$ интеграл $I_1$ описывает $q$-аналог функции Мейера,
изучавшийся Слэтер в монографии \cite{sla:generalized}.

При $|q|>1$, уравнение $\Delta(y+1)=R(q^y)\Delta(y)$ имеет следующее
общее решение
$$
\Delta(y)= \frac{\prod_{j=1}^n(t_jq^{y-1};q^{-1})_\infty
\prod_{j=1}^m(w_j^{-1}q^{-y};q^{-1})_\infty}{\prod_{j=n+1}^r
(t_j^{-1}q^{-y};q^{-1})_\infty \prod_{j=m+1}^s(w_jq^{y-1};q^{-1})_\infty}
\, x^{y}\varphi(y),
$$
т.е. мы имеем эффективную замену $q\to q^{-1}$ и переопределение параметров
в (\ref{q-mej}).

Напомним, что $q=e^{2\pi i \sigma}$. Параметр $\sigma$
дает второй масштаб, который может использоваться для наложения
естественного ограничения на вид $\Delta(y)$. Функция $q^y$ периодична
при сдвигах $y\to y+\sigma^{-1}$ и, поэтому, функция (\ref{q-mej})
удовлетворяет уравнению
$\Delta(y+\sigma^{-1})/\Delta(y)=x^{1/\sigma}\varphi(y+\sigma^{-1})/
\varphi(y)$. Мы можем определить $\varphi(y)$ исходя из дополнительного
требования
$$
\varphi(y+\sigma^{-1})=\tilde R(e^{2\pi i y})\varphi(y),
$$
где $\tilde R$ обозначает другую рациональную функцию $y$. В согласии с условием
периодичности $\varphi(y+1)=\varphi(y)$, имеем
$$
\tilde R(e^{2\pi i y})=
\frac{\prod_{j=1}^{n'}(1-\tilde t_je^{-2\pi iy})
\prod_{j=n'+1}^{r'}(1-\tilde t_j^{-1}e^{2\pi iy})}
{\prod_{j=1}^{m'}(1-\tilde w_j^{-1}e^{2\pi iy})
\prod_{j=m'+1}^{s'}(1-\tilde w_je^{-2\pi iy})},
$$
где $\tilde t_j$ и $\tilde w_j$ обозначают новые произвольные
параметры. Отметим, что мы не можем умножить функцию $\tilde R$ на выражения
вида $\rho e^{2\pi i k y}, k\in\mathbb{Z}, \rho\in\mathbb{C},$ если они отличны
от единицы, поскольку тогда условие периодичности $\varphi(y)$
будет нарушено. При $|q|<1$, общее мероморфное решение разностного
уравнения для $\varphi(y)$ есть
\begin{equation}
\varphi(y)= \frac{\prod_{j=m'+1}^{s'}(\tilde w_je^{-2\pi iy};
\tilde q)_\infty \prod_{j=n'+1}^{r'}(\tilde q\tilde t_j^{-1}e^{2\pi iy};
\tilde q)_\infty} {\prod_{j=1}^{m'} (\tilde q\tilde w_j^{-1}e^{2\pi iy};
\tilde q)_\infty \prod_{j=1}^{n'}(\tilde t_je^{-2\pi iy};
\tilde q)_\infty}\, \tilde \varphi(y),
\label{tq-mej}\end{equation}
где $\tilde q=e^{-2\pi i/\sigma}$ обозначает модулярно преобразованный параметр $q$.
Действительно, при $\mbox{Im}(\sigma)>0$ мы имеем  $\mbox{Im}(\sigma^{-1})<0$,
так что функция (\ref{tq-mej}) хорошо определена. Функция $\tilde\varphi(y)$
в выражении (\ref{tq-mej}) обозначает произвольную эллиптическую функцию с периодами
$1$ и $\sigma^{-1}$. Она однозначно характеризуется
положением полюсов и нулей в фундаментальном параллелограмме периодов и имеет
$2k-1$ свободных параметров, где $k\in\N$ обозначает порядок $\tilde \varphi(y)$.
Таким образом, пространство решений уже не такое широкое: оно становится
конечномерным (в смысле числа свободных параметров).

Рассмотрим теперь случай $|q|=1$. Обозначим $\sigma=\omega_1/\omega_2$,
$u\equiv y\omega_1$ и будем считать что $\mbox{Re}(\sigma)>0$.
Теперь мы можем выбрать $t_j, \tilde t_j$ и другие параметры специальным
образом, таким что бесконечные произведения $(t_jq^y;q)_\infty$,
$(\tilde t_je^{-2\pi i y};\tilde q)_\infty$, и т.д. в выражениях (\ref{q-mej}) и
(\ref{tq-mej}) комбинируются в функцию
$S(u+g_j;\omega_1,\omega_2)$ для некоторых $g_j$, где
\begin{equation}
S(u;\omega_1,\omega_2) = \frac{(e^{2\pi i u/\omega_2}; q)_\infty}
{(e^{2\pi iu/\omega_1}\tilde q; \tilde q)_\infty},
\label{2d-sin}\end{equation}
есть функция двойного синуса хорошо определенная в пределе $|q|\to 1$.
Действительно, можно проверить, что нули и полюсы функции
(\ref{2d-sin}) совпадают с нулями и полюсами отношения двойных
гамма-функций Барнса $\Gamma_2(\omega_1+\omega_2-u;\mathbf{\omega})/\Gamma_2(u;
\mathbf{\omega})$, задающем мероморфную функцию $u$
при $\omega_1/\omega_2>0$. Теория многократной гамма-функции Барнса изложена в
работах \cite{bar:theory,bar:multiple} и не рассматривается здесь.

Поскольку $\sigma$ действительно, то при несоизмеримости этого параметра с 1
имеем $\tilde\varphi(y)=const$ (т.е. функция $\Delta(y)$ определяется вполне
уникально). Свойства функции двойного синуса и некоторые ее приложения описаны,
например, в работах \cite{kur:multiple,jim-miw:quantum,fkv:strongly,
kls:unitary,nis-uen:integral,rui:first,rui:generalized,volk}.
В частности, интегралы, введенные Джимбо и Мивой в статье \cite{jim-miw:quantum}
как решения некоторых $q$-разностных уравнений
при $|q|=1$, послужили первыми примерами $q$-гипергеометрических интегралов
для $q$ на единичной окружности. Концепция модулярного дубля Фаддеева
для квантовых групп \cite{fad:discrete,fad:modular} так же связана с
функцией (\ref{2d-sin}).

Таким образом, мир $q$-аналогов функции Мейера оказался
значительно богаче стандартного гипергеометрического.
Введение дополнительного уравнения, использующего сдвиги по $\sigma^{-1}$,
внесло новые нетривиальные элементы в структуру интегралов и сузило
функциональную свободу в определении мероморфной функции $\Delta(y)$
к эллиптической функции $\tilde\varphi(y)$, содержащей конечное число
свободных параметров.

Обратимся теперь к однократным эллиптическим гипергеометрическим интегралам.
Мы берем две комплексные переменные $q, p\in\C$, удовлетворяющие
ограничениям $|q|,|p|<1$. Ключевая тета-функция Якоби имеет вид
$$
\theta(z;p)=(z;p)_\infty(pz^{-1};p)_\infty,
$$
где $(a;p)_\infty=\prod_{j=0}^\infty(1-ap^j)$.
Ее основные свойства симметрии отражены в преобразованиях
$$
\theta(z^{-1};p)=\theta(pz;p)=-z^{-1}\theta(z;p).
$$
Стандартная эллиптическая гамма-функция определяется с помощью
двойного бесконечного произведения \cite{rui:first}:
\begin{equation}
\Gamma(z;q,p) = \prod_{j,k=0}^\infty\frac{1-z^{-1}q^{j+1}p^{k+1}}
{1-zq^jp^k}.
\label{ell-gamma}\end{equation}
Она удовлетворяет двум разностным уравнениям первого порядка
$$
\Gamma(qz;q,p)=\theta(z;p)\Gamma(z;q,p),\qquad
\Gamma(pz;q,p)=\theta(z;q)\Gamma(z;q,p).
$$
Удобно опускать базисные переменные $q$ и $p$ в обозначении $\Gamma(z;q,p)$
и использовать соглашения
$$
\Gamma(tz^\pm)=\Gamma(tz,tz^{-1}),\quad
\Gamma(z^{\pm 2})=\Gamma(z^2,z^{-2}),
\quad \Gamma(t_1,\ldots,t_m)\equiv\prod_{r=1}^m \Gamma(t_r;q,p).
$$
Аналогично, для тета-функций мы подразумеваем
$$
\theta(t_1,\ldots,t_m;p)=\theta(t_1;p)\cdots\theta(t_m;p),
\quad \theta(tz^\pm;p)=\theta(tz,tz^{-1};p).
$$

Удовлетворительное с современной точки зрения описание эллиптической
гамма-функ\-ции как независимой фундаментальной специальной функции
впервые было дано в работе Рюйсенаарса \cite{rui:first}.
Раннее функции такого типа были рассмотрены Джексоном
\cite{jac:basic}, но его результаты не привлекли к себе внимания.
В несколько неявном виде функция \re{ell-gamma} появилась в знаменитой
статье Бакстера по решению восьмивершинной модели \cite{bax:partition}
и более поздних исследованиях точно решаемых моделей статистической
механики. После рассмотрения \cite{rui:first}, систематическое
исследование функции \re{ell-gamma}  было продолжено Фельдером и
Варченко \cite{fel-var:elliptic}. Другой тип эллиптической
гамма-функции, который будет описан ниже (см. также главу 5),
был найден автором в работе \cite{spi:integrals}.

Общая эллиптическая функция порядка $r+1$ может быть представлена в виде
\cite{whi-wat:course}:
\be
h(y)=e^\gamma\prod_{j=0}^r\frac{\theta_1(u_j+y;\sigma,\tau)}
{\theta_1(v_j+y;\sigma, \tau)}=e^\gamma
\frac{\theta(t_0q^y,\ldots,t_rq^y ;p)}
{\theta(w_0q^y,\ldots,w_rq^y;p)},
\label{n=1} \ee
где $p=e^{2\pi i\tau }, \; \mbox{Im}(\tau)>0, \; q=e^{2\pi i\sigma}$.
Параметр $\gamma\in\C$ произволен, но $t_i\equiv q^{u_i}$ и $w_i\equiv q^{v_i}$
удовлетворяют условию балансировки
\begin{equation}
\sum_{i=0}^r(u_i-v_i)=0, \quad \mbox{или}\quad
\prod_{i=0}^rt_i=\prod_{i=0}^rw_i,
\label{balance}\end{equation}
гарантирующему двоякую периодичность мероморфной функции $h(y)$:
$$
h(y+\sigma^{-1})=h(y),\qquad   h(y+\tau\sigma^{-1})=h(y).
$$
Заметим, что при $\tau=\sigma$ (что требует $\mbox{Im}(\sigma)> 0$), функция $h(y)$
задает явный вид множителя $\tilde\varphi(y)$ в выражении (\ref{tq-mej}).

Для того, чтобы найти подынтегральную функцию $\Delta(y)$, необходимо
решить разностное уравнение первого порядка
\begin{equation}
\Delta(y+1)=h(y)\Delta(y)
\label{1-eq}\end{equation}
в классе мероморфных функций. Теория уравнений такого типа
была разработана достаточно давно (см., например, \cite{bar:linear}).
Поскольку $h(y)$ факторизуется на отношение произведений тета-функций,
достаточно найти мероморфное решение уравнения
\begin{equation}
f(y+1)=\theta(q^y;p)f(y),
\label{ell-g-eq}\end{equation}
приводящее к различным эллиптическим гамма-функциям.
Простейшая функция такого типа (\ref{ell-gamma}) определяется из уравнения
(\ref{ell-g-eq}) только с точностью до периодической функции
$\varphi(y+1)=\varphi(y)$ и, более того, она требует
$\mbox{Im}(\sigma)>0$ (или $|q|<1$), что не подразумевалось в функции (\ref{n=1}).

Введем переменную $z=q^y$, так что сдвиг $y\to y+1$ становится эквивалентным
умножению $z\to qz$. Тогда общее решение уравнения (\ref{1-eq}) выглядит так:
\begin{equation}
\Delta(y)=\prod_{j=0}^r\frac{\Gamma(t_jz;q,p)}
{\Gamma(w_jz;q,p)} e^{\gamma y}\varphi(y),
\label{DE}\end{equation}
где подразумевается справедливость условия балансировки (\ref{balance})
и $\varphi(y)$ обозначает произвольную периодическую функцию,
$\varphi(y+1)=\varphi(y)$. Используя формулы отражения
\be
\Gamma(pz, qz^{-1};q,p)=\Gamma(qz,pz^{-1};q,p)
=\Gamma(pqz,z^{-1};q,p)=\Gamma(z,pqz^{-1};q,p)=1,
\label{refl-eq}\ee
в выражении (\ref{DE}) можно заменить некоторые эллиптические гамма-функции с аргументами,
содержащими $z$, на функции с аргументами $\propto z^{-1}$. После этого,
$\Delta(y)$ будет похожей на подынтегральные функции в определении
функции Мейера и ее $q$-аналогов, однако в эллиптическом случае
это не приводит к каким-либо обобщениям, поскольку правая часть в
равенствах (\ref{refl-eq}) тривиальна.

В области $\mbox{Im}(\sigma)<0$, т.е. $|q|>1$,
общее решение (\ref{1-eq}) может быть записано в виде
\begin{equation}
\Delta(y)=\prod_{j=0}^r\frac{\Gamma(w_jq^{y-1};q^{-1},p)}
{\Gamma(t_jq^{y-1};q^{-1},p)} e^{\gamma y}\varphi(y).
\label{D'}\end{equation}
Эффективно мы имеем простое переопределение параметров и замену $q\to q^{-1}$
в эллиптических гамма-функциях в (\ref{DE}) (ср. с определением этой
функции при $|q|>1$, данном в работе \cite{fel-var:elliptic}).

Положим $\varphi(y)=1$. Тогда функция (\ref{DE}) удовлетворяет двум простым
разностным уравнениям первого порядка:
\begin{eqnarray} \label{2-eq}
&&\Delta(y+\sigma^{-1})=e^{\gamma/\sigma}\Delta(y),  \\
&&\Delta(y+\tau\sigma^{-1})=e^{\gamma\tau/\sigma}
\prod_{j=0}^r \frac{\theta(t_jq^y;q)}{\theta(w_jq^y;q)}\Delta(y).
\label{3-eq}\end{eqnarray}
Предположим, что $1, \sigma^{-1}, \tau\sigma^{-1}$ попарно несоизмеримы.
Тогда система трех уравнений (\ref{1-eq}), (\ref{2-eq}), и
(\ref{3-eq}) однозначно определяет $\Delta(y)$ с точностью до численного
множителя. Так же как и в $q$-гипергеометрическом случае, мы можем обобщить уравнения
(\ref{2-eq}) и (\ref{3-eq}), использовать их как естественные средства для фиксации
функциональной свободы в $\Delta(y)$ и получить таким способом качественно
различные эллиптические гипергеометрические интегралы.

Отношение $\Delta(y+\tau\sigma^{-1})/\Delta(y)$ в (\ref{3-eq})
оказывается эллиптической функцией с периодами 1 и $\sigma^{-1}$.
Поэтому, естественно потребовать чтобы отношение $\Delta(y+\sigma^{-1})/\Delta(y)$
так же было эллиптической функцией, периоды которой по симметрии равны
1 и $\tau\sigma^{-1}$.

\begin{theorem}
Предположим, что $\Delta(y)$ удовлетворяет уравнению (\ref{1-eq}) и
что 1, $\sigma^{-1}$, $\tau\sigma^{-1}$ попарно несоизмеримы.
Обозначим $\tilde q=e^{-2\pi i/\sigma}$ и ${r} =e^{2\pi i\tau/\sigma}$
и для простоты ограничимся случаем $\mbox{Im}(\sigma)>0$ (т.е., $|q|<1$).
Если отношение $\Delta(y+\tau\sigma^{-1})/\Delta(y)$ равно эллиптической функции
с периодами 1 и $\tau\sigma^{-1}$, то для
$\mbox{Im}(\tau/\sigma)>0$ наиболее общая мероморфная функция $\Delta(y)$ имеет вид
\begin{equation}\label{delta-ell}
\Delta(y)=\prod_{j=0}^r\frac{\Gamma(t_jq^y;q,p)}{\Gamma(w_jq^y;q,p)}
\prod_{j=0}^{n}\frac{\Gamma(\tilde t_je^{-2\pi iy};\tilde q,{r} )}
{\Gamma(\tilde w_je^{-2\pi iy};\tilde q,{r} )}\,e^{\gamma y+\delta},
\end{equation}
где $\prod_{j=0}^rt_jw_j^{-1}=\prod_{j=0}^{n}\tilde t_j\tilde w_j^{-1}=1$.
При $\mbox{Im}(\tau/\sigma)<0$, имеем
\begin{equation}\label{delta-ell'}
\Delta(y)=\prod_{j=0}^r\frac{\Gamma(t_jq^y;q,p)}{\Gamma(w_jq^y;q,p)}
\prod_{j=0}^{n}\frac{\Gamma(\tilde w_je^{-2\pi iy}{r} ^{-1};
\tilde q,{r} ^{-1})} {\Gamma(\tilde t_je^{-2\pi iy}{r} ^{-1};
\tilde q,{r} ^{-1})}\,e^{\gamma y+\delta}.
\end{equation}
\end{theorem}
{\bf Доказательство.}
Прежде всего отметим, что при $\mbox{Im}(\sigma)>0$ автоматически
выполняется неравенство $\mbox{Im}(\sigma^{-1})<0$, т.е. $|\tilde q|<1.$
Поэтому, при $\mbox{Im}(\tau/\sigma)>0$ функция
$\Gamma(z;\tilde q,{r} )$ хорошо определена.

Функция $\Delta(y)$ (\ref{DE}) дает общее решение уравнения
(\ref{1-eq}). Предположим, что отношение $\Delta(y+\sigma^{-1})/\Delta(y)$
задает эллиптическую функцию порядка $n+1$ с периодами 1 и $\tau\sigma^{-1}$.
При $\mbox{Im}(\tau/\sigma)>0$, это требование эквивалентно следующему
уравнению на $\varphi(y)$:
\begin{equation}
\frac{\varphi(y+\sigma^{-1})}{\varphi(y)}=
\prod_{j=0}^{n}\frac{\theta(\tilde t_je^{-2\pi iy}; {r} )}
{\theta(\tilde w_je^{-2\pi iy}; {r} )},
\label{phi-eq}\end{equation}
где $\prod_{j=0}^{n}\tilde t_j\tilde w_j^{-1}=1$.
Заметим, что правую часть этого уравнения нельзя домножать даже
на произвольную константу отличную от единицы так как это нарушит
условие что его решение периодично, $\varphi(y+1)=\varphi(y)$.

Мероморфное решение (\ref{phi-eq}) имеет вид
\begin{equation}
\varphi(y)=\prod_{j=0}^{n}\frac{\Gamma(\tilde t_je^{-2\pi iy};\tilde q,
{r} )}{\Gamma(\tilde w_je^{-2\pi iy};\tilde q,{r} )}\,
\tilde \varphi(y),
\label{phi-form}\end{equation}
где $\tilde \varphi(y)$ обозначает эллиптическую функцию с периодами
1 и $\sigma^{-1}$, т.е.
$$
\tilde\varphi(y)=\prod_{j=1}^m\frac{\theta(a_je^{-2\pi iy};\tilde q)}
{\theta(b_je^{-2\pi iy};\tilde q)}\, e^\delta=\prod_{j=1}^m\frac{\Gamma(a_je^{-2\pi iy}
{r} ,b_je^{-2\pi iy};\tilde q,{r} )}
{\Gamma(a_je^{-2\pi iy},b_je^{-2\pi iy}{r} ;\tilde q,{r} )}\, e^\delta,
$$
где $\prod_{j=1}^ma_jb_j^{-1}=1$. В следствие такого представления можно поглотить
функцию $\tilde\varphi(y)$ в отношение эллиптических гамма-функций
в (\ref{phi-form}) заменив $n\to n+2m$ и отождествив
$\tilde t_k={r}  a_k, \tilde w_k=a_k$ при $k=n+1,\ldots, n+m$
и $\tilde t_k=\tilde b_k, \tilde w_k={r}  \tilde b_k$ при
$k=n+m+1,\ldots,n+2m$. Поскольку $n$, $\tilde t_j$, $\tilde w_j$
произвольны, без потери общности можно положить $\tilde\varphi(y)=e^\delta=const$,
что дает требуемое выражение (\ref{delta-ell}).

Функция (\ref{delta-ell}) удовлетворяет уравнениям
\begin{eqnarray} \label{2-eq'}
&&\Delta(y+\sigma^{-1})=e^{\gamma/\sigma}
\prod_{j=0}^{n}\frac{\theta(\tilde t_je^{-2\pi iy};{r} )}
{\theta(\tilde w_je^{-2\pi i y};{r} )}\, \Delta(y),  \\
&&\Delta(y+\tau\sigma^{-1})=e^{\gamma\tau/\sigma}
\prod_{j=0}^r \frac{\theta(t_jq^y;q)}{\theta(w_jq^y;q)}
\prod_{j=0}^{n}\frac{\theta(\tilde w_je^{-2\pi iy}{r} ^{-1};\tilde q)}
{\theta(\tilde t_je^{-2\pi iy}{r} ^{-1};\tilde q)}\, \Delta(y).
\label{3-eq'}\end{eqnarray}
Эллиптические функции, определенные произведениями $\prod_{j=0}^r$ и
$\prod_{j=0}^{n}$ в (\ref{3-eq'}), выглядят по разному
хотя обе имеют периоды 1 и $\sigma^{-1}$. Их формы связаны друг с другом
модулярным преобразованием $\sigma\to -1/\sigma$
для соответствующих тета-функций.

Рассмотрим теперь область $\mbox{Im}(\tau/\sigma)<0$. Уравнения
(\ref{1-eq}) и (\ref{3-eq'}) хорошо определены в этом случае. Они
могут быть использованы для определения $\Delta(y)$, и можно проверить,
что, действительно, функция (\ref{delta-ell'}) дает их общее решение.
Уравнение (\ref{2-eq'}) теперь заменяется на
\begin{equation}
\frac{\Delta(y+\sigma^{-1})}{\Delta(y)}= e^{\gamma/\sigma} \prod_{j=0}^{n}
\frac{\theta(\tilde w_je^{-2\pi iy};{r} ^{-1})}
{\theta(\tilde t_je^{-2\pi iy}; {r} ^{-1})},
\label{phi-eq'}\end{equation}
т.е. ${r} $ в (\ref{phi-eq}) заменяется на ${r} ^{-1}$,
и параметры $\tilde t_j, \tilde w_j$ заменяются на
${r} ^{-1}\tilde w_j$ и ${r} ^{-1}\tilde t_j$, соответственно.
Используя (\ref{D'}), легко построить функцию $\Delta(y)$, удовлетворяющую
(\ref{1-eq}), (\ref{2-eq'}) и (\ref{3-eq'}), и при $|q|>1$.
\hfill{Q.E.D.}

\smallskip

Для того, чтобы иметь возможность работать с $q$, лежащей на единичной
окружности, $|q|=1$, необходима другая эллиптическая гамма-функция --
некоторый аналог функции двойного синуса (\ref{2d-sin}). Обозначим
\be
q=e^{2\pi i\frac{\omega_1}{\omega_2}}, \qquad
\tilde q =e^{-2\pi i\frac{\omega_2}{\omega_1}},
\qquad  p=e^{2\pi i\frac{\omega_3}{\omega_2}}, \qquad
r =e^{2\pi i\frac{\omega_3}{\omega_1}},
\label{ell-bases}\ee
где $\omega_i\in\C$ обозначают произвольные в общем случае несоизмеримые
параметры. Предположим, что $\omega_1/\omega_2>0$ и
$\mbox{Im}(\omega_3/\omega_2)>0$ (т.е., $|p|<1$). Тогда имеем
$\mbox{Im}(\omega_3/\omega_1)=(\omega_2/\omega_1)\mbox{Im}(\omega_3/\omega_2)>0$,
т.е. $|r |<1$ автоматически. Поэтому при анализе уравнения (\ref{1-eq})
с $|q|=1$ необходимо подразумевать, что $|p|, | r |<1$.

\begin{definition}
Пусть $|q|,|p|, |r |<1$. Тогда модифицированная эллиптическая
гамма-функция определяется формулой
\begin{equation}
G(u;\mathbf{\omega})= \prod_{j,k=0}^\infty
\frac{(1-e^{-2\pi i\frac{u}{\omega_2}}q^{j+1}p^{k+1})
(1-e^{2\pi i \frac{u}{\omega_1}}{\tilde q}^{j+1}{r }^k)}
{(1-e^{2\pi i \frac{u}{\omega_2}}q^jp^k)
(1-e^{-2\pi i \frac{u}{\omega_1}}{\tilde q}^j{r }^{k+1})}.
\label{ell-d}\end{equation}
\end{definition}

В пределе $p\to 0$, взятом таким образом чтобы одновременно иметь
$ r \to 0$, получаем $G(u;\omega_1,\omega_2,\omega_3)\to
S^{-1}(u;\omega_1,\omega_2)$, где $S(u;\mathbf{\omega})$
есть функция двойного синуса, определенная равенством (\ref{2d-sin}).

Как показано в работах \cite{spi:integrals,die-spi:unit} (см. также
первый параграф главы 5), функция (\ref{ell-d}) хорошо определена
для действительных $\omega_1, \omega_2$ с $\omega_1/\omega_2>0$
(и любом комплексном $\omega_3$), как и в случае функции двойного
синуса. Поэтому естественно ожидать, что $G(u;\mathbf{\omega})$
позволит расширить известные $|q|=1$ результаты в теории
$q$-гипергеометрических функций на эллиптический уровень.
В частности, интересно было бы найти эллиптические расширения принципа
модулярного удвоения Фаддеева $q$-деформированных алгебр
\cite{fad:modular,kls:unitary} и некоммутитивной гипергеометрии
Волкова \cite{volk}.

В результате, при $|q|=1$ мы получим решение $\Delta(y)$ уравнения
(\ref{1-eq}) простой заменой $\Gamma(q^y;q,p)$ в (\ref{DE})
на $G(y\omega_1;\mathbf{\omega})$. Отметим, что область базовых
переменных $|p|=|q|=1$ не правомочна при рассмотрении эллиптических
функций. В некотором смысле область действительных значений всех
$\omega_{1,2,3}$ достижима только на уровне чистых многократных
гамма-функций Барнса \cite{bar:multiple}.

\section{Тета и эллиптический аналоги функции Мейера}

Интегралы с ядром вида (\ref{DE}) могут рассматриваться
как эллиптические расширения функций Мейера специального вида.
Общие тета-функциональные аналоги функции Мейера
возникают в случае, если $h(y)$ квазипериодическая функция, соответствующая
четвертому случаю в определении тета-гипергеометрических интегралов
(\ref{quasi-per}) с $n=1$ .

Обозначим $P_3(y)=\sum_{i=1}^3\alpha_iy^i$ произвольный полином третьей степени
со свойством $P_3(0)=0$. Тогда функция, определенная интегралом
\begin{equation}
G_r^s\left({\mathbf{t} \atop \mathbf{w}};\mathbf{\alpha};q,p\right)
=\int_C\frac{\prod_{j=0}^{s}\Gamma(t_jq^y;q,p)}
{\prod_{k=0}^r\Gamma(w_kq^y;q,p)} e^{P_3(y)}dy,
\label{e-meijer}\end{equation}
где $C$ некоторый контур на комплексной плоскости гарантирующий сходимость
интеграла, называется тета-аналогом функции Мейера. Отметим, что
на целые числа $r, s$ и параметры $t_j, w_k\in\C$ нет никаких ограничений.

Подынтегральная функция $\Delta(y)$ в (\ref{e-meijer}) удовлетворяет уравнению
\begin{equation}
\frac{\Delta(y+1)}{\Delta(y)}=h(y)=e^{P_2(y)}
\frac{\theta(t_0q^y,\ldots,t_sq^y ;p)}
{\theta(w_0q^y,\ldots,w_rq^y;p)},
\label{theta-int}\end{equation}
где $P_2(y)=P_3(y+1)-P_3(y)$ является полиномом $y$ второй степени.
Из рассмотрения, проведенного в \cite{spi:theta} и изложенного в
предыдущей главе следует, что $h(y)$ такого вида задает наиболее общую
мероморфную функцию $y$, для которой $1/h(y)$ так же мероморфна и
выполняются условия квазипериодичности
\begin{equation}
h(y+\sigma^{-1})=e^{ay+b}h(y),\qquad
h(y+\tau\sigma^{-1})=e^{cy+d} h(y)
\label{qua}\end{equation}
для некоторых параметров $a,b,c,d$. Функцию $h(y)$ можно также интерпретировать
как общую мероморфную модулярную форму Якоби в смысле Эйхлера-Загира
\cite{eic-zag:theory}.

Однако, интеграл (\ref{e-meijer}) не является наиболее общим интегралом,
приводящим к (\ref{theta-int}). Используя подходящие модификации
выражений (\ref{delta-ell}) и (\ref{delta-ell'}) и заменив $y$ на
$y/\omega_1$, получаем общий тета-аналог функции Мейера.

\begin{definition}
Пусть базисные переменные (\ref{ell-bases}) удовлетворяют
условиям $|q|, |p|<1$. Тогда, при $| r |<1$ интеграл
\begin{eqnarray}\nonumber
\lefteqn{G_{rm}^{sn}\left({\mathbf{t},\mathbf{\tilde t}\atop \mathbf{w},
\mathbf{\tilde w}};\mathbf{\alpha};\mathbf{\omega}\right) } &&
\\ &&
=\int_C\frac{\prod_{k=0}^{s}\Gamma(t_ke^{\frac{2\pi iy}{\omega_2}};q,p)
\prod_{j=0}^{n}\Gamma(\tilde t_je^{-\frac{2\pi iy}{\omega_1}};
\tilde q, r )}{\prod_{k=0}^r\Gamma(w_ke^{\frac{2\pi iy}{\omega_2}};q,p)
\prod_{j=0}^m\Gamma(\tilde w_je^{-\frac{2\pi iy}{\omega_1}};
\tilde q, r )}\, e^{P_3(y)}dy
\label{gen-int}\end{eqnarray}
называется общим тета-гипергеометрическим интегралом одной переменной.
При $| r |>1$, полагаем
\begin{eqnarray}\nonumber
\lefteqn{G_{rm}^{sn}\left({\mathbf{t},\mathbf{\tilde t}\atop \mathbf{w},
\mathbf{\tilde w}};\mathbf{\alpha};\mathbf{\omega}\right) }&&
\\ &&
=\int_C\frac{\prod_{k=0}^{s}\Gamma(t_ke^{\frac{2\pi iy}{\omega_2}};q,p)
\prod_{j=0}^{m}\Gamma(\tilde w_je^{-\frac{2\pi iy}{\omega_1}} r ^{-1};
\tilde q, r ^{-1})}{\prod_{k=0}^r
\Gamma(w_ke^{\frac{2\pi iy}{\omega_2}};q,p)
\prod_{j=0}^n\Gamma(\tilde t_je^{-\frac{2\pi iy}{\omega_1}} r ^{-1};
\tilde q, r ^{-1})}\, e^{P_3(y)}dy.
\label{gen-int'}\end{eqnarray}
Контур $C$, параметры $r, s, n, m \in\mathbb{N}$ и
$t_j,\tilde t_j, w_k, \tilde w_k\in\C$ произвольны, единственным ограничением
служит только условие существования интегралов \re{gen-int} и \re{gen-int'}.
\end{definition}

Подынтегральные функции (\ref{gen-int}) и (\ref{gen-int'}) удовлетворяют
уравнениям $\Delta(y+\omega_i)/\Delta(y)=h_i(y)$, $i=1,2,3,$ с квазипериодическими
функциями $h_i$: $h_i(y+\omega_k)=e^{a_{ik} y+b_{ik}}h(y)$,
$i\neq k$, где $a_{ik}, b_{ik}$ некоторые константы, связанные с параметрами
$\mathbf{t}, \mathbf{\tilde t},\mathbf{w},\mathbf{\tilde w},
\mathbf{\alpha}$ и $\mathbf{\omega}$. Интеграл (\ref{gen-int'})
определяется из условия, что он имеет те же самые
функции $h_1(y)$ и $h_3(y)$, что и (\ref{gen-int}). При специальном
выборе параметров $\mathbf{t}, \mathbf{\tilde t},\mathbf{w},
\mathbf{\tilde w},\mathbf{\alpha}$, в пределе $|p|, | r |\to 0$ или
$|p|, | r |^{-1}\to 0$ функция
$G_{rm}^{sn}\left({\mathbf{t},\mathbf{\tilde t}\atop \mathbf{w},
\mathbf{\tilde w}};\mathbf{\alpha};\mathbf{\omega}\right)$ сводится
к общему $q$-гипергеометрическому интегралу, рассмотренному
в предыдущем параграфе, см. выражения (\ref{q-mej}) и (\ref{tq-mej}).

Общий тета-гипергеометрический ряд одной переменной имеет вид:
\be
{_{s+1}E_r}\left({t_0,\ldots, t_{s} \atop w_1,\ldots,w_r};
\mathbf{\alpha}; q,p\right)
= \sum_{n=0}^\infty \frac{\theta(t_0,t_1,\ldots,t_{s};p;q)_n}
{\theta(q,w_1,\ldots,w_r;p;q)_n}\, e^{P_3(n)}.
\label{_rE_s}\ee
Это наиболее общий ряд, удовлетворяющий определению тета-гипергеометрических
функций в виде однократного ряда. Случай, рассмотренный в предыдущей главе,
соответствует выбору $\alpha_3=0$ и его расширение на (\ref{_rE_s}) вполне
естественно. Отношение $c_{n+1}/c_n$ для этого ряда равно функции $h(n)$
(\ref{theta-int}) с $w_0=q$ и $y=n$, что не является случайным фактом.
Рассмотрим последовательность полюсов подынтегральной функции в
(\ref{e-meijer}) расположенных в точках $y=y_0+n, n\in\mathbb{N}$, для некоторого
$y_0$. Обозначим  $\kappa c_n, c_0=1,$ вычеты этих полюсов, т.е. при
$y\to y_0+n$ имеем $\Delta(y)\to \kappa c_n/(y-y_0-n) +O(1)$.
Теперь несложно заметить, что
$$
\lim_{y\to y_0+n}\frac{\Delta(y+1)}{\Delta(y)}=\frac{c_{n+1}}{c_n}=
\lim_{y\to y_0+n} h(y)=h(y_0+n).
$$
В частности это означает, что суммы вычетов в интеграле (\ref{e-meijer}),
появляющихся после подходящих деформаций контура $C$, образуют
тета-гипергеометрические ряды вида \re{_rE_s}) с определенным выбором
параметров.

По аналогии с классификацией введенной для рядов, будем называть интеграл
(\ref{e-meijer}) эллиптическим гипергеометрическим или сбалансированным
тета-гипергеометрическим интегралом, если
\begin{equation}
s=r, \quad \alpha_3=\alpha_2=0, \quad \prod_{j=0}^rt_j=\prod_{j=0}^r w_j.
\label{theta-balance}\end{equation}
В этом случае $h(y)$ в (\ref{theta-int}) становится эллиптической функцией $y$.

Если $h(y)$ эллиптична и по $y$ и по всем параметрам $u_j, v_j$
(напомним, что $t_j=q^{u_j}, w_j=q^{v_j}$), то функция (\ref{e-meijer}) определяет
{\em полностью эллиптический гипергеометрический} интеграл.
Как показано в предыдущей главе, в дополнение к сбалансированности,
такое свойство требует наложения условий
\begin{equation}\label{well-poised-2i}
t_jw_j=\rho=const, \qquad j=0,\ldots,r,
\end{equation}
известных как условия вполне уравновешенности для $q$-гипергеометрических
рядов \cite{gas-rah:basic}.  Подынтегральная функция $\Delta(y)$
в этом случае приобретает вид
\begin{equation}
\Delta(y)=\prod_{j=0}^{m-1}\frac{\Gamma(t_jz;q,p)}
{\Gamma(\rho t_j^{-1}z;q,p)}\; \frac{\Gamma(\rho^{\frac{m+1}{2}}
\prod_{j=0}^{m-1}t_j^{-1}\, z;q,p)}{\Gamma(\rho^{\frac{1-m}{2}}
\prod_{j=0}^{m-1}t_j\, z;q,p)}\,e^{\gamma y},
\label{wp-int}\end{equation}
где мы обозначили $z=q^y$ и $\gamma=\alpha_1$.
Параметр $\rho$ излишен, он может быть исключен из рассмотрения масштабными
преобразованиями $t_i\to \rho^{1/2}t_i, z\to \rho^{-1/2}z$, но он удерживается
для удобства. Заметим, что без потери общности один из параметров
в (\ref{e-meijer}) может быть устранен сдвигом $y$.

Интеграл (\ref{e-meijer}) называется {\em совершенно уравновешенным}, если,
в дополнение к условиям (\ref{well-poised-2i}), восемь параметров $t_i$
фиксированы следующим образом:
\begin{eqnarray}
(t_{m-8},\ldots,t_{m-1})=(\pm(pq)^{1/2}, \pm q^{1/2}p, \pm p^{1/2}q,\pm pq).
\label{vwp}\end{eqnarray}
Эти ограничения приводят к удвоению аргумента эллиптической гамма-функции
\begin{equation}
\prod_{j=m-8}^{m-1}\Gamma(t_jz;q,p)=\Gamma(pqz^{2};q,p)=\frac{1}{\Gamma(z^{-2};q,p)}
\label{doubling}\end{equation}
и подынтегральная функция совершенно уравновешенного сбалансированного
тета-ги\-пер\-гео\-ме\-три\-чес\-ко\-го интеграла принимает вид
\begin{equation}
\Delta(y)=\prod_{j=0}^{m-9}\frac{\Gamma(t_jz;q,p)}
{\Gamma(\rho t_j^{-1}z;q,p)}\, \frac{\Gamma(\rho^{\frac{m+1}{2}}
p^{-6}q^{-6}\prod_{j=0}^{m-9}t_j^{-1}z;q,p)\,e^{\gamma y} }
{\Gamma(z^{-2},(\rho/pq)^2 z^2, \rho^{\frac{1-m}{2}}p^6q^6
\prod_{j=0}^{m-9}t_j\, z;q,p)}.
\label{vwp-int}\end{equation}

После наложения условий (\ref{vwp}), параметр $\rho$ больше не излишен
и его выбор играет важную роль. Если мы положим $\rho=pq,$ то $\Delta(y)$
принимает особенно симметричный вид
\begin{equation}
\Delta(y)=\frac{\prod_{j=0}^{m-9}\Gamma(t_jz,t_jz^{-1};q,p)}
{\Gamma(z^2,z^{-2}, Az, Az^{-1};q,p)}\,e^{\gamma y},
\label{vwp-int-fin}\end{equation}
где $A=(pq)^{\frac{13-m}{2}} \prod_{j=0}^{m-9}t_j$.
Случаи четного и нечетного значений $m$ качественно отличаются друг от друга.
Выбор $m=13$ (содержащий в себе $m=9$ и $11$ случаи)  играет важную роль в
дальнейших рассмотрениях. Он соответствует самой простой форме $\Delta(y)$
и при $\gamma =0$ приводит к эллиптическому бета-интегралу, вычисленному
в \cite{spi:elliptic,spi:umn} и описываемому в следующем параграфе.

Что касается многомерных интегралов гипергеометрического типа, общий вид
$\Delta(\mathbf{y})$ в обычном и $q$-гипергеометрических случаях
может быть выведен из теоремы Оре-Сато для рядов (см., например, ее
подробное обсуждение в \cite{ggr:general}). Общая структура многократных
тета и эллиптических гипергеометрических рядов и интегралов еще не
установлена и это одна из важных открытых проблем в этой области.
Интересные тета-гипергеометрические функции появляются в
модели Фельдера-Варченко \cite{fel-var:qkz}, связанной с теорией уравнений
Книжника-Замолодчикова над эллиптическими кривыми.
В дальнейшем мы будем рассматривать только эллиптические гипергеометрические
интегралы, связанные с одномерным эллиптическим бета-интегралом и его
многомерными расширениями на корневые системы  $A_n$ и $C_n$.

\section{Одномерная эллиптическая бета-функция}

Возьмем переменную $z\in\C$ и пять комплексных параметров $t_m, m=1, \dots,5.$
Предположим, что $|q|,|p|<1$ и скомпонуем ядро эллиптического гипергеометрического
интеграла
\begin{equation}
\rho(z,t_1,\ldots,t_5)=\frac{\prod_{m=1}^5\Gamma(t_mz^\pm,At_m^{-1})}
{\Gamma(z^{\pm 2},Az^\pm)\prod_{1\leq m<s\leq 5}\Gamma(t_mt_s)},
\label{kernel}\end{equation}
где $A=\prod_{m=1}^5t_m$. Эта функция имеет последовательности полюсов,
сходящихся к $z=0$ по точкам
$$
\P=\{t_mq^jp^k,\, A^{-1}q^{j+1}p^{k+1}
\}_{m=1,\ldots, 5,\, j,k \in\N}
$$
и уходящих к бесконечности вдоль их партнеров $\P^{-1}$, получающихся отражением
$z\to 1/z$. Обозначим $C$ некоторый гладкий контур на комплексной плоскости,
ориентированный против часовой стрелки и разделяющий множества точек $\P$ и
$\P^{-1}$. Например, при $|t_m|<1,\; |pq|<|A|,$
контур $C$ может совпадать с единичной окружностью $\mathbb{T}$.

\begin{theorem}
\begin{equation}\label{ell-int}
\int_C\rho(z,t_1,\ldots,t_5)\frac{dz}{z}=\frac{4\pi i}
{(q;q)_\infty(p;p)_\infty}.
\end{equation}
\end{theorem}
{\bf Доказательство.}
Первый шаг в вычислении этого интеграла состоит в проверке
справедливости следующего $q$-разностного уравнения для ядра:
\begin{equation}
\rho(z,qt_1,t_2,\ldots,t_5)-\rho(z,t_1,\ldots,t_5)
=g(q^{-1}z,t_1,\ldots,t_5)-g(z,t_1,\ldots,t_5),
\label{eqn}\end{equation}
где
\begin{equation}
g(z,t_1,\ldots,t_5)=\rho(z,t_1,\ldots,t_5)
\frac{\prod_{j=1}^5\theta(t_jz;p)}{\prod_{j=2}^5\theta(t_1t_j;p)}
\frac{\theta(t_1A;p)}{\theta(z^2,Az;p)}\frac{t_1}{z}.
\label{g}\end{equation}
В последовательных пределах $p\to 0$ и затем $t_5\to 0$,
равенство \re{eqn} сводится к уравнению, использованному в работе \cite{wz:invent}
для доказательства интеграла Аски-Вильсона \cite{ask-wil:some}.
После деления уравнения (\ref{eqn}) на $\rho(z,t_1,\ldots,t_5)$ оно принимает
форму
\begin{eqnarray}\nonumber
\lefteqn{\frac{\theta(t_1z,t_1z^{-1};p)}{\theta(Az,Az^{-1};p)}
\prod_{j=2}^5\frac{\theta(At_j^{-1};p)}{\theta(t_1t_j;p)}-1 } &&
\\ &&
=\frac{t_1\theta(t_1A;p)}{z\theta(z^2;p)\prod_{j=2}^5\theta(t_1t_j;p)}
\left(\frac{z^4\prod_{j=1}^5\theta(t_jz^{-1};p)}{\theta(Az^{-1};p)} -
\frac{\prod_{j=1}^5\theta(t_jz;p)}{\theta(Az;p)}\right).
\label{eqn-exp}\end{eqnarray}

Обе стороны этого равенства представляют собой эллиптические функции переменной
$\log z$ (т.е. они инвариантны относительно преобразования $z\to pz$, что легко
проверяется) с одинаковыми полюсами и их вычетами. Например,
$$
\lim_{z\to A}\theta(Az^{-1};p)\, (\mbox{левая часть}) =
\frac{\theta(t_1A,t_1A^{-1};p)}{\theta(A^2;p)}
\prod_{j=2}^5\frac{\theta(At_j^{-1};p)}{\theta(t_1t_j;p)}
$$
и точно такой же результат имеет место для правой части.
Поэтому разность левой и правой частей равенства (\ref{eqn-exp})
является эллиптической функцией без полюсов, т.е. константой. Эта
константа должна быть равна нулю, так как при $z=t_1$ справедливость
(\ref{eqn-exp}) очевидна.

Теперь мы проинтегрируем уравнение (\ref{eqn}) по переменной $z\in C$
и получим
\begin{equation}
I(qt_1,t_2,\ldots,t_5)-I(t_1,\ldots,t_5)=
\left(\int_{q^{-1}C}-\int_C\right)g(z,t_1,\ldots,t_5)
\frac{dz}{z},
\label{int-eqn}\end{equation}
где $I(t_1,\ldots,t_5)=\int_C\rho(z,t_1,\ldots,t_5)dz/z$ и
$q^{-1}C$ обозначает контур $C$, растянутый множителем $q^{-1}$ по отношению
к точке $z=0$. Функция (\ref{g}) имеет сходящиеся к $z=0$ последовательности
полюсов при $z=t_mq^jp^k, A^{-1}q^jp^{k+1}$ и уходящие на бесконечность при
$z=t_m^{-1}q^{-1-j}p^{-k},
Aq^{-j-1}p^{-k-1}$ для $m=1,\ldots,5$ и $j,k\in\mathbb{N}$.
Взяв $C=\mathbb{T}$, мы видим, что при $|t_m|<1$ и $|p|<|A|$
кольцо $1\leq |z|\leq |q|^{-1}$ не содержит полюсов. Поэтому
мы можем деформировать $q^{-1}\mathbb{T}$ к $\mathbb{T}$
в (\ref{int-eqn}) и получить ноль в правой части этого равенства.
Это дает $I(qt_1,t_2,\ldots,t_5)=I(t_1,\ldots,t_5)$.
Наложив более жесткие ограничения $|p|,|q|<|A|$, из перестановочной симметрии
по переменным $p$ и $q$, получаем
$I(pt_1,t_2,\ldots,t_5)=I(t_1,\ldots,t_5)$. Дальнейшие преобразования
растяжения $t_1\to q^{\pm1}t_1$ и $t_1\to p^{\pm1}t_1$ могут
быть произведены только если они не выводят параметр $t_1$ из кольца
аналитичности функции $I(t_1,\ldots,t_5)$.

Пусть, временно, $p$ и $q$ действительны, $p<q$ и $p^n\neq q^k$ для любых
$n,k\in\mathbb{N}$. Предположим так же, что аргументы комплексных параметров
$t_m^{\pm1},\, m=1,\ldots,5,$ и $A^{\pm 1}$ попарно отличаются. Возьмем
теперь в качестве $C$ контур, охватывающий $\P$ и два отрезка
$c_1=[t_1, t_1p^2]$, $c_2=[(pq/A)p^{-2},pq/A]$ и исключающий их $z\to 1/z$ партнеры.
Теперь будем делать преобразования сжатия $t_1\to t_1q^k$, $k=1,2,\ldots,$
до тех пор, пока $t_1q^k$ не войдет в интервал $[t_1p,t_1p^2]$.
После этого произведем растяжение $t_1\to t_1p^{-1}$, которое не выводит параметры
$t_1$ и $pq/A$ из $c_1$ и $c_2$, соответственно. Таким образом, мы получаем
равенство $I(q^jp^{-k}t_1,t_2,\ldots,t_5)=I(t_1,\ldots,t_5)$
при таких $j,k\in\N$, что $q^jp^{-k}\in [1,p]$. Поскольку такое множество  точек
плотно, то мы приходим к выводу, что $I$ не зависит от $t_1$ и, по симметрии, от
всех $t_m$. Поэтому $I$ есть постоянная, зависящая только от $p$ и $q$.
Ее значение, даваемое выражением в правой части (\ref{ell-int}),
находится с помощью анализа структуры вычетов полюсов подынтегральной функции.
Для этого необходимо деформировать контур $C$ так, чтобы он пересек полюса
в точках $z=t_1^{\pm 1}$, учесть соответствующие вычеты и взять предел
$t_2\to 1/t_1$ (при котором никакие другие полюса не пересекают  $C$),
подобно тому как это было сделано в работе \cite{die-spi:elliptic}.
Более подробное рассмотрение этой процедуры с более общими выводами будет
произведено в конце этой главы.

После доказательства формулы интегрирования (\ref{ell-int}) в ограниченной
области параметров, мы можем произвести аналитическое продолжение по параметрам
в область, допустимую выбранным контуром $C$.
\hfill{Q.E.D.}

\begin{remark}
Интеграл (\ref{ell-int}) был впервые вычислен в работе \cite{spi:elliptic} с помощью
другого уравнения по параметрам с использованием $_2\psi_2$ формулы суммирования Бэйли
\cite{bai:well,bai:some} и аналитического
продолжения по дискретным значениям четырех параметров. Соответствующий метод
являлся эллиптическим обобщением метода Аски \cite{ask:element,ask:beta},
предложенного для доказательства интеграла Рахмана
\cite{rah:integral} (см. также \cite{nas-rah:projection})
получающегося из \re{ell-int} в пределе $p\to 0$. Приведенное
доказательство, представленное в работе \cite{spi:short},
проще и оно не использует никаких тождеств для $q$-гипергеометрических рядов
или интегралов, а аналитическое продолжение используется в минимальной форме.
\end{remark}

\section{Це\-поч\-ки преобразований для эллиптических гипергеометрических
интегралов}

Цепочки Бэйли предоставляют мощные средства для порождения бесконечных
последовательностей формул суммирования или преобразования для рядов
гипергеометрического типа \cite{and:bailey,war:50}.
В предыдущей главе описано применение этой техники
к эллиптическим гипергеометрическим рядам.
При работе над статьей \cite{spi:bailey}, автор пришел к
принципиальному заключению, что должны существовать цепочки Бэйли
для интегралов, но первые попытки построения простых примеров таких
цепочек не увенчались успехом. Нетривиальное преобразование симметрии
типа Бэйли для пары эллиптических гипергеометрических интегралов
было сконструировано в \cite{spi:integrals}. Оказалось, что этот
результат дает средства необходимые для соответствующего обобщения
эллиптической цепочки Бэйли для рядов.
В данном параграфе, основанном на работе \cite{spi:tree}, сконструированы
два примера цепочек Бэйли для интегралов с помощью эллиптического
бета-интеграла. При совместном использовании
эти цепочки образуют двоичное дерево тождеств для
эллиптических гипергеометрических интегралов.

Основное концептуальное определение вводит понятие
пар Бэйли для интегралов.

\begin{definition}
Две функции $\alpha(z,t)$ и $\beta(z,t)$, $z, t\in \mathbb{C}$,
по определению образуют интегральную эллиптическую пару Бэйли
по отношению к параметру $t$ если они удовлетворяют следующему равенству
\begin{equation}
\beta(w,t)=\kappa
\int_{\mathbb{T}}\Gamma(tw^\pm z^\pm)\, \alpha(z,t)\frac{dz}{z},
\qquad \kappa\equiv\frac{(p;p)_\infty (q;q)_\infty}{4\pi i}.
\label{bint-pair}\end{equation}
\end{definition}

Предположим, что мы нашли частный пример интегральной эллиптической
па\-ры Бэй\-ли. Тогда существует способ построения бесконечной
последовательности таких пар исходя из заданной пары,
определяемый следующей интегральной леммой Бэйли.

\begin{theorem}\label{ibl-1}
Для всякой пары функций $\alpha(z,t)$ и $\beta(z,t)$,
образующей интегральную эллиптическую пару Бэйли по отношению
к $t$, новые функции, определенные для $w\in\mathbb{T}$ как
\begin{equation}
\alpha'(w,st)=\frac{\Gamma(tuw^\pm)}
{\Gamma(ts^2uw^\pm)}\,\alpha(w,t)
\label{a'}\end{equation}
и
\begin{equation}
\beta'(w,st)=\kappa
\frac{\Gamma(t^2s^2,t^2suw^\pm)}{\Gamma(s^2,t^2,suw^\pm)}
\int_{\mathbb{T}} \frac{\Gamma(sw^\pm x^\pm,ux^\pm)}
{\Gamma(x^{\pm 2},t^2s^2ux^\pm)}\,\beta(x,t)\frac{dx}{x}
\label{b'}\end{equation}
образуют интегральную эллиптическую пару Бэйли по отношению
к параметру $st$. Параметры $t,s,u\in\C$ удовлетворяют
ограничениям $|t|,|s|,|u|<1, |pq|<|t^2s^2u|$.
\end{theorem}
{\bf Доказательство.}
Подставим определение (\ref{bint-pair}) в соотношение (\ref{b'})
\begin{eqnarray}\nonumber
\lefteqn{\beta'(w,st)=\kappa^2
\frac{\Gamma(t^2s^2,t^2suw^\pm)}{\Gamma(s^2,t^2,suw^\pm)} } &&
\\ && \makebox[4em]{}
\times \int_{\mathbb{T}^2}
\frac{\Gamma(sw^\pm x^\pm,tz^\pm x^\pm,ux^\pm)}
{\Gamma(x^{\pm 2},t^2s^2ux^\pm)}\,\alpha(z,t)\frac{dz}{z}\frac{dx}{x}.
\label{proof}\end{eqnarray}
Поскольку интегрируемая функция ограничена на $\mathbb{T}$, мы можем
изменить порядок интегрирований (т.е. проинтегрировать сначала по $x$)
и применить формулу (\ref{ell-int}). Это дает
\begin{eqnarray}\nonumber
&& \beta'(w,st)=\kappa
\int_{\mathbb{T}}\frac{\Gamma(stw^\pm z^\pm,tuz^\pm)}
{\Gamma(ts^2uz^\pm)}\; \alpha(z,t)\frac{dz}{z}
\\ && \makebox[4em]{}
= \kappa
\int_{\mathbb{T}}\Gamma(stw^\pm z^\pm)\; \alpha'(z,st)\frac{dz}{z},
\label{st}\end{eqnarray}
то есть функции $\alpha'(w,st)$ и $\beta'(w,st)$ образуют
интегральную пару Бэйли по отношению к параметру $st$.
\hfill{Q.E.D.}

\smallskip

Если мы подставим в интеграл (\ref{ell-int}) соотношения
$t_3=tw, t_4=tw^{-1}$, тогда нетрудно найти начальную
интегральную эллиптическую пару Бэйли:
\begin{eqnarray}\label{a_1}
&& \alpha(z,t)=\frac{\prod_{r=0}^2\Gamma(t_rz^\pm)}
{\Gamma(z^{\pm 2}, t^2t_0t_1t_2z^\pm)},
\\ &&
\beta(w,t)=\Gamma(t^2)\prod_{0\leq r<j\leq 2}\frac{\Gamma(t_rt_j)}
{\Gamma(t^2t_rt_j)}\,
\frac{\prod_{r=0}^2\Gamma(tt_rw^\pm)}
{\Gamma(tt_0t_1t_2w^\pm)},
\label{b_1}\end{eqnarray}
где $|t|, |t_r|<1, |pq|<|t^2t_0t_1t_2|$.
Подставляя эти выражения в цепочечные правила (\ref{a'}) и
(\ref{b'}), мы получаем
\begin{eqnarray*}
&& \alpha'(z,st)=\frac{\Gamma(tuz^\pm)\prod_{r=0}^2\Gamma(t_rz^\pm)}
{\Gamma(ts^2uz^\pm,z^{\pm2},t^2t_0t_1t_2z^\pm)},
\\
&& \beta'(w,st)=\kappa \frac{\Gamma(t^2s^2,t^2suw^\pm)}
{\Gamma(s^2,suw^\pm)}
\prod_{0\leq r<j\leq 2}\frac{\Gamma(t_rt_j)}{\Gamma(t^2t_rt_j)}
\\
&& \makebox[4em]{} \times
\int_{\mathbb{T}}
\frac{\Gamma(sw^\pm x^\pm,ux^\pm)\prod_{r=0}^2\Gamma(tt_rx^\pm)}
{\Gamma(x^{\pm 2},t^2s^2ux^\pm,tt_0t_1t_2x^\pm)}\frac{dx}{x}.
\end{eqnarray*}
Будучи подставленными в основное генерирующее соотношение
(\ref{st}), эти функции порождают преобразование симметрии
для двух эллиптических гипергеометрических интегралов:
\begin{eqnarray}\nonumber
&& \prod_{j=0}^{2}\frac{\Gamma(Bt_j^{-1})}{\Gamma(t^2Bt_j^{-1})}
\int_{\mathbb{T}}\frac{\prod_{j=1}^{3}
\Gamma(tt_jz^\pm,s_jz^\pm)}
{\Gamma(z^{\pm2},t^2Sz^\pm,tBz^\pm)}\frac{dz}{z}
\\ && \makebox[2em]{}
=\prod_{j=0}^{2}\frac{\Gamma(Ss_j^{-1})}{\Gamma(t^2Ss_j^{-1})}
\int_{\mathbb{T}}\frac{\prod_{j=1}^{3}
\Gamma(ts_jz^\pm,t_jz^\pm)}
{\Gamma(z^{\pm2},t^2Bz^\pm,tSz^\pm)}\frac{dz}{z},
\label{an-trans}\end{eqnarray}
где $s_0=u,s_1=sw,s_2=sw^{-1}$ и
$B=\prod_{j=0}^{2}t_j, S=\prod_{j=0}^{2}s_j$.
Тождество (\ref{an-trans}) было установлено в работе \cite{spi:integrals}
и интегральная эллиптическая цепочка Бэйли (\ref{bint-pair})-(\ref{b'})
построена в результате обобщения соответствующего способа вывода этого равенства.

Итеративное применение правил (\ref{a'}) и (\ref{b'})
порождает бесконечную цепочку тождеств для эллиптических
гипергеометрических интегралов.  Результат $m$-кратной
итерации имеет вид
\begin{eqnarray}\label{a-m}
&& \alpha^{(m)}(x,\prod_{l=1}^ms_l\, t)=
\prod_{k=1}^m\frac{\Gamma(t\prod_{l=1}^{k-1}s_l\, u_kx^\pm)}
{\Gamma(t\prod_{l=1}^{k-1}s_l\, s_k^2u_kx^\pm)} \alpha^{(0)}(x,t),
\\ \label{b-m}  &&
\beta^{(m)}(x_{m+1},\prod_{l=1}^ms_l\, t)=\kappa^m
\prod_{k=1}^m\frac{\Gamma(t^2\prod_{l=1}^{k}s_l^2)}
{\Gamma(s_k^2,t^2\prod_{l=1}^{k-1}s_l^2)}
\\ && \makebox[2em]{} \times
\int_{\mathbb{T}^m} \prod_{k=1}^m\frac
{\Gamma(t^2\prod_{l=1}^{k-1}s_l^2\, s_ku_kx_{k+1}^\pm,
s_kx_{k+1}^\pm x_k^\pm,u_kx_k^\pm)}
{\Gamma(s_ku_kx_{k+1}^\pm,x_k^{\pm 2}, t^2\prod_{l=1}^ks_l^2\,
u_kx_k^\pm)} \beta^{(0)}(x_1,t)\prod_{l=1}^m\frac{dx_l}{x_l}.
\nonumber\end{eqnarray}

Используя (\ref{a_1}) и (\ref{b_1}) в качестве начальных функций
$\alpha^{(0)}(x,t)$ и $\beta^{(0)}(x,t)$ и подставляя
(\ref{a-m}) и (\ref{b-m}) в (\ref{bint-pair}), мы получаем тождество
\begin{eqnarray}\nonumber
\lefteqn{ \kappa^{m-1}\prod_{k=1}^m
\frac{\Gamma(t^2\prod_{l=1}^{k}s_l^2)}
{\Gamma(s_k^2,t^2\prod_{l=1}^{k-1}s_l^2)}
\int_{\mathbb{T}^m}\frac{\prod_{r=0}^2\Gamma(tt_rx_1^\pm)}
{\Gamma(tt_0t_1t_2x_1^\pm)}  } &&
\\ \nonumber && \makebox[2em]{} \times
\prod_{k=1}^m\frac
{\Gamma(t^2\prod_{l=1}^{k-1}s_l^2\, s_ku_kx_{k+1}^\pm,
s_kx_{k+1}^\pm x_k^\pm,u_kx_k^\pm)}
{\Gamma(s_ku_kx_{k+1}^\pm,x_k^{\pm 2}, t^2\prod_{l=1}^ks_l^2\,
u_kx_k^\pm)} \prod_{l=1}^m\frac{dx_l}{x_l}
\\ \nonumber &&
=\frac{1}{\Gamma(t^2)}\prod_{0\leq r<j\leq 2}\frac{\Gamma(t^2t_rt_j)}
{\Gamma(t_rt_j)} \int_{\mathbb{T}}
\frac{\Gamma(t\prod_{l=1}^ms_l\, x_{m+1}^\pm x^\pm)
\prod_{r=0}^2\Gamma(t_rx^\pm)}{\Gamma(x^{\pm 2},t^2t_0t_1t_2x^\pm)}
\\ && \makebox[2em]{} \times
\prod_{k=1}^m \frac{\Gamma(t\prod_{l=1}^{k-1}s_l\, u_kx^\pm)}
{\Gamma(t\prod_{l=1}^{k-1}s_l\, s_k^2u_kx^\pm)}\frac{dx}{x},
\label{id-seq}\end{eqnarray}
которое может быть переписано как (\ref{an-trans}) для $m=1$.
В некотором смысле, соотношение (\ref{id-seq}) представляет собой
обобщение высокого уровня тождеств типа Роджерса-Рамануджана
для многократных рядов, которые были выведены в статье \cite{and:multiple}.
Интегральные аналоги некоторых из тождеств Роджерса-Рамануджана
были рассмотрены в работе \cite{gis:variants}. Публикации по исследованиям
цепочек Бэйли для рядов достаточно многочисленны, но насколько известно
автору интегральные аналоги этих цепочек совершенно не рассматривались
в литературе. Некоторые примеры однократных преобразований
типа Бэйли для интегралов были известны в течение долгого времени.
Фактически, все изменения порядков интегрирований в многократных
интегралах, использованных в исследованиях \cite{and:short,den-gus:beta,gus:some2,
gus-rak:beta,die-spi:selberg,spi:integrals},
могут быть интерпретированы как интегральные обобщения изменений
порядка суммирований для рядов \cite{and:bailey,war:50}.
Общая идея использования изменения порядка интегрирований для
получения полезной информации об интегралах восходит к методам
вывода обычного бета-интеграла, предложенных Пуассоном и Якоби \cite{aar:special}.

Известно, что для определенных пар Бэйли существует несколько различных
цепочечных правил подстановки (или лемм Бэйли). В этом случае
совокупность получающихся тождеств имеет более сложную структуру,
которая называется решеткой \cite{aab:bailey} или деревом Бэйли
\cite{and:bailey}. Оказывается, что в нашем случае можно расширить
интегральную цепочку Бэйли, описанную в Теореме 1, до структуры дерева
благодаря дополнительной интегральной лемме Бэйли.

\begin{theorem}\label{ibl-2}
Для любых двух функций $\alpha(z,t)$ и $\beta(z,t)$
образующих интегральную эллиптическую пару Бэйли по отношению к параметру $t$,
$|t|<1$, новые функции, определенные для $w\in\mathbb{T}$ как
\begin{eqnarray}
\alpha'(w,t)=\kappa \frac{\Gamma(s^2t^2,uw^\pm)}
{\Gamma(s^2,t^2,w^{\pm2},t^2s^2uw^{\pm})}
\int_{\mathbb{T}} \frac{\Gamma(t^2sux^\pm,sw^\pm x^\pm)}
{\Gamma(sux^\pm)}\,\alpha(x,st)\frac{dx}{x}
\label{a2}\end{eqnarray}
и
\begin{equation}
\beta'(w,t)=\frac{\Gamma(tuw^\pm)}
{\Gamma(ts^2uw^\pm)}\,\beta(w,st),
\label{b2}\end{equation}
также образуют интегральную эллиптическую пару Бэйли
по отношению к $t$. Здесь $s$ и $u$ являются произвольными комплексными
параметрами (они не связаны с $s,u$, фигурирующими в Теореме 1),
удовлетворяющими ограничениям $|s|,|u|<1, |pq|<|t^2s^2u|$.
\end{theorem}
{\bf Доказательство.}
Проверка справедливости этого утверждения проводится так же как и в предыдущей
теореме. Для этого необходимо заменить $\alpha(w,t)$ и $\beta(w,t)$ в
равенстве (\ref{bint-pair}) штрихованными выражениями, подставить
в него соотношение (\ref{a2}), изменить порядок интегрирований и применить
формулу (\ref{ell-int}). Мы опускаем детали рассмотрения в виду их простоты.
\hfill{Q.E.D.}

\smallskip

В работе \cite{spi:bailey} была найдена только одна лемма Бэйли для
эллиптических гипергеометрических рядов. Теорема (\ref{ibl-1}) была получена
по некоторой эвристической аналогии с ней. Естественно ожидать, что
существует некоторый предел, который позволил бы перейти от
теоремы (\ref{ibl-1}) к соответствующей цепочке Бэйли для рядов. Аналогично,
теорема (\ref{ibl-2}) также должна иметь некоторый аналог для рядов, который скорее
всего совпадет со второй леммой Бэйли \re{a2-ser}, описанной в предыдущей
главе. Соответствующее эллиптическое дерево Бэйли для рядов отличается
от дерева преобразований для $q$-ги\-пер\-гео\-ме\-три\-чес\-ких рядов
работ \cite{and:bailey,and-ber:bailey}.

В качестве иллюстрации, мы опишем одну частную последовательность
тождеств, порожденных сконструированным деревом Бэйли. Для этого
возьмем функции (\ref{a_1}) и (\ref{b_1}) и применим к ним преобразования
(\ref{a2}) и (\ref{b2}) с параметрами $s=s_1,u=u_1$.
К результирующим $\alpha'(w,t)$ и $\beta'(w,t)$ мы применяем
преобразования (\ref{a'}) и (\ref{b'}) с параметрами $s=s_2,u=u_2$.
В итоге мы получаем
\begin{eqnarray}\nonumber
&& \alpha''(w,s_2t)=\kappa % \frac{(p;p)_\infty (q;q)_\infty}{4\pi i}
\frac{\Gamma(t^2s_1^2,tu_2w^\pm,u_1w^\pm)}
{\Gamma(s_1^2,t^2,w^{\pm2},ts_2^2u_2w^\pm,t^2s_1^2u_1w^\pm)}
\\ && \makebox[5em]{} \times
\int_{\mathbb{T}} \frac{\Gamma(t^2s_1u_1x^\pm,s_1w^\pm x^\pm)
\prod_{r=0}^2\Gamma(t_rx^\pm)}{\Gamma(x^{\pm2},s_1u_1x^\pm,
t^2s_1^2t_0t_1t_2x^\pm)}\frac{dx}{x},
\\ \nonumber &&
\beta''(w,s_2t)=\kappa % \frac{(p;p)_\infty (q;q)_\infty}{4\pi i}
\frac{\Gamma(t^2s_2^2,t^2s_1^2,t^2s_2u_2w^\pm)}
{\Gamma(s_2^2,t^2,s_2u_2w^\pm)}\prod_{0\leq r<j\leq 2}
\frac{\Gamma(t_rt_j)}{\Gamma(t^2s_1^2t_rt_j)}
\\ && \makebox[5em]{} \times
\int_{\mathbb{T}}\frac{\Gamma(s_2w^\pm x^\pm,tu_1x^\pm,u_2x^\pm)
\prod_{r=0}^2\Gamma(ts_1t_r x^\pm)}{\Gamma(x^{\pm2},t^2s_2^2u_2x^\pm,
ts_1^2u_1x^\pm,ts_1t_0t_1t_2x^\pm)}\frac{dx}{x}.
\label{ab''}\end{eqnarray}
Подставляя эти функции в основное соотношение
(\ref{bint-pair}) с параметром $t$ замененным на $s_2t$,
мы получаем равенство для двух интегралов
\begin{eqnarray} \nonumber
&& \int_{\mathbb{T}}\frac{\Gamma(s_2w^\pm x^\pm,tu_1x^\pm,u_2x^\pm)
\prod_{r=0}^2\Gamma(ts_1t_r x^\pm)}{\Gamma(x^{\pm2},t^2s_2^2u_2x^\pm,
ts_1^2u_1x^\pm,ts_1t_0t_1t_2x^\pm)}\frac{dx}{x}
\\ && \nonumber \makebox[4em]{}
=\kappa \prod_{0\leq r<j\leq 2}\frac{\Gamma(t^2s_1^2t_rt_j)}{\Gamma(t_rt_j)}
\frac{\Gamma(s_2^2,s_2u_2w^\pm)}{\Gamma(s_1^2,t^2s_2^2,t^2s_2u_2w^\pm)}
\\ \nonumber && \makebox[5em]{}\times
\int_{\mathbb{T}^2}\Gamma(s_1x_2^\pm x_1^\pm)
\frac{\Gamma(ts_2w^\pm x_2^\pm,tu_2x_2^\pm,u_1x_2^\pm)}
{\Gamma(x_2^{\pm2},ts_2^2u_2x_2^\pm,t^2s_1^2u_1x_2^\pm)}
\\ && \makebox[5em]{}\times
\frac{\Gamma(t^2s_1u_1x_1^\pm)\prod_{r=0}^2
\Gamma(t_rx_1^\pm)}
{\Gamma(x_1^{\pm2},s_1u_1x_1^\pm,t^2s_1^2t_0t_1t_2x_1^\pm)}
\frac{dx_1}{x_1}\frac{dx_2}{x_2}.
\label{ident1}\end{eqnarray}
Это тождество было получено при условиях, что все параметры лежат внутри
единичной окружности и что $|t^2s_1^2t_0t_1t_2|, |t^2s_1^2u_2|, |t^2s_2^2u_2|
>|pq|$. Однако, благодаря аналитичности, оно остается справедливым для
более широкой области значений параметров при условии, что мы заменим
$\mathbb{T}$ контурами интегрирования, которые охватывают последовательности
полюсов интегрируемой функции, сходящиеся к нулю, и исключают другие полюсы.

Нетрудно найти результат $m$-кратной итерации преобразований (\ref{a2}) и
(\ref{b2}) с различными параметрами $u_k,s_k$:
\begin{eqnarray}\nonumber
&& \alpha^{(m)}(x_{m+1},t)=\kappa^m\prod_{k=1}^m
\frac{\Gamma(\prod_{l=k}^ms_l^2t^2)}
{\Gamma(s_k^2,\prod_{l=k+1}^ms_l^2t^2)}
\\  && \makebox[4em]{} \times
\int_{\mathbb{T}^m}\prod_{k=1}^m\frac{\Gamma(u_kx_{k+1}^\pm,
\prod_{l=k+1}^ms_l^2t^2s_ku_kx_k^\pm,s_kx_{k+1}^\pm x_k^\pm)}
{\Gamma(x_{k+1}^{\pm2},\prod_{l=k}^ms_l^2t^2u_kx_{k+1}^\pm,s_ku_kx_k^\pm)}
\nonumber \\  && \makebox[6em]{} \times
\alpha^{(0)}\left(x_1,\prod_{l=1}^ms_l\, t\right)
\frac{dx_1}{x_1}\cdots\frac{dx_m}{x_m},
\label{a_m}\\
&& \beta^{(m)}(w,t)=\prod_{k=1}^m\frac{\Gamma(t\prod_{l=k+1}^ms_lu_kw^\pm)}
{\Gamma(t\prod_{l=k+1}^ms_ls_k^2u_kw^\pm)}\,
\beta^{(0)}\left(w,\prod_{l=1}^ms_l\, t\right).
\label{b_m}\end{eqnarray}
Применение к этим функциям преобразований (\ref{a'}) и (\ref{b'}) с
параметрами $s=s_{m+1}, u=u_{m+1}$ дает:
\begin{eqnarray}\nonumber
&& \alpha'(x_{m+1},s_{m+1}t)=\frac{\Gamma(tu_{m+1}x_{m+1}^\pm)}
{\Gamma(ts_{m+1}^2u_{m+1}x_{m+1}^\pm)}\, \alpha^{(m)}(x_{m+1},t),
\\ &&
\beta'(w,s_{m+1}t)=\kappa\frac{\Gamma(t^2s_{m+1}^2,t^2s_{m+1}u_{m+1}w^\pm)}
{\Gamma(s_{m+1}^2,t^2,s_{m+1}u_{m+1}w^\pm)}
\nonumber \\ && \makebox[4em]{} \times
\int_{\mathbb{T}} \frac{\Gamma(s_{m+1}w^\pm x^\pm,u_{m+1}x^\pm)}
{\Gamma(x^{\pm2},t^2s_{m+1}^2u_{m+1}x^\pm)}\,\beta^{(m)}(x,t)\frac{dx}{x}
\nonumber \\ && \makebox[4em]{}
=\kappa\int_{\mathbb{T}}\Gamma(s_{m+1}tw^\pm x_{m+1}^\pm)
\,\alpha'(x_{m+1},s_{m+1}t)\frac{dx_{m+1}}{x_{m+1}}.
\label{ident4}\end{eqnarray}
Мы подставляем в последнее соотношение выражения (\ref{a_m}) и
(\ref{b_m}) с функциями (\ref{a_1}) и (\ref{b_1}) выступающими
в качестве $\alpha^{(0)}(w,t)$ и $\beta^{(0)}(w,t)$.
В результате мы получаем тождество
\begin{eqnarray}\nonumber
\lefteqn{ \int_{\mathbb{T}}\frac{\Gamma(s_{m+1}w^\pm x^\pm,u_{m+1}x^\pm)}
{\Gamma(x^{\pm2},t^2s_{m+1}^2u_{m+1}x^\pm)}
\prod_{k=1}^m\frac{\Gamma(t\prod_{l=k+1}^ms_lu_kx^\pm)}
{\Gamma(t\prod_{l=k+1}^ms_ls_k^2u_kx^\pm)} } &&
\\ &&
\times \frac{\prod_{r=0}^2\Gamma(t\prod_{k=1}^ms_k t_rx^\pm)}
{\Gamma(t\prod_{k=1}^ms_kt_0t_1t_2x^\pm)}\frac{dx}{x}
=\frac{\kappa^m\Gamma(s_{m+1}^2,t^2,s_{m+1}u_{m+1}w^\pm)}
{\Gamma(t^2s_{m+1}^2,t^2\prod_{k=1}^ms_k^2,t^2s_{m+1}u_{m+1}w^\pm)}
\nonumber \\ &&
\times \prod_{k=1}^m\frac{\Gamma(\prod_{l=k}^ms_l^2t^2)}
{\Gamma(s_k^2,\prod_{l=k+1}^ms_l^2t^2)}\prod_{0\leq r<j\leq 2}
\frac{\Gamma(t^2\prod_{k=1}^ms_k^2 t_rt_j)}{\Gamma(t_rt_j)}\,
\nonumber \\ &&
\times \int_{\mathbb{T}^{m+1}}
\frac{\Gamma(ts_{m+1}w^\pm x_{m+1}^\pm,tu_{m+1}x_{m+1}^\pm)
\prod_{r=0}^2\Gamma(t_rx_1^\pm)}{\Gamma(ts_{m+1}^2u_{m+1}x_{m+1}^\pm,
x_1^{\pm 2}, t^2\prod_{k=1}^ms_k^2t_0t_1t_2x_1^\pm)}
\nonumber \\ &&\makebox[1em]{} \times
\prod_{k=1}^m\frac{\Gamma(u_kx_{k+1}^\pm,
\prod_{l=k+1}^ms_l^2t^2s_ku_kx_k^\pm,s_kx_{k+1}^\pm x_k^\pm)}
{\Gamma(x_{k+1}^{\pm2},\prod_{l=k}^ms_l^2t^2u_kx_{k+1}^\pm,s_ku_kx_k^\pm)}
\frac{dx_1}{x_1}\cdots\frac{dx_{m+1}}{x_{m+1}}.
\label{identfin}\end{eqnarray}
Для $m=1$ это равенство сводится к соотношению (\ref{ident1}).

Как показано в работе \cite{spi-war:inversions}, две интегральные леммы Бэйли,
описанные выше, связаны друг с другом обращением интегрального
оператора, фигурирующего в определении пар Бэйли
(данные результаты не вошли в настоящую диссертацию ввиду позднего
завершения работы над ними). Более того,
в этой же работе были предложены обобщения интегрального преобразования,
стоящего за парами Бэйли, на системы корней $A_n$ и $C_n$.
Так же как и в одномерном случае эти преобразования приводят к
бесконечным последовательностям тождеств для эллиптических
гипергеометрических интегралов следствия которых еще предстоит
изучить. В частности, представляет интерес прояснение возможной
связи обобщенных интегральных преобразований Бэйли с преобразованиями
симметрии для интегралов разных размерностей, описанных в статьях
\cite{rai:trans} и \cite{tar-var:duality}.

\section{Многомерные эллиптические бета-интегралы для $C_n$ корневой системы}

\subsection{Многопараметрический $C_n$ интеграл (тип I)}

Многопараметрический эллиптический бета-интеграл для
$C_n$ системы корней был предложен ван Диехеном и автором в
\cite{die-spi:selberg}. Для его описания необходимо $n$ переменных
$z=(z_1,\ldots,z_n)\in\C^n$ и $2n+3$ комплексных параметра
$t=(t_1,\ldots,t_{2n+3})$. Составим подынтегральная функцию
\begin{equation}
\rho(z,t;C_n)=\prod_{1\leq i<j\leq n}\frac{1}{\Gamma(z_i^\pm z_j^\pm)}
\prod_{i=1}^n\frac{\prod_{m=1}^{2n+3}\Gamma(t_mz_i^\pm)}
{\Gamma(z_i^{\pm 2},Az_i^\pm)}\frac{\prod_{m=1}^{2n+3}
\Gamma(At_m^{-1})}{\prod_{1\leq m<s\leq 2n+3}\Gamma(t_mt_s)},
\label{kernel-C}\end{equation}
где $A=\prod_{m=1}^{2n+3}t_m$. Пусть
$$
\P=\{t_mq^jp^k,\, A^{-1}q^{j+1}p^{k+1}
\}_{m=1,\ldots, 2n+3,\, j,k \in\N}
$$
обозначает множество точек на комплексной плоскости, в которых находятся
последовательности полюсов (\ref{kernel-C}) для всех $z_i$, сходящиеся к нулю,
а $C$ обозначает положительно ориентированный контур, разделяющий $\P$ и $\P^{-1}$;
$dz/z\equiv \prod_{i=1}^ndz_i/z_i$.

\begin{theorem}
\begin{equation}\label{ell-int-C}
\int_{C^n}\rho(z,t;C_n)\frac{dz}{z}=\frac{2^nn!(2\pi i)^n}
{(q;q)_\infty^n(p;p)_\infty^n}.
\end{equation}
\end{theorem}
{\bf Доказательство.}
Ядро интеграла $\rho(z,t;C_n)$ удовлетворяет $q$-разностному уравнению, аналогичному (\ref{eqn}),
\begin{eqnarray} \nonumber
&& \rho(z,qt_1,t_2,\ldots,t_{2n+3};C_n)-\rho(z,t;C_n)
\\ && \makebox[4em]{}
=\sum_{i=1}^n\left(g_i(z_1,...,q^{-1}z_i,\ldots,z_n,t)-g_i(z,t)\right),
\label{eqn-C}\end{eqnarray}
где
\begin{equation}
g_i(z,t)=\rho(z,t;C_n)\prod_{j=1,\neq i}^n\frac{\theta(t_1z_j^\pm;p)}
{\theta(z_iz_j^\pm;p)}\frac{\prod_{j=1}^{2n+3}\theta(t_jz_i;p)}
{\prod_{j=2}^{2n+3}\theta(t_1t_j;p)}
\frac{\theta(t_1A;p)}{\theta(z_i^2,Az_i;p)}\frac{t_1}{z_i}.
\label{g-C}\end{equation}
Разделив уравнение (\ref{eqn-C}) на $\rho(z,t;C_n)$, мы получаем
\begin{eqnarray}\nonumber
\lefteqn{\prod_{i=1}^n\frac{\theta(t_1z_i^\pm;p)}{\theta(Az_i^\pm;p)}
\prod_{j=2}^{2n+3}\frac{\theta(At_j^{-1};p)}{\theta(t_1t_j;p)}-1
=\frac{t_1\theta(t_1A;p)}{\prod_{j=2}^{2n+3}\theta(t_1t_j;p)}
\sum_{i=1}^n\frac{1}{z_i\theta(z_i^2;p)}
}&&
\\ &&
\times
\prod_{j=1,\neq i}^n\frac{\theta(t_1z_j^\pm;p)}{\theta(z_iz_j^\pm;p)}
\left(z_i^{2n+2}\frac{\prod_{j=1}^{2n+3}\theta(t_jz_i^{-1};p)}
{\theta(Az_i^{-1};p)} - \frac{\prod_{j=1}^{2n+3}\theta(t_jz_i;p)}
{\theta(Az_i;p)}\right).
\label{eqn-exp-C}\end{eqnarray}
Обе стороны этого равенства инвариантны относительно преобразования
$z_1\to pz_1$
и имеют одинаковые полюсы (сингулярности в точках
$z_1=z_j,z_j^{-1}, j=2,\ldots, n,$ и $z_1=\pm p^{k/2},\ k\in\N$
в правой части сокращаются) и их вычеты. Например,
$$
\lim_{z_1\to A}\theta(Az_1^{-1};p)\, (\mbox{обе стороны}) =
\frac{\theta(t_1A^\pm;p)}{\theta(A^2;p)}\prod_{j=1,\neq i}^n\frac{
\theta(t_1z_i^\pm;p)}{\theta(Az_i^\pm;p)}
\prod_{j=2}^{2n+3}\frac{\theta(At_j^{-1};p)}{\theta(t_1t_j;p)}.
$$
Поэтому, функции в двух сторонах равенства (\ref{eqn-exp-C})
отличаются только аддитивной константой, не зависящей от $z_1$, которая равна
нулю поскольку при $z_1=t_1$ равенство (\ref{eqn-exp-C}) легко проверяется.

Интегрируя (\ref{eqn-C}) по переменным $z\in C^n$, получаем
\begin{equation}
I(qt_1,t_2,\ldots,t_{2n+3})-I(t)=\sum_{i=1}^n
\left(\int_{C^{i-1}\times(q^{-1}C)\times C^{n-i}}-\int_{C^n}\right)g_i(z,t)\frac{dz}{z},
\label{int-eqn-C}\end{equation}
где $I(t)=\int_{C^n}\rho(z,t;C_n)dz/z$ и
$q^{-1}C$ обозначает контур $C$ растянутый относительно нулевой точки.

Полюсы функции (\ref{g-C}) по переменным $z_i$ сходятся к нулю по точкам
$z_i=t_mq^jp^k,$  $A^{-1}q^{j}p^{k+1}$ и уходят на бесконечность
при $z_i=t_m^{-1}q^{-1-j}p^{-k},Aq^{-j}p^{-k-1}$, где $m=1,\ldots,2n+3$,
$j,k\in\mathbb{N}$. При $|t_m|<1$ и $|p|<|A|$ область
$1\leq |z_i|\leq |q|^{-1}$ не содержит полюсов, так что можно
положить $C=\T$, деформировать $q^{-1}\mathbb{T}$ к $\mathbb{T}$
в (\ref{int-eqn-C}) и получить равенство $I(qt_1,t_2,\ldots,t_{2n+3})=I(t)$.

Повторяя почти дословно процедуру аналитического продолжения,
использовавшуюся в $n=1$ случае, получаем, что $I$ есть постоянная,
зависящая только от $p$ и $q$. Ее значение находится с помощью анализа
структуры вычетов, проведенного в  \cite{die-spi:elliptic},
и оно равно правой части (\ref{ell-int-C}). Пример такого анализа подробно рассмотрен
в конце этой главы.
\hfill{Q.E.D.}

\begin{remark}
Формула (\ref{ell-int-C}) была предложена в работе \cite{die-spi:elliptic},
с нетривиальными частично подтверждающими аргументами. Она была полностью
доказана в \cite{rai:trans} после редукции интеграла к вычислению
довольно сложных последовательностей детерминантов на плотных множествах
значений параметров. Приведенное доказательство \cite{spi:short}
значительно короче и основывается на полностью элементарных рассуждениях.
\end{remark}

\subsection{Эллиптический аналог интеграла Сельберга (тип II)}

Возьмем шесть комплексных параметров $t$ и $t_m$, $m=0,\dots,4$ и введем
ядро многомерного $C_n$ (или, как его еще классифицируют, $BC_n$) интеграла
\begin{eqnarray}
\Delta_n^{II} ({\bf z};p,q)&=&\frac{1}{(2\pi i)^n}\prod_{1\leq j<k\leq n}
\frac{\Gamma(tz_jz_k,tz_jz_k^{-1},tz_j^{-1}z_k,tz_j^{-1}z_k^{-1};p,q)}
{\Gamma(z_jz_k,z_jz_k^{-1},z_j^{-1}z_k,z_j^{-1}z_k^{-1};p,q)} \nonumber \\
&& \times \prod_{j=1}^n\frac{\prod_{r=0}^4\Gamma(t_rz_j,t_rz_j^{-1};p,q)}
{\Gamma(z_j^2,z_j^{-2},
  B z_j, B z_j^{-1} ;p,q)} ,
\label{esintB}
\end{eqnarray}
где $B\equiv t^{2n-2}\prod_{s=0}^4t_s$.

\begin{theorem}\label{B:thm}
Пусть $|p|,|q|, |t|, |t_r| <1$ (где $r=0,\ldots , 4$)
и $|pq|<|t^{2n-2}\prod_{s=0}^4t_s|$. Тогда эллиптический аналог интеграла Сельберга
имеет вид
\begin{eqnarray}
\lefteqn{\int_{\T^n}\Delta_n^{II} ({\bf z};p,q)\frac{dz_1}{z_1}
\cdots \frac{dz_n}{z_n}=} && \nonumber \\
&& \frac{2^n n!}{(p;p)_\infty^n(q;q)_\infty^n}
\prod_{j=1}^n\frac{\Gamma(t^j;p,q)}{\Gamma(t;p,q)}
\frac{\prod_{0\leq r<s\leq 4}\Gamma(t^{j-1}t_rt_s;p,q)}
{\prod_{r=0}^4\Gamma(t^{1-j}t_r^{-1}B;p,q)} .\label{SintB}
\end{eqnarray}
\end{theorem}
{\bf Доказательство.}
Основная идея состоит в рассмотрении составного интеграла
\begin{eqnarray}
&& \frac{1}{(2\pi i)^{2n-1}}\int_{\T^n}\int_{\T^{n-1}}
\prod_{1\leq j<k\leq n}
\Gamma^{-1}(z_jz_k,z_jz_k^{-1},z_j^{-1}z_k,z_j^{-1}z_k^{-1};p,q) \nonumber \\
&& \times \prod_{j=1}^n\frac{\prod_{r=0}^{4}\Gamma(t_rz_j,t_rz_j^{-1};p,q)}
{\Gamma(z_j^2,z_j^{-2},
   z_jt^{n-1}\prod_{0\leq s\leq 4}t_s,
   z_j^{-1}t^{n-1}\prod_{0\leq s\leq 4}t_s;p,q)} \nonumber \\
&&\times \prod_{\stackrel{1\leq j\leq n}{1\leq k\leq n-1}}
\Gamma (t^{1/2}z_jw_k,t^{1/2}z_jw_k^{-1},
t^{1/2}z_j^{-1}w_k,t^{1/2}z_j^{-1}w_k^{-1};p,q) \nonumber \\
&&\times
\prod_{1\leq j<k\leq n-1}
\Gamma^{-1}(w_jw_k,w_jw_k^{-1},w_j^{-1}w_k,w_j^{-1}w_k^{-1};p,q)\nonumber \\
&& \times \prod_{j=1}^{n-1}
\frac{\Gamma(w_jt^{n-3/2}\prod_{0\leq s\leq 4}t_s,
   w_j^{-1}t^{n-3/2}\prod_{0\leq s\leq 4}t_s  ;p,q)}
{\Gamma(w_j^2,w_j^{-2},
   w_jt^{2n-3/2}\prod_{0\leq s\leq 4}t_s,
   w_j^{-1}t^{2n-3/2}\prod_{0\leq s\leq 4}t_s;p,q)} \nonumber\\
&&\makebox[1em]{}
\times \frac{dw_1}{w_1}\cdots\frac{dw_{n-1}}{w_{n-1}}
\frac{dz_1}{z_1}\cdots\frac{dz_n}{z_n}\label{compint}
\end{eqnarray}
с параметрами, удовлетворяющими указанным выше ограничениям.
Обозначим интеграл в левой части \re{SintB} как
$I_n^{II}(t,t_r;p,q)$. Интегрируя в составном интеграле \re{compint} по
$w$-циклам с помощью формулы \re{ell-int-C}, получаем
\begin{equation}\label{compa}
\frac{2^{n-1}(n-1)!}{(p;p)_\infty^{n-1} (q;q)_\infty^{n-1}}
\frac{\Gamma^n (t;p,q)}{\Gamma (t^n;p,q)}
I_n^{II}(t,t_r;p,q).
\end{equation}
Аналогично, интегрируя  в \re{compint} по переменным $z_k$ с помощью
той же формулы \re{ell-int-C}, получаем
\begin{equation} \label{compb}
\frac{2^n n!}{(p;p)_\infty^n (q;q)_\infty^n}
\frac{\Gamma^{n-1} (t;p,q) \prod_{0\leq r< s\leq 4} \Gamma (t_rt_s;p,q)}
     {\prod_{r=0}^4 \Gamma (t^{n-1}t_r^{-1}\prod_{s=0}^4 t_s;p,q)}
I_{n-1}^{II}(t,t^{1/2}t_r;p,q).
\end{equation}
Сравнивая выражения \re{compa} и \re{compb}, получаем рекуррентное соотношение
для  $I_n^{II}$ по размерности $n$:
\begin{eqnarray}\label{recB}
\lefteqn{I_n^{II}(t,t_r;p,q)=} && \\
&& \frac{2n}{(p;p)_\infty (q;q)_\infty}
\frac{\Gamma(t^n;p,q)}{\Gamma(t;p,q)}
\frac{\prod_{0\leq r<s\leq 4}\Gamma(t_rt_s;p,q)}
{\prod_{r=0}^4\Gamma(t^{n-1}t_r^{-1}\prod_{s=0}^4t_s;p,q)}
I_{n-1}^{II}(t,t^{1/2}t_r;p,q) .\nonumber
\end{eqnarray}
Итерируя это соотношение начиная с известного значения
$I_n$ при $n=1$ \re{ell-int}, получаем
\begin{equation}\label{Bconstant}
I_n^{II}(t,t_r;p,q)=
\frac{2^n n!}{(p;p)_\infty^n(q;q)_\infty^n}
\prod_{j=1}^n\frac{\Gamma(t^j;p,q)}{\Gamma(t;p,q)}
\frac{\prod_{0\leq r<s\leq 4}\Gamma(t^{j-1}t_rt_s;p,q)}
{\prod_{r=0}^4\Gamma(t^{1-j}t_r^{-1}B;p,q)},
\end{equation}
что и требовалось показать.
\hfill{Q.E.D.}

\begin{remark}
Изложенная техника доказательства была использована Густафсоном для
доказательства многократного $q$-бета интеграла, соответствующего $p=0$
в \re{SintB}, \cite{gus:some2}. Последующее устранение
одного из параметров $t_4\to 0$ приводит к $q$-аналогу интеграла
Сельберга \cite{gus:generalization}. Аналогичная техника использовалась
несколько ранее Андерсоном \cite{and:short} для доказательства самого
интеграла Сельберга \cite{sel:bemerkninger}:
\begin{eqnarray}
\lefteqn{\int_0^1\cdots\int_0^1
\prod_{1\leq j\leq n} x_j^{\alpha-1} (1-x_j)^{\beta-1}
\prod_{1\leq j<k\leq n}|x_j-x_k|^{2\gamma}\; dx_1\cdots dx_n } && \nonumber \\
&& =\prod_{1\leq j\leq n}
\frac{\Gamma(\alpha+(j-1)\gamma) \Gamma(\beta+(j-1)\gamma) \Gamma (1+j\gamma)}
     {\Gamma (\alpha+\beta+(n+j-2)\gamma) \Gamma (1+\gamma)}, \label{Sint}
\end{eqnarray}
где $\text{Re}(\alpha ), \text{Re} (\beta) >0$ и $\text{Re}(\gamma ) >$
$- \min (1/n,$ $\text{Re} (\alpha )/(n-1),$  $\text{Re} (\beta)/(n-1) )$,
а $\Gamma(u)$ обозначает обычную гамма-функцию Эйлера.
Благодаря указанной связи, формула \re{SintB} описывает эллиптическое
расширение интеграла Сельберга. Другие примеры многократных
гипергеометрических интегралов, связанных с интегралом Сельбергом,
могут быть найдены в работах \cite{var,ner1,ner3,tar-var:selberg}
и приведенных в них ссылках, но эллиптические расширения этих интегралов
пока еще не найдены.
\end{remark}

\begin{remark}
Хорошо известно, что точно вычисляемые контурные интегралы по единичной
окружности от мероморфных функций могут быть интерпретированы как тождества
для постоянных членов разложения в ряды Фурье \cite{aar:special}.
В качестве примера укажем на доказательство в работе \cite{bre-zei:proof}
$q$-тождества Дайсона, предложенного Эндрюсом. Другой тип таких
тождеств, связанных с обычными и аффинными системами корней,
был предложен в виде набора гипотез Макдональдом \cite{mac:constant}
и Моррисом \cite{mor:constant}. Часть из них связана с многомерными
полиномами Якоби \cite{aom:jacobi}, Макдональда \cite{mac:constant}
и Корнвиндера \cite{kor:pol}. Для системы корней $BC_n$ соответствующее
тождество было доказано в \cite{kad:proof,gus:generalization},
для исключительных корневых систем см.
\cite{gar-gon:macdonald's,hab:q-conjecture}. В этом ракурсе, интеграл
\re{SintB} может интерпретироваться  как эллиптическое
тождество Макдональда-Морриса для постоянных членов \cite{die-spi:elliptic}.
\end{remark}

\subsection{$C_n$ интеграл, связанный с детерминантом Варнаара (тип III)}

Для комплексных параметров $t, t_{1,2,3}, x_i$, $i=1,\dots, n,$
определим функцию
\begin{eqnarray}\nonumber
&& \Delta^{III}(\mathbf{z};C_n) = \frac{1}{(2\pi i)^n}
\prod_{1\leq i<j\leq n}z_j\theta(z_iz_j^{-1},z_i^{-1}z_j^{-1};p)
\\ && \makebox[2em]{} \times
\prod_{i=1}^n\prod_{\nu=\pm 1}\frac{\Gamma(z_i^\nu x_i,z_i^\nu t_1,
z_i^\nu t_2,z_i^\nu t_3,z_i^\nu t/x_i)}
{\Gamma(z_i^{2\nu},z_i^\nu A)},
\label{delta-c3}\end{eqnarray}
где $A=tt_1t_2t_3q^{n-1}$.

\begin{theorem}
Наложим следующие ограничения на параметры: $|x_i|,$ $|t_k|<1,$
$|t|<|x_i|$, где $i=1,\ldots, n$, $k=1,2,3$, и $|pq|<|A|$. Тогда
\begin{eqnarray}\nonumber
&& \int_{\mathbb{T}^n}\Delta^{III}(\mathbf{z};C_n)
\frac{dz_1}{z_1}\cdots\frac{dz_n}{z_n}
=\frac{2^n}{(p;p)_\infty^n(q;q)_\infty^n}
\prod_{1\leq i<j\leq n}x_j\theta\left(x_i/x_j,t/x_ix_j;p\right)
\\ && \makebox[4em]{} \times
\Gamma^n(t)\prod_{i=1}^n\left(\frac{\prod_{1\leq r<s\leq 3}
\Gamma(t_rt_sq^{i-1})}{\Gamma(A/x_i,Ax_i/t)} \prod_{k=1}^3\frac{
\Gamma(x_it_k,tt_k/x_i)}{\Gamma(Aq^{1-i}/t_k)}\right).
\label{c3}\end{eqnarray}
\end{theorem}
{\bf Доказательство.}
Рассмотрим детерминант
\begin{eqnarray}\nonumber
&& \det_{1\leq i,j\leq n}\left(\int_\mathbb{T}
\Delta_E(z,x_i,t_1q^{n-j},t_2q^{j-1},t_3,tx_i^{-1})\frac{dz}{z}\right)
\\  && \makebox[2em]{}
=\frac{1}{(2\pi i)^n}\int_{\mathbb{T}^n}\frac{dz_1}{z_1}\cdots
\frac{dz_n}{z_n}\prod_{i=1}^n G_i(z_i)G_i(z_i^{-1})\; D(\mathbf{z}),
\label{det1}\end{eqnarray}
где $\Delta_E$ обозначает ядро эллиптического бета-интеграла \re{ell-int}:
\begin{equation}\label{weight}
\Delta_E(z,\mathbf{t}) = \frac{1}{2\pi i}
\frac{\prod_{m=0}^4\Gamma(zt_m, z^{-1}t_m; q,p)}
{\Gamma(z^2,z^{-2}, zA, z^{-1}A;q,p)}.
\end{equation}
Выражение, стоящее в правой части (\ref{det1}), возникает в результате
выноса знаков интегрирования за знак детерминанта, что приводит к указанному
многократному интегралу с
\begin{eqnarray*}
&& G_i(z_i)=\frac{\Gamma(z_ix_i,z_it_1,z_it_2,z_it_3,z_it/x_i)}
{\Gamma(z_i^2,z_iA)},
\\
&& D(\mathbf{z})=\det_{1\leq i,j\leq n}\left(\theta(z_it_1,
z_i^{-1}t_1;p;q)_{n-j}\theta(z_it_2,z_i^{-1}t_2;p;q)_{j-1}\right).
\end{eqnarray*}
Детерминант $D(\mathbf{z})$ может быть переписан в виде
\be
D(\mathbf{z})=\frac{\prod_{i=1}^n \theta(z_it_2,z_i^{-1}t_2;p;q)_{n-1}}
{t_2^{2\binom{n}{2}} q^{4\binom{n}{3}} }
\det_{1\leq i,j \leq n}
\left(\frac{\theta(z_it_1,z_i^{-1}t_1;p;q)_{n-j}}
{\theta(q^{2-n}/z_it_2,q^{2-n}z_i/t_2;p;q)_{n-j}}\right).
\nonumber \ee
В работе \cite{war:summation}, Варнаар вычислил эллиптическое обобщение
детерминанта Кратенталера  \cite{kra:major}:
\begin{eqnarray}
&& \det_{1 \leq i,j \leq n} \left(
\frac{ \theta(aX_i,ac/X_i;p;q)_{n-j} }
{ \theta(bX_i,bc/X_i;p;q)_{n-j} } \right)
= a^{\binom{n}{2}}q^{\binom{n}{3}}
\prod_{1\leq i<j\leq n}X_j\theta(X_iX_j^{-1},cX_i^{-1}X_j^{-1};p)
\nonumber \\ && \makebox[6em]{}
\times \prod_{i=1}^n \frac{\theta(b/a,abcq^{2n-2i};p;q)_{i-1}}
{\theta(bX_i,bc/X_i;p;q)_{n-1}}.
\label{e-kratt} \end{eqnarray}
Используя это тождество при $X_i=z_i$, $a=t_1, b=q^{2-n}/t_2$, и $c=1$,
находим
\begin{eqnarray}\nonumber
&& D(\mathbf{z})=(t_1t_2^2)^{\binom{n}{2}}q^{3\binom{n}{3}}
\prod_{1\leq i<j\leq n}z_j\theta(z_iz_j^{-1},z_i^{-1}z_j^{-1};p)
\\ && \makebox[4em]{} \times
\prod_{i=1}^n \theta(q^{2-n}/t_1t_2,t_1q^{n+2-2i}/t_2;p;q)_{i-1}.
\label{det-aux}\end{eqnarray}
В результате, детерминант (\ref{det1}) дает выражение пропорциональное
левой части рассматриваемого $C_n$-интеграла (\ref{c3}).

Подставим теперь результат вычисления эллиптического бета-интеграла
(\ref{ell-int}) в детерминант (\ref{det1}). Это дает
\begin{eqnarray}\nonumber
&& \mbox{(\ref{det1})} = \frac{2^n}{(q;q)_\infty^n(p;p)_\infty^n}
\\ && \makebox[1em]{}  \times
\prod_{i=1}^n\frac{\Gamma(x_it_1,t_1t/x_i,x_it_2,t_2t/x_i,
x_it_3,t,t_1t_2q^{n-1},t_1t_3q^{n-i},t_2t_3q^{i-1}, t_3t/x_i)}
{\Gamma(A/x_i,A/t_3,x_iA/t,Aq^{1-i}/t_2,Aq^{i-n}/t_1)}
\nonumber \\ && \makebox[1em]{}
\times \det_{1\leq i,j\leq n}\left(\theta(x_it_1,t_1t/x_i;p;q)_{n-j}
\theta(x_it_2,t_2t/x_i;p;q)_{j-1}\right).
\label{rhs}\end{eqnarray}
Благодаря равенству (\ref{e-kratt}), детерминант в последней строчке принимает вид
$$
(t_1t_2^2t)^{\binom{n}{2}}q^{3\binom{n}{3}}
\prod_{1\leq i<j\leq n}x_j\theta\left(x_i/x_j,t/x_ix_j;p\right)
\prod_{i=1}^n \theta(q^{2-n}/t_1t_2t,t_1q^{n+2-2i}/t_2;p;q)_{i-1}.
$$
Приравнивая получающееся выражение в (\ref{rhs}) правой части равенства
(\ref{det1}), получаем требуемое значение интеграла (\ref{c3}).
\hfill{Q.E.D.}

\smallskip

Этот $C_n$-интеграл несимметричен по $p$ и $q$, в отличие от всех
остальных интегралов рассмотренных выше.
Метод, использованный в его вычислении,
является естественным развитием техники детерминантных формул
для многократных
базисных гипергеометрических рядов, см. например, статью
Густафсона-Кратенталера \cite{gus-kra:determinant}, развитую Шлоссером
\cite{sch:summation,sch:nonterminating} и Варнааром \cite{war:summation}.
Аналогичные рассмотрения для вычисления некоторых интегралов
гипергеометрического типа проведены в работах Варченко \cite{var}
и Тарасова-Варченко \cite{tar-var:geometry}.
Интегралы, напрямую связанные с детерминантами, мы условно классифицируем
как интегралы третьего типа.

\section{Интегралы для $A_n$ корневой системы}

\subsection{Многопараметрический $A_n$ интеграл (тип I)}

Опишем теперь многопараметрический эллиптический бета-интеграл для
корневой системы $A_n$, предложенный автором в статье \cite{spi:integrals}.
Возьмем $z=(z_1,\ldots,z_n)\in\C^n$, определим переменную $z_{n+1}$
через соотношение $\prod_{i=1}^{n+1}z_i=1$, и введем  $2n+3$
комплексных параметра $t=(t_1,\ldots,t_{n+1})$ и $s=(s_1,\ldots,s_{n+2})$.
Ядро $A_n$-интеграла имеет вид
\begin{eqnarray}\nonumber
&& \rho(z,t,s;A_n)=
\prod_{i=1}^{n+1}\frac{\prod_{m=1}^{n+1}\Gamma(t_mz_i^{-1})
\prod_{j=1}^{n+2}\Gamma(s_jz_i)\, \Gamma(St_i)}
{\Gamma(TSz_i)\prod_{j=1}^{n+2}\Gamma(t_is_j)} % } &&
\\ && \makebox[4em]{}
\times \prod_{1\leq i<j\leq n+1}\frac{1}{\Gamma(z_iz_j^{-1},z_i^{-1}z_j)}
\frac{1}{\Gamma(T)}\prod_{j=1}^{n+2}\frac{\Gamma(STs_j^{-1})}{\Gamma(Ss_j^{-1})},
\label{kernel-A}\end{eqnarray}
где $T=\prod_{m=1}^{n+1}t_m$ и $S=\prod_{j=1}^{n+2}s_j$.
Эта функция имеет полюсы в точках
$$
z_i=\{t_m q^jp^k,\, (TS)^{-1}q^{j+1}p^{k+1}\}, \; i=1,\ldots, n,\quad
z_{n+1}^{-1}=z_1\cdots z_n= \{s_lq^{j}p^{k}\},
$$
с $m=1,\ldots, n+1,\,l=1,\ldots, n+2,\, j,k \in\N,$ сходящихся
к нулю, и точках
$$
z_i=\{s_l^{-1} q^{-j}p^{-k}\},\; i=1,\ldots, n,\quad
z_{n+1}^{-1}= \{t_m^{-1}q^{-j}p^{-k}, TSq^{-j-1}p^{-k-1}\},
$$
уходящих к бесконечности.

\begin{theorem}
Предположим, что $|t_m|, |s_l|<1$ и $|pq|<|TS|$. Тогда
\begin{equation}\label{ell-int-A}
\int_{\T^n}\rho(z,t,s;A_n)\frac{dz}{z}=\frac{(n+1)!(2\pi i)^n}
{(q;q)_\infty^n(p;p)_\infty^n}.
\end{equation}
\end{theorem}
{\bf Доказательство.}
Ядро $A_n$-интеграла  удовлетворяет $q$-разностному уравнению
аналогичному уравнению (\ref{eqn})
\begin{eqnarray}\nonumber
&& \rho(z,qt_1,t_2,\ldots,t_{n+1},s;A_n)-\rho(z,t,s;A_n)
\\ && \makebox[2em]{}
=\sum_{i=1}^n\left(g_i(z_1,\ldots,q^{-1}z_i,\ldots,z_n,t,s)-g_i(z,t,s)\right),
\label{eqn-A}\end{eqnarray}
где
\begin{equation}
\frac{g_i(z,t,s)}{\rho(z,t,s;A_n)}=
\prod_{j=1,\neq i}^{n+1}\frac{\theta(t_1z_j^{-1};p)}{\theta(z_iz_j^{-1};p)}
\prod_{j=1}^{n+2}\frac{\theta(z_is_j;p)}{\theta(t_1s_j;p)}
\frac{\theta(z_iT/t_1,TSt_1;p)}{\theta(T,TSz_i;p)}\frac{t_1}{z_i}.
\label{g-A}\end{equation}
Используя соотношение
\begin{eqnarray}\nonumber
&& \frac{\rho(\ldots,q^{-1}z_i,\ldots;A_n)}
{\rho(\ldots, z_i,\ldots;A_n)}=
\prod_{m=1}^{n+2}\frac{\theta(s_mz_{n+1};p)}{\theta(q^{-1}s_mz_i;p)}
\prod_{j=1}^{n+1}\frac{\theta(t_jz_i^{-1};p)}{\theta(q^{-1}t_jz_{n+1}^{-1};p)}
\frac{\theta(q^{-1}TSz_i;p)}{\theta(TSz_{n+1};p)}
\\ && \makebox[2em]{} \times
\prod_{j=1,\neq i}^n \frac{\theta(q^{-1}z_iz_j^{-1},q^{-1}z_jz_{n+1}^{-1};p)}
{\theta(z_i^{-1}z_j,z_j^{-1}z_{n+1};p)}
\frac{\theta(q^{-2}z_iz_{n+1}^{-1};p)z_i^2}
{\theta(z_iz_{n+1}^{-1};p)qz_{n+1}^2},
\label{eqn-A-zi}\end{eqnarray}
можно привести уравнение (\ref{eqn-A}) к виду
\begin{eqnarray}\nonumber
\lefteqn{\prod_{i=1}^{n+1}\frac{\theta(t_1z_i^{-1};p)}{\theta(TSz_i;p)}
\frac{\theta(t_1S;p)}{\theta(T;p)}\prod_{j=1}^{n+2}
\frac{\theta(TSs_j^{-1};p)}{\theta(t_1s_j;p)}-1 } &&
\\ \nonumber &&
=\frac{t_1\theta(t_1TS;p)}{\theta(T;p)\prod_{j=1}^{n+2}\theta(t_1s_j;p)}
\sum_{i=1}^n\frac{1}{z_i\theta(z_iz_{n+1}^{-1};p)}
\prod_{j=1,\neq i}^n\frac{\theta(t_1z_j^{-1};p)}{\theta(z_iz_j^{-1};p)}
\\ \nonumber && \makebox[2em]{} \times
\Biggl(\frac{z_i^{n+1}}{z_{n+1}^{n+1}}
\prod_{j=1}^{n+2}\theta(s_jz_{n+1};p)
\frac{\prod_{j=1}^{n+1}\theta(t_jz_i^{-1};p)\theta(q^{-1}z_iTt_1^{-1};p)}
{\prod_{j=2}^{n+1}\theta(q^{-1}t_jz_{n+1}^{-1};p)\theta(TSz_{n+1};p)}
\\ && \makebox[3em]{}
\times \prod_{j=1,\neq i}^n\frac{\theta(q^{-1}z_jz_{n+1}^{-1};p)}
{\theta(z_jz_{n+1}^{-1};p)} - \prod_{j=1}^{n+2}\theta(s_jz_i;p)
\frac{\theta(t_1z_{n+1}^{-1},z_iTt_1^{-1};p)}{\theta(TSz_i;p)} \Biggr).
\label{eqn-exp-A}\end{eqnarray}
Обе стороны этого равенства инвариантны относительно преобразования
$z_1\to pz_1$ и имеют одинаковые полюсы по переменной $z_1$
и их вычеты. Действительно, нетрудно проверить совпадение вычетов
полюса $z_1=(TS)^{-1}$ и мы пропускаем их рассмотрение.
Рассмотрение вычетов для $z_{n+1}^{-1}=TS$ приводит к равенству
\begin{eqnarray}\nonumber
&& \sum_{i=1}^n \theta(q^{-1}z_iTt_1^{-1};p)
\prod_{j=1,\neq i}^n\frac{\theta(TSz_jq^{-1};p)}
{\theta(z_i^{-1}z_j;p)}
\prod_{j=2}^{n+1}\theta(t_jz_i^{-1};p)
\\ &&\makebox[2em]{}
=\theta(t_1^{-1}S^{-1};p)\prod_{j=2}^{n+1}\theta(q^{-1}t_jTS;p).
\label{A-function}\end{eqnarray}
Если домножить обе стороны на
$\prod_{1\leq i<j\leq n } z_j\theta(z_i/z_j;p),$
то они становятся антисимметричными голоморфными тета-функциями
переменных $z_1,\dots,z_{n-1}$ обладающими $A_{n-1}$ симметрией и
равными нулю при $z_{n}=z_i, i<n$. Но любая такая функция
пропорциональна указанному множителю. Поэтому левая часть
(\ref{A-function}) не должна зависеть от $z_i$. Полагая
$z_1=q(TS)^{-1}$, получаем правую часть равенства.

Сумма вычетов полюсов в точках $z_{n+1}=t_kq^{-1},\, k=2,\ldots,n+1,$
пропорциональна выражению
\begin{equation}
\sum_{i=1}^n \frac{\prod_{j=2,\neq k}^{n+1}\theta(t_j/z_i;p)}
{\prod_{j=1,\neq i}^n \theta(z_j/z_i;p)}\theta(qt_1/z_iT;p),
\label{zn+1=tj/q}\end{equation}
которое есть $TS=q/t_k$ подслучай (\ref{A-function}), и поэтому оно равно нулю.

Полюсы $z_1=z_j, z_j^{-1},\, j=2,\ldots, n,$ в правой части  равенства
(\ref{eqn-exp-A}) сокращают друг друга. Наконец, сравнение вычетов
полюсов $z_1=z_{n+1}$ приводит к равенству
\begin{equation}
\sum_{i=1}^n \theta(z_iT/qt_1;p)
\prod_{j=1,\neq i}^n\frac{\theta(z_j/qz_1;p)}{\theta(z_j/z_i;p)}
\prod_{j=2}^{n+1}\frac{\theta(t_j/z_i;p)}
{\theta(t_j/qz_1;p)}=\theta(z_1T/t_1;p).
\label{z1=zn+1}\end{equation}
После домножения обеих сторон на
$\prod_{1\leq i<j\leq n } z_j \theta(z_i/z_j;p)$  получаем
голоморфные $A_{n-2}$ антисимметричные функции
$z_2,\ldots,z_{n-1}$ ($z_1$ рассматривается как параметр),
равные нулю при $z_j=z_1,z_n$ и поэтому они равны друг другу
с точностью до множителя, не зависящего от $z_2,\dots,z_{n-1}$. Для
$z_2=qz_1$, равенство (\ref{z1=zn+1}) верно и поэтому оно справедливо
и в общем случае.

Из перечисленных свойств следует, что разница функций с двух сторон
(\ref{eqn-exp-A}) не зависит от $z_1$. Полагая $z_1=t_1$, мы видим, что
эта разница равна нулю и равенство (\ref{eqn-exp-A}) верно в общем случае.

Интегрируя равенство (\ref{eqn-A}) по переменным $z\in \T^n$, получаем
\begin{equation}
I(qt_1,t_2,\ldots,t_{n+1},s)-I(t,s)=\sum_{i=1}^n
\left(\int_{\T^{i-1}\times(q^{-1}\T)\times \T^{n-i}}-\int_{\T^n}\right)
g_i(z,t,s)\frac{dz}{z},
\label{int-eqn-A}\end{equation}
где $I(t,s)=\int_{\T^n}\rho(z,t,s;A_n)dz/z$.
Полюсы функции (\ref{g-A}) при $n>1$ находятся в точках
\begin{eqnarray*}
&& z_i=t_mq^jp^k, \; m=1,\ldots,n+1, \quad (TS)^{-1}q^{j}p^{k+1},\quad \\
&& z_{n+1}^{-1}=s_mq^jp^j,\; m=1,\ldots,n+2,
\end{eqnarray*}
сходящихся к нулю, и
\begin{eqnarray*}
&& z_i=s_m^{-1}q^{-j-1}p^{-k},\; m=1,\ldots,n+2,\quad \\
&& z_{n+1}^{-1}=t_1^{-1}q^{-j-1}p^{-k},\quad t_m^{-1}q^{-j}p^{-k},\;
m=2,\ldots,n+1,\quad TSq^{-j-1}p^{-k-1},
\end{eqnarray*}
уходящих к бесконечности при $j,k\in\mathbb{N}$.
Приняв ограничение $|t_m|<|q|,\, m=2,\ldots,n+1,$ и $|p|<|TS|$, мы видим,
что область $1\leq |z_i|\leq |q|^{-1}$ не содержит полюсов и мы
можем деформировать $q^{-1}\T$ к $\T$ в (\ref{int-eqn-A}), что приводит
к нулю в правой части. Таким образом, мы получаем равенство
$I(qt_1,t_2,\ldots,t_{n+1},s)=I(t,s)$.

Подходящая деформация контура интегрирования позволяет доказать также
инвариантность относительно растяжения $t_1\to pt_1$ и итерировать
определенным образом $q,p$-растяжения $t_1$. Процедура аналитического
продолжения, аналогичная $C_n$-случаю, показывает, что  $I(t,s)$ не
зависит от $t_m$.  Другой способ добиться такого же результата состоит
в разложении интеграла \re{ell-int-A} в ряд по параметру $p$, каждый
член которого оказывается аналитической функцией $t_1$ при $|t_1|<1$
инвариантной относительно преобразования $t_1\to qt_1$. Поскольку точка
накопления последовательности $t_1q^j,\, j\in \N,$ лежит в области аналитичности
каждого члена ряда по отдельности (только условие сходимости всего ряда в
целом накладывает ограничения на значения $|t_1|$ снизу), все они не зависят
от $t_1$ и, таким образом, от всех $t_m$.

Интеграл (\ref{ell-int-A}) может быть переписан в
$t_m\leftrightarrow s_m$ симметричном виде
\begin{eqnarray}\nonumber
\lefteqn{ \int_{\T^n} \frac{\prod_{i=1}^{n+1}\prod_{m=1}^{n+2}
\Gamma(t_mz_i^{-1},s_jz_i)}
{\prod_{1\leq i<j\leq n+1}\Gamma(z_iz_j^{-1},z_i^{-1}z_j)}\frac{dz}{z}
} && \\ &&
= \frac{(n+1)!(2\pi i)^n}{(q;q)_\infty^n(p;p)_\infty^n}
\prod_{j=1}^{n+2}\Gamma(Ss_j^{-1},At_j^{-1})
\prod_{j,k=1}^{n+2}\Gamma(t_js_k).
\label{ell-int-A-sym}\end{eqnarray}
где $S=\prod_{j=1}^{n+2}s_j, A=\prod_{j=1}^{n+2}t_j,
AS=pq$. Поэтому $I(t,s)$ так же не зависит и от $s_m$.
Явный вид постоянной $I=I(q,p)$, имеющей вид правой части (\ref{ell-int-A}),
находится с помощью анализа вычетов, проведенного автором в
работе \cite{spi:integrals}.
\hfill{Q.E.D.}

\begin{remark}
Точная формула интегрирования (\ref{ell-int-A}) была предложена автором в работе
\cite{spi:integrals}, где был представлен ряд аргументов частично подтверждающих ее справедливость.
Первое полное доказательство было получено Рэйнсом в \cite{rai:trans} с помощью
той же самой детерминантной техники, использовавшейся для вывода $C_n$ эллиптического
бета-интеграла типа I. Приведенное доказательство найдено автором в работе
\cite{spi:short} и оно значительно проще.
\end{remark}

\subsection{Обобщение смешанных  $A_n$ и $C_n$ интегралов Густафсона (тип II)}

Исходя из формулы \re{ell-int-A} в этом и следующем параграфах мы выведем
несколько новых различных эллиптических бета-интегралов для системы
корней $A_n$, смешанных с $C_n$-интегралами. Введем обозначение
\begin{eqnarray}\nonumber
&& \Delta^{II}({\bf z};A_n)=
\frac{1}{(2\pi i)^n}\prod_{1\leq i<j\leq n+1}
\frac{\Gamma(tz_iz_j,sz_i^{-1}z_j^{-1})}
{\Gamma(z_iz_j^{-1},z_i^{-1}z_j)}
\\ && \makebox[4em]{} \times
\prod_{j=1}^{n+1}\frac{\Gamma(t_1z_j,t_2z_j,t_3z_j,t_4z_j^{-1},t_5z_j^{-1})}
{\Gamma(z_j(ts)^{n-1}\prod_{j=1}^5 t_j)}.
\label{a2-int}\end{eqnarray}

\begin{theorem}
Из $A_n$ и $C_n$ многократных эллиптических бета-интегралов
(\ref{ell-int-A}) и (\ref{ell-int-C}) вытекает две новые точные формулы
интегрирования. Для нечетных $n=2m-1$, мы имеем
\begin{eqnarray}\nonumber
&& \int_{\mathbb{T}^n}\Delta^{II}({\bf z};A_n)
\frac{dz_1}{z_1}\ldots \frac{dz_n}{z_n}
= \frac{ (n+1)! }{ (q;q)_\infty^n (p;p)_\infty^n }
\\ && \makebox[4em]{} \times
\frac{ \Gamma(t^m,s^m,s^{m-1}t_4t_5)\prod_{1\leq i<j\leq 3}
\Gamma(t^{m-1}t_it_j) }
{ \prod_{k=4}^5\Gamma(t^{2m-2}s^{m-1}t_1t_2t_3t_k) }
\nonumber \\ && \makebox[4em]{} \times
\prod_{j=1}^m \frac{ \prod_{i=1}^3\prod_{k=4}^5
\Gamma((ts)^{j-1}t_it_k) }
{ \prod_{1\leq i<\ell\leq 3} \Gamma((ts)^{m+j-2}t_it_\ell t_4t_5) }
\nonumber \\ && \makebox[4em]{} \times
\prod_{j=1}^{m-1} \frac{ \Gamma((ts)^j,t^js^{j-1}t_4t_5)
\prod_{1\leq i<\ell\leq 3}\Gamma(t^{j-1}s^jt_it_\ell) }
{ \prod_{k=4}^5\Gamma(t^{m+j-2}s^{m+j-1}t_1t_2t_3 t_k) }.
\label{int-an-odd}\end{eqnarray}

Для четных $n=2m$, мы имеем
\begin{eqnarray}\nonumber
&&\int_{\mathbb{T}^n}\Delta^{II}({\bf z};A_n)
\frac{dz_1}{z_1}\ldots \frac{dz_n}{z_n}
= \frac{ (n+1)! }{ (q;q)_\infty^n (p;p)_\infty^n }
\\ && \makebox[4em]{} \times
\frac{ \prod_{i=1}^3\Gamma(t^mt_i)\prod_{k=4}^5\Gamma(s^mt_k)\,
\Gamma(t^{m-1}t_1t_2t_3) }
{ \Gamma(t^{2m-1}s^{m-1}\prod_{i=1}^5t_i,t^{2m-1}s^m t_1t_2t_3) }
\nonumber \\ && \makebox[4em]{} \times
\prod_{j=1}^m \frac{ \Gamma((ts)^j,t^js^{j-1}t_4t_5)
\prod_{i=1}^3\prod_{k=4}^5 \Gamma((ts)^{j-1}t_it_k) }
{ \prod_{k=4}^5\Gamma(t^{m+j-2}s^{m+j-1}t_k^{-1}\prod_{i=1}^5t_i)}
\nonumber \\ && \makebox[4em]{} \times
\prod_{j=1}^m \prod_{1\leq i<\ell\leq3 }
\frac{\Gamma(t^{j-1}s^jt_it_\ell)}
{\Gamma((ts)^{m+j-1}t_it_\ell t_4t_5)}.
\label{int-an-even}\end{eqnarray}
\end{theorem}
{\bf Доказательство.}
Структура доказательств повторяет процедуру, использованную
Густафсоном в \cite{gus:some2} для доказательства $p=0$ случая интегралов
(\ref{int-an-odd}) и (\ref{int-an-even}). Начнем со случая нечетного $n=2m-1$
и рассмотрим следующий $(4m-1)$-кратный интеграл:
\begin{eqnarray}\nonumber
\lefteqn{ \int_{\mathbb{T}^{4m-1}} \frac{\prod_{i=1}^{2m}\prod_{j=1}^m
\Gamma(t^{1/2}z_iw_j,t^{1/2}z_iw_j^{-1},
s^{1/2}z_i^{-1}x_j,s^{1/2}z_i^{-1}x_j^{-1})}
{\prod_{i,j=1;\: i\neq j}^{2m}\Gamma(z_iz_j^{-1}) \prod_{\nu=\pm1}
\prod_{1\leq i<j\leq m}\Gamma(w_i^\nu w_j^\nu,w_i^\nu w_j^{-\nu},
x_i^\nu x_j^\nu,x_i^\nu x_j^{-\nu})} }&&
\\ && \times
\prod_{i=1}^{2m}\frac{\Gamma(z_i(ts)^{m-2}\prod_{k=1}^5 t_k)}
{\Gamma(z_i(ts)^{2m-2}\prod_{k=1}^5 t_k)}
\prod_{k=1}^{2m-1}\frac{dz_k}{z_k}\prod_{j=1}^m\Biggl(\frac{dw_j}{w_j}
\frac{dx_j}{x_j} \label{n=2m-1} \\ && \times
\prod_{\nu=\pm1} \Biggl( \frac{\prod_{k=1}^3\Gamma(t^{-1/2}t_kw_j^\nu)
\Gamma(x_j^\nu t^{m-2}s^{-1/2}t_1t_2t_3)
\prod_{k=4}^5\Gamma(s^{-1/2}t_kx_j^\nu)}
{\Gamma(w_j^\nu t^{m-3/2}t_1t_2t_3,x_j^\nu t^{m-2}s^{m-3/2}\prod_{k=1}^5t_k,
w_j^{2\nu},x_j^{2\nu})} \Biggr)\Biggr),
\nonumber\end{eqnarray}
где $\prod_{i=1}^{2m}z_j=1$.
Используя $C_n$ эллиптический бета-интеграл типа $I$ (\ref{ell-int-C}),
мы сначала вычисляем в (\ref{n=2m-1}) интегралы по переменным
$w_j, j=1,\ldots,m,$ а затем по $x_j, j=1,\ldots,m$. Результирующее
выражение равно левой стороне (\ref{int-an-odd}) с точностью до множителя
\begin{eqnarray*}
&& (2\pi i)^{4m-1}\frac{2^{2m}(m!)^2}{(p;p)_\infty^{2m}(q;q)_\infty^{2m}}
\frac{\Gamma(s^{-1}t_4t_5)}{\Gamma(s^{m-1}t_4t_5)}
\\ && \times
\prod_{1\leq i<k\leq 3}\frac{\Gamma(t^{-1}t_it_k)}{\Gamma(t^{m-1}t_it_k)}
\prod_{k=4,5}\frac{\Gamma(t^{m-2}s^{-1}t_k^{-1}\prod_{i=1}^5t_k)}
{\Gamma(t^{m-2}s^{m-1}t_k^{-1}\prod_{i=1}^5t_k)}.
\end{eqnarray*}

В этой двухшаговой схеме необходимо наложить ограничения на параметры:
$$
|t|<1, \quad |t_{1,2,3}|<|t|^{1/2}, \quad |pq|<|t^{m-3/2}t_1t_2t_3|
$$
и
$$
|s|<1,\quad |t_{4,5}|<|s|^{1/2},\quad
|pq|<|t^{m-2}s^{m-3/2}\prod_{k=1}^5t_k|,
$$
соответственно. Однако, результирующее выражение может быть аналитически
продолжено в область $|t_k|<1, \, k=1,\ldots,5$, $|pq|<|(ts)^{2m-2}
\prod_{k=1}^5t_k|$ без изменения значений интеграла.

Поскольку подынтегральное выражение в (\ref{n=2m-1}) ограничено
на единичной окружности, можно поменять порядок интегрирований и
взять сначала интегралы по $z_i,\: i=1,\ldots,2m-1,$ используя
$A_n$-формулу (\ref{ell-int-A}). После этого мы применяем
формулу (\ref{ell-int-C}) для того, чтобы вычислить интегралы по
переменным $x_j, j=1,\ldots,m$. Наконец, мы используем исконный эллиптический
аналог интеграла Сельберга (\ref{SintB}) для взятия интегралов
по $w_j,\: j=1,\ldots,m.$ Все это приводит к выражению
\begin{eqnarray*}
&& \frac{(2\pi i)^{4m-1}(2m)! (2^mm!)^2}{((p;p)_\infty(q;q)_\infty)^{4m-1}}
\frac{\Gamma(t^m,s^m,t_4t_5s^{-1})} {\Gamma((ts)^m,t^ms^{m-1}t_4t_5)}
\prod_{k=4}^5\frac{\Gamma(t^{m-2}s^{-1}t_k^{-1}\prod_{i=1}^5t_i)}
{\Gamma(t^{2m-2}s^{m-1}t_k^{-1}\prod_{i=1}^5t_i)}
\\ && \times
\prod_{j=1}^m\frac{\Gamma((ts)^j,t^js^{j-1}t_4t_5)
\prod_{1\leq i<k\leq 3}\Gamma(t^{j-2}s^{j-1}t_it_k)
\prod_{k=4}^5 \prod_{i=1}^3\Gamma((ts)^{j-1}t_kt_i)}
{\prod_{k=1}^3\Gamma((ts)^{m+j-2}t_k^{-1}\prod_{i=1}^5t_i)
\prod_{k=4}^5 \Gamma(t^{m+j-3}s^{m+j-2}t_k^{-1}\prod_{i=1}^5t_i)}.
\end{eqnarray*}
Сравнивая два полученных выражения, получаем необходимую формулу интегрирования
при нечетных $n=2m-1$.

Для доказательства формулы (\ref{int-an-even}), мы рассмотрим $4m$-кратный интеграл
\begin{eqnarray}\nonumber
&& \int_{\mathbb{T}^{4m}} \frac{\prod_{i=1}^{2m+1}\prod_{j=1}^m
\Gamma(t^{1/2}z_iw_j,t^{1/2}z_iw_j^{-1},
s^{1/2}z_i^{-1}x_j,s^{1/2}z_i^{-1}x_j^{-1})}
{\prod_{i,j=1;\: i\neq j}^{2m+1}\Gamma(z_iz_j^{-1}) \prod_{\nu=\pm1}
\prod_{1\leq i<j\leq m}\Gamma(w_i^\nu w_j^\nu,w_i^\nu w_j^{-\nu},
x_i^\nu x_j^\nu,x_i^\nu x_j^{-\nu})}
\\ && \makebox[2em]{}\times
\prod_{i=1}^{2m+1}\frac{\Gamma(t_3z_i,s^{m-1}t_4t_5z_i,
t^{m-1}t_1t_2z_i^{-1})}{\Gamma(z_i(ts)^{2m-1}\prod_{k=1}^5 t_k)}
\prod_{j=1}^m\Biggl( \frac{dw_j}{w_j}\frac{dx_j}{x_j}
\label{n=2m} \\ && \makebox[2em]{}\times
\prod_{\nu=\pm1} \Biggl(\frac{\prod_{k=1}^2\Gamma(t^{-1/2}t_kw_j^\nu,
s^{-1/2}t_{k+3}x_j^\nu)}{\Gamma(t^{m-1/2}t_1t_2w_j^\nu,
s^{m-1/2}t_4t_5x_j^\nu, w_j^{2\nu},x_j^{2\nu})} \Biggr)\Biggr)
\prod_{k=1}^{2m}\frac{dz_k}{z_k},
\nonumber\end{eqnarray}
где $\prod_{i=1}^{2m+1}z_j=1$. Повторяя тот же самый трюк как и в случае
нечетного $n$ (то есть интегрируя последовательно по переменным $w_j$ и
$x_j$, а затем меняя порядок интегрирований в этом выражении),
мы получаем равенство (\ref{int-an-even}).
\hfill{Q.E.D.}

\subsection{Обобщение интегралов Густафсона-Ракха (тип II)}

Аналогичным образом можно построить эллиптические аналоги $A_n$
$q$-бета интегралов Густафсона и Ракха \cite{gus-rak:beta}.
Для этого нам понадобятся $A_n$ и $C_n$ эллиптические бета-интегралы первого типа
(\ref{ell-int-A}) и (\ref{ell-int-C}).

\begin{theorem}
Обозначим
\begin{eqnarray}\nonumber
&& \Delta^{II}({\bf z};A_n)=
\frac{1}{(2\pi i)^n}\prod_{1\leq i<j\leq n+1}
\frac{\Gamma(tz_iz_j)}{\Gamma(z_iz_j^{-1},z_i^{-1}z_j) }
\\ && \makebox[4em]{} \times
\frac{\prod_{i=1}^{n+1}\left(\prod_{k=1}^{n+1}\Gamma(t_kz_i^{-1})
\prod_{k=n+2}^{n+4}\Gamma(tt_kz_i)\right)}
{\prod_{j=1}^{n+1}\Gamma(Az_j^{-1}) },
\label{a3}\end{eqnarray}
где $A=t^{n+2}\prod_{i=1}^{n+4}t_i$ и $\prod_{j=1}^{n+1}z_j=1$.
Справедливы следующие формулы интегрирования: для нечетных $n=2l-1$,
\begin{eqnarray}
&& \int_{\mathbb{T}^n}\Delta^{III}({\bf z};A_n)
\frac{dz_1}{z_1}\ldots\frac{dz_n}{z_n}
= \frac{ (n+1)! }{ (q;q)_\infty^n (p;p)_\infty^n }
\frac{\Gamma(t^l,\prod_{k=1}^{2l}t_k)}
{\Gamma(t^l\prod_{k=1}^{2l}t_k)}
\label{int-an-o} \\ && \makebox[1em]{} \times
\frac{ \prod_{i=1}^{2l}\prod_{j=2l+1}^{2l+3}\Gamma(tt_it_j)
\prod_{1\leq i<j\leq 2l} \Gamma(tt_it_j)
\prod_{2l+1\leq i<j\leq 2l+3}\Gamma(t^{l+1}t_it_j)}
{\prod_{i=1}^{2l}\Gamma(t^{2l+1}t_i^{-1}\prod_{k=1}^{2l+3}t_k)
\prod_{i=2l+1}^{2l+3}\Gamma(t^{l+1}t_i^{-1}\prod_{k=1}^{2l+3}t_k) }.
\nonumber\end{eqnarray}
и для четных $n=2l$,
\begin{eqnarray}\nonumber
&&\int_{\mathbb{T}^n}\Delta^{III}({\bf z};A_n)
\frac{dz_1}{z_1}\ldots\frac{dz_n}{z_n}
= \frac{ (n+1)! }{ (q;q)_\infty^n (p;p)_\infty^n }
\frac{\Gamma(\prod_{k=1}^{2l+1}t_k,t^{l+2}\prod_{k=2l+2}^{2l+4}t_k)}
{\Gamma(t^{l+2}\prod_{k=1}^{2l+4}t_k)}
\\ && \makebox[1em]{} \times
\frac{ \prod_{i=1}^{2l+1}\prod_{j=2l+2}^{2l+4}\Gamma(tt_it_j)
\prod_{1\leq i<j\leq 2l+1} \Gamma(tt_it_j)
\prod_{i=2l+2}^{2l+4}\Gamma(t^{l+1}t_i)}
{\prod_{i=1}^{2l+1}\Gamma(t^{2l+2}t_i^{-1}\prod_{k=1}^{2l+4}t_k)
\prod_{i=2l+2}^{2l+4}\Gamma(t^{l+1}t_i\prod_{k=1}^{2l+1}t_k) },
\label{int-an-e} \end{eqnarray}
где подразумеваются ограничения на параметры такие, что все
последовательности полюсов сходящихся к нулю находятся внутри $\T$.
\end{theorem}
{\bf Доказательство.}
В соответствии с процедурой, использованной в работе \cite{gus-rak:beta},
рассмотрим $(3l-1)$-кратный интеграл
\begin{eqnarray*}
&& \int_{\mathbb{T}^{3l-1}}
\frac{\prod_{i=1}^{2l}\prod_{j=1}^l\Gamma(t^{1/2}z_iw_j,t^{1/2}z_iw_j^{-1})
\prod_{i=0}^{2l}\prod_{j=1}^{2l}\Gamma(t_iz_j^{-1})}
{\prod_{i,j=1; i\neq j}^{2l}\Gamma(z_iz_j^{-1})
\prod_{j=1}^{2l} \Gamma(t^l\prod_{i=0}^{2l}t_iz_j^{-1})}
\\ && \makebox[1em]{} \times
\prod_{\nu=\pm1}\prod_{1\leq i<j\leq l}\Gamma^{-1}
(w_i^\nu w_j^{\nu},w_i^\nu w_j^{-\nu})
\prod_{j=1}^l
\frac{\prod_{k=2l+1}^{2l+3}\Gamma(t^{1/2}t_kw_j^\nu)}
{\Gamma(w_j^{2\nu},t^{l+3/2}\prod_{k=2l+1}^{2l+3}t_kw_j^\nu)}
\\ && \makebox[1em]{} \times
\frac{dw_1}{w_1}\ldots\frac{dw_l}{w_l}\frac{dz_1}{z_1}
\ldots\frac{dz_{2l-1}}{z_{2l-1}},
\end{eqnarray*}
где $\prod_{i=1}^{2l}z_i=1$ и $t_0=t^{l+1}\prod_{k=2l+1}^{2l+3}t_k$.
Интегрирование по переменным $w_j$ с помощью формулы (\ref{ell-int-C}),
мы получаем левую сторону (\ref{int-an-o}) с точностью до некоторого
множителя. Меняя порядок интегрирования, мы сначала интегрируем по $z_i$
с помощью (\ref{ell-int-A}) (где необходимо заменить $z_k$ на $z_k^{-1}$),
а затем по $w_j$ с помощью (\ref{ell-int-C}). Приравнивая два выражения,
мы получаем формулу (\ref{int-an-o}).

Аналогичным образом, в случае четного $n=2l$ мы рассмотрим $(3l+1)$-кратный
интеграл
\begin{eqnarray*}
&& \int_{\mathbb{T}^{3l+1}}
\frac{\prod_{i=1}^{2l+1}\left(\prod_{j=1}^{l+1}
\Gamma(t^{1/2}z_iw_j,t^{1/2}z_iw_j^{-1})
\prod_{j=1}^{2l+1}\Gamma(t_iz_j^{-1})\right)}
{\prod_{i,j=1; i\neq j}^{2l+1}\Gamma(z_iz_j^{-1})
\prod_{j=1}^{2l+1} \Gamma(t^{l+1}\prod_{i=1}^{2l+1}t_iz_j)}
\\ && \makebox[1em]{} \times
\prod_{\nu=\pm1}\prod_{1\leq i<j\leq l+1}\Gamma^{-1}
(w_i^\nu w_j^{\nu},w_i^\nu w_j^{-\nu})
\prod_{j=1}^{l+1}
\frac{\prod_{k=2l+1}^{2l+5}\Gamma(t^{1/2}t_kw_j^\nu)}
{\Gamma(w_j^{2\nu},t^{l+5/2}\prod_{k=2l+1}^{2l+5}t_kw_j^\nu)}
\\ && \makebox[1em]{} \times
\frac{dw_1}{w_1}\ldots\frac{dw_{l+1}}{w_{l+1}}\frac{dz_1}{z_1}
\ldots\frac{dz_{2l}}{z_{2l}},
\end{eqnarray*}
где $\prod_{i=1}^{2l+1}z_i=1$ и $t_{2l+5}=t^l\prod_{k=1}^{2l+1}t_k$.
Повторяя тот же самый трюк как и в предыдущем случае, получаем
необходимую формулу интегрирования (\ref{int-an-e}).
\hfill{Q.E.D.}

\section{Дополнительный $(2n+3)$-параметрический эллиптический
$A_n$ бета-интеграл (тип I)}

В работе \cite{spi-war:inversions}, Варнаар и автор построили
еще один многомерный эллиптический бета-интеграл, связанный
с $A_n$ системой корней. Поскольку он оказывается связанным с
гипергеометрическими рядами, классифицированными
как $D_n$ ряды \cite{B,ros:elliptic,Sch1}, формально его можно
отнести также к этой системе корней.
Обозначив переменные интегрирования как $z=(z_1,\ldots,z_{n})\in\C^{n}$,
$z_{n+1}=\prod_{i=1}^{n}z_i^{-1}$, и введя параметры
$t=(t_1,\ldots,t_{n+3})\in \C^{n+3}$ и $s=(s_1,\ldots,s_{n})\in\C^{n}$,
мы определим ядро интеграла
\begin{eqnarray}\nonumber
&& \rho(z,t,s)=\prod_{1\leq i<j\leq n+1}\frac{\Gamma(Dz_i^{-1}z_j^{-1})}
{\Gamma(z_iz_j^{-1},z_i^{-1}z_j)}
\prod_{i=1}^{n+1}\frac{\prod_{m=1}^{n}\Gamma(s_mz_i)
\prod_{j=1}^{n+3}\Gamma(t_jz_i^{-1})}
{\prod_{m=1}^{n}\Gamma(Ds_mz_i^{-1})}
\\ && \makebox[6em]{}
\times \prod_{j=1}^{n+3}\prod_{m=1}^n\frac{\Gamma(Ds_mt_j^{-1})}
{\Gamma(s_mt_j)}\prod_{1\leq j<k\leq n+3}
\frac{1}{\Gamma(Dt_j^{-1}t_k^{-1})},
\label{kernel-D}\end{eqnarray}
где $D=\prod_{j=1}^{n+3}t_j$.
Эта мероморфная функция имеет полюсы в точках
$$
z_i=\{t_l q^jp^k\}, \, i=1,\ldots,n, \quad
z_{n+1}^{-1}=z_1\cdots z_n= \{s_mq^{j}p^{k},\, (Ds_m)^{-1}q^{j+1}p^{k+1}\},
$$
с $l=1,\ldots, n+3,\,m=1,\ldots, n,\, j,k \in\N,$ и их $z\to 1/z$ партнеры
$$
z_i=\{s_m^{-1} q^{-j}p^{-k}, Ds_mq^{-j-1}p^{-k-1}\},\quad i=1,\ldots, n,\qquad
z_{n+1}^{-1}= \{t_l^{-1}q^{-j}p^{-k}\}.
$$
Помимо них существуют специальные последовательности полюсов, появляющиеся из
ус\-ло\-вий  $z_iz_j=\{Dq^{a}p^{b}\}_{a,b\in\N}$ при $i,j=1,\ldots,n+1$.

Нам потребуется два конечно-разностных уравнения на ядро интеграла.
Растяжение параметра $s_1$ на $q$ дает
\begin{equation}
\frac{\rho(qs_1)}{\rho(s_1)}=
\prod_{i=1}^{n+1}\frac{\theta(s_1z_i;p)}{\theta(Ds_1z_i^{-1};p)}
\prod_{j=1}^{n+3}\frac{\theta(Ds_1t_j^{-1};p)}{\theta(s_1t_j;p)}.
\label{eq-D-s1}\end{equation}
Аналогично, для переменных $z_i$ получаем
\begin{eqnarray}\nonumber
\lefteqn{\frac{\rho(q^{-1}z_i)}{\rho(z_i)}=
\prod_{m=1}^n\frac{\theta(s_mz_{n+1},q^{-1}Ds_mz_{n+1}^{-1};p)}
{\theta(q^{-1}s_mz_i,Ds_mz_i^{-1};p)}
\prod_{j=1}^{n+3}\frac{\theta(t_jz_i^{-1};p)}{\theta(q^{-1}t_jz_{n+1}^{-1};p)}
}&& \\ && \times
\prod_{j=1,\neq i}^n \frac{\theta(q^{-1}z_iz_j^{-1},q^{-1}z_jz_{n+1}^{-1},
Dz_i^{-1}z_j^{-1};p)}{\theta(z_i^{-1}z_j,z_j^{-1}z_{n+1},
q^{-1}Dz_j^{-1}z_{n+1}^{-1};p)}
\frac{\theta(q^{-2}z_iz_{n+1}^{-1};p)z_i^2}
{\theta(z_iz_{n+1}^{-1};p)qz_{n+1}^2}.
\label{eq-D-zi}\end{eqnarray}

Теперь сформулируем основное утверждение.

\begin{theorem}
Предположим, что $|t_l|, |s_m|<1$ и $|pq|<|Ds_m|$. Тогда,
\begin{equation}\label{ell-int-D}
\int_{\T^n}\rho(z,t,s)\frac{dz}{z}=\frac{(n+1)!(2\pi i)^n}
{(q;q)_\infty^n(p;p)_\infty^n}.
\end{equation}
\end{theorem}
{\bf Доказательство.} %\begin{proof}
В качестве первого шага, установим смешанное $q$-разностное уравнение
для ядра:
\begin{eqnarray}\nonumber
&& \rho(z,t,qs_1,s_2,\ldots,s_n)-\rho(z,t,s)
\\ && \makebox[2em]{}
=\sum_{i=1}^n\left(g_i(z_1,...,q^{-1}z_i,\ldots,z_n,t,s)-g_i(z,t,s)\right),
\label{eqn-D-s}\end{eqnarray}
где мероморфные функции $g_i$ имеют вид
\begin{eqnarray}\nonumber &&
\frac{g_i(z,t,s)}{\rho(z,t,s)}=
s_1z_{n+1}\prod_{m=2}^n
\frac{\theta(Ds_m/qz_i,s_mz_i;p)}{\theta(Ds_m/z_{n+1},s_mz_{n+1}/q;p)}
\frac{\theta(Ds_1^2;p)\prod_{j=1}^n\theta(s_1z_j;p)}
{\theta(z_i/z_{n+1},Ds_1/z_{n+1};p)}
\\&& \makebox[2em]{} \times
\prod_{j=1}^{n+3}\frac{\theta(t_j/z_{n+1};p)}{\theta(t_js_1;p)}
\prod_{j=1,\neq i}^n\frac{\theta(z_{n+1}/qz_j,D/z_jz_{n+1};p)}
{\theta(z_i/z_j,z_j/z_{n+1},D/qz_iz_j;p)}.
\label{g/rho}\end{eqnarray}
Подставляя это выражение в (\ref{eqn-D-s}), получаем
\begin{eqnarray}\nonumber
\lefteqn{
\prod_{i=1}^{n+1}\frac{\theta(s_1z_i;p)}{\theta(Ds_1/z_i;p)}
\prod_{j=1}^{n+3}\frac{\theta(Ds_1/t_j;p)}
{\theta(s_1t_j;p)} -1 } &&
\\ \nonumber &&
=\frac{s_1z_{n+1}\theta(Ds_1^2;p)}{\prod_{j=1}^{n+3}\theta(t_js_1;p)}
\sum_{i=1}^n\frac{\prod_{j=1,\neq i}^n\theta(s_1z_j;p)}
{\theta(z_i/z_{n+1};p)}
\\ \nonumber && \makebox[1em]{} \times
\Biggl(\frac{z_i^2}{z_{n+1}^2}
\frac{\theta(s_1z_{n+1};p)}{\theta(Ds_1/z_i;p)}
\frac{\prod_{j=1}^{n+3}\theta(t_j/z_i;p)}
{\prod_{j=1,\neq i}^n \theta(z_j/z_i;p)}
\\ \nonumber && \makebox[1em]{}
-\frac{\theta(s_1z_i;p)}{\theta(Ds_1/z_{n+1};p)}
\frac{\prod_{j=1}^{n+3}\theta(t_j/z_{n+1};p)}
{\prod_{j=1,\neq i}^n\theta(z_i/z_j;p)}
\\  && \times
\prod_{m=2}^n\frac{\theta(Ds_m/qz_i,s_mz_i;p)}
{\theta(Ds_m/z_{n+1},s_mz_{n+1}/q;p)}
\prod_{j=1,\neq i}^n\frac{\theta(z_{n+1}/qz_j,D/z_jz_{n+1};p)}
{\theta(z_j/z_{n+1},D/qz_iz_j;p)} \Biggr).
\label{eqn-exp-D-s}\end{eqnarray}
Обе стороны этого равенства являются эллиптическими функциями
$\log z_1$ с полюсами при $z_1=Ds_1$ и $z_{n+1}=Ds_1$ (в обеих
сторонах) и $z_{n+1}=\{Ds_k,q/s_k,z_j\}_{j,k=1,\ldots,n}$,
$z_1=\{z_k,q/Dz_k\}_{k=2,\dots,n}$ (в правой части).

Совпадение вычетов полюсов $z_1=Ds_1, \, j=1,\dots,n,$
легко проверяется и мы пропускаем их. Сравнивая вычеты в полюсах
$z_{n+1}=Ds_1$, получаем равенство
\begin{eqnarray}\nonumber
&& \sum_{i=1}^n\prod_{m=2}^n\theta(Ds_m/qz_i,s_mz_i;p)
\prod_{j=1,\neq i}^n\frac{\theta(Ds_1/qz_j,1/s_1z_j;p)}
{\theta(D/qz_iz_j,z_i/z_j;p)}
\\ && \makebox[4em]{}
=\prod_{m=2}^n\theta(s_m/s_1,Ds_1s_m/q;p).
\label{res-zn+1}\end{eqnarray}
Умножая обе стороны на
$ \prod_{1\leq i<j\leq n}z_j\theta(D/qz_iz_j,z_i/z_j;p), $
получаем две голоморфные тета-функции $z_1,\ldots,z_n$
антисимметричные по этим переменным и равные нулю при $z_kz_l=q/D$,
для $k,l=1,\ldots,n$. Любая такая функция равна
указанному множителю с точностью до фактора, не зависящего
от  $z_j$. При $z_1=1/s_1$, равенство (\ref{res-zn+1}) очевидно
и, поэтому, оно верно  в общем случае.

Складывая вычеты при $z_{n+1}=Ds_k$ или $z_{n+1}=q/s_k$ для
некоторого фиксированного $k=2,\ldots,n,$ мы получаем выражение
$$
\sum_{i=1}^n\prod_{m=2}^n\theta(Ds_m/qz_i,s_mz_i;p)
\prod_{j=1,\neq i}^n\frac{\theta(Ds_k/qz_j,1/s_kz_j;p)}
{\theta(D/qz_iz_j,z_i/z_j;p)},
$$
равное нулю, поскольку оно не зависит от $z_j$ (аналогично
предыдущему случаю) и зануляется при $z_i=1/s_k$.

Сокращение полюсов при $z_1=z_j,j=2,\dots,n,$ легко проверяется.
При $z_{n+1}=z_k,\, k=1,\ldots,n,$ сумма вычетов приводит к равенству
\begin{eqnarray}\nonumber
&& \sum_{i=1}^n\prod_{m=2}^n\theta(Ds_m/qz_i,s_mz_i;p)
\prod_{j=1,\neq i}^n\frac{\theta(z_k/qz_j,D/z_kz_j;p)}
{\theta(D/qz_iz_j,z_i/z_j;p)}
\\ && \makebox[4em]{}
=\prod_{m=2}^n\theta(s_mz_k/q,Ds_m/z_k;p),
\label{zn+1=zk}\end{eqnarray}
являющемуся $s_1=q/z_k$ подслучаем (\ref{res-zn+1}).
Наконец, имеется только два полюса в правой части равенства
(\ref{eqn-exp-D-s}) при $z_1=D/qz_k$ для некоторого
фиксированного $k=2,\dots,n,$ и они сокращают друг друга.

Из перечисленных свойств следует, что разность выражений правой и
левой частей равенства (\ref{eqn-exp-D-s}) есть эллиптическая
функция $\log z_1$ без полюсов, то есть константа. Последняя
равна нулю поскольку при $z_1=1/s_1$ справедливость
(\ref{eqn-exp-D-s}) очевидна.

Интегрируя (\ref{eqn-D-s}) по переменным $z\in \T^n$, получаем
\begin{equation}
I(t,qs_1,s_2,\dots,s_n)-I(t,s)=\sum_{i=1}^n
\left(\int_{\T^{i-1}\times(q^{-1}\T)\times \T^{n-i}}-\int_{\T^n}\right)
g_i(z,t,s)\frac{dz}{z},
\label{int-eqn-D-s}\end{equation}
где $I(t,s)=\int_{\T^n}\rho(z,t,s)dz/z$.
Функции $g_i(z,t,s)$ имеют полюсы по $z_i$ в точках
\begin{eqnarray*}
&& z_i=\{t_kq^ap^b, \quad  Dz_l^{-1}q^{a-1}p^b, \quad
s_1^{-1}q^{-a-1}p^{-b}, \quad s_m^{-1}q^{-a-1}p^{-b},
\\ && \makebox[2em]{}
Ds_1q^{-a-1}p^{-b-1}, \quad Ds_mq^{-a-2}p^{-b-1}\},
\\
&& z_{n+1}^{-1}=\{s_1q^ap^b, \quad s_mq^{a+1}p^b, \quad
D^{-1}s_1^{-1}q^ap^{b+1},
\\ && \makebox[2em]{}
D^{-1}s_m^{-1}q^ap^{b+1},\quad
t_k^{-1}q^{-a-1}p^{-b},\quad z_lD^{-1}q^{-a-1}p^{-b}\},
\end{eqnarray*}
где $m=2,\dots,n,\, k=1,\dots,n+3,\, l=1,\ldots, n, l\neq i,\, a,b\in\N$.
При $|s_j|, |t_k|<1$ для всех $j,k$ и $\max\{|p|,|q|\}<|Ds_m|$,
$m=1,\ldots, n,$ функции $g_i(z,t,s)$ не имеют полюсов в области
$1\leq |z_i|\leq 1/|q|$ и правая часть (\ref{int-eqn-D-s})
равна нулю.

Разлагая рассматриваемый интеграл по $p$ и совершая почленно
итеративное растяжение $s_1\to qs_1$ с последующим
аналитическим продолжением, получаем, что $I$ не зависит от
$s_1$. По симметрии, $I$ не зависит от всех $s_m$, т.е.
$I(t,s)=f(t,q,p)$ зависит только от параметров $t_k$ и
базовых переменных $p,q$.

Для того, чтобы найти $f(t,q,p)$ необходимо применить анализ
вычетов. А именно, надо растянуть $t_1,\ldots, t_n$ из области
$|t_k|<1$ в кольцо $1<|t_k|<1/|q|$. В результате этого конечное число
полюсов пересекает $\T$ и мы получаем интеграл по $\T^n$
и сумму вычетов пересеченных полюсов. Взяв предел
$s_k\to 1/t_k, k=1,\ldots,n,$ находим, что
$f(t,q,p)$ имеет вид указанный в правой части (\ref{ell-int-D}).
\hfill{Q.E.D.}

\section{Вывод формул суммирования для эллиптических гипергеометрических
 рядов с помощью вычетов}

Рассмотрим эллиптический аналог бета-интеграла Сельберга \re{SintB}.
Соответствующая подынтегральная
функция $\Delta_n(\mathbf{z};p,q)$ имеет полюсы по $z_j$ внутри единичного
круга $\T$ в точках $\{ t_rp^lq^m \}_{l,m\in \mathbb{N}}$ ($r=0,\ldots ,4$),
$\{ z_k^{\pm 1}t p^lq^m \}_{l,m\in \mathbb{N}}$, и
$\{ p^{l+1}q^{m+1}t^{2-2n}\prod_{s=0}^4t_s^{-1}\}_{l,m\in \mathbb{N}}$.
Более того, благодаря $z_j\rightarrow z_j^{-1}$ симметрии, полюсы снаружи
$\T$ связаны с внутренними простой инверсией.

Переместим теперь параметр $t_0$ из области $|t_0|<1$ в $|t_0|>1$. При этом конечное
число полюсов пересечет единичную окружность. Конкретнее, возьмем
$q^{-N}<|t_0|<q^{-N-1}$ для некоторого целого $N\in\mathbb{N}$ и $0<q<1$.
Тогда, при достаточно маленьком действительном $p$ ($p<1/|t_0|$)
полюсы $z_j=t_0p^lq^m $ перемещаются за $\T$ при $l=0$ и $m=0,\ldots ,N$.
Соответственно, полюсы $z_j=1/t_0p^lq^m$ перемещаются во внутрь $\T$.
Опишем формулу интегрирования с учетом вычетов, появляющихся при такой процедуре.
Аналогичные формулы для вычетов многомерных интегралов приведены в
работах \cite{sto:basic,hec-opd:yang}.

\begin{theorem}\label{residue:thm}
Пусть $0<q,t<1$ и
$$
\Delta_n ({\bf z};p,q)=\frac{1}{(2\pi i)^n}\prod_{1\leq j<k\leq n}
\frac{\Gamma(tz_j^\pm z_k^\pm)}{\Gamma(z_j^\pm z_k^\pm)}
\prod_{j=1}^n\frac{\prod_{r=0}^4\Gamma(t_rz_j^\pm)}
{\Gamma(z_j^{\pm 2},t^{2n-2}\prod_{s=0}^4t_s z_j^\pm)}
$$
с параметрами $t_0,\ldots ,t_4$ общего положения такими, что
аргументы $\{ \arg (t_r), \arg (t_r^{-1}) \mid r=0,\ldots ,4 \}$
независимы над $\Z$
и $t_r^{-1}\prod_{s=0}^4t_s\not\in [1,+\infty [$ для $r=0,\ldots ,4$.
Предположим далее, что $|t_0|>1$ и $|t_r|<1$ при $r=1,\ldots ,4$,
а также $0<p<\min (|t_0|^{-1},q^{-1}t^{2n-2}\prod_{s=0}^4|t_s|)$.
Тогда имеем
\begin{eqnarray}\label{ell-res}
\lefteqn{\int_{C^n} \Delta_n(\mathbf{z};p,q)
\frac{dz_1}{z_1}\cdots\frac{dz_n}{z_n}=} && \\
&&
\sum_{m=0}^n 2^mm! \binom{n}{m}
\sum_{\stackrel{0\leq\lambda_1\leq\cdots \leq \lambda_m}
               {|\tau_mq^{\lambda_m}| >1}}
\int_{\T^{n-m}} \mu_m(\lambda,\mathbf{z};p;q)
\frac{dz_1}{z_1}\cdots\frac{dz_{n-m}}{z_{n-m}}, \nonumber
\end{eqnarray}
где $\tau_j=t_0t^{j-1}$, $j=1,\ldots ,n$,
\begin{equation*}
\mu_m(\lambda,\mathbf{z};p;q) =
\kappa_m \nu_m(\lambda;p;q)
\delta_{m,n-m}(\lambda ,\mathbf{z})\Delta_{n-m}(\mathbf{z};p,q),
\end{equation*}
\begin{eqnarray*}
\kappa_m &=&
\prod_{1\leq j<k\leq m}
\frac{\Gamma(t\tau_k\tau_j^{-1},t\tau_k^{-1}\tau_j^{-1};p,q)}
{\Gamma(\tau_k\tau_j^{-1},\tau_k^{-1}\tau_j^{-1};p,q)}
\\
&& \times \prod_{j=1}^m
\frac{\prod_{r=1}^4\Gamma(t_r\tau_j,t_r\tau_j^{-1};p,q)}
{(q;q)_\infty(p;p)_\infty
\Gamma(\tau_j^{-2},\tau_j^{-1}t^{2n-2}\prod_{s=0}^4t_s,
       \tau_j t^{2n-2}\prod_{s=0}^4t_s ;p,q)},
\end{eqnarray*}
\begin{eqnarray*}
&& \nu_m(\lambda;p;q)=
q^{\sum_{j=1}^m\lambda_j} t^{2 \sum_{j=1}^m (n-j)\lambda_j}
\prod_{1\leq j<k\leq m}
\Biggl(
\frac{\theta(\tau_k\tau_jq^{\lambda_k+\lambda_j},
\tau_k\tau_j^{-1}q^{\lambda_k-\lambda_j};p)}
{\theta(\tau_k\tau_j,\tau_k\tau_j^{-1};p)} \\
&& \makebox[8em]{}\times
\frac{\theta(t\tau_k\tau_j;p;q)_{\lambda_k+\lambda_j}}
     {\theta(qt^{-1}\tau_k\tau_j;p;q)_{\lambda_k+\lambda_j}}
\frac{\theta(t\tau_k\tau_j^{-1};p;q)_{\lambda_k-\lambda_j}}
     {\theta(qt^{-1}\tau_k\tau_j^{-1};p;q)_{\lambda_k-\lambda_j}}
\Biggr) \\
&&\makebox[3em]{} \times \prod_{j=1}^m
\Biggl( \frac{\theta(\tau_j^2q^{2\lambda_j};p)}{\theta(\tau_j^2;p)}
\prod_{r=0}^5 \frac{\theta(t_r\tau_j;p;q)_{\lambda_j}}
     {\theta(qt_r^{-1}\tau_j;p;q)_{\lambda_j}}\Biggr)  ,
\end{eqnarray*}
и
\begin{equation*}
\delta_{m,n-m}(\lambda, \mathbf{z}) =
\prod_{\begin{subarray}{c}1\leq j\leq m\\ 1\leq k\leq n-m\end{subarray}}
\frac{\Gamma(t\tau_jq^{\lambda_j}z_k,t\tau_jq^{\lambda_j}z_k^{-1},
t\tau_j^{-1}q^{-\lambda_j}z_k,t\tau_j^{-1}q^{-\lambda_j}z_k^{-1};p,q)}
{\Gamma(\tau_jq^{\lambda_j}z_k,\tau_jq^{\lambda_j}z_k^{-1},
\tau_j^{-1}q^{-\lambda_j}z_k,\tau_j^{-1}q^{-\lambda_j}z_k^{-1};p,q)}.
\end{equation*}
Здесь $\T$ обозначает единичную окружность ориентированную против часовой стрелки
и контур интегрирования $C\subset\mathbb{C}$ есть гладкая кривая Жордана с той
же ориентацией, причем (i) каждый луч, исходящий из нулевой точки пересекает
$C$ только один раз, (ii) $C^{-1}:=\{ z\in\mathbb{C} \mid z^{-1}\in C\} = C$,
и (iii) $C$ разделяет полюсы $z_j$ при
$\{ t_rp^lq^m \}_{l,m\in \mathbb{N}}$, $r=0,\ldots ,4$,
(находящиеся внутри $C$) от их $z_j\to1/z_j$ партнеров (снаружи $C$).
Более того, $t_5$ (в $\nu_m(\lambda;p;q)$) определяется через
$q$, $t$ и $t_0,\ldots ,t_4$ с помощью  условия балансировки
$q^{-1}t^{2n-2}\prod_{r=0}^5 t_r =1$.
(Отметим, что эти условия гарантируют, что $C$ также разделяет оставшиеся полюсы
при $z_j$ равных $\{ z_k^{\pm 1}t p^lq^m \}_{l,m\in \mathbb{N}}$,
$\{ z_k^{\pm 1} p^{l+1}q^{m+1} \}_{l,m\in \mathbb{N}}$, и $\{ p^{l+1}q^{m+1}
t^{2-2n}\prod_{s=0}^4t_s^{-1}\}_{l,m\in \mathbb{N}}$ от их $z_j\to 1/z_j$ партнеров.)
\end{theorem}
{\bf Доказательство.}
Доказательство основывается на формуле для вычетов многомерных интегралов типа
Аски-Вильсона, полученной Стокманом \cite{sto:basic}. Пусть
\begin{eqnarray}\label{AWweight}
\bar{\Delta}_n(\mathbf{z};q)&=& \frac{1}{(2\pi i)^n}
\prod_{1\leq j<k\leq n}
\frac{(z_jz_k,z_jz_k^{-1},z_j^{-1}z_k,z_j^{-1}z_k^{-1};q)_\infty}
{(tz_jz_k,tz_jz_k^{-1},tz_j^{-1}z_k,tz_j^{-1}z_k^{-1};q)_\infty} \\
&& \times
\prod_{1\leq j\leq n}\frac{(z_j^2,z_j^{-2};q)_\infty}
{(t_0z_j,t_0z_j^{-1};q)_\infty} ,\nonumber
\end{eqnarray}
с $0<q,t<1$ и $t_0\in \mathbb{C}\setminus [1,+\infty [$.
В работе  \cite{sto:basic} доказана формула:
\begin{eqnarray}\label{stokres:form}
\lefteqn{\int_{C^n} \bar{\Delta}_n(\mathbf{z};q) f(\mathbf{z})
\frac{dz_1}{z_1}\cdots\frac{dz_n}{z_n}=} && \\
&&
\sum_{m=0}^n 2^m m! \binom{n}{m}
\sum_{\stackrel{0\leq\lambda_1\leq\cdots \leq \lambda_m}
               {|\tau_mq^{\lambda_m}| >1}}
\int_{\T^{n-m}} \bar{\mu}_m(\lambda ,\mathbf{z};q) f_m(\lambda ,z)
\frac{dz_1}{z_1}\cdots\frac{dz_{n-m}}{z_{n-m}}, \nonumber
\end{eqnarray}
где $\tau_j=t_0t^{j-1}$,
\begin{eqnarray*}
\bar{\mu}_m (\lambda ,\mathbf{z};q) & = &
\bar{\kappa}_{m}\, \bar{\nu}_m(\lambda ;q )
\bar{\delta}_{m,n-m}(\lambda ,\mathbf{z})
\bar{\Delta}_{n-m}(\mathbf{z};q) ,\\
f_m(\lambda ,\mathbf{z}) &=&
f(\tau_1 q^{\lambda_1},\ldots ,\tau_m q^{\lambda_m},z_1,\ldots ,z_{n-m}), \\
\bar{\kappa}_{m} &=&
\prod_{1\leq j<k\leq m}
\frac{(\tau_k\tau_j^{-1},\tau_k^{-1}\tau_j^{-1};q)_\infty}
     {(t\tau_k\tau_j^{-1},t\tau_k^{-1}\tau_j^{-1};q)_\infty}
\prod_{1\leq j\leq m}
\frac{(\tau_j^{-2};q)_\infty }{(q;q)_\infty } ,
\end{eqnarray*}
\begin{eqnarray*}
\bar{\nu}_m(\lambda ;q) &=& q^{-\sum_{j=1}^m
   \lambda_j(3\lambda_j+1)/2}
(-t_0^{-4})^{\sum_{j=1}^m\lambda_j}
t^{-5 \sum_{j=1}^m (j-1)\lambda_{j}} \\
&& \makebox[3em]{} \times \prod_{1\leq j<k\leq m}
\Biggl(
\frac{1-\tau_k\tau_jq^{\lambda_k+\lambda_j}}{1-\tau_k\tau_j}
\frac{1-\tau_k\tau_j^{-1}q^{\lambda_k-\lambda_j}}{1-\tau_k\tau_j^{-1}}
\Biggr) \\
&& \makebox[10em]{}\times
\frac{(t\tau_k\tau_j;q)_{\lambda_k+\lambda_j}}
     {(qt^{-1}\tau_k\tau_j;q)_{\lambda_k+\lambda_j}}
\frac{(t\tau_k\tau_j^{-1};q)_{\lambda_k-\lambda_j}}
     {(qt^{-1}\tau_k\tau_j^{-1};q)_{\lambda_k-\lambda_j}} \\
&&\makebox[3em]{} \times \prod_{1\leq j\leq m}
\Biggl( \frac{1-\tau_j^2q^{2\lambda_j}}{1-\tau_j^2}\Biggr)
\frac{(t_0\tau_j;q)_{\lambda_j}}
     {(qt_0^{-1}\tau_j;q)_{\lambda_j}}  ,
\end{eqnarray*}
и
\begin{equation*}
\bar{\delta}_{m,n-m}(\lambda , \mathbf{z})  =\!\!\!
\prod_{\begin{subarray}{c} 1\leq j\leq m\\[0.3ex]
                           1\leq k\leq n-m        \end{subarray}}
\frac{(\tau_jq^{\lambda_j}z_k,\tau_jq^{\lambda_j}z_k^{-1},
\tau_j^{-1}q^{-\lambda_j}z_k,\tau_j^{-1}q^{-\lambda_j}z_k^{-1};q)_\infty}
{(t\tau_jq^{\lambda_j}z_k,t\tau_jq^{\lambda_j}z_k^{-1},
t\tau_j^{-1}q^{-\lambda_j}z_k,t\tau_j^{-1}q^{-\lambda_j}z_k^{-1};q)_\infty},
\end{equation*}
где контур $C$ удовлетворяет условиям указанным в формулировке нашей теоремы
при $p=0$.
Функция $f(\mathbf{z})$ обозначает произвольную комплексную симметрическую функцию
$z_1,\ldots ,z_n$, инвариантную относительно $z\to z^{-1}$ преобразований
и голоморфную по $z_j$ в области, заметаемой контуром $C$ при деформации к $\T$
вдоль радиальных лучей.

Очевидно, что $S_n\ltimes\mathbb{Z}_2^n$ симметричная функция
\begin{equation}\label{f}
f(\mathbf{z}) =
\frac{\Delta_n (\mathbf{z};p,q)}{\bar{\Delta}_n(\mathbf{z};q)}
\end{equation}
голоморфна по $z_j$ внутри требуемой области. Действительно, полюсы
$\Delta_n (\mathbf{z};p,q)$ при $z_j^{\pm 1}=t_0q^l$ с $|t_0q^l|>1$
компенсируются соответствующими полюсами $\bar{\Delta}_n(\mathbf{z};q)$.
Подставляя $f(\mathbf{z})$ \re{f} в \re{stokres:form}, мы получаем
утверждение теоремы \ref{residue:thm} благодаря пределу
\begin{equation*}
\lim_{y_m\rightarrow \tau_mq^{\lambda_m}}\cdots
\lim_{y_1\rightarrow \tau_1q^{\lambda_1}}
f(y_1,\ldots ,y_m,z_1,\ldots ,z_{n-m})
= \frac{\mu_m(\lambda,\mathbf{z};p,q)}{\bar{\mu}_m(\lambda ,\mathbf{z};q)}.
\end{equation*}
Это соотношение следует индуктивно из равенства
\begin{equation}
\lim_{z_{n-m+1}\rightarrow \tau_{m}q^{\lambda_{m}}}
 \frac{\mu_{m-1}(\lambda,\mathbf{z};p,q)}
      {\bar{\mu}_{m-1}(\lambda ,\mathbf{z};q)}
= \frac{\mu_{m}(\lambda,\mathbf{z};p,q)}{\bar{\mu}_{m}
(\lambda ,\mathbf{z};q)},
\end{equation}
проверяемого прямым расчетом, и перестановочной симметрии.
\hfill{Q.E.D.}

\begin{remark} При $n=1$, формула для вычетов теоремы \ref{residue:thm}
значительно упрощается
\begin{eqnarray}
\lefteqn{\frac{1}{2\pi i}\int_C
\frac{\prod_{r=0}^4\Gamma(zt_r, z^{-1}t_r; p,q)}
{\Gamma(z^2,z^{-2}, z\prod_{r=0}^4t_r, z^{-1}\prod_{r=0}^4t_r;p,q)}
\frac{d z}{z} = } && \\
&&
\frac{1}{2\pi i}\int_T
\frac{\prod_{r=0}^4\Gamma(zt_r, z^{-1}t_r; p,q)}
{\Gamma(z^2,z^{-2}, z\prod_{r=0}^4t_r, z^{-1}\prod_{r=0}^4t_r;p,q)}
\frac{d z}{z}
+2 \kappa \sum_{\stackrel{\lambda \geq 0}{|t_0q^\lambda |>1}}
\nu (\lambda ;p;q) ,\nonumber
\end{eqnarray}
где
\begin{eqnarray}
\kappa &=&
\frac{\prod_{r=1}^4\Gamma(t_rt_0,t_rt_0^{-1};p,q)}
{(q;q)_\infty(p;p)_\infty
\Gamma(t_0^{-2},t_0^{-1}\prod_{r=0}^4t_r, t_0 \prod_{r=0}^4t_r;p,q)}, \\
\nu (\lambda ;p;q) &=& q^\lambda\,
 \frac{\theta(t_0^2q^{2\lambda};p)}{\theta(t_0^2;p)}
\prod_{r=0}^5 \frac{\theta(t_rt_0;p;q)_{\lambda}}
     {\theta(qt_r^{-1}t_0;p;q)_{\lambda}}  ,
\end{eqnarray}
и $q^{-1}\prod_{r=0}^5t_r =1$.
Эта формула определяет эллиптическое расширение известного анализа вычетов
для $q$-бета интегралов Аски-Вильсона и Рахмана \cite{ask-wil:some,rah:biorthogonality}.
\end{remark}

С помощью эллиптического бета-интеграла \re{SintB} можно вывести многомерное
обобщение формулы суммирования Френкеля-Тураева для соответствующей системы корней
\re{multi-1}.

\begin{corollary}\label{wftsum:cor}
Пусть $N\in\mathbb{N}$ и параметры удовлетворяют ограничениям
\begin{equation}\label{pc}
\begin{array}{rl}
q^{-1}t^{2n-2}\prod_{r=0}^5t_r=1 & \mbox{(условие балансировки)}, \\
q^Nt^{n-1}t_0t_4=1 & \mbox{(условие обрыва)}.
\end{array}
\end{equation}
Тогда справедлива формула суммирования
\begin{equation}\label{sum-war}
\sum_{0\leq \lambda_1\leq \lambda_2\leq \cdots \leq \lambda_n\leq N}
\nu_n (\lambda ;p;q) = {\cal N}_n(p;q)
\end{equation}
как тождество для мероморфных функций по параметрам,
где $\nu_n(\lambda;p;q)$ указана в теореме \ref{residue:thm} и
\begin{eqnarray}
{\cal N}_n(p;q)
 &=&
\prod_{1\leq j< k\leq n}
\frac{\theta (q\tau_k\tau_j,qt\hat{\tau}_k^{-1}\hat{\tau}_j^{-1};p;q)_N}
  {\theta (qt^{-1}\tau_k\tau_j,q\hat{\tau}_k^{-1}\hat{\tau}_j^{-1};p;q)_N} \\
&& \times \prod_{1\leq j\leq n}
\frac{\theta (q\tau_j^2;p;q)_N\prod_{r=1}^3\theta(q\hat{t}_r\hat{\tau}_j^{-1};p;q)_N}
     {\theta (q\hat{\tau}_j^{-2};p;q)_N \prod_{r=1}^3\theta(qt_r^{-1}\tau_j;p;q)_N}
\nonumber \\
&= &\prod_{j=1}^n
\frac{\theta(qt^{n+j-2}t_0^2;p;q)_N
      \prod_{1\leq r <s \leq 3} \theta(qt^{1-j} t_r^{-1}t_s^{-1};p;q)_N}
      {\theta(qt^{2-n-j}\prod_{r=0}^3t_r^{-1};p;q)_N
      \prod_{r=1}^3 \theta(qt^{j-1}t_0t_r^{-1};p;q)_N}. \nonumber
\end{eqnarray}
Здесь использовались параметры
$\hat{t}_0=(t_0t_1t_2t_3)^{1/2}$,
$\hat{t}_1=(t_0t_1t_2^{-1}t_3^{-1})^{1/2}$,
$\hat{t}_2=(t_0t_1^{-1}t_2t_3^{-1})^{1/2}$,
$\hat{t}_3=(t_0t_1^{-1}t_2^{-1}t_3)^{1/2}$ и обозначения
$\tau_j=t_0t^{j-1}$, $\hat{\tau}_j=\hat{t}_0t^{j-1}$,
$j=1,\ldots ,n$.
\end{corollary}
{\bf Доказательство.}
Выберем параметры в соответствии с требованиями теоремы \ref{residue:thm}
так что $t^{1-n}q^{-N}< |t_0|<t^{1-n}q^{-N-1}$ для некоторого $N\in\mathbb{N}$.
Разделив формулу \re{ell-res} на $2^n n! \kappa_{n}$, и устремив $t_4$ к
$t_0^{-1}t^{1-n}q^{-N}$, мы получим указанную формулу. Для того чтобы показать
это, необходимо переписать постоянную $\kappa_{m}$ в виде
\begin{equation*}\label{km}
\kappa_{m}=\frac{1}{(q;q)_\infty^m (p;p)_\infty^m}
\prod_{j=1}^m
\frac{\Gamma(t^j;p,q)
\prod_{r=1}^4\Gamma (t^{j-1}t_0t_r,t^{1-j}t_0^{-1}t_r;p,q)}
{\Gamma(t,t^{2-m-j}t_0^{-2},t^{2n+j-3}t_0^2
\prod_{r=1}^4t_r,t^{2n-j-1}\prod_{r=1}^4t_r;p,q)},
\end{equation*}
сократив общие множители в числителе и знаменателе. Из этого выражения
видно, что отношение $\kappa_{m}/\kappa_{n}$ с $m<n$ стремится к нулю
при $t_4\to t_0^{-1}t^{1-n}q^{-N}$. Коэффициент $\kappa_{n}$ содержит полюс
при $t_4=t_0^{-1}t^{1-n}q^{-N}$ благодаря множителю $\Gamma (t^{n-1}t_0t_4;p,q)$,
который отсутствует в $\kappa_m$ при $m<n$. Поскольку интеграл
$$
\int_{\T^{n-m}} \Delta_{n-m} (\mathbf{z};p,q)
\delta_{m,n-m}(\lambda ,\mathbf{z})
\frac{d z_1}{z_1}\cdots \frac{d z_{n-m}}{z_{n-m}}
$$
и коэффициенты $\nu_m(\lambda ; p; q )$
остаются ограниченными при $t_4\rightarrow t_0^{-1}t^{1-n}q^{-N}$,
можно заключить, что члены с $m<n$ в сумме правой части формулы для вычетов
стремятся к нулю. Член с $m=n$ приводит к сумме
$\sum_{0\leq\lambda_1\leq\cdots\leq\lambda_n\leq N} \nu_n(\lambda ;p;q)$,
равной ряду в левой части равенства \re{sum-war}.

После деления правой части на $2^n n! \kappa_{n}$, получаем
\begin{equation}
\prod_{j=1}^n
\frac{\Gamma (t^{2-n-j}t_0^{-2},t^{2n+j-3}t_0^2\prod_{r=1}^4t_r;p,q)
       \prod_{1\leq r<s\leq 4} \Gamma (t^{j-1}t_r t_s;p,q)}
      {\prod_{r=1}^4
\Gamma (t^{1-j}t_0^{-1}t_r,t^{n+j-2}t_r^{-1}\prod_{s=1}^4t_s;p,q)} ,
\end{equation}
что в пределе $t_4\rightarrow t_0^{-1}t^{1-n}q^{-N}$ равно правой части
указанной формулы. Полученный результат может быть аналитически продолжен
до значений параметров $q$, $t$ и $t_r$, $r=0,\ldots ,5$, ограниченных
только условиями \re{pc}.
\hfill{Q.E.D.}

\smallskip

Данная формула суммирования была впервые предложена Варнааром в работе
\cite{war:summation}. Ее $p=0$ аналог был доказан ван Диехеном и автором в
статье \cite{die-spi:elliptic} с помощью указанной выше процедуры
(ее чисто гипергеометрическое вырождение приведено в обзоре \cite{die-spi:review}),
примененной к установленному ранее Густафсоном
$p=0$ предельному значению интеграла \cite{gus:some2}
(другое доказательство было дано в работе \cite{rai:pols}). Впервые эта формула была
полностью доказана Розенгреном по индукции в статье \cite{ros:proof}. После
доказательства эллиптических бета-интегралов в работах \cite{rai:trans,spi:short}
указанная выше теорема получила свой полноценный статус. Отметим, что до анализа
эллиптических рядов, была известна только формула суммирования такого типа
на уровне многомерного $_6\varphi_5$ ряда \cite{die:certain} (которая, в свою
очередь, связана с многомерными формулами суммирования работ
\cite{aom:elliptic,ito:theta}). А самый общий анализ вычетов, проведенный
Стокманом в \cite{sto:basic}, соответствовал выводу многомерной
${}_6\varphi_5$ формулы суммирования \cite{die:certain} с помощью вычетов в
Густафсоновском аналоге
многомерного интеграла Аски-Вильсона \cite{gus:generalization,kad:proof}.
Другие многомерные аналоги ${}_8\varphi_7$ формулы суммирования Джексона
можно найти в работах
\cite{mil:multiple,mil:multidimensional,mil-lil:consequences,sch:summation}.

В качестве открытой проблемы мы приведем гипотезу, выдвинутую в статье \cite{spi:integrals},
о формуле суммирования, которая должна быть связана с анализом вычетов
в $A_n$ эллиптических бета-интегралах (\ref{int-an-o})-(\ref{int-an-e}).
Эта гипотеза определяет эллиптический аналог теоремы 1.2 работы \cite{gus-rak:beta}.

\begin{conjecture}
Предположим, что $N$ положительное целое и $\prod_{k=1}^n t_k=q^{-N}$.
Тогда
\begin{eqnarray}\nonumber
&&
\sum_{\stackrel{\lambda_k=0,\ldots,N}{\lambda_1+\ldots+\lambda_n=N}}
\frac{
\prod_{1\leq i<j\leq n}\theta(tt_it_j)_{\lambda_i+\lambda_j}
\prod_{i=1}^n\prod_{j=n+1}^{n+3}\theta(tt_it_j)_{\lambda_i}
\prod_{i,j=1}^n\theta(t_it_j^{-1})_{-\lambda_j}
}{
\prod_{i,j=1;i\neq j}^n\theta(t_it_j^{-1})_{\lambda_i-\lambda_j}
\prod_{j=1}^n\theta(t^{n+1}t_j^{-1}\prod_{k=1}^{n+3}t_k)_{-\lambda_j}
} \\ && \makebox[2em]{}
= \left\{ \begin{aligned} % \begin{array}{ll}
\frac{\theta(1)_{-N}}{\theta(t^{n/2})_{-N}
\prod_{n+1\leq i<j\leq n+3}\theta(t^{(n+2)/2}t_it_j)_{-N}},
& \quad n \quad \mbox{четно,}  \\
\frac{\theta(1)_{-N}}{\prod_{i=n+1}^{n+3}\theta(t^{(n+1)/2}t_i)_{-N}
\theta(t^{(n+3)/2}\prod_{i=n+1}^{n+3}t_i)_{-N} },
& \quad n \quad \mbox{нечетно,}  \\
\end{aligned} \right.
\label{new-sums} \end{eqnarray}
где $\theta(a)_\lambda\equiv\theta(a;p;q)_\lambda$.
\end{conjecture}

Некоторое свидетельство в пользу формулы (\ref{new-sums})
дает следующая теорема.

\begin{theorem}
Обозначим $t=q^g, t_i=q^{g_i},i=1,\ldots,n+3$ (так что
$\sum_{j=1}^{n}g_j+N=0$). Ряд
$\sum_{\mathbf{\lambda}}c(\mathbf{\lambda})$,
стоящий в правой части (\ref{new-sums}) является полностью эллиптическим
гипергеометрическим рядом, то есть отношение соседних членов ряда
$$
h_k(\mathbf{\lambda})=\frac{c(\lambda_1,\ldots,\lambda_k+1,\ldots,
\lambda_n)}{c(\lambda_1,\ldots,\lambda_n)}
$$
есть эллиптические функции всех свободных переменных в множестве
$(\lambda_1,\ldots,\lambda_n,$ $g,$ $g_1,$ $\ldots,$ $ g_{n+3})$.
Более того, $h_k(\mathbf{\lambda})$ являются $PSL(2,\mathbb{Z})$
модулярно инвариантными функциями. Отношения выражений в левой и правой частях
(\ref{new-sums}) так же являются эллиптическими функциями $g$ и $n+2$
свободных параметров из множества $(g_1,\ldots,g_{n+3})$ и эти отношения так
же модулярно инвариантны.
\end{theorem}

Доказательство этого утверждения достаточно длинное, но оно сводится к элементарным
выкладкам по проверке сокращения простых множителей, появляющихся при соответствующих
преобразованиях. По структуре оно не отличается от процедуры, описанной в предыдущей
главе и работах \cite{spi:theta,spi:modularity} для других рядов. Поскольку параболические
формы с весом меньшим 12 не существуют, так же как и в работах
\cite{die-spi:elliptic,spi:modularity} из этой теоремы следует справедливость
(\ref{new-sums}) для первых 11 членов разложения по малому параметру $\sigma$.
Естественно ожидать, что эллиптические бета-интегралы (\ref{int-an-odd}) и (\ref{int-an-even}) так же
приводят к некоторым нетривиальным многомерным $A_n$ формулам суммирования,
аналогичным (\ref{new-sums}). Как показано в \cite{spi:integrals},
эллиптическая $A_n$-формула суммирования Милна \re{A_n} вытекает из анализа вычетов в
$A_n$ эллиптическом бета-интеграле типа I \re{ell-int-A}. В работе \cite{die-spi:modular}
таким образом выведена формула \re{multi-2}  из многопараметрического
$C_n$ эллиптического бета-интеграла.
Аналогично, эллиптическое обобщение $D_n$-формулы суммирования Бхатнагара-Шлоссера \re{D_n}
вытекает из формулы интегрирования \re{ell-int-D}, открытой Варнааром и автором \cite{spi-war:inversions}.
Представляет интерес вывод тождеств для эллиптических гипергеометрических
функций при $q$ равном примитивным корням единицы, $q^n=1$. В частности,
с помощью анализа вычетов можно вывести тождества для эллиптических
гипергеометрических рядов, обобщающие результаты работы \cite{SZAW}
(которые охватывают хорошо известные суммы Гаусса).

%% file: CHAPTER4.TEX
\chapter[Биортогональные функции]{Биортогональные функции}

В этой главе мы опишем новый класс биортогональных функций, определенных
с помощью эллиптических гипергеометрических рядов. Для удобства изложения
мы дадим вывод основных свойств этих функций в мультипликативной системе
обозначений.

\section{Эллиптическое гипергеометрическое уравнение}

Предположим, что $|p|<1$ и $|q|<1$ и положим
$$
\kappa=\frac{\pqinf}{4\pi i}.
$$
Возьмем восемь комплексных параметров $t_1,\ldots, t_8\in\C$, удовлетворяющих
ограничению
$$
\prod_{k=1}^8t_k=p^2q^2,
$$
обозначим $t=(t_1,\ldots,t_8)$, и определим $V$-функцию
\ba\label{V-fn}
&& V(t)\equiv V(t_1,\ldots,t_8;p,q)
=\kappa\int_C \Delta(z,t)\frac{dz}{z}, \quad
\\ \lab{kern} && \makebox[3em]{}
\Delta(z,t)=\frac{\prod_{k=1}^8\Gamma(t_kz^\pm)}{\Gamma(z^{\pm2})},
\qquad \Gamma(z)\equiv \Gamma(z;q,p),
\ea
где $C$ обозначает любой гладкий контур, ориентированный против часовой стрелки и
разделяющий последовательности полюсов подынтегральной функции в
\eqref{V-fn}, сходящихся к точке $z=0$, от их $z\to 1/z$ отраженных
партнеров, уходящих на бесконечность.

Как уже было показано в третьей главе
$$
V(t_1,\ldots,t_7,pq/t_7)=\prod_{1\leq j<k\leq 6}\Gamma(t_jt_k),
$$
где $\prod_{k=1}^6t_k=pq$, что совпадает с эллиптическим бета-интегралом.
С помощью этой точной формулы интегрирования, в работе
\cite{spi:integrals} было доказано следующее нетривиальное преобразование
симметрии для эллиптической гипергеометрической функции $V(t)$.

\begin{theorem}
Возьмем произвольные параметры $c,a_{1,2,3},b_{1,2,3}\in\C$
в общем положении и положим $A=a_1a_2a_3$, $B=b_1b_2b_3$.
Тогда,
\begin{eqnarray}\nonumber
&& \prod_{j=1}^3\frac{\Gamma(A/a_j)}{\Gamma(c^2A/a_j)}
\int_{C}\frac{\prod_{j=1}^3\Gamma(ca_jz^\pm,b_jz^\pm)}
{\Gamma(z^{\pm2},cAz^\pm,c^2Bz^\pm)}\frac{dz}{z}
\\ && \makebox[4em]{}
= \prod_{j=1}^3\frac{\Gamma(B/b_j)}{\Gamma(c^2B/b_j)}
\int_{C}\frac{\prod_{j=1}^3\Gamma(a_jz^\pm,cb_jz^\pm)}
{\Gamma(z^{\pm2},c^2Az^\pm,cBz^\pm)}\frac{dz}{z},
\label{symm-tr}\end{eqnarray}
где контур $C$ разделяет сходящиеся к $z=0$ и уходящие на бесконечность
последовательности полюсов подынтегральных функций в обоих
сторонах равенства (существование такого контура -- единственное условие,
ограничивающее выбор параметров).
\end{theorem}

Это утверждение легко выводится рассмотрением составного интеграла
\be
\kappa \int_{C^{2}}
\frac{\prod_{j=1}^{3}\Gamma(a_jz^\pm,b_jw^\pm)\;
\Gamma(cz^\pm w^\pm)}{\Gamma(z^{\pm2},w^{\pm2},c^{2}Az^\pm,c^2Bw_k^\pm)}
\frac{dz}{z}\frac{dw}{w}
\nonumber\ee
и двумя возможными путями снятия одного из интегрирований: либо интегрированием
по $w$, либо по $z$ с помощью эллиптического бета-интеграла. При этом перестановка
порядков интегрирования разрешена ввиду ограниченности подынтегральной функции
на контуре интегрирования.

На уровне рядов, тождество \eqref{symm-tr} соответствует эллиптическому
обобщению четырех\-членного преобразования Бэйли для необрывающихся
совершенно уравновешенных $_{10}\varphi_9$ рядов, хотя его и невозможно
записать в виде тождества для комбинации бесконечных
эллиптических гипергеометрических рядов из-за проблем со сходимостью.
Оно сыграло ключевую роль
при выводе интегральных аналогов леммы Бэйли в работе \cite{spi:tree}, позволивших
построить бесконечное двоичное дерево тождеств для многократных эллиптических
гипергеометрических интегралов различных кратностей.

После обозначения
$$
t_{1,2,3}= ca_{1,2,3},\quad t_4=\frac{pq}{cA},\quad
t_{5,6,7}=b_{1,2,3}, \quad t_8=\frac{pq}{c^2B},
$$
равенство \eqref{symm-tr} может быть переписано в более симметричной форме
\begin{equation}
V(t)=\prod_{1\leq i<j\leq 4}\Gamma(t_it_j,t_{i+4}t_{j+4})\, V(s),
\label{symm}\end{equation}
где
$$
s_{1,2,3,4}=\epsilon^{-1} t_{1,2,3,4},\quad s_{5,6,7,8}=\epsilon t_{5,6,7,8},
$$
и
$$
\epsilon=\sqrt{\frac{t_1t_2t_3t_4}{pq}}=\sqrt{\frac{pq}{t_5t_6t_7t_8}}.
$$
Что касается явного вида контура интегрирования $C$, фигурирующего в обеих
частях \eqref{symm-tr}, при $\sqrt{|pq|}<|t_k|<1$ мы можем взять $C=\T$
(единичная окружность).

Благодаря перестановочной симметрии в левой части, мы имеем
 $C^4_8=70$ различных преобразований такого типа. Два из них
являются инволюциями (или отражениями), а именно, если мы повторно применим
то же самое преобразование к правой части \eqref{symm}
с параметрами $s_{1,2,3,4}$, играющими роль $t_{1,2,3,4}$
(или $t_{5,6,7,8}$), то вернемся к функции $V(t)$.

Применение $2C^1_4C^3_4=32$ преобразований, для которых
роль $t_{1,2,3,4}$ исполняется параметрами $s_{2,3,4,5}$ или $s_{1,6,7,8}$
(и другими аналогичными комбинациями параметров), приводит к простой перестановке
параметров в \eqref{symm}. Однако, преобразования с параметрами  $s_{3,4,5,6}$,
играющими роль $t_{1,2,3,4}$ (имеется $C^2_4C^2_4=36$
таких преобразований),  дает
\begin{equation}\label{double-tr}
V(t)=\prod_{1\leq i<j\leq 4}\Gamma(t_it_j,t_{i+4}t_{j+4},
s_{i+2}s_{j+2},s_{i+6}s_{j+6})\, V(f),
\end{equation}
где для удобства считается, что $s_{k+8}\equiv s_k$, $k=1,\ldots,8$, и
$$
f_{1,2,7,8}=\rho s_{1,2,7,8},\quad f_{3,4,5,6}=\rho^{-1}s_{3,4,5,6},
\quad\rho=\sqrt{\frac{s_3s_4s_5s_6}{pq}}=\sqrt{\frac{t_3t_4t_5t_6}{pq}}.
$$
Переставим в этом преобразовании пару $t_3,t_4$ с $t_5, t_6$. Это приводит к
следующему соотношению
\begin{equation}\label{double-tr-perm}
V(t)=\prod_{i,j=1}^4 \Gamma(t_it_{j+4})\, V\left(\frac{\sqrt{T}}{t_1},\ldots,
\frac{\sqrt{T}}{t_4},\frac{\sqrt{U}}{t_5},\ldots,\frac{\sqrt{U}}{t_8}\right),
\end{equation}
где $T=t_1t_2t_3t_4,\, U=t_5t_6t_7t_8$.

Приравнивая правые части соотношений \eqref{symm} и \eqref{double-tr-perm}
и выражая параметры $t_1,\ldots,t_8$  через $s_1,\dots, s_8$, мы приходим  к
равенству
 \begin{equation}
V(s)=\prod_{1\leq i<j\leq 8}\Gamma(s_is_j)\, V\left(\frac{\sqrt{pq}}{s}\right),
\label{symm-2}\end{equation}
где $\sqrt{pq}/s=(\sqrt{pq}/s_1,\ldots,\sqrt{pq}/s_8)$.
Тождества \eqref{double-tr-perm} и \eqref{symm-2} были доказаны Рэйнсом в
 \cite{rai:trans} весьма сложным образом, но, как мы показали,  эти
тождества принадлежат группе преобразованием симметрии, порожденной преобразованием
найденным ранее автором в работе  \cite{spi:integrals}.
Этот результат позволяет предположить, что
аналогичные преобразования симметрии для многомерных эллиптических гипергеометрических
интегралов, построенные в  \cite{rai:trans}, так же можно будет вывести
достаточно простым способом из многомерного обобщения на корневые системы
интегральных пар Бэйли работы \cite{spi:tree}.
Отметим, что указанные преобразования симметрии связаны с группой отражений в
решетке корневой системы $E_7$, частично описанной в работах \cite{kmnoy}
и \cite{rai:trans}.

Выведем теперь соотношения сопряжения для $V$-функции.
Возьмем соотношение Римана
\be\label{ident}
\theta(xw^\pm,yz^\pm;p)-\theta(xz^\pm,yw^\pm;p)=
yw^{-1}\theta(xy^\pm,wz^\pm;p)
\ee
и положим в нем $y=t_1, w=t_2, $ и $x=q^{-1}t_8$. Тогда оно приводит к следующему
$q$-разностному уравнению для подынтегральной функции \re{kern}:
\ba
&& \Delta(z,qt_1,t_2,\ldots,t_7,q^{-1}t_8)
-\frac{\theta(t_1t_2^\pm;p)}{\theta(q^{-1}t_8t_2^\pm;p)} \Delta(z,t)
\\ && \makebox[2em]{}
=\frac{t_1}{t_2}\frac{\theta(q^{-1}t_8t_1^\pm;p)}{\theta(q^{-1}t_8t_2^\pm;p)}
\Delta(z, t_1,qt_2,t_3,\ldots, t_7,q^{-1}t_8).
\ea
Интегрируя это соотношение по $z$ вдоль контура $C$ (выбор которого в этом
случае вообще говоря произволен), не пересекающего сингулярности
$\Delta(z,t)$, мы получаем первое соотношение сопряжения для $V$-функции
\ba\nonumber
&& t_2\theta(q^{-1}t_8t_2^\pm;p)V(qt_1,q^{-1}t_8)-
t_1\theta(q^{-1}t_8t_1^\pm;p)V(qt_2,q^{-1}t_8)
\\ && \makebox[8em]{}
=t_2\theta(t_1t_2^\pm;p)V(t).
\lab{cont-1}\ea
Очевидно, что мы имеем целую совокупность таких соотношений, получающихся
перестановками параметров. Одно из них использовалось автором в самом
первом доказательстве эллиптического бета-интеграла \cite{spi:elliptic,spi:umn}.

Вторая совокупность соотношений сопряжения возникает после подстановки
преобразований симметрии \re{symm} и \re{symm-2} в \re{cont-1} и это
накладывает ограничения на выбор контура интегрирования $C$.
В частности, \re{symm-2} приводит к равенству
\ba\nonumber
&& t_1\theta\left(\frac{t_2}{qt_8};p\right)\prod_{k=3}^7
\theta\left(\frac{t_1t_k}{q};p\right) V(q^{-1}t_1,qt_8)
\\ \nonumber && \makebox[2em]{}
- t_2\theta\left(\frac{t_1}{qt_8};p\right)\prod_{k=3}^7
\theta\left(\frac{t_2t_k}{q};p\right) V(q^{-1}t_2,qt_8)
\\  && \makebox[6em]{}
=t_1\theta\left(\frac{t_2}{t_1};p\right) \prod_{k=3}^7\theta(t_kt_8;p)
\, V(t).
\lab{cont-2}\ea
Комбинируя равенства \re{cont-1} и \re{cont-2}, мы получаем третий
набор соотношений сопряжения
\ba\nonumber
&& \frac{\theta(t_1/qt_8,t_1t_8,t_8/t_1;p)}
{\theta(t_1/t_2,t_2/qt_1;p)}\prod_{k=3}^7
\theta(t_2t_k/q;p)\left( U(t)-U(qt_1,q^{-1}t_2)\right)
\\ \nonumber && \makebox[2em]{}
+ \frac{\theta(t_2/qt_8,t_2t_8,t_8/t_2;p)}
{\theta(t_2/t_1,t_1/qt_2;p)}\prod_{k=3}^7
\theta(t_1t_k/q;p)\left( U(t)-U(q^{-1}t_1,qt_2)\right)
\\  && \makebox[8em]{}
=\theta(t_1t_2/q;p) \prod_{k=3}^7\theta(t_kt_8;p)\, U(t),
\lab{cont-3}\ea
где
$$
U(t)=\frac{V(t)}{\prod_{j=1}^7\Gamma(t_jt_8,t_j/t_8)}.
$$
Если подставить в этом равенстве параметры $t_1=az, t_2=a/z$ и заменить
функцию $U(t)$ на неизвестную функцию $f(z)$, то
мы получим $q$-разностное уравнение второго порядка, которое мы называем
эллиптическим гипергеометрическим уравнением
\ba\nonumber
&& \frac{\theta(az/qt_8,azt_8,t_8/az;p)}
{\theta(z^2,1/qz^2;p)}\prod_{k=3}^7
\theta(at_k/qz;p)\left( f(z)-f(qz)\right)
\\ \nonumber && \makebox[2em]{}
+ \frac{\theta(a/qzt_8,at_8/z,zt_8/a;p)}
{\theta(1/z^2,z^2/q;p)}\prod_{k=3}^7
\theta(azt_k/q;p)\left( f(z)-f(q^{-1}z\right)
\\  && \makebox[8em]{}
=\theta(a^2/q;p) \prod_{k=3}^7\theta(t_kt_8;p)\, f(z),
\lab{ehe}\ea
где $a^2\prod_{k=3}^8t_k=p^2q^2$.
Одно решение этого уравнения нам уже известно
$$
f(z)=\frac{V(az,a/z,t_3,\ldots,t_8)}{\Gamma(azt_8,az/t_8,at_8/z,a/zt_8)}.
$$
Необходимо отметить, что при подходящих ограничениях на параметры функции $V(t)$
и контура интегрирования $C$ она становиться полностью симметричной
по своим аргументам и, поэтому, она удовлетворяет
$C^2_8=28$ уравнениям такого типа. Помимо этого для каждого из этих уравнений
можно переставлять параметры в коэффициентах уравнений, например, $t_8$ с $t_k,
k=3,\ldots,7$, в \eqref{cont-3}.  Более того, все эти уравнения симметричны по $q$ и $p$,
что удваивает число уравнений. По всей видимости существует только одно решение,
удовлетворяющее всем этим уравнениям.
Однако, если ограничиться только одним
разностным уравнением второго порядка \eqref{ehe}, то тогда легко найти
линейно-независимые решения. Например, ряд из них получается простым
домножением $z$ или любого из параметров $f(z)$ на $p$ (при этом нарушается
симметрия между $p$ и $q$ и перестановочная симметрия между параметрами).
Другие решения получаются заменой стандартных
эллиптических гамма-функций на их модифицированную версию, хорошо определенную
при $|q|=1$ (см. главу 5). При этом работает несколько другая параметризация
переменных, не использующая только базисные переменные $p$ и $q$.
Другой возможный способ построения независимых решений состоит в модификации
контура интегрирования $C$ в определении $V(t)$ на такой, что $f(z)$ по
прежнему задает решение \eqref{ehe} по аналитическому продолжению,
но равенства \eqref{symm}, \eqref{double-tr-perm} или \eqref{symm-2}
уже не выполняются. Соответствующие решения были бы $p,q$-симметричны, но
не обладали бы $E_7$-симметрией. Эта и другие возможности построения решений
эллиптического гипергеометрического уравнения (т.е., описание эллиптического аналога
24 Куммеровских решений обычного гипергеометрического уравнения) требуют
тщательного анализа. Следующим шагом является построение трехчленных
соотношений между этими независимыми решениями с определенными
коэффициентами - эта задача будет рассмотрена отдельно от настоящей диссертации.

Рассмотрим совершенно уравновешенный эллиптический гипергеометрический ряд
\be
_{r+1}V_r(t_0;t_1,\ldots,t_{r-4};q,p;x)
\equiv \sum_{n=0}^\infty \frac{\theta(t_0q^{2n};p)}{\theta(t_0;p)}
\prod_{m=0}^{r-4}\frac{\theta(t_m;p;q)_n}{\theta(qt_0t_m^{-1};p;q)_n}
\, (qx)^n.
\label{vwp-series}\ee
Общее условие балансировки (\ref{theta-balance}) для этого ряда имеет вид
$$
\prod_{m=1}^{r-4}t_m= \pm t_0^{(r-5)/2}q^{(r-7)/2}.
$$
В этом общем случае
сбалансированный $_{r+1}V_r$ ряд инвариантен при независимых сдвигах $t_0\to p^2t_0$
и $t_j\to pt_j, \, j=1,\ldots, r-5,$ и отношение последовательных членов ряда
$h(n)=c_{n+1}/c_n$ есть эллиптическая функция $n$. Для четных $r+1$ можно добиться
инвариантности относительно замены $t_0\to pt_0$, т.е. условия полной эллиптичности
с равными периодами. Для этого необходимо в условии
балансировки выбрать положительный знак.
При $p\to 0$, ряд $_{r+1}V_r$ сводится к совершенно уравновешенному
$q$-гипергеометрическому ряду $_{r-1}W_{r-2}$ от аргумента $qx$
(в обозначениях монографии \cite{gas-rah:basic}).

Сбалансированный $_{12}V_{11}$ ряд с $x=1$ играет важную роль
в приложениях. Например, эллиптические решения уравнения Янга-Бакстера,
полученные Дате и др. в работах \cite{djkmo:exactly1,djkmo:exactly2}, выражаются
через такой ряд при специальном выборе параметров \cite{fre-tur:elliptic}.
Во всех случаях рассмотренных ниже, $_{12}V_{11}$ ряды имеют $x=1$;
поэтому мы опускаем зависимость от этого аргумента в обозначениях.

Рассмотрим функцию ${\cal E}(\mathbf{t}) \equiv
{_{12}}V_{11}(t_0;t_1,\ldots,t_7;q,p),$
где $\prod_{m=1}^7t_m=t_0^3q^2$, и будем считать, что этот ряд обрывается
благодаря условию $t_m=q^{-n},\, n\in\mathbb{N},$ для некоторого $m$.
В работах \cite{spi-zhe:spectral,spi-zhe:classical} (см. также главу 1),
были установлены два соотношения сопряжения для ${\cal E}(\mathbf{t})$:
\begin{eqnarray}\label{1_con}
\lefteqn{{\cal E}(\mathbf{t}) - {\cal E}(t_0;t_1,\ldots,t_5,q^{-1}t_6,
qt_7) } &&  \\ &&
=\frac{\theta(qt_0,q^2t_0,qt_7/t_6,t_6t_7/qt_0;p)}
{\theta(qt_0/t_6,q^2t_0/t_6,t_0/t_7,t_7/qt_0;p)}
\prod_{r=1}^5\frac{\theta(t_r;p)}{\theta(qt_0/t_r;p)}\,
{\cal E}(q^2t_0;qt_1,\ldots,qt_5,t_6,qt_7),
 \nonumber\\  \nonumber
\lefteqn{\frac{\theta(t_7;p)}{\theta(t_6/qt_0,t_6/q^2t_0,t_6/t_7;p)}
\prod_{r=1}^5 \theta(t_rt_6/qt_0;p)\,
{\cal E}(q^2t_0;qt_1,\ldots,qt_5,t_6,qt_7) } &&
\\ && \makebox[4em]{}
+\frac{\theta(t_6;p)}{\theta(t_7/qt_0,t_7/q^2t_0,t_7/t_6;p)}
\prod_{r=1}^5 \theta(t_rt_7/qt_0;p)\,
{\cal E}(q^2t_0;qt_1,\ldots,qt_6,t_7)
\nonumber \\ && \makebox[8em]{}
=\frac{1}{\theta(qt_0,q^2t_0;p)}
\prod_{r=1}^5\theta(qt_0/t_r;p)\, {\cal E}(\mathbf{t}).
\label{2_con} \end{eqnarray}
Для доказательства первого из них достаточно применить тождество для тета-функций
\eqref{ident} к разности рядов стоящих в левой части \eqref{1_con}.
Второе равенство вытекает из первого после применения специальной комбинации
двух преобразований Бэйли
для $_{12}V_{11}$ рядов (см. подробный вывод в статье \cite{spi-zhe:theory}).

Другой способ вывода этих соотношений состоит в применении анализа вычетов
в интегральных соотношениях сопряжения \eqref{cont-1} и \eqref{cont-2}.
Для этого достаточно взять $|t_j|<1$ и $C=\T$, затем растянуть $t_8$ на
область $|t_8|>1$ с соответствующей деформацией $C$, чтобы ни один полюс
подынтегральной функции в $V(t)$ не пересек контура интегрирования.
После этого необходимо деформировать $C$ назад к $\T$ и просуммировать вычеты
всех пересеченных полюсов. Это приводит к тому, что разница между интегралами
по контурам $C$ и $\T$ выглядит как конечная сумма членов $_{12}V_{11}$ ряда.
В пределе $t_mt_8\to q^{-n}$ для некоторого фиксированного $m$
вычеты начинают расходиться, но после простой
перенормировки мы получаем, что интегральная часть зануляется, а функция
$V(t)$ эффективно заменяется на обрывающийся $_{12}V_{11}$ ряд. После
замены $t_8\to t_0$ и некоторых других простых переобозначений параметров,
мы получим указанные соотношения сопряжения для рядов.

Комбинация \re{1_con} и \re{2_con} дает
\begin{eqnarray}\nonumber
\lefteqn{
\frac{\theta(t_7,t_0/t_7,qt_0/t_7;p)}{\theta(qt_7/t_6,t_7/t_6;p)}
\prod_{r=1}^5\theta(qt_0/t_6t_r;p)\left(
{\cal E}(t_0;t_1,\ldots,t_5,q^{-1}t_6,qt_7) -
{\cal E}(\mathbf{t}) \right) } &&
\\  \nonumber
&& +\frac{\theta(t_6,t_0/t_6,qt_0/t_6;p)}
{\theta(qt_6/t_7,t_6/t_7;p)}\prod_{r=1}^5\theta(qt_0/t_7t_r;p)
\left({\cal E}(t_0;t_1,\ldots,t_5,qt_6,q^{-1}t_7)-
{\cal E}(\mathbf{t})\right) \\
&& \makebox[4em]{}
+\theta(qt_0/t_6t_7;p)\prod_{r=1}^5\theta(t_r;p)\,
{\cal E}(\mathbf{t})=0.
\label{3_con}\end{eqnarray}
При $p\to 0$, указанные три уравнения сводятся к соотношениям сопряжения
для обрывающегося совершенно уравновешенного сбалансированного
$_{10}\varphi_9$ ряда, полученных Гуптой и Массоном \cite{gup-mas:contiguous}.
Аналогичные соотношения на уровне $_8\varphi_7$ рядов были получены ранее
Исмаилом и Рахманом в работе \cite{ism-rah:associated}.

\section{Обобщенная задача на собственные значения и трехчленное рекуррентное соотношение}

Изменим параметризацию ${\cal E}$-функции и рассмотрим соотношение
(\ref{3_con}) для следующей функции:
\begin{equation}
R_{n}(z;q,p)\equiv {_{12}V_{11}}\left(\frac{t_3}{t_4};
\frac{q}{t_0t_4},\frac{q}{t_1t_4},\frac{q}{t_2t_4},t_3z,\frac{t_3}{z},q^{-n},
\frac{Aq^{n-1}}{t_4};q,p\right),
\label{R_n}\end{equation}
где $A=\prod_{r=0}^4 t_r$. Мы заменяем параметр $t_0$
в (\ref{3_con}) на $t_3/t_4$; переменные $t_1,t_2,$ и $t_3$
на $q/t_0t_4,q/t_1t_4,$ и $q/t_2t_4$; переменные $t_4$ и $t_5$ на
$q^{-n}$ и $Aq^{n-1}/t_4$; и параметры $t_6$ и $t_7$ на $t_3z$ и $t_3/z$,
соответственно. В результате мы видим, что $R_{n}(z;q,p)$ задает
специальное решение следующего конечно-разностного уравнения:
\begin{equation}
{\cal D}_\mu f(z)=0,\qquad
{\cal D}_\mu=V_\mu(z)(T-1)+V_\mu(z^{-1})(T^{-1}-1)+\kappa_\mu,
\label{eig}\end{equation}
где $T$ обозначает оператор $q$-растяжения, $Tf(z)=f(qz)$, и
\begin{eqnarray}\label{V}
&& V_\mu(z)=\theta\left(\frac{t_4}{q\mu z},\frac{A\mu}{q^2z},\frac{t_4z}{q};p\right)
\frac{\prod_{r=0}^4\theta(t_rz;p)}{\theta(z^2,qz^2;p)},
\\
&& \kappa_\mu=\theta\left(\frac{A\mu}{qt_4},\mu^{-1};p\right)
\prod_{r=0}^3\theta\left(\frac{t_rt_4}{q};p\right).
\label{kappa}\end{eqnarray}
Функции $f(z)=R_n(z;q,p)$ определяют решение (\ref{eig})
при $\mu=q^{n}, n\in\mathbb{N}.$

Уравнение (\ref{eig}) выглядит как нестандартная задача на собственные
значения со ``спектральным параметром" $\mu$. Действительно,
оно может быть переписано как обобщенная задача на собственные значения
\begin{equation}\label{gevp}
 {\cal D}_\eta f(z)=\lambda {\cal D}_\xi f(z),
\end{equation}
где спектральный параметр $\lambda$ имеет вид
\begin{equation}
\lambda=\frac{\theta\left(\frac{\mu A\eta}{qt_4},
\frac{\mu }{\eta};p\right)}{\theta\left(\frac{\mu A\xi}{qt_4},
\frac{\mu }{\xi};p\right)}
\label{lambda}\end{equation}
и операторы ${\cal D}_{\xi}$, ${\cal D}_{\eta}$ получаются
из ${\cal D}_{\mu}$ после замены $\mu$ на произвольные калибровочные
параметры $\xi,\eta \in\mathbb{C},\, \xi\neq\eta p^k,
qt_4p^k/A\eta,$ $k\in\mathbb{Z}$. Применение тождества для тета-функций
(\ref{ident}) к операторной части уравнения (\ref{gevp}) дает
$$
{\cal D}_\eta-\lambda {\cal D}_\xi
=\frac{\theta(A\eta\xi/qt_4,\xi/\eta;p)}
{\theta(A\mu \xi/qt_4, \xi/\mu ;p)}\, {\cal D}_{\mu}.
$$
Как видно, калибровочные параметры $\xi,\eta$ полностью выпадают
из уравнения, определяющего $f(z)$, ${\cal D}_\mu f(z)=0$.

Трехчленное рекуррентное соотношение для функций $R_n(z;q,p)$ было
получено в работах \cite{spi-zhe:spectral,spi-zhe:classical}. Оно появляется из
формулы (\ref{3_con}) после замены $t_6$ на $q^{-n}$ и
$t_7$ на $Aq^{n-1}/t_4$, а также подстановок $t_1\to q/t_0t_4,
t_2\to q/t_1t_4, t_3\to q/t_2t_4$, $t_4\to t_3z$,
$t_5\to t_3/z$. После некоторых выкладок, эта формула приобретает вид
\begin{eqnarray}\nonumber
&& (\gamma(z)-\alpha_{n+1})B(Aq^{n-1}/t_4)\left(R_{n+1}(z;q,p)-R_n(z;q,p)
\right) \\ \nonumber && \makebox[2em]{}
+(\gamma(z)-\beta_{n-1})B(q^{-n})
\left(R_{n-1}(z;q,p)-R_n(z;q,p)\right)
\\ && \makebox[4em]{}
+\delta\left(\gamma(z)-\gamma(t_3)\right)R_n(z;q,p)=0,
\label{ttr-ell}\end{eqnarray}
где
\begin{eqnarray}\label{B}
&& B(x)=\frac{\theta\left(x,\frac{t_3}{t_4x},
\frac{qt_3}{t_4x},\frac{qx}{t_0t_1},\frac{qx}{t_0t_2},
\frac{qx}{t_1t_2},\frac{q^2\eta x}{A},\frac{q^2x}{A\eta};p
\right)}{\theta\left(\frac{qt_4x^2}{A},\frac{q^2t_4x^2}{A};p\right)},
\\ \label{gamma}
&& \delta=\theta\left(\frac{q^2t_3}{A},\frac{q}{t_0t_4},
\frac{q}{t_1t_4},\frac{q}{t_2t_4},t_3\eta,\frac{t_3}{\eta};p\right),
\\ &&
\gamma(z)=\frac{\theta(z\xi,z/\xi;p)}{\theta(z\eta,z/\eta;p)},
\label{ab} \\ &&
\alpha_n=\gamma(q^n/t_4),\qquad \beta_n=\gamma(q^{n-1}A).
\end{eqnarray}
Здесь  $\xi$ и $\eta \neq \xi p^k,\xi^{-1}p^k,$ $k\in\mathbb{Z},$
произвольные калибровочные параметры (они не связаны с $\xi,\eta$
в разностном уравнении, но мы используем те же самые обозначения).
Подставляя выражения (\ref{B})-(\ref{ab}) в уравнение (\ref{ttr-ell}) и применяя тождество
(\ref{ident}), можно увидеть, что вспомогательные калибровочные параметры
$\xi, \eta$ полностью выпадают из результирующего рекуррентного соотношения.

Поскольку $B(q^{-n})=0$ при $n=0$, неизвестная величина $R_{-1}$
не входит в (\ref{ttr-ell}) при $n=0$. Можно сказать, что функции $R_n(z;q,p)$
генерируются трехчленным рекуррентным соотношением (\ref{ttr-ell})
при начальных условиях $R_{-1}=0, R_0=1$. Все рекуррентные коэффициенты в
(\ref{ttr-ell}) зависят линейно от переменной $\gamma(z)$, которая поглощает всю
$z$-зависимость. Поэтому, $R_n(z;q,p)$ есть рациональные функции $\gamma(z)$,
где $n$ означает степень полиномов по $\gamma(z)$ в
числителе и знаменателе $R_n$. Более того, полюсы этих функций находятся в
точках $\gamma(z)=\alpha_1,\ldots,\alpha_n$.

При специальном выборе одного из параметров и дискретизации значений
$z$, функции $R_n(z;q,p)$ определяют эллиптическое расширение Вильсоновских
конечномерных $_9F_8$ и $_{10}\varphi_9$ рациональных функций
\cite{wil:hypergeometric,wil:orthogonal}.
Они были получены в работе \cite{spi-zhe:spectral} с помощью автомодельных
редукций нелинейных интегрируемых цепочек с дискретным временем (см. главу I).
Дискретные аналоги уравнений (\ref{eig}), (\ref{gevp}) для этого семейства
дискретных биортогональных рациональных функций были получены в \cite{spi-zhe:gevp}.
Конечно-разностное уравнение для $_{10}\varphi_9$ функции, появляющейся
из $R_n(z;q,p)$ в пределе $p\to 0$, изучалось Рахманом и Сусловым
\cite{rah-sus:classical}. Общие трехчленные рекуррентные соотношения типа
(\ref{ttr-ell}) рассматривались в \cite{zhe:bio} и в другой форме, связанной
с $R_{II}$-цепными дробями, в статье \cite{ism-mas:general}. Общие ортогональные
рациональные функции изучались в монографии \cite{Bult}.

Известно, что решения обобщенной задачи на собственные значения
образуют биортогональную систему функций \cite{Wilk}. Покажем, что эллиптический
бета-интеграл (\ref{ell-int}) служит мерой биортогональности для решений
уравнения (\ref{gevp}). Рассмотрим скалярное произведение
\begin{equation}
\int_C\Delta_E(z;\mathbf{t}) \Psi(z)\left(
{\cal D}_\eta-\lambda {\cal D}_\xi\right)\Phi(z)\frac{dz}{z},
\label{scalar}\end{equation}
где $\Delta_E(z;\mathbf{t})$ обозначает ядро интеграла (\ref{ell-int}),
а $\Phi(z), \Psi(z)$ некоторые комплексные функции. Выражение (\ref{scalar})
может быть переписано в виде
\begin{eqnarray}\nonumber
&& \int_{C}\Delta_E(z;\mathbf{t})(\kappa_\eta-V_\eta(z)-V_\eta(z^{-1}))
\Psi(z)\Phi(z)\frac{dz}{z}
\\ \nonumber && \makebox[2em]{}
+\int_{C_-}\Delta_E(q^{-1}z;\mathbf{t}) V_\eta(q^{-1}z)
\Psi(q^{-1}z)\Phi(z)\frac{dz}{z}
\\ && \makebox[2em]{}
+ \int_{C_+}\Delta_E(qz;\mathbf{t}) V_\eta(q^{-1}z^{-1})
\Psi(qz)\Phi(z)\frac{dz}{z} - \lambda \{\eta\to\xi\},
\label{scalar'}\end{eqnarray}
где $\{\eta\to\xi\}$ означает предыдущее выражение с переменной $\eta$
замененной на $\xi$. Контуры интегрирования $C_\pm$ получаются из
$C$ с помощью растяжений $z\to q^{\pm 1} z$. Предположим, что полюсы
$\Delta_E(z;\mathbf{t})$ и сингулярности функций $\Phi(z), \Psi(q^{\pm 1}z)$
не лежат в области, заметаемой контурами $C_\pm$ при их деформации к $C$.
Тогда (\ref{scalar'}) принимает вид
\begin{equation}
\int_C\Delta_E(z;\mathbf{t}) \Phi(z)\left(
{\cal D}_\eta^T-\lambda {\cal D}_\xi^T\right)\Psi(z)\frac{dz}{z},
\label{scalar-conj}\end{equation}
где сопряженный (или транспонированный) оператор ${\cal D}_\xi^T$ имеет вид
\begin{eqnarray} \nonumber
&& {\cal D}_\xi^T=\frac{\Delta_E(qz;\mathbf{t})}{\Delta_E(z;\mathbf{t})}
V_\xi(q^{-1}z^{-1})T+\frac{\Delta_E(q^{-1}z;\mathbf{t})}
{\Delta_E(z;\mathbf{t})}V_\xi(q^{-1}z)T^{-1}
\\ && \makebox[4em]{}
-V_\xi(z)-V_\xi(z^{-1})+\kappa_\xi.
\label{op-conj}\end{eqnarray}

Предположим, что $\Phi_\lambda(z)$ есть решение уравнения
$\left({\cal D}_\eta-\lambda {\cal D}_\xi\right)\Phi(z)=0$
и $\Psi_{\lambda'}(z)$ есть решение сопряженного уравнения
$\left({\cal D}_\eta^T-\lambda'{\cal D}_\xi^T\right)\Psi(z)=0$
для некоторых $\lambda$ и $\lambda'$. Обе эти функции могут быть умножены
на произвольные функции $f(z)$, удовлетворяющие условию периодичности
в логарифмическом масштабе, $f(qz)=f(z)$. После замены $\Phi(z)$ и $\Psi(z)$ в
(\ref{scalar}) в (\ref{scalar-conj}) на $\Phi_\lambda(z)$ и
$\Psi_{\lambda'}(z)$, эти выражения становятся равными нулю.
В частности, выражение (\ref{scalar-conj}) приводит к соотношению
\begin{eqnarray}\nonumber
&& \int_C\Delta_E(z;\mathbf{t}) \Phi_\lambda(z)\left(
{\cal D}_\eta^T-\lambda {\cal D}_\xi^T\right)\Psi_{\lambda'}(z)
\frac{dz}{z}
\\ && \makebox[4em]{}
=(\lambda'-\lambda)\int_C\Delta_E(z;\mathbf{t})
\Phi_\lambda(z){\cal D}_\xi^T\Psi_{\lambda'}(z)\frac{dz}{z}=0,
\label{biort-rel}\end{eqnarray}
показывающему, что при $\lambda'\neq\lambda$ функция $\Phi_\lambda(z)$
ортогональна ${\cal D}_\xi^T\Psi_{\lambda'}(z).$

Построим функцию $g(z)$, такую что
$$
g^{-1}(z)\left({\cal D}_\eta^T-\lambda {\cal D}_\xi^T\right)g(z)
={\cal D}_\eta-\lambda {\cal D}_\xi.
$$
После подстановки известных выражений для
$\Delta_E(qz;\mathbf{t})/\Delta_E(z;\mathbf{t})$ и спектрального параметра
(см. $\lambda$ в (\ref{lambda})), мы получаем уравнение для $g(z)$:
$$
g(qz)=\frac{\theta\left(\frac{q}{t_4z},\frac{q\mu z}{t_4},
Az,\frac{A\mu}{q^2z};p\right)}{\theta\left(\frac{q^2z}{t_4},
\frac{\mu}{t_4 z},\frac{A}{qz},\frac{A\mu z}{q};p\right)}\, g(z),
$$
которое легко решается (заметим, что мы подразумеваем $|q|<1$):
\begin{equation}
g(z)=\frac{\Gamma(\frac{q\mu z}{t_4},\frac{\mu q}{t_4 z},
Az,\frac{A}{z};q,p)}{\Gamma(\frac{q^2z}{t_4},\frac{q^2}{t_4 z},
\frac{A\mu z}{q}, \frac{A\mu}{qz};q,p)}.
\label{g-function}\end{equation}
Здесь мы пренебрегли уже упоминавшимся произвольным множителем $f(qz)=f(z)$.
В результате, мы находим прямую связь между
$\Phi_\lambda(z)$ и $\Psi_\lambda(z)$:
$\Psi_\lambda(z)=g(z)\Phi_\lambda(z)$, где $g(z)$ так же зависит от $\lambda$.

Обозначим $\lambda_n\equiv \lambda|_{\mu=q^n}$ и
$$
g_n(z)\equiv g(z)|_{\mu=q^n}=
\frac{\theta\left(\frac{q^2z}{t_4},\frac{q^2}{t_4z};p;q\right)_{n-1}}
{\theta\left(Az,\frac{A}{z};p;q\right)_{n-1}}.
$$
Подставляя $\Phi_{\lambda_n}(z)=R_n(z;q,p)$ в полученное соотношения
биортогональности, мы видим, что функции $R_n(z;q,p)$ формально
ортогональны ${\cal D}_\xi^T g_m(z)R_m(z;q,p)$ при $n\neq m$.

Общие рассуждения работ \cite{spi-zhe:theory,zhe:bio}, представленные в первой главе,
показывают что $R_n(z;q,p)$, будучи рациональными функциями $\gamma(z)$ с полюсами
при $\gamma(z)=\alpha_1,\ldots,\alpha_n$, ортогональны другим рациональным
функциям $\gamma(z)$, которые мы обозначим $T_m(z;q,p)$, с полюсами
при $\gamma(z)=\beta_1,\ldots,\beta_{n}$. Выбор
$\alpha_n,\beta_n$ и других рекуррентных коэффициентов в (\ref{ttr-ell})
определяют $R_n$ и $T_n$ единственным образом, так что перестановка всех $\alpha_n$
с $\beta_n$ переставляет $R_n$ и $T_n$. В нашем случае видно, что параметры
$\beta_n$ получаются из $\alpha_n=\gamma(q^n/t_4)$ после замены $t_4$ на $pq/A$.
Эквивалентно, эта замена превращает $\beta_n$ в $\alpha_n$.
Важным фактом является инвариантность весовой функции $\Delta_E(z,\mathbf{t})$
относительно этого преобразования. Поэтому, мы можем получить $T_n$ из
$R_n$ простой $t_4\to pq/A$  инволюцией, что дает
\begin{equation}
T_n(z;q,p)={_{12}V_{11}}\left(\frac{At_3}{q};\frac{A}{t_0},\frac{A}{t_1},
\frac{A}{t_2},t_3z,\frac{t_3}{z},q^{-n},\frac{Aq^{n-1}}{t_4};q,p\right),
\label{T_n}\end{equation}
где зависимость от $p$ в параметрах выпадает благодаря полной эллиптичности
данного $_{12}V_{11}$ ряда. Сравнивая с предыдущим рассмотрением, мы получаем
$$
{\cal D}_\xi^T g_n(z)R_n(z;q,p) = \rho_n T_n(z;q,p)
$$
для некоторой константы пропорциональности $\rho_n$, не существенной для наших
дальнейших рассмотрений.

Таким образом, операторный формализм приводит к следующему формальному соотношению
биортогональности:
\begin{equation}
\int_C T_n(z;q,p)R_m(z;q,p)\Delta_E(z,\mathbf{t})\frac{dz}{z}=\tilde h_n
\delta_{nm}
\label{formal-ort}\end{equation}
для некоторых констант $\tilde h_n$. Предположим, что $C=\mathbb{T}$
и $|t_r|<1, |qp|<|A|$. Часть полюсов функций $R_n,T_m$
сокращается с нулями  $\Delta_E(z;\mathbf{t})$. Остающиеся полюсы приближаются
к $\mathbb{T}$ по мере возрастания  целых $n,m$. Начиная с некоторых
значений $n$ и $m$, контур $\mathbb{T}$ уже не удовлетворяет условиям,
использованным в (\ref{scalar-conj}), и он должен быть деформирован.

На первый взгляд, соотношение (\ref{formal-ort}) должно оставаться верным и при
умножении $R_n(z;q,p)$ или $T_n(z;q,p)$ на произвольную функцию $f(z)$,
удовлетворяющую равенству $f(qz)=f(z)$. Однако такая нетривиальная $f(z)$
должна содержать сингулярности, которые пересекаются контурами $C_\pm$ при их
деформации к $C$ (иначе, $f(z)=const$). Поэтому влияние таких дополнительных
множителей должно быть изучено более тщательно. Более того, только при весьма
специальных $f(z)$ нормировочные константы $\tilde h_n$ могут быть точно
вычислены.

Весовая функция $\Delta_E(z,\mathbf{t})$ симметрична по $q$ и $p$,
в то время как ни $R_n(z;q,p)$ ни $T_n(z;q,p)$ не обладают таким свойством.
Можно попытаться восстановить эту симметрию, пользуясь свободой в множителе
$f(z)=f(qz)$. Положим $f(z)=R_k(z;p,q)$, $k\in\mathbb{N}$, то есть
воспользуемся самими функциями $R_n(z;q,p)$ с переставленными базисными
переменными $q$ и $p$. Тогда произведение
$$
R_{nk}(z)\equiv R_n(z;q,p)R_k(z;p,q)
$$
удовлетворяет двум обобщенным спектральным задачам: (\ref{eig}) и
$p$-разностному уравнению, получающемуся из него перестановкой
$q$ и $p$. Для (\ref{eig}) мы должны иметь $\mu=q^n$, а для его
партнера $\mu=p^k$. Функция (\ref{lambda}) не меняется
при замене $\mu\to p\mu$. Поэтому, выбор
$\mu=q^np^k,\, n,k\in \mathbb{N},$
дает ``спектр" для обеих обобщенных спектральных задач.
Первый множитель в $R_{nk}(z)$ есть рациональная функция
$\gamma(z;p)$ (мы указываем зависимость от $p$ явным образом),
но второй множитель является уже рациональной функцией $\gamma(z;q)$.
Поэтому, для $q$ и $p$ общего положения необходимо рассматривать
функции $R_{nk}(z)$ не как рациональные функции некоторого аргумента,
но как мероморфные функции $z$.

Точно таким же образом, условие обрыва ряда $t_6=q^{-n}$ в (\ref{R_n})
может быть заменено на $t_6=q^{-n}p^{-k}$, обрывающее одновременно
и $_{12}V_{11}$ ряд $R_k(z;p,q)$. Свойство полной эллиптичности
сбалансированного $_{r+1}V_{r}(t_0;$ $t_1,\ldots,$ $t_{r-4};q,p)$
ряда играет ключевую роль в этом месте: любой параметр
$t_1,\ldots,t_{r-5}$ может быть умножен на $p$ в произвольной
целой степени без изменений.

\section{Доказательство соотношения двухиндексной би\-ор\-то\-го\-наль\-нос\-ти}

Благодаря удвоению числа уравнений на собственные значения,
функции $R_{nk}(z)$ удовлетворяют довольно необычному условию
биортогональности (\ref{ort2}), характерному для функций {\em двух}
независимых переменных.
Для ясности изложения дадим явный вид $R_{nk}(z)$ и их сопряженных партнеров:
\begin{eqnarray} \nonumber
&&R_{nm}(z)={_{12}V_{11}}\left(\frac{t_3}{t_4};
\frac{q}{t_0t_4},\frac{q}{t_1t_4},\frac{q}{t_2t_4},t_3z,\frac{t_3}{z},q^{-n},
\frac{Aq^{n-1}}{t_4};q,p\right) \\  \label{R_nm}
&&\makebox[3em]{}\times {_{12}V_{11}}\left(\frac{t_3}{t_4};
\frac{p}{t_0t_4},\frac{p}{t_1t_4},\frac{p}{t_2t_4},t_3z,\frac{t_3}{z},p^{-m},
\frac{Ap^{m-1}}{t_4};p,q\right), \\ \nonumber
&&T_{nm}(z)={_{12}V_{11}}\left(\frac{At_3}{q};\frac{A}{t_0},\frac{A}{t_1},
\frac{A}{t_2},t_3z,\frac{t_3}{z},q^{-n},\frac{Aq^{n-1}}{t_4};q,p\right)
\\ \label{T_nm}
&& \makebox[3em]{}\times {_{12}V_{11}}\left(\frac{At_3}{p};\frac{A}{t_0},\frac{A}{t_1},
\frac{A}{t_2},t_3z,\frac{t_3}{z},p^{-m},\frac{Ap^{m-1}}{t_4};p,q\right),
\end{eqnarray}
где $n,m\in \mathbb{N}$. Очевидно, что эти функции симметричны при перестановке
$p$ и $q$ одновременно с $n$ и $m$.
Условие балансировки уже использовалось в этих рядах для того чтобы выразить
один из параметров через другие.

При $m\neq 0$ и фиксированных параметрах $q, t_r, r=0,\ldots,4,$
предел $p\to 0$ не определен ни для $R_{nm}(z)$ ни для $T_{nm}(z)$.
Причина этого кроется в квазипериодичности функции $\theta(z;p)$, которая
не позволяет определить $z\to 0$ и $z\to\infty$ пределы.

Теперь мы рассмотрим следующий интеграл:
\begin{equation}\label{ort-int}
J_{mn,kl}\equiv\int_{C_{mn,kl}}T_{nl}(z)R_{mk}(z)
\Delta_E(z,\mathbf{t}) \frac{d z}{z},
\end{equation}
где $\Delta_E(z,\mathbf{t})$ обозначает весовую функцию (\ref{weight})
и $C_{mn,kl}$ некоторый контур интегрирования. Мы хотим показать, что
специальный выбор $C_{mn,kl}$ (описанный ниже) приводит к формуле
\begin{equation}\label{ort}
J_{mn,kl}=h_{nl}\delta_{mn}\delta_{kl},
\end{equation}
где $h_{nl}$ некоторые нормировочные константы.

Эллиптический аналог преобразования Бэйли для обрывающегося $_{10}\varphi_9$
ряда, найденный Френкелем и Тураевым в работе \cite{fre-tur:elliptic},
имеет вид (альтернативное доказательство находится с помощью цепочки Бэйли -- оно
подробно описано в статье \cite{spi:bailey}):
\begin{eqnarray}\nonumber
{_{12}V_{11}}(t_0;t_1,\dots,t_7;q,p) &=&
\frac{\theta(qt_0,qs_0/t_4,qs_0/t_5,qt_0/t_4t_5;p;q)_N}
{\theta(qs_0,qt_0/t_4,qt_0/t_5,qs_0/t_4t_5;p;q)_N}   \\
&& \times {_{12}V_{11}}(s_0;s_1,\dots,s_7;q,p),
\label{ft-b}\end{eqnarray}
где $\prod_{m=1}^7t_m=t_0^3q^2$, $t_6=q^{-N}$, $N\in\mathbb{N}$,
\begin{equation}
s_0=\frac{qt_0^2}{t_1t_2t_3},\quad s_1=\frac{s_0t_1}{t_0},\quad
s_2=\frac{s_0t_2}{t_0},\quad s_3=\frac{s_0t_3}{t_0},
\label{s-param}\end{equation}
и $s_4, s_5, s_6, s_7$ равны произвольной перестановке
$t_4, t_5, t_6, t_7$. Используя это тождество,
можно переписать $R_{mk}(z)$ так
\begin{eqnarray} \nonumber
&& R_{mk}(z) = \kappa_m(p;q)\kappa_k(q;p) \\
&& \makebox[3em]{}  \times
{_{12}V_{11}}\left(\frac{t_0}{t_4};\frac{q}{t_1t_4},\frac{q}{t_2t_4},
\frac{q}{t_3t_4},t_0z,\frac{t_0}{z},q^{-m},\frac{Aq^{m-1}}{t_4};q,p\right)
\nonumber \\
&& \makebox[3em]{} \times
{_{12}V_{11}}\left(\frac{t_0}{t_4};\frac{p}{t_1t_4},\frac{p}{t_2t_4},
\frac{p}{t_3t_4},t_0z,\frac{t_0}{z},p^{-k},\frac{Ap^{k-1}}{t_4};p,q\right),
\label{R-transf}\end{eqnarray}
где
$$ \kappa_m(p;q)=\frac{\theta(qt_3/t_4,t_0t_1,t_0t_2,A/qt_0;p;q)_m}
{\theta(qt_0/t_4,t_1t_3,t_2t_3,A/qt_3;p;q)_m}. $$

Подстановка явного вида рядов  \re{T_nm} и
\re{R-transf} в определение $J_{mn,kl}$ (\ref{ort-int}) дает
\begin{eqnarray}\nonumber
J_{mn,kl}&=& \kappa_m(p;q)\kappa_{k}(q;p)
\sum_{r=0}^n\sum_{r'=0}^l\sum_{s=0}^{m}\sum_{s'=0}^k q^{r+s}p^{r'+s'}
\\ \nonumber
&& \times \frac{\theta(At_3q^{2r-1},t_0q^{2s}/t_4;p)}{\theta(At_3/q,t_0/t_4;p)}
\frac{\theta(At_3p^{2r'-1},t_0p^{2s'}/t_4;q)}{\theta(At_3/p,t_0/t_4;q)}
\\  \nonumber
&& \times \frac{\theta(At_3/q,A/t_0,A/t_1,A/t_2,q^{-n},Aq^{n-1}/t_4;p;q)_r}
{\theta(q,t_0t_3,t_1t_3,t_2t_3,At_3q^n,t_3t_4q^{1-n};p;q)_r}
\\ \nonumber
&& \times \frac{\theta(At_3/p,A/t_0,A/t_1,A/t_2,p^{-l},Ap^{l-1}/t_4;q;p)_{r'}}
{\theta(p,t_0t_3,t_1t_3,t_2t_3,At_3p^l,t_3t_4p^{1-l};q;p)_{r'}}
\\ \nonumber
&& \times\frac{\theta(t_0/t_4,q/t_1t_4,q/t_2t_4,q/t_3t_4,q^{-m},
Aq^{m-1}/t_4;p;q)_s}
{\theta(q,t_0t_1,t_0t_2,t_0t_3,t_0q^{m+1}/t_4,t_0q^{2-m}/A;p;q)_s}
\\ \nonumber
&& \times\frac{\theta(t_0/t_4,p/t_1t_4,p/t_2t_4,p/t_3t_4,p^{-k},
Ap^{k-1}/t_4;q;p)_{s'}}
{\theta(p,t_0t_1,t_0t_2,t_0t_3,t_0p^{k+1}/t_4,t_0p^{2-k}/A;q;p)_{s'}}\,
I_{rs,r's'},
\end{eqnarray}
где
\begin{eqnarray}\nonumber
&I_{rs,r's'}& = \int_{C_{mn,kl}}\Delta_E(z,\mathbf{t})
\frac{\theta(zt_3,z^{-1}t_3;p;q)_r}{\theta(zA,z^{-1}A;p;q)_r}
\frac{\theta(zt_0,z^{-1}t_0;p;q)_s}
{\theta(zqt_4^{-1},z^{-1}qt_4^{-1};p;q)_s} \\
&& \times
\frac{\theta(zt_3,z^{-1}t_3;q;p)_{r'}}{\theta(zA,z^{-1}A;q;p)_{r'}}
\frac{\theta(zt_0,z^{-1}t_0;q;p)_{s'}}
{\theta(zpt_4^{-1},z^{-1}pt_4^{-1};q;p)_{s'}} \frac{d z}{z}.
\label{int3}\end{eqnarray}
Введем обозначения
$$
\tilde t_0=t_0q^sp^{s'},\quad \tilde t_1=t_1,\quad  \tilde t_2=t_2,
\quad \tilde t_3=t_3q^rp^{r'}, \quad \tilde t_4=t_4q^{-s}p^{-s'},
$$
так что
$$
\tilde A=\prod_{m=0}^4\tilde t_m=Aq^rp^{r'}.
$$
Тогда, пользуясь преобразованиями
$$
\theta(z;p;q)_l=(-z)^lq^{l(l-1)/2}\theta(z^{-1}q^{-l+1};p;q)_l
=\frac{(-z)^lq^{l(l-1)/2}}{\theta(qz^{-1};p;q)_{-l}},
$$
мы можем переписать интеграл (\ref{int3}) в виде
\begin{equation}\label{int4}
I_{rs,r's'}=\left(\frac{t_0t_4}{pq}\right)^{2ss'}
\left(\frac{t_3}{A}\right)^{2rr'}\frac{t_4^{2(s+s')}}{q^{s(s+1)}p^{s'(s'+1)}}
\int_{C_{mn,kl}}
\Delta_E(z,\mathbf{\tilde t}) \frac{d z}{z}.
\end{equation}
Но интеграл в правой части (\ref{int4}) совпадает с эллиптическим бета-интегралом
(\ref{ell-int}), при условии, что мы выберем $C_{mn,kl}=\mathbb{T}$ и наложим
ограничения $|\tilde t_m|<1$, $|pq|<|\tilde A|$. Однако значения целых
$r,s,r',s'\in\mathbb{N}$ не ограничены. Начиная с некоторого достаточно большого
значения, мы получим либо $|\tilde A|=|q^rp^{r'}A|<|pq|$ либо
$|\tilde t_4|=|q^{-s}p^{-s'}t_4|>1$. Именно эти соображения и определяют критерии
выбора контура $C_{mn,kl}$: мы выберем его из условия, чтобы формула
(\ref{ell-int}) оставалась применимой при любых $r,s,r',s'$.
Точнее, пусть контур $C_{mn,kl}$ будет деформацией $\mathbb{T}$, такой что
он разделяет полюсы при $z=t_{0,1,2,3}p^aq^{b},$ $t_4p^{a-k}q^{b-m},$
и $A^{-1}p^{a+1-l}q^{b+1-n},$ $a,b\in\mathbb{N}$, лежащие внутри
$C_{mn,kl}$ и сходящиеся к нулю, от их $z\to z^{-1}$ партнеров,
уходящих к бесконечности. Индексы $m,n,k,l$ в обозначении контура
показывает, что форма $C_{nm,kl}$ зависит от индексов функций $T_{nl}(z)$ и $R_{mk}(z)$.

Для такого контура $C_{mn,kl}$, мы получаем
\begin{eqnarray}\nonumber
I_{rs,r's'} &=&
\left(\frac{t_0t_4}{pq}\right)^{2ss'}
\left(\frac{t_3}{A}\right)^{2rr'}\frac{t_4^{2(s+s')}}{q^{s(s+1)}p^{s'(s'+1)}}
{\cal N}_E(\mathbf{\tilde t}) \\  \nonumber
&=&{\cal N}_E(\mathbf{t})
\frac{\theta(t_1t_3,t_2t_3,t_3t_4;p;q)_r}
{\theta(A/t_0,A/t_1,A/t_2;p;q)_r}  \\  \nonumber
&&\times \frac{\theta(t_0t_3;p;q)_{r+s}}{\theta(A/t_4;p;q)_{r+s}}
\frac{\theta(t_0t_1,t_0t_2,q^{1-r}t_0/A;p;q)_s}
{\theta(q/t_1t_4,q/t_2t_4,q^{1-r}/t_3t_4;p;q)_s}
\\  \nonumber
&& \times
\frac{\theta(t_1t_3,t_2t_3,t_3t_4;q;p)_{r'}}
{\theta(A/t_0,A/t_1,A/t_2;q;p)_{r'}}  \\  \nonumber
&&\times \frac{\theta(t_0t_3;q;p)_{r'+s'}}{\theta(A/t_4;q;p)_{r'+s'}}
\frac{\theta(t_0t_1,t_0t_2,p^{1-r'}t_0/A;q;p)_{s'}}
{\theta(p/t_1t_4,p/t_2t_4,p^{1-r'}/t_3t_4;q;p)_{s'}},
\end{eqnarray}
где
\begin{equation}\label{result}
{\cal N}_E(\mathbf{t})=\frac{2\prod_{0\leq m<s\leq 4} \Gamma(t_mt_s;q,p)}
{(q;q)_\infty(p;p)_\infty\prod_{m=0}^4\Gamma(At_m^{-1};q,p)}.
\end{equation}

В результате этих манипуляций величина  $J_{mn,kl}$
расщепляется на произведение двух двукратных рядов, каждый из которых
зависит только от индексов $m, n$ и $k,l$ по отдельности. После применения
соотношения $(a;p;q)_{r+s}=(aq^r;p;q)_s(a;p;q)_r$ и различных упрощений,
можно записать
$$
J_{mn,kl}={\cal N}_E(\mathbf{t})J_{mn}(p;q)J_{kl}(q;p),
$$
где
\begin{eqnarray}\nonumber
J_{mn}(p;q) &=& \kappa_m(p;q) \sum_{r=0}^nq^r
\frac{\theta(At_3q^{2r-1};p)}{\theta(At_3/q;p)}
\frac{\theta(At_3/q,q^{-n},Aq^{n-1}/t_4,t_3t_4;p;q)_r}
{\theta(q,At_3q^n,t_3t_4q^{1-n},A/t_4;p;q)_r}
\\  \label{8E7-aux}
&& \times {_{10}V_9}\left(\frac{t_0}{t_4};\frac{q}{t_3t_4},
t_0t_3q^r,\frac{q^{1-r}t_0}{A},\frac{Aq^{m-1}}{t_4},q^{-m};q,p\right).
\end{eqnarray}

Ограничение $t_2t_3=qt_0$, наложенное на соотношение (\ref{ft-b}),
превращает $_{12}V_{11}$-ряд в левой части в обрывающийся $_{10}V_9$
ряд, в то время как в ряде в правой части остается только один член.
В результате получаем эллиптический аналог формулы суммирования Джексона:
\be
{_{10}V_9}(t_0;t_1,t_4, t_5,t_6,t_7;q,p)
= \frac{\theta (qt_0,qt_0/t_1t_4,qt_0/t_1t_5,qt_0/t_4t_5;p;q)_N}
{\theta(qt_0/t_1t_4t_5,qt_0/t_1,qt_0/t_4,qt_0/t_5;p;q)_N},
\label{ft-s}\ee
где $t_1t_4t_5t_6t_7 =qt_0^2$ и $t_6=q^{-N}$, $N\in \mathbb{N}$.
Применение этой формулы к $_{10}V_9$ ряду в \re{8E7-aux} дает
$$
{_{10}V_9}(\ldots)=\frac{\theta(qt_0/t_4,t_1t_2,At_3q^{r-1},q^{-r};p;q)_m}
{\theta(t_0t_3,A/qt_0,Aq^r/t_4,q^{1-r}/t_3t_4;p;q)_m}.
$$
Очевидно, что это выражение зануляется при $m>r$. Это означает, что $J_{mn}=0$
при $m>n$. При $m\leq n$, мы получаем
\begin{eqnarray}\nonumber
J_{mn}(p;q)=\kappa_m(p;q)
\frac{\theta(At_3;p;q)_{2m}}{\theta(A/t_4;p;q)_{2m}}
\frac{\theta(Aq^{n-1}/t_4,qt_0/t_4,t_1t_2,q^{-n};p;q)_m}
{\theta(t_3t_4q^{1-n},t_0t_3,A/qt_0,At_3q^n;p;q)_m}
\\  \nonumber
\times (t_3t_4)^m\, {_8V_7}(At_3q^{2m-1};t_3t_4,Aq^{n+m-1}/t_4,q^{m-n};q,p).
\end{eqnarray}
Применяя формулу (\ref{ft-s}) к этому $_8V_7$ ряду, получаем
$$
{_8V_7}(\ldots)=\frac{\theta(At_3q^{2m},q^{m-n+1},Aq^{m+n}/t_4,qt_3t_4;p;q)_{n-m}}
{\theta(q,Aq^{2m}/t_4,t_3t_4q^{m-n+1},At_3q^{m+n};p;q)_{n-m}},
$$
что равно нулю при $n>m$ благодаря множителю
$\theta(q^{m-n+1};p;q)_{n-m}$. В результате получаем, что
$J_{mn}(p;q)$ $=h_n(p;q)\delta_{mn}$, где нормировочные константы имеют вид
\begin{equation}\label{norm}
h_n(p;q)=\frac{\theta(A/qt_4;p)
\theta(q,qt_3/t_4,t_0t_1,t_0t_2,t_1t_2,At_3;p;q)_nq^{-n}}
{\theta(Aq^{2n-1}/t_4;p)
\theta(1/t_3t_4,t_0t_3,t_1t_3,t_2t_3,A/qt_3,A/qt_4;p;q)_n}.
\end{equation}

Тот факт, что $J_{mn}=0$ при $n\neq m$ приводит к требуемой биортогональности
(\ref{ort}). Суммируем полученные результаты в форме теоремы, анонсированной
первый раз в заметке \cite{spi:special} (для совпадения формул необходимо в
\cite{spi:special} перейти к обозначениям для рядов предложенным в
работах \cite{spi:theta,spi:bailey} и переставить параметры $t_3$ и $t_4$).

\begin{theorem}
Пусть $C_{mn,kl}$ обозначает контур, ориентированной против часовой
стрел\-ки и отделяющий точки $z=\{ t_{0,1,2,3}p^aq^b,$ $t_4p^{a-k}q^{b-m},$
$A^{-1}p^{a+1-l}q^{b+1-n}\}_{a,b\in\mathbb{N}}$
от точек с отраженными ($z\to z^{-1}$) координатами.
Тогда $R_{mk}(z)$ и $T_{nl}(z)$ удовлетворяют следующему соотношению
биортогональности
\begin{equation}\label{ort2}
\int_{C_{mn,kl}}T_{nl}(z)R_{mk}(z)
\Delta_E(z,{\bf t}) \frac{d z}{z}=
h_{nl}{\cal N}_E({\bf t})\delta_{mn}\delta_{kl},
\end{equation}
где нормировочные константы $h_{nl}$ имеют следующий явный вид
\begin{eqnarray} \nonumber
&& h_{nl}= \frac{\theta(A/qt_4;p)
\theta(q,qt_3/t_4,t_0t_1,t_0t_2,t_1t_2,At_3;p;q)_nq^{-n}}
{\theta(Aq^{2n}/qt_4;p)
\theta(1/t_3t_4,t_0t_3,t_1t_3,t_2t_3,A/qt_3,A/qt_4;p;q)_n} \\
&&\makebox[2em]{}\times  \frac{\theta(A/pt_4;q)
\theta(p,pt_3/t_4,t_0t_1,t_0t_2,t_1t_2,At_3;q;p)_lp^{-l}}
{\theta(Ap^{2l}/pt_4;q)
\theta(1/t_3t_4,t_0t_3,t_1t_3,t_2t_3,A/pt_3,A/pt_4;q;p)_l}.
\label{norm2}\end{eqnarray}
\end{theorem}

Как ясно из определений (\ref{R_n}) и (\ref{T_n}), мы имеем
$R_m(z;q,p)$ $= R_{m0}(z)$ и $T_n(z;q,p)= T_{n0}(z)$.

\begin{corollary}
Функции $R_m(z;q,p)$ и $T_n(z;q,p)$ удовлетворяют условию биортогональности
\begin{equation}\label{R_mT_n-ort}
\int_{C_{mn}}T_{n}(z;q,p)R_{m}(z;q,p)
\Delta_E(z,{\bf t}) \frac{d z}{z}= h_{n}{\cal N}_E({\bf t})\delta_{mn},
\end{equation}
где константы $h_n$ фиксированы в (\ref{norm}),
а контур $C_{mn}$ окружает полюсы в точках $z=\{ t_{0,1,2,3}q^ap^b,
t_4p^aq^{b-m},A^{-1}p^{a+1}q^{b+1-n}\}_{a,b\in\mathbb{N}}$ и отделяет их
от полюсов с отраженными $z\to z^{-1}$ координатами.
\end{corollary}

Биортогональные рациональные функции $R_m(z;q,p)$ и $T_n(z;q,p)$
описывают эллиптическое расширение множества непрерывных $_{10}\varphi_9$
функций работ \cite{rah:integral,rah:biorthogonality}, к которым они сводятся
в пределе $p\to 0$. Соответственно, формула (\ref{R_mT_n-ort}) в этом
пределе сводится к условию биортогональности Рахмана. Последующее
вырождение приводит к биортогональным $_9F_8$ функциям, обобщающим
полиномы Вильсона \cite{wil:hypergeometric} (см. также
недавние связанные исследования \cite{ner2}).

Биортогональные функции, описанные выше, привлекли к себе
существенное внимание экспертов по специальным функциям.
Ряд исследований, последовавших за работами
\cite{spi:special,spi:integrals,spi-zhe:spectral},
был посвящен теоретико-групповой интерпретации этих функций.
Возможные алгебраические структуры, скрытые за самодуальными
обобщенными задачами на собственные значения для матриц Якоби,
обсуждались на уровне классических Пуассоновских алгебр в
\cite{spi-zhe:poisson}. При этом понадобились обобщения классического
аналога известной алгебры Склянина \cite{skl:alg}. Анализ
следствий этого подхода на квантовом уровне еще не рассматривался.
Другой подход основан на описании $R_m(z;q,p)$ как функций перекрытий
базисов, появляющихся в некоторых спектральных задачах.
Хотя у авторов \cite{spi-zhe:poisson} и были предположения,
что для этого необходимо рассматривать решения
обобщенных задач на собственные значения для генераторов алгебры Склянина,
ключевое наблюдение о справедливости такого подхода было сделано
Рэйнсом \cite{rai:abelian}. Он нашел правильную комбинацию Склянинских
операторов, нулевые моды которых порождают необходимые базисные
(интерполирующие) функции, а Розенгрен в
\cite{ros:elementary,ros:sklyanin} детально описал как их перекрытия
порождают конечномерные биортогональные функции \cite{spi-zhe:spectral},
для которых формула суммирования Френкеля-Тураева определяет нормировку дискретной
весовой функции.

В работах ван Диехена и автора \cite{die-spi:elliptic,die-spi:selberg} и в
независимых исследованиях автора \cite{spi:special,spi:theta,spi:integrals}
была сформулирована программа по построению специальных функций многих
переменных, связанных с многомерными эллиптическими бета-интегралами,
которые должны описывать меру в соотношениях биортогональности для них.
В частности, автором проводились поиски таких многомерных обобщений
в связи с моделями типа Калоджеро (см. об этом ниже). Первый пример
многомерного обобщения биортогональных функций, построенных в работах
автора \cite{spi:special,spi:integrals} и описанных в этом параграфе,
на корневую систему $BC_n$ был найден Рэйнсом в \cite{rai:abelian}.
При этом важную роль сыграла работа Окунькова по интерполирующим
полиномам Макдональда \cite{oku} (ранее Окуньковым и Ольшанским были
построены аналогичные обобщения полиномов Джека \cite{oo}).
Другой тип многомерных функций гипергеометрического типа над
эллиптическими кривыми построен в работе \cite{fel-var:mac}.

\section{Интегральное представление для симметричного
произведения двух $_{12}V_{11}$ рядов}

Построим интегральное представление для произведения двух $_{12}V_{11}$
рядов со специфическим выбором параметров. Для этого мы используем эллиптическое
обобщение техники, предложенной в статье \cite{rah:integral}, для вывода
интегрального представления для обрывающегося совершенно
уравновешенного $_{10}\varphi_9$ ряда.

\begin{theorem}
Возьмем пару целых чисел $m,n\in\N$ и обозначим $C_{mn}$ контур, ориентированный
против часовой стрелки и разделяющий точки  $z=\{t_kp^aq^b,$ $A^{-1}q^{b+1-m}p^{a+1-n}\}$
для любых $a,b \in\mathbb{N}$ от их $z\to z^{-1}$ партнеров. При этих условиях
справедливо следующее интегральное представление для произведения двух
обрывающихся $_{12}V_{11}$ рядов с аргументом $x=1$:
\begin{eqnarray}\nonumber
&&{_{12}V_{11}}\left(\frac{At_0}{q};\alpha,t_0t_1,t_0t_2,t_0t_3,t_0t_4,
q^{-m},\frac{A^2q^{m-1}}{\alpha};q,p\right) \\  \nonumber
&& \times {_{12}V_{11}}
\left(\frac{At_0}{p};\beta,t_0t_1,t_0t_2,t_0t_3,t_0t_4,
p^{-n},\frac{A^2p^{n-1}}{\beta};p,q\right) \\  \nonumber
&& \makebox[2em]{}
=\frac{1}{{\cal N}_E({\bf t})}
\frac{\theta(At_0,\frac{A}{t_0};p;q)_m\theta(At_0,\frac{A}{t_0};q;p)_n}
{\theta(\frac{A}{\alpha t_0},\frac{At_0}{\alpha};p;q)_m
\theta(\frac{A}{\beta t_0},\frac{At_0}{\beta};q;p)_n}  \\
&& \makebox[3em]{}
\times \int_{C_{mn}}
\Delta_E(z,{\bf t})\frac{\theta(\frac{Az}{\alpha},\frac{A}{\alpha z};p;q)_m
\theta(\frac{Az}{\beta},\frac{A}{\beta z};q;p)_n }
{\theta(Az,\frac{A}{z};p;q)_m\theta(Az,\frac{A}{z};q;p)_n}\frac{d z}{z},
\label{intrep}\end{eqnarray}
где $\alpha$ и $\beta$ произвольные комплексные параметры.
\end{theorem}
{\bf Доказательство.}
При указанных условиях справедливы соотношения
\begin{eqnarray}\nonumber
\lefteqn{
\int_{C_{ij}}\Delta_E(z,\mathbf{t})
\frac{\theta(zt_0,z^{-1}t_0;p;q)_i\theta(zt_0,z^{-1}t_0;q;p)_j}
{\theta(zA,z^{-1}A;p;q)_i\theta(zA,z^{-1}A;q;p)_j}\frac{dz}{z}
} && \\ &&
= \left(\frac{t_0}{A}\right)^{2ij}
{\cal N}_E(t_0q^ip^j,t_1,\ldots,t_4)
\nonumber \\  &&
=\frac{\theta(t_0t_1,\ldots,t_0t_4;p;q)_i\theta(t_0t_1,\ldots,t_0t_4;q;p)_j}
{\theta(A/t_1,\ldots,A/t_4;p;q)_i\theta(A/t_1,\ldots,A/t_4;q;p)_j}
\, {\cal N}_E(\mathbf{t}).
\label{identity}\end{eqnarray}
Умножая (\ref{identity}) на выражение
\begin{eqnarray*}
  && q^i\frac{\theta(At_0q^{2i-1};p)}{\theta(At_0/q;p)}
\frac{\theta(At_0/q,\alpha,q^{-m},A^2q^{m-1}/\alpha;p;q)_i}
     {\theta(q,At_0/\alpha,At_0q^m,\alpha q^{1-m}t_0/A;p;q)_i}
\\
 && \times p^j\frac{\theta(At_0p^{2j-1};q)}{\theta(At_0/p;q)}
\frac{\theta(At_0/p,\beta,p^{-n},A^2p^{n-1}/\beta;q;p)_j}
     {\theta(p,At_0/\beta,At_0p^n,\beta p^{1-n}t_0/A;q;p)_j},
\end{eqnarray*}
где $\alpha$ и $\beta$ произвольные комплексные параметры, и суммируя по $i$ от
0 до $m$ и по $j$ от 0 до $n$, получаем равенство
\begin{eqnarray*}
&&\int_{C_{mn}} \Delta_E(z,\mathbf{t})\,
{_{10}V_9}(At_0/q;t_0z,t_0z^{-1},\alpha,q^{-m},A^2q^{m-1}/\alpha;q,p)
\\ && \makebox[4em]{}
\times {_{10}V_9}(At_0/p;t_0z,t_0z^{-1},\beta,p^{-n},A^2p^{n-1}/\alpha;p,q)
\,\frac{dz}{z}  \\
&&\makebox[2em]{} ={_{12}V_{11}}(At_0/q;\alpha,q^{-m},A^2q^{m-1}/\alpha,
t_0t_1,\ldots,t_0t_4;q,p) \\
&& \makebox[3em]{} \times
{_{12}V_{11}}(At_0/p;\beta,p^{-n},A^2p^{n-1}/\beta,t_0t_1,\ldots,t_0t_4;
p,q)\, {\cal N}_E(\mathbf{t}).
\end{eqnarray*}
Применение формулы суммирования Френкеля-Тураева к $_{10}V_9$ рядам, стоящим под знаком интеграла,
приводит к (\ref{intrep}).
\hfill{Q.E.D.}

\smallskip

При $n=0$, мы получаем интегральное представление для одного обрывающегося
$_{12}V_{11}$ ряда, который может быть редуцирован к $_{10}\varphi_9$
$q$-гипергеометрическому ряду работы \cite{rah:integral} в пределе $p\to 0$.
Однако, при $n\neq 0$ предел $p\to 0$ хорошо не определен и формула
(\ref{intrep}) существует только на эллиптическом уровне.

\section{Обрывающаяся цепная дробь}

В данном параграфе мы опишем самую общую обрывающуюся цепную дробь
гипергеометрического типа известную в настоящее время. Она была получена
в работах \cite{spi-zhe:gevp,spi-zhe:theory}. Эта дробь является эллиптическим
обобщением $_{10}\varphi_9$ семейства обрывающихся цепных дробей построенных
Массоном \cite{Masson} и Гуптой и Массоном \cite{gup-mas:watson,gup-mas:contiguous}.
Эти дроби включают в себя $q$-гипергеометрическую дробь Ватсона \cite{Watson},
являющуюся $q$-аналогом знаменитой цепной дроби Рамануджана (номер 40)
связанной со специальным случаем совершенно уравновешенного 2-сбалансированного
${_9}F_8$ гипергеометрического ряда \cite{Berndt,Ramanu}.

Предположим, что мы имеем трехчленное рекуррентное соотношение
\begin{equation}
\psi_{n+1}= \xi_n \psi_n + \eta_n \psi_{n-1}, \quad n\in\mathbb{N},
\label{ttr-cf}\end{equation}
для некоторых несингулярных коэффициентов $\xi_n$  и $\eta_n$.
Обозначим $U_n$ и $V_n$ две последовательности,  удовлетворяющие \re{ttr-cf}
с начальными условиями $U_0=0, \; U_1 = 1$ и $V_0 = 1, \;
V_1 = \xi_0.$ Известно, что отношение $U_n/V_n$ равно цепной дроби \cite{JT}
\begin{equation}
\frac{U_n}{V_n}={1\over\displaystyle \xi_0 + {\strut \eta_1
\over\displaystyle \xi_1 +{\strut \eta_2 \over {\xi_2 + \dots
{\atop \displaystyle +\frac{\eta_{n-1}}{\xi_{n-1}} }}}}},
\qquad n=1,2,\dots.
\label{cf-chain} \end{equation}
В случае ортогональных полиномов (определяемых последовательностью $V_n$),
$\xi_n$ линейны по аргументу полиномов $z$, а $\eta_n$ не зависят от $z$.
Если $\eta_n(z)$ квадратичны по $z$, а $\xi_n(z)$ линейны по $z$, то мы получаем
цепные дроби, названные в \cite{ism-mas:general} дробями $R_{II}$-типа
(они еще известны как оскуляторные цепные дроби \cite{Wu}).
Как уже упоминалось, такие рекуррентные соотношения порождают биортогональные
рациональные функции. Положим
\begin{equation}
\xi_n(z) = \rho_n z-v_n, \qquad \eta_n(z) = - u_n(z-\alpha_n)(z-\beta_{n-1}),
\label{rec_z} \end{equation}
где $\rho_n, v_n, u_n, \alpha_n, \beta_n$ некоторые последовательности
чисел. В результате $V_n(z)$ $=P_n(z)$ становятся полиномами $n$-й степени
с начальными условиями $P_0=1, \; P_1(z) = \rho_0 z-v_0$.
Полиномы $U_n(z) = P^{(1)}_{n-1}(z)$, называемые ассоциированными полиномами,
имеют степень $n-1$. Они удовлетворяют рекуррентным соотношениям
\begin{equation}
P^{(1)}_{n}(z) + (v_n -\rho_n z) P^{(1)}_{n-1}(z) +
u_n(z-\alpha_n)(z-\beta_{n-1}) P^{(1)}_{n-2}(z) =0
\label{rec_R} \end{equation}
с начальными условиями $P^{(1)}_{-1}=0, \; P^{(1)}_0(z)=1$.

Подставляя рекуррентные коэффициенты \re{rec_z} в \re{cf-chain}, получаем
\begin{equation}
F_N(z)\equiv
\frac{P^{(1)}_{N}(z)}{P_{N+1}(z)}={1\over\displaystyle \rho_0 z-v_0 -
{\strut u_1(z-\alpha_1)(z-\beta_0) \over\displaystyle \rho_1 z-v_1 -
\dots {\atop \strut \displaystyle -\frac{u_N(z-\alpha_N)(z-\beta_{N-1})}
{\rho_N z-v_N} }}}.
\label{chain_P} \end{equation}
Поскольку $F_N(z)$ есть рациональная функция $z$, мы можем разложить ее
в частные дроби:
\begin{equation}
F_N(z)=\sum_{s=0}^N \frac{g_s}{z-z_s},
\label{par_frac} \end{equation}
где $z_s, \: s=0,1,\dots,N,$ обозначает нули полинома $P_{N+1}(z)$ и
\begin{equation}
g_s = \frac{P^{(1)}_N(z_s)}{P_{N+1}'(z_s)}.
\label{g_s} \end{equation}
Мы подразумеваем, что $P_{N+1}(z)$ имеет только простые нули, то есть
$z_s \ne z_s'$ при $s \ne s'$.

Любые два решения $U_n, V_n$ рекуррентного соотношения \re{ttr-cf}
удовлетворяют Вронскианному соотношению
$$
U_{n+1}V_n-U_nV_{n+1}=(-1)^n\eta_1\cdots\eta_n (U_1V_0-U_0V_1),
$$
которое в нашем случае дает
\begin{equation}
P_n(z) P^{(1)}_n(z) - P_{n+1}(z)P^{(1)}_{n-1}(z) = h_n A_n(z) \tilde B_n(z),
\label{Wron} \end{equation}
где $h_n = u_1 u_2 \cdots u_n$ и
$$
A_n(z)=\prod_{i=1}^n(z-\alpha_i), \qquad
\tilde B_n(z)=\prod_{i=1}^n(z-\beta_{i-1}).
$$
Взяв $n=N$ и $z=z_s, \; s=0,1,\dots,N,$ в \re{Wron}, можно выразить
$P^{(1)}_N(z_s)$ через $P_N(z_s),$ $h_n,$  $A_n(z_s)$ и
$\tilde B_N(z_s)$. Это приводит к следующему удобному для вычислений
выражению для $g_s$:
\begin{equation}
g_s = \frac{h_N A_N(z_s)\tilde B_N(z_s)}{P_{N+1}'(z_s) P_N(z_s)}.
\label{ex_g} \end{equation}

Рассмотрим теперь семейство рациональных функций $R_n(z(u))$,
определенных в первой главе и предыдущих параграфах этой главы,
\begin{eqnarray}
&&R_n(z(u))= {_{12}}v_{11}(1-x_1;1+x_0-d_3,1+x_0-d_4,1+x_0-d_5, \nonumber \\
&& \makebox[5em]{}
1+u+x_0-e_1,1-u+x_0-e_2,1-x_2+n,-n),
\label{R_e} \end{eqnarray}
выражающихся через эллиптические гипергеометрические ряды
(в мультипликативной системе обозначений).
Вычислим $g_s$ для соответствующих полиномов
\begin{eqnarray}
&& P_n(z)=\kappa_n A_n(z)R_n(z),
\label{pols}\\
&& \kappa_n=G_{n-1}\cdots G_1G_0=\frac{[1-x_2,-x_0]_n
\prod_{i=1}^5[1-x_2+d_i]_n}{[1-x_2]_{2n}[2-x_1]_n}.
\nonumber\end{eqnarray}
с $d_3=x_2-N-1$. В этом случае $u_{N+1}=0$ и цепная дробь \re{chain_P}
обрывается автоматически. Прежде всего отметим, что
\begin{equation}
P'_{N+1}(z_s)=\gamma_{N+1}(z_s-z_0)\cdots(z_s-z_{s-1})(z_s-z_{s+1})
\cdots(z_s-z_{N}).
\label{P'} \end{equation}
Это выражение может быть вычислено с помощью соотношения
\begin{equation}
z_s-z_k=\mu_s\frac{[k-s,2+k+s+2x_0-d_1-d_2]}
{[k+x_0-d_1+1,k+x_0-d_2+1]},
\label{zz} \end{equation}
$$
\mu_s=\frac{[d_2-e_1,e_1-d_1]}{[s+1+x_0-d_1,s+1+x_0-d_2]}.
$$
В результате,
\begin{eqnarray*}
&& P_{N+1}'(z_s)=\gamma_{N+1}\mu_s^N
\frac{[s+1+x_0-d_1,s+1+x_0-d_2]}{[2s+2+2x_0-d_1-d_2]}
\\ && \makebox[4em]{}
\times\frac{[s+2+2x_0-d_1-d_2]_{N+1}[1]_N[1]_s}
{[1+x_0-d_1,1+x_0-d_2]_{N+1}[-N]_s}.
\end{eqnarray*}
Полином $\tilde B_n(z_s)$ легко строится
$$
\tilde B_n(z_s)
=\mu_s^n\frac{[s+2+x_0-d_1-d_2,-s-x_0]_n}{[1-d_1,1-d_2]_n}.
$$
Подставляя $d_3=x_2-N-1$ и
\begin{equation}
z= z_s\equiv z(u_s),\quad u_s= s+x_0+1-e_2, \quad s=0,1,\dots, N,
\label{xis}\end{equation}
в \re{R_e} при $n=N$, мы находим
\begin{eqnarray}\nonumber
&& R_N(z_s)={_{10}}v_9(1-x_1; 1+x_0-d_4,1+x_0-d_5,
\\ && \makebox[4em]{}
N-x_2+1,s+2+2x_0-d_1-d_2,-s).
\label{RN} \end{eqnarray}
Этот ряд суммируется с помощью формулы Френкеля-Тураева, что дает
\begin{equation}
R_N(z_s)= \frac{[2-x_1,d_4+d_5-x_0-x_2,-N+d_4,-N+d_5]_s}
{[1+d_4-x_2,1+d_5-x_2,1+x_0-N,d_4+d_5-1-N-x_0]_s}.
\label{RN-11} \end{equation}

Приняв во внимание, что $P_N(z_s)/A_N(z_s)=\kappa_NR_N(z_s)$, и
подставив все необходимые выражения в \re{ex_g}, мы находим
\begin{eqnarray}
&& g_s= \frac{t_N [2s+2+2x_0-d_1-d_2]}{[s+1+x_0-d_1,s+1+x_0-d_2]}
\label{gs}\\ && \makebox[4em]{}
\times \frac{[x_0+1,2+2x_0-d_1-d_2,1+d_4-x_2]_s}
{[1,2+x_0-d_1-d_2,3+2x_0-d_1-d_2+N]_s}
\nonumber \\ && \makebox[4em]{}
\times \frac{[1+d_5-x_2,1+x_0+x_2-d_1-d_2,-N]_s}
{[2-x_1,d_4-N,d_5-N]_s},
\nonumber \\ &&
t_N= \frac{[1-d_4,1-d_5,2-x_1,2+x_0-d_1-d_2]_N}
{[2-x_2,1+x_0-d_4,1+x_0-d_5,2+2x_0-d_1-d_2]_{N+1}}.
\nonumber \end{eqnarray}

Теперь мы можем вычислить саму цепную дробь:
\begin{eqnarray*}
&& F_N(z(u))=\sum_{s=0}^{N} \frac{g_s}{z(u)-z_s}
=\frac{[u+d_2-e_1,u+d_1-e_1]}{[d_2-e_1,e_1-d_1]}
\\ && \makebox[4em]{} \times
\sum_{s=0}^N g_s\frac{[s+1+x_0-d_1,s+1+x_0-d_2]}
{[s+1+u+x_0-e_1,s+1-u+x_0-e_2]}
\\ && % \makebox[4em]{}
=\frac{t_N[2+2x_0-d_1-d_2,u+d_2-e_1,u+d_1-e_1]}
{[d_2-e_1,d_1-e_1,u+1+x_0-e_1,u-1-x_0+e_2]}\:
{_{12}v_{11}}(u_0;u_1, \ldots, u_7;\sigma,\tau;1),
\end{eqnarray*}
где
\begin{eqnarray*}
&& u_0=2+2x_0-d_1-d_2, \quad u_1=1+u+x_0-e_1, \quad u_2=1-u+x_0-e_2, \quad
u_4=1+x_0,
\\  &&
u_3=1+x_0+x_2-d_1-d_2, \quad u_5=1+d_4-x_2, \quad u_6=1+d_5-x_2,
\quad  u_7=-N.
\end{eqnarray*}
Применим теперь к этому $_{12}v_{11}$ ряду эллиптическое преобразование Бэйли
и получим
\begin{eqnarray} \nonumber
&& F_N(z(u))=
K_N\frac{[u+d_2-e_1,u+d_1-e_1]}{[u+1+x_0-e_1,u-1-x_0+e_2]}\:
{_{12}v_{11}}(2-x_1;1,1+x_0,
\\ && \makebox[2em]{}
1+u-x_2+e_2,1-u-x_2+e_1, 1+d_4-x_2,1+d_5-x_2,-N;\sigma,\tau;1),
\label{cont-fr} \\ &&
K_N=t_N\frac{[2+2x_0-d_1-d_2][3+2x_0-d_1-d_2]_N}
{[d_2-e_1,d_1-e_1][3-x_1]_N}
\nonumber \\ && \makebox[4em]{}
\times \frac{[2-d_4+x_0,2-d_5+x_0,x_2-N-1]_N}
{[d_4-N,d_5-N,d_1+d_2-x_0-N-1]_N}
\nonumber \\ && \makebox[4em]{}
=\frac{[2-x_1]}{[2+N-x_1,1+x_0-d_4,1+x_0-d_5,d_2-e_1,d_1-e_1]}.
\nonumber\end{eqnarray}
Этот результат был анонсирован в работе \cite{spi-zhe:gevp} и опубликован в
статье \cite{spi-zhe:theory}. Придавая параметрам $_{12}v_{11}$ ряда специальные
значения, так же как это делалось в случае  $_{10}\varphi_9$ ряда в
\cite{gup-mas:contiguous}, можно его редуцировать к $_{10}v_9$ рядам,
суммировать их и выражать соответствующие цепные дроби в виде отношения
произведений тета-функций. Данное эллиптическое расширение обрывающейся
цепной дроби Рамануджана-Ватсона-Гупты-Массона представляет собой самую
общую явную цепную дробь, найденную к настоящему времени.

\section{Связь с квантовыми моделями типа Ка\-ло\-дже\-ро}

Модели типа Калоджеро (или Калоджеро-Сазерленда-Мозера) описывают
полностью интегрируемые одномерные многочастичные квантовые системы с рациональными,
тригонометрическими (или гиперболическими) и эллиптическими потенциалами.
Они нашли многочисленные приложения в
теоретической физике и их исследованию посвящено большое количество работ.
В работах Ольшанецкого и Переломова установлена глубокая связь между
такими моделями и классическими системами корней (см., например,
систематический обзор \cite{ols-per:quantum} и ссылки, приведенные в нем).
Одна из наиболее общих моделей с эллиптическим взаимодействием для $BC_n$
системы корней была предложена Иноземцевым \cite{ino:lax}. Мы ограничимся
кратким описанием релятивистского обобщения этой модели, предложенного
ван Диехеном \cite{die:integrability,die:difference} и подробно
исследованного Комори и Хиками \cite{kom,kom-hik:quantum}.
Эта модель является также многопараметрическим обобщением модели
Рюйсенаарса \cite{rui:complete}. Оказалось, что одночастичный сектор
этой модели связан с с эллиптическими гипергеометрическими функциями.

Гамильтониан Иноземцева имеет вид
\be
H=-\frac{1}{2}\sum_{i=1}^n\frac{\partial^2}{\partial x_i^2}+\nu\sum_{1\leq i<j\leq n}
\left({\mathcal P}(x_i-x_j)+{\mathcal P}(x_i+x_j)\right)
+\sum_{i=1}^n\sum_{k=0}^3\nu_k{\mathcal P}(x_i+\omega_k),
\lab{ino}\ee
где ${\mathcal P}(x)$ функция Вейерштрасса c полупериодами $\omega_{1,2}$
и $\omega_0=0,\, \omega_3=\omega_1+\omega_2$. Он содержит пять свободных
параметров $\nu, \nu_0,\ldots,\nu_3$. Релятивистcкое или конечно-разностное
обобщение этого оператора содержит большее число параметров и имеет вид
(в системе обозначений работы \cite{kom})
\be
H=\sum_{j=1}^n\left( a_j({\bf u})\tau_j+a_j(-{\bf u})\tau_j^{-1}\right)
+b({\bf u}),
\lab{hik}\ee
где $\tau_j$ являются операторами сдвига
$$
\tau_j^{\pm 1} f(u_1,\ldots,u_n)=f(u_1,\ldots, u_j\pm \kappa,\ldots,u_n).
$$
Потенциалы $a_j$ и $b$ имеют вид
\begin{eqnarray*}
&& a_j({\bf u})=\prod_{k=1\atop k\neq j}^n\frac{\theta_1(u_j-u_k-\mu)
\theta_1(u_j+u_k-\mu)}{\theta_1(u_j-u_k)\theta_1(u_j+u_k)}
\prod_{r=0}^3\frac{\theta_r(u_j-\mu_r)\theta_r(u_j+\kappa/2-\mu_r')}
{\theta_r(u_j)\theta_r(u_j+\kappa/2)},
\\
&& b({\bf u})=\left(\frac{\pi}{\theta'_1(0)}\right)^2
\frac{2}{\theta_1(\mu)\theta_1(\kappa+\mu)}
\sum_{r=0}^3 \left(\prod_{k=0}^3\theta_k(\mu_{\pi_rk}+\kappa/2)
\theta_k(\mu'_{\pi_rk})\right)
\\ && \makebox[10em]{}
\times \prod_{j=1}^n\frac{\theta_r(u_j-\kappa/2-\mu)\theta_r(-u_j-\kappa/2-\mu)}
{\theta_r(u_j-\kappa/2)\theta_r(-u_j-\kappa/2)},
\end{eqnarray*}
где используются перестановки $\pi_1=\mbox{id}$, $\pi_1=(12)(03)$,
$\pi_3=(13)(02)$ и $\pi_0=(01)(23)$. Тета-функции Якоби $\theta_r(u)=\theta_r(u|\tau)$
имеют стандартный вид
\begin{eqnarray*}
&& \theta_1(u)=-i\sum_{n\in{\Z }} (-1)^n e^{i\pi\tau(n+1/2)^2+2\pi iu(n+1/2)}
=ip^{1/8}(p;p)_\infty z^{-1/2}\theta(z;p),
\\ && \theta_2(u)=\sum_{n\in{\Z }} e^{i\pi\tau(n+1/2)^2+2\pi iu(n+1/2)}
=p^{1/8}(p;p)_\infty z^{-1/2}\theta(-z;p),
\\&& \theta_3(u)=\sum_{n\in{\Z }}e^{i\pi\tau n^2+2\pi iu n}
=(p;p)_\infty \theta(-p^{1/2}z;p),
\\&& \theta_0(u)=\sum_{n\in{\Z }} (-1)^n e^{i\pi\tau n^2+2\pi iu n}
=(p;p)_\infty \theta(p^{1/2}z;p),
\\ && \makebox[2em]{}
z=e^{2\pi iu}, \qquad p=e^{2\pi i\tau}, \qquad
\theta_1'(0)=2\pi p^{1/8}(p;p)_\infty^3.
\end{eqnarray*}

Воспользуемся мультипликативной системой обозначений
\begin{eqnarray*}
&& q=e^{2\pi i\kappa}, \quad z_j=e^{2\pi i u_j}, \quad t=e^{-2\pi i\mu},\quad
\\&& t_r=\epsilon_r e^{-2\pi i\mu_r}, \quad
t_{r+4}=\epsilon_r e^{\pi i\kappa-2\pi i\mu_r'},\; r=0,1,2,3,
\end{eqnarray*}
где $\epsilon_1=1, \epsilon_2=-1, \epsilon_3=-p^{1/2}, \epsilon_0=p^{1/2}$,
и наложим ограничение на параметры $\prod_{r=0}^7t_r=p^2q^2$.
Тогда одночастичный гамильтониан этой системы может быть переписан в виде
$H=\rho L$, где $\rho=\exp\{\pi i(\mu_1+\mu_1'+\mu_2+\mu_2')\}$ и,
с точностью до аддитивной постоянной,
\begin{eqnarray}\nonumber
&& L= A(z)(T-1)+A(z^{-1})(T^{-1}-1), \qquad Tf(z)=f(qz),
\\ && \makebox[4em]{}
A(z)=\frac{\prod_{r=0}^7\theta(t_rz;p)}{\theta(z^2,qz^2;p)}.
\lab{dkh-ham}\end{eqnarray}
Благодаря ограничению, наложенному на $t_r$, мы имеем $A(pz)=A(z)$,
т.е. функция $A(z)$ эллиптична по $\log z$.

Сравним стандартную задачу на собственные значения для оператора $L$,
$L\psi=\lambda\psi$,  с $q$-разностным уравнением для биортогональных функций
\re{eig}. Легко увидеть, что они совпадают при следующих ограничениях,
наложенных на параметры $t_5, t_6, t_7$ и $\lambda$:
\begin{eqnarray}\nonumber
&& t_5=\frac{t_4}{q},\qquad t_6=\frac{pq^2}{\mu A}, \qquad t_7=\frac{pq\mu}{t_4},
\\ &&
\lambda=-\theta\left(\frac{A\mu}{qt_4},\frac{1}{\mu};p\right)
\prod_{r=0}^3\theta\left(\frac{t_rt_4}{q};p\right),
\lab{par-restr}\end{eqnarray}
где $A=\prod_{r=0}^5t_r$. Таким образом, $_{12}V_{11}$
эллиптический гипергеометрический ряд описывает собственные функции для
выделенного собственного значения $\lambda$ в спектральной
задаче для одночастичного сектора разностной модели типа Калоджеро со
специально подобранными значениями параметров (аналогичное наблюдение
было сделано Комори на основе работы \cite{spi-zhe:gevp}). С точки зрения
разностных уравнений всего имеется три ограничения: условие балансировки
$\prod_{r=0}^7t_r=p^2q^2$, связь между двумя параметрами $t_4=qt_5$ и
фиксация собственного значения $\lambda$ как явной функции параметров.
Квантование $\mu=q^n$ существенно для решения в виде
$_{12}V_{11}$-ряда, а для функции $V(t)$ такого ограничения уже нет.

Однако имеется более глубокая связь между моделью \re{hik} и
эллиптическими гипергеометрическими функциями. Гамильтониан $L$
(без ограничений \re{par-restr}) формально самосопряжен относительно
свертки
$$
\langle\psi|\chi\rangle=\frac{(p;p)_\infty(q;q)_\infty}{4\pi i}
\int_{\T}\omega(z)\psi(z)\chi(z)\frac{dz}{z},
\quad \omega(z)=\frac{\prod_{r=0}^7\Gamma(t_rz,t_rz^{-1};q,p)}
{\Gamma(z^2,z^{-2};q,p)},
$$
и $\psi(z)=1$ есть его собственная функция с собственным значением
$\lambda=0$, $L\psi=0$.
Нормировка этой собственной функции (при некоторых ограничениях на абсолютные
значения параметров) в точности совпадает с функцией $V(t)$,
обладающей $E_7$ группой преобразований симметрии.

Как показано в работе \cite{kmnoy}, эллиптический гипергеометрический
ряд $_{12}V_{11}$ дает нетривиальное мероморфное решение эллиптического
уравнения Пенлеве, найденного Сакаи \cite{sak} и включающего
в качестве предельных случаев все шесть уравнений Пенлеве и их обычных
разностных и $q$-разностных аналогов. (В статье \cite{kmnoy} использовалось
старое обозначение для этого ряда $_{10}E_9$.) Поскольку это решение
выводилось с помощью редукции к эллиптическому гипергеометрическому
уравнению, то наша функция $V(t)$ дает новое мероморфное решение
уравнения Сакаи.  Следуя работе Манина \cite{man:painleve},
Левин и Ольшанецкий \cite{LO} ввели понятие Пенлеве-Калоджеро соответствия,
согласно которому можно связать с каждым уравнением Пенлеве некоторую
модель типа Калоджеро. Однако это соответствие не было распространено
на дискретные уравнения Пенлеве. Автор предполагает, что в рамках расширенного
трактования этой концепции должна существовать связь между уравнением Сакаи и моделью
ван Диехена \cite{die:integrability,die:difference} ранга $n=1$ и
что в обоих случаях работает одна и та же исключительная
группа симметрий $E_8$.

В многомерном случае, при выполнении условия балансировки
$t^{2n-2}\prod_{m=0}^7t_m=p^2q^2$ гамильтониан \re{hik} эквивалентен
оператору
\begin{eqnarray}\nonumber
&& L= \sum_{j=1}^n\left(A_j(z)(T_j-1)+A_j(z^{-1})(T_j^{-1}-1)\right),
\quad T_jf(z)=f(\ldots,z_{j-1},qz_j,z_{j+1},\ldots),
\\ && \makebox[4em]{}
A_j(z)=\frac{\prod_{r=0}^7\theta(t_rz_j;p)}{\theta(z_j^2,qz_j^2;p)}
\prod_{k=1 \atop k\neq j}^n \frac{\theta(tz_jz_k^\pm;p)}
{\theta(z_jz_k^\pm;p)}.
\lab{multi-op}\end{eqnarray}
При подходящих ограничениях на параметры, этот оператор формально
самосопряжен относительно свертки
\be
\langle\psi|\chi\rangle = \kappa_n^C \int_{\T^n}\omega(z)\psi(z)\chi(z)
\frac{dz}{z},\qquad \frac{dz}{z}\equiv\frac{dz_1}{z_1}\cdots
\frac{dz_n}{z_n},\quad \kappa_n^C=\frac{(p;p)_\infty^n(q;q)_\infty^n}
{(2\pi i)^n2^nn!},
\lab{inner-pr}\ee
с весовой функцией, которая определяет многомерный аналог функции
$V(t)$ для $BC_n$ системы корней:
$$
V({\bf t};BC_n)=\kappa_n^C\int_{\T^n}\omega(z)\frac{dz}{z}=\kappa_n^C
\int_{\T^n}\prod_{1\leq j<k\leq n}\frac{\Gamma(tz_j^\pm z_k^\pm;q,p)}
{\Gamma(z_j^\pm z_k^\pm;q,p)}\prod_{j=1}^n
\frac{\prod_{r=0}^7\Gamma(t_rz_j^\pm;q,p)}
{\Gamma(z_j^{\pm2};q,p)}\frac{dz}{z}.
$$
Поскольку $\psi(z)=1$ является $\lambda=0$ собственной функцией оператора
$L$, то ее норма совпадает с многомерной эллиптической гипергеометрической
функцией, $||1||^2=V({\bf t};BC_n)$. В работе \cite{rai:trans} был найден ряд нетривиальных
преобразований симметрии для функции $V({\bf t};BC_n)$.

С точки зрения эллиптических гипергеометрических интегралов,  можно
построить аналоги функции $V({\bf t};BC_n)$ для всех известных
многомерных эллиптических бета-ин\-те\-гра\-лов и сконструировать
операторы, для которых они будут определять нормировку
простейших собственных функций. Например, определим весовую функцию
связанную с многопараметрическим $C_n$ бета-интегралом
$$
\omega(z)=\prod_{1\leq j<k\leq n}\frac{1}{\Gamma(z_j^\pm z_k^\pm;q,p)}
\prod_{j=1}^n\frac{\prod_{m=1}^{2n+6}\Gamma(t_mz_j^\pm;q,p)}
{\Gamma(z_j^{\pm2};q,p)},
$$
где $\prod_{m=1}^{2n+6}t_m=p^2q^2$. Тогда скалярное произведение
\re{inner-pr} приводит к  оператору
\begin{eqnarray}\nonumber &&
L= \sum_{j=1}^n\left(A_j(z)(T_j-1)+A_j(z^{-1})(T_j^{-1}-1)\right),
\\ && \makebox[2em]{}
A_j(z)=\frac{\prod_{r=1}^{2n+6}\theta(t_rz_j;p)}{\theta(z_j^2,qz_j^2;p)}
\prod_{k=1 \atop k\neq j}^n \frac{1}{\theta(z_jz_k^\pm;p)},
\lab{multi}\end{eqnarray}
который формально самосопряжен относительно него,
$\langle\psi|L \chi\rangle=\langle L\psi|\chi\rangle$.
Для ранга $n=1$ эта модель совпадает с \re{dkh-ham}, однако
неизвестно является ли она полностью интегрируемой.
Точно таким же образом можно рассмотреть и другие эллиптические
бета-интегралы, описанные в предыдущей главе---для всех них
довольно легко найти аналоги оператора \re{multi} и все они
совпадают в ранге $n=1$. К сожалению, эти новые предположительно
интегрируемые системы вырождены---в той или иной форме константа
связи парного взаимодействия частиц этих моделей является фиксированной.

Описанная в этом параграфе связь между эллиптической интегрируемой
системой типа Калоджеро и эллиптическими гипергеометрическими функциями
заслуживает более детального изучения. В частности, было бы интересно
прояснить возможную связь многомерного обобщения биортогональных
функций этой главы, предложенного Рэйнсом \cite{rai:abelian},
и моделью ван Диехена.
Другая важная проблема состоит в выяснении действительно ли являются
интегрируемыми многочастичные модели, определяемые многомерными аналогами
функции $V(t)$ для других известных эллиптических бета-интегралов --- трех
многопараметрических $C_n$ и $A_n$ интегралов типа I и двух
$A_n$ интегралов типа II. Если эта гипотеза подтвердится, то
необходимо будет найти их обобщения на случаи без условий балансировки по
аналогии с гамильтонианом \re{hik} и развить $R$-матричный подход
\cite{kom,kom-hik:quantum} к этим расширенным квантовым моделям.

%% file: CHAPTER5.TEX
\chapter[Эллиптические гипергеометрические функции с $|q|=1$]
{Эллиптические гипергеометрические функции с $|q|=1$}

Все эллиптические бета-интегралы обсуждавшиеся до сих пор определены
только при базисных параметрах, лежащих внутри диска единичного
радиуса, $|q|,|p|<1$. Эллиптические гипергеометрические ряды, обсуждавшиеся
во второй главе требуют $|p|<1$ и не накладывают ограничений на $q$.
Эллиптические гипергеометрические интегралы, хорошо определенные при
$|q|=1$, были введены в работе \cite{spi:integrals} и они были частично описаны
в третьей главе. Для их определения необходима модифицированная эллиптическая
гамма-функция. В этой главе мы описываем эллиптические бета-интегралы
на корневых системах $A_n$ и $C_n$, хорошо определенные при $|p|<1$ и $ |q|\leq 1$.

\section{Модифицированная эллиптическая гамма-функ\-ция}

Удобно использовать экспоненциальную параметризацию базисных параметров:
\begin{eqnarray}\nonumber
&& q=e^{2\pi i\frac{\omega_1}{\omega_2}}, \qquad p=e^{2\pi
i\frac{\omega_3}{\omega_2}},\qquad r=e^{2\pi
i\frac{\omega_3}{\omega_1}},
\\
&& \tilde q =e^{-2\pi i\frac{\omega_2}{\omega_1}}, \qquad \tilde
p=e^{-2\pi i\frac{\omega_2}{\omega_3}},\qquad \tilde r=e^{-2\pi
i\frac{\omega_1}{\omega_3}},
\label{ell-b}\end{eqnarray}
где $\omega_{1,2,3}$ некоторые комплексные переменные. Для
$|q|, |p|, |r|<1$, в работе \cite{spi:integrals} автором была определена
модифицированная эллиптическая гамма-функция
\begin{equation}
G(u;\boldsymbol{\omega})= \prod_{j,k=0}^\infty \frac{(1-e^{-2\pi
i\frac{u}{\omega_2}}q^{j+1}p^{k+1}) (1-e^{2\pi i
\frac{u}{\omega_1}}{\tilde q}^{j+1}{r }^k)} {(1-e^{2\pi i
\frac{u}{\omega_2}}q^jp^k) (1-e^{-2\pi i
\frac{u}{\omega_1}}{\tilde q}^j{r}^{k+1})}.
\label{ell-d2}\end{equation}
Она удовлетворяет трем линейным разностным уравнениям первого порядка
\begin{eqnarray}
&& G(u+\omega_1;\boldsymbol{\omega})=\theta(e^{2\pi
i\frac{u}{\omega_2}};p) G(u;\boldsymbol{\omega}),
\label{ell-1eq1} \\
&& G(u+\omega_2;\boldsymbol{\omega})=\theta(e^{2\pi
i\frac{u}{\omega_1}};r) G(u;\boldsymbol{\omega}),
\label{ell-2eq2} \\
&& G(u+\omega_3;\boldsymbol{\omega})=e^{-\pi iB_{2,2}(u;\mathbf{\omega})}
G(u;\boldsymbol{\omega}),
\label{ell-3eq-r}\end{eqnarray}
где
$$
B_{2,2}(u;\boldsymbol{\omega})=\frac{u^2}{\omega_1\omega_2}
-\frac{u}{\omega_1}-\frac{u}{\omega_2}+
\frac{\omega_1}{6\omega_2}+\frac{\omega_2}{6\omega_1}+\frac{1}{2}.
$$
В статье \cite{die-spi:unit} было найдено другое представление для этой
функции, связанное с модулярными преобразованиями для эллиптической гамма-функции
\cite{fel-var:elliptic},
\begin{equation}
G(u;\boldsymbol{\omega})=e^{-\pi iP(u)}\Gamma(e^{-2\pi
i\frac{u}{\omega_3}}; \tilde r, \tilde p),
\label{gamma-tr}\end{equation}
где $P(u)$ есть полином третьей степени
\begin{equation}
P(u)=\frac{1}{3\omega_1\omega_2\omega_3}
\left(u-\frac{1}{2}\sum_{n=1}^3\omega_n\right)
\left(u^2-u\sum_{n=1}^3\omega_n+\frac{\omega_1\omega_2\omega_3}{2}
\sum_{n=1}^3\frac{1}{\omega_n}\right).
\label{p(u)}\end{equation}
Из него легко увидеть, что, в отличие от $\Gamma(z;q,p)$,
модифицированная эллиптическая гамма-функция хорошо определена
при $\omega_1/\omega_2>0$ (или $|q|=1$). Доказательство соотношения
\re{gamma-tr} совершенно элементарное. Необходимо просто проверить,
что выражение в правой части в \re{gamma-tr} удовлетворяет тем же
самым трем уравнениям \re{ell-1eq1}-\re{ell-3eq-r}, что достигается
с помощью модулярных преобразований для тета-функции. Поскольку
не бывает нетривиальных трижды периодических функций с несоизмеримыми
периодами, то правая часть \re{gamma-tr} пропорциональна $G(u;\boldsymbol{\omega})$.
Коэффициент пропорциональности равен единице поскольку обе функции
равны единице при выборе $u=(\omega_1+\omega_2+\omega_3)/2$.

В пределе $p,r\to 0$, функция $G(u;\boldsymbol{\omega})$ становится
обратной к функции двойного синуса
\begin{equation}
\lim_{p,r\to 0} \frac{1}{G(u;\boldsymbol{\omega})}
=S(u;\omega_1,\omega_2) = \frac{(e^{2\pi i u/\omega_2}; q)_\infty}
{(e^{2\pi iu/\omega_1}\tilde q; \tilde q)_\infty}.
\label{2d-sin2}\end{equation}

Мы используем так же обозначения
\begin{eqnarray*}
&&
G(u_1,\ldots,u_m)\equiv G(u_1;\boldsymbol{\omega})\cdots
G(u_1;\boldsymbol{\omega}),
\\ &&
S(u_1,\ldots,u_m)\equiv S(u_1;\boldsymbol{\omega})\cdots
S(u_1;\boldsymbol{\omega}).
\end{eqnarray*}

При $\mbox{Im}(\omega_{1}/\omega_2)>0$ (или $|q|<1$) функция
двойного синуса может быть определена бесконечным произведением
(\ref{2d-sin2}). Однако, она остается хорошо определенной мероморфной функцией
$u$ и в области $\omega_1/\omega_2>0$ (или $|q|=1$) \cite{kls:unitary}.
Ее нули расположены в точках $u=-\omega_1\N-\omega_2\N$
и полюсы в $u=\omega_1(1+\N)+\omega_2 (1+\N)$.
Асимптотически,
\begin{eqnarray}\label{asymp1}
&& \lim_{\mbox{Im}(\frac{u}{\omega_1}),\; \mbox{Im}(\frac{u}{\omega_2})\to
+\infty}S(u;\boldsymbol{\omega})=1,
\\ &&
\lim_{\mbox{Im}(\frac{u}{\omega_1}),\; \mbox{Im}(\frac{u}{\omega_2})\to
-\infty} e^{\pi iB_{2,2}(u;\boldsymbol{\omega})}
S(u;\boldsymbol{\omega}) =1 .
\label{asymp2}\end{eqnarray}

\section{Модифицированные эллиптические бета-ин\-те\-гра\-лы с $|q|\leq 1$}

\subsection{Одномерный интеграл}

Обратимся к эллиптическим бета-интегралам. Формула \re{ell-int} справедлива
при $|q|,|p|<1$ и $|t_n|<1$, $n=0, \dots,4$,  $ |pq|<|A|$, где $A\equiv\prod_{n=0}^4t_n$.
Следующая теорема описывает модифицированный эллиптический бета-интеграл,
который, в отличие от \re{ell-int}, хорошо определен при $|q|=1$.

\begin{theorem}\label{circle-int:thm}
Пусть $\mbox{Im} (\omega_1/\omega_2)\geq 0$ и $\mbox{Im}
(\omega_3/\omega_1), \mbox{Im} (\omega_3/\omega_2)>0$. Возьмем пять
комплексных параметров $g_n$, $n=0,\ldots,4$, удовлетворяющих ограничениям
$$
\mbox{Im}(g_n/\omega_3)<0, \qquad
\mbox{Im}((\mathcal{A}-\omega_1-\omega_2)/\omega_3)>0,
$$
с $\mathcal{A}\equiv \sum_{n=0}^4g_n$. Тогда справедлива следующая формула
интегрирования
\begin{equation}
\int_{-\omega_3/2}^{\omega_3/2} \frac{\prod_{n=0}^4 G(g_n\pm
u;\boldsymbol{\omega})} {G(\pm 2u, \mathcal{A}\pm
u;\boldsymbol{\omega})} \frac{du}{\omega_2} = \kappa\, \frac{
\prod_{0\leq n<m\leq 4}G(g_n+g_m;\boldsymbol{\omega})}
{\prod_{n=0}^4G(\mathcal{A}-g_n;\boldsymbol{\omega})},
\label{circle-int}\end{equation}
где
\begin{equation}
\kappa= \frac{-2(\tilde q;\tilde q)_\infty}
{(q;q)_\infty(p;p)_\infty(r;r)_\infty}. \label{kappa1}
\end{equation}
Здесь интеграл берется вдоль сегмента прямой линии, соединяющей
точку $-\omega_3/2$ с $\omega_3/2$, а также использовано компактное
обозначение $G(a\pm b;\boldsymbol{\omega})\equiv G(a+b,a-b;\boldsymbol{\omega})$.
\end{theorem}
{\bf Доказательство.}
Начнем с подстановки соотношения (\ref{gamma-tr}) в левую часть (\ref{circle-int}).
Это дает
\begin{equation}
e^{\pi ia}\int_{-\omega_3/2}^{\omega_3/2}\frac{\prod_{n=0}^4
\Gamma(e^{-2\pi i\frac{g_n\pm u}{\omega_3}};\tilde r,\tilde p)}
{\Gamma(e^{\mp 4\pi i\frac{u}{\omega_3}}, e^{-2\pi
i\frac{\mathcal{A}\pm u}{\omega_3}};\tilde r,\tilde p)}
\frac{du}{\omega_2}, \label{cir-int}\end{equation}
где
\begin{eqnarray}\nonumber
&&
a=\frac{2}{3\omega_1\omega_2\omega_3}\left(\mathcal{A}^3-\sum_{n=0}^4g_n^3
\right) -\frac{\sum_{m=1}^3\omega_m}{\omega_1\omega_2\omega_3}
\left(\mathcal{A}^2-\sum_{n=0}^4g_n^2 \right)
\\  && \makebox[4em]{}
+\frac{1}{2}\left(\sum_{m=1}^3\omega_m\right)
\left(\sum_{m=1}^3\frac{1}{\omega_m}\right).
\nonumber\end{eqnarray}
Ограничения на параметры позволяют использовать формулу (\ref{ell-int})
с заменами
\begin{equation*}
z\to e^{2\pi i\frac{u}{\omega_3}}, \quad  t_n\to e^{-2\pi i\frac{g_n
}{\omega_3}},\quad p\to e^{-2\pi i\frac{\omega_1}{\omega_3}},\quad
q\to e^{-2\pi i\frac{\omega_2}{\omega_3}} ,
\end{equation*}
которые дают для (\ref{cir-int})
\be\nonumber
 \frac{2\omega_3\omega_2^{-1} e^{\pi ia}}{(\tilde r; \tilde r)_\infty
(\tilde p;\tilde p)_\infty}\frac{\prod_{0\leq n<m\leq 4}
\Gamma(e^{-2\pi i\frac{g_n+g_m}{\omega_3}};\tilde r,\tilde p)}
{\prod_{n=0}^4\Gamma(e^{-2\pi
i\frac{\mathcal{A}-g_n}{\omega_3}};\tilde r,\tilde p)}
= \kappa\, \frac{\prod_{0\leq n<m\leq 4}
G(g_n+g_m;\boldsymbol{\omega})} {\prod_{n=0}^4
G(\mathcal{A}-g_n;\boldsymbol{\omega})}, % && \nonumber
\ee
с
$$
\kappa= \frac{2\omega_3e^{\frac{\pi i}{12}(\sum_{m=1}^3\omega_m)
(\sum_{m=1}^3\omega_m^{-1}) }}
{\omega_2(\tilde r; \tilde r)_\infty(\tilde p;\tilde p)_\infty}.
$$
Применяя модулярное преобразование
\begin{equation}
e^{-\frac{\pi i}{12\tau}} \left(e^{-2\pi i/\tau};e^{-2\pi
i/\tau}\right)_\infty =(-i\tau)^{1/2}e^{\frac{\pi i\tau}{12}}
\left(e^{2\pi i\tau};e^{2\pi i\tau}\right)_\infty
\label{ded}\end{equation}
к бесконечным произведениям, входящим в $\kappa$, мы получаем
$$
\kappa= -2\sqrt{\frac{\omega_1}{i\omega_2}}
\frac{e^{\frac{\pi i}{12}(\frac{\omega_1}{\omega_2}
+\frac{\omega_2}{\omega_1}) }}{(r;r)_\infty (p;p)_\infty}.
$$
Использовав (\ref{ded}) еще раз, мы можем заменить экспоненциальный
множитель на отношение бесконечных произведений, что приводит к
искомому виду $\kappa$ в (\ref{kappa1}).
\hfill{Q.E.D.}

\smallskip

Рассмотрим теперь формальный предел $p,r\to 0$ в интеграле
(\ref{circle-int}). Для этого мы фиксируем квазипериоды
$\omega_{1,2}$, такие что $\mbox{Im}(\omega_1/\omega_2)\geq 0$ и
$\mbox{Re}(\omega_1/\omega_2)> 0$, а так же мы берем
$\omega_3=it\omega_2$ с $t>0$. При $t\to +\infty$ формула
интегрирования теоремы \ref{circle-int:thm} формально вырождается в равенство
\begin{equation}
\int_{\mathbb{L}}\frac{S(\pm 2u, \mathcal{A}\pm
u;\boldsymbol{\omega})} {\prod_{n=0}^4 S(g_n\pm
u;\boldsymbol{\omega})}\frac{du}{\omega_2}= -2\frac{(\tilde
q;\tilde q)_\infty}{(q;q)_\infty}
\frac{\prod_{n=0}^4S(\mathcal{A}-g_n;\boldsymbol{\omega})} {
\prod_{0\leq n<m\leq 4}S(g_n+g_m;\boldsymbol{\omega})},
\label{red1}\end{equation}
где интегрирование ведется по линии $\mathbb{L}\equiv i\omega_2\mathbb{R}$,
а параметры подвержены ограничениям $\mbox{Re}(g_n/\omega_2)>0$ и
$\mbox{Re}((\mathcal{A}-\omega_1)/\omega_2)<1$.
Этот интеграл был получен аналогичным формальным пределом
из стандартного эллиптического бета-интеграла \re{ell-int} и
строго доказан независимым способом Стокманом в работе \cite{sto:hyperbolic}.

\begin{remark}
Обобщенная гамма-функция $[z;\tau]_\infty$, использованная в
\cite{sto:hyperbolic}, совпадает с функцией двойного синуса
(\ref{2d-sin2}) после отождествления $z=(\omega_2-u)/\omega_2$ и
$\tau=-\omega_2/\omega_1$.
\end{remark}

\begin{remark}
Для пояснения связи между интегралами \re{ell-int} и  \re{circle-int},
напомним что контурный интеграл называется эллиптическим гипергеометрическим
интегралом если для некоторого параметра $\omega_1$
отношение $\Delta(u+\omega_1)/\Delta(u)$ для подынтегральной функции
определяет эллиптическую функцию $u$ с некоторыми периодами $\omega_2$ и $\omega_3$.
После замены переменных
$$
z=e^{2\pi i\frac{u}{\omega_2}}, \qquad t_n=e^{2\pi i\frac{g_n}{\omega_2}},
\qquad A=e^{2\pi i\frac{\mathcal{A}}{\omega_2}},
$$
интеграл (\ref{ell-int}) переходит в аддитивную форму
\begin{equation}\label{ell-int-add}
\int_{-\omega_2/2}^{\omega_2/2}\Delta(u)\frac{du}{\omega_2}=
\frac{2\prod_{0\leq n<m\leq 4} \Gamma(e^{2\pi
i\frac{g_n+g_m}{\omega_2}};q,p)}
{(q;q)_\infty(p;p)_\infty\prod_{n=0}^4\Gamma(e^{2\pi
i\frac{\mathcal{A}-g_n}{\omega_2}};q,p)},
\end{equation}
с подынтегральной функцией
\begin{equation}
\Delta(u)=\frac{\prod_{n=0}^4\Gamma(e^{2\pi i \frac{g_n\pm
u}{\omega_2}};q,p)} {\Gamma(e^{\pm 4\pi
i\frac{u}{\omega_2}},e^{2\pi i\frac{\mathcal{A}\pm u}{\omega_2}};q,p)},
\label{int}\end{equation}
где $\Gamma(e^{a\pm b};q,p)\equiv\Gamma(e^{a+b},e^{a-b};q,p)$.
Легко увидеть, что
\begin{equation}
\frac{\Delta(u+\omega_1)}{\Delta(u)}= e^{2\pi
i\frac{\omega_1}{\omega_2}}\frac{\theta(e^{4\pi
i\frac{u+\omega_1}{\omega_2}};p)\theta(e^{2\pi
i\frac{u+\omega_1-\mathcal{A}}{\omega_2}};p)} {\theta(e^{4\pi
i\frac{u}{\omega_2}};p)\theta(e^{2\pi i\frac{u+\mathcal{A}}{\omega_2}};p)}
\prod_{n=0}^4\frac{\theta(e^{2\pi i\frac{u+g_n}{\omega_2}};p)}
{\theta(e^{2\pi i\frac{u+\omega_1-g_n}{\omega_2}};p)}
\label{int-eq}\end{equation}
есть эллиптическая функция $u$ с периодами $\omega_2$ и $\omega_3$.
С помощью разностного уравнения \re{ell-1eq1} для модифицированной гамма-функции,
нетрудно проверить, что подынтегральная функция в \re{circle-int} так же
дает решение уравнения \re{int-eq}. Следовательно, оба интеграла \re{ell-int} и
\re{circle-int} удовлетворяют взятому определению с одним и тем же уравнением
\re{int-eq}. В то время как интеграл \re{ell-int} связан с решением
\re{int-eq}, определенном только при $|q|<1$ (или при $|q|>1$,
связанном с инверсией $q\to q^{-1}$ в эллиптической гамма-функции
\cite{spi:integrals}), модифицированный интеграл \re{circle-int} соответствует
решению, которое может быть аналитически продолжено из области $|q|<1$
на единичную окружность $|q|=1$. Именно в этом смысле мы будем называть
интеграл \re{circle-int} эллиптическим бета-интегралом на единичной
окружности $|q|=1$.
\end{remark}

\subsection{Многомерные интегралы}

Следующая теорема описывает многомерный аналог формулы интегрирования
теоремы \ref{circle-int:thm}, являющийся модифицированной версией
$C_n$ эллиптического бета-интеграла работы \cite{die-spi:elliptic}
(см. формулу \re{SintB} в третьей главе).

\begin{theorem}\label{mcircle-int:thm}
Пусть $\mbox{Im} (\omega_1/\omega_2)\geq 0$ и $\mbox{Im}
(\omega_3/\omega_1), \mbox{Im} (\omega_3/\omega_2) >0$. Возьмем шесть
комплексных параметров $g$, $g_m$, $m=0,\ldots,4$, удовлетворяющих ограничениям
$$
\mbox{Im}(g/\omega_3),\mbox{Im}(g_n/\omega_3)<0, \qquad
\mbox{Im}((\mathcal{B}-\omega_1-\omega_2)/\omega_3)>0,
$$
с $\mathcal{B}\equiv (2n-2)g+\sum_{m=0}^4g_m$. Тогда
\begin{eqnarray}\nonumber
\int_{-\frac{\omega_3}{2}}^{\frac{\omega_3}{2}}\cdots
\int_{-\frac{\omega_3}{2}}^{\frac{\omega_3}{2}} \prod_{1\leq
j<k\leq n} \frac{G(g\pm u_j \pm u_k;\boldsymbol{\omega})} {G(\pm
u_j \pm u_k;\boldsymbol{\omega})} \prod_{j=1}^n\frac{
\prod_{m=0}^4 G(g_m \pm u_j;\boldsymbol{\omega}) } { G(\pm 2u_j,
\mathcal{B}\pm u_j;\boldsymbol{\omega})}\frac{du_1}{\omega_2}
\cdots \frac{du_n}{\omega_2}  && \nonumber \\  = \kappa^n
n!\prod_{j=1}^n\frac{G(jg;\boldsymbol{\omega})}{G(g;\boldsymbol{\omega})}
\frac{\prod_{0\leq k<m\leq 4}G((j-1)g+g_k+g_m;\boldsymbol{\omega})}
{\prod_{m=0}^4G((1-j)g+\mathcal{B}-g_m;\boldsymbol{\omega})},
\makebox[2em]{} &&
\label{ds-circle}\end{eqnarray}
с $\kappa$ описанной в (\ref{kappa1}) и $G(c\pm a \pm b;\boldsymbol{\omega})\equiv
G(c+a+b,c+a-b,c-a+b,c-a-b;\boldsymbol{\omega})$.
\end{theorem}
{\bf Доказательство.}
Аналогично случаю теоремы \ref{circle-int:thm}, после подстановки
(\ref{gamma-tr}) в левую часть (\ref{ds-circle}) и применения
многократного бета-интеграла (\ref{SintB}), мы получаем некоторую
комбинацию эллиптических гамма-функций. После этого необходимо выразить
получающийся ответ через модифицированную гамма-функцию $G(u;\boldsymbol{\omega})$,
что дает требуемую правую часть с точностью до некоторого экспоненциального
множителя. Для того, чтобы убедиться в правильности коэффициента пропорциональности
$\kappa^n n!$, достаточно проверить, что зависимость от $g$ в указанном экспоненциальном
множителе сокращается. Поскольку при $g\to 0$ интеграл \re{ds-circle} сводится
к $n$-й степени интеграла \re{circle-int}, этот коэффициент устанавливается однозначно.
\hfill{Q.E.D.}

\smallskip

Перейдем теперь к описанию многопараметрических модифицированных эллиптических
бета-интегралов на  корневых системах $C_n$ и $A_n$.

\begin{theorem}
Для корневой системы $C_n$, мы берем $n$ переменных $u=(u_1,\ldots,u_n)$ $\in\C^n$,
$2n+3$ комплексных параметра $g=(g_1,\ldots,g_{2n+3})$, и определяем ядро
\begin{eqnarray}\nonumber
&& \rho(u,g;C_n)=
\prod_{i=1}^n\frac{\prod_{m=1}^{2n+3}G(g_m\pm u_i)}
{G(\pm 2u_i,\mathcal{A}\pm u_i)}\frac{\prod_{m=1}^{2n+3}
G(\mathcal{A}-g_m)}{\prod_{1\leq m<s\leq 2n+3}G(g_m+g_s)}
\\ && \makebox[6em]{} \times
\prod_{1\leq i<j\leq n}\frac{1}{G(\pm u_i\pm u_j)},
\label{kernel-C-unit}\end{eqnarray}
где $\mathcal{A}=\sum_{m=1}^{2n+3}g_m$, модифицированная эллиптическая гамма-функция
 $G(u)\equiv G(u;\boldsymbol{\omega})$, и
$$
G(c\pm a \pm b)\equiv G(c+a+b,c+a-b,c-a+b,c-a-b).
$$
Предположим теперь, что $\mbox{Im} (\omega_1/\omega_2)\geq 0$ и $\mbox{Im}
(\omega_3/\omega_1), \mbox{Im} (\omega_3/\omega_2) >0$, и
$$
\mbox{Im}(g_m/\omega_3)<0,\; m=1,\ldots,2n+3, \qquad
\mbox{Im}((\mathcal{A}-\omega_1-\omega_2)/\omega_3)>0.
$$
Тогда
\begin{equation}
\int_{-\frac{\omega_3}{2}}^{\frac{\omega_3}{2}}\cdots
\int_{-\frac{\omega_3}{2}}^{\frac{\omega_3}{2}} \rho(u,g;C_n)
\frac{du_1}{\omega_2}\cdots \frac{du_n}{\omega_2}
= \kappa^n 2^n n!
\label{C_n-circle}\end{equation}
с
\begin{equation}
\kappa= \frac{-(\tilde q;\tilde q)_\infty}
{(q;q)_\infty(p;p)_\infty(r;r)_\infty}
\label{kappa2}\end{equation}
и интегрированием вдоль прямой линии, соединяющей точки $-\omega_3/2$
и $\omega_3/2$.
\end{theorem}
{\bf Доказательство.}
После подстановки выражения (\ref{gamma-tr}) в (\ref{kernel-C-unit}) и некоторых
прямых вычислений, мы получаем
\begin{eqnarray}\nonumber
&& \rho(u,g;C_n)=e^{-\pi i nP(0)}
\prod_{1\leq i<j\leq n}\frac{1}{\Gamma(z_i^\pm z_j^\pm;\tilde r,\tilde p)}
\prod_{i=1}^n\frac{\prod_{m=1}^{2n+3}\Gamma(t_mz_i^\pm;\tilde r,\tilde p)}
{\Gamma(z_i^{\pm 2},Az_i^\pm;\tilde r,\tilde p)}
\\ && \makebox[6em]{} \times
\frac{\prod_{m=1}^{2n+3}\Gamma(At_m^{-1};\tilde r,\tilde p)}
{\prod_{1\leq m<s\leq 2n+3}\Gamma(t_mt_s;\tilde r,\tilde p)},
\label{rho-C'}\end{eqnarray}
где $z_j=e^{2\pi i u_j/\omega_3},$ $t_m=e^{-2\pi i g_m/\omega_3}$, и
$A=\prod_{m=1}^{2n+3}t_m$. Поскольку
$
du_j=\omega_3dz_j/2\pi i z_j,
$
интегралы по $u_j$ эквивалентны интегралам по $z_j\in\T$. Применяя $C_n$
формулу интегрирования (\ref{ell-int-C}) с $C=\T$ (что позволено благодаря
взятым ограничениям $|t_m|<1, |pq|<|A|$) и $q,p$ замененными на $\tilde r,\tilde p$,
мы получаем, что выражение в левой части (\ref{C_n-circle}) равно
$\kappa^n 2^n n!,$ где
$$
\kappa=\frac{\omega_3e^{\frac{\pi i}{12}(\sum_{m=1}^3\omega_m)
(\sum_{m=1}^3\omega_m^{-1}) }}
{\omega_2(\tilde r; \tilde r)_\infty(\tilde p;\tilde p)_\infty}.
$$
Используя формулу модулярного преобразования для $\eta$-функции Дедекинда
$$
e^{-\frac{\pi i}{12\tau}} \left(e^{-2\pi i/\tau};e^{-2\pi
i/\tau}\right)_\infty =(-i\tau)^{1/2}e^{\frac{\pi i\tau}{12}}
\left(e^{2\pi i\tau};e^{2\pi i\tau}\right)_\infty ,
$$
так же как и в одномерном случае получаем $\kappa$ в виде (\ref{kappa2}).
\hfill{Q.E.D.}

\begin{theorem}
Для простейшего интеграла на $A_n$ корневой системе, мы берем переменные
интегрирования $u=(u_1,\ldots,u_n)$ $\in\C^n$
и комплексные параметры $g=(g_1,\ldots,g_{n+1})$, $h=(h_1,\ldots,h_{n+2})$.
Обозначив $u_{n+1}=-\sum_{k=1}^n u_k$, определим ядро
\begin{eqnarray}\nonumber
&& \rho(u,g,h;A_n)=
\prod_{i=1}^{n+1}\frac{\prod_{m=1}^{n+1}G(g_m-u_i)
\prod_{j=1}^{n+2}G(h_j+u_i)\, G(H+g_i)}
{G(F+H+u_i)\prod_{j=1}^{n+2}G(g_i+h_j)}
\\ && \makebox[4em]{}
\times \prod_{1\leq i<j\leq n+1}\frac{1}{G(u_i-u_j,u_j-u_i)}
\frac{1}{G(F)}\prod_{j=1}^{n+2}\frac{G(F+H-h_j)}{G(H-h_j)},
\label{kernel-A-unit}\end{eqnarray}
где $F=\sum_{m=1}^{n+1}g_m$ и $H=\sum_{j=1}^{n+2}h_j$.
Предположим, что
\begin{eqnarray*}
&& \mbox{Im}\, \frac{\omega_1}{\omega_2}\geq 0,\quad
\mbox{Im}\, \frac{\omega_3}{\omega_1} >0,\quad
\mbox{Im}\, \frac{\omega_3}{\omega_2} >0, \\
&& \mbox{Im}\, \frac{g_m}{\omega_3}<0,\quad
\mbox{Im}\, \frac{h_j}{\omega_3}<0,\quad
\mbox{Im}\, \frac{F+H-\omega_1-\omega_2}{\omega_3}>0
\end{eqnarray*}
для $m=1,\ldots,n+1,\, j=1,\ldots,n+2$. Тогда имеем
\begin{equation}
\int_{-\frac{\omega_3}{2}}^{\frac{\omega_3}{2}}\cdots
\int_{-\frac{\omega_3}{2}}^{\frac{\omega_3}{2}} \rho(u,g,h;A_n)
\frac{du_1}{\omega_2}\cdots \frac{du_n}{\omega_2}
= \kappa^n (n+1)!
\label{A_n-circle}\end{equation}
с тем же самым $\kappa$ как и в предыдущей теореме.
\end{theorem}
{\bf Доказательство.}
Подставляя (\ref{gamma-tr}) в (\ref{kernel-A-unit}), получаем
\begin{eqnarray}\nonumber
e^{\pi i nP(0)}\lefteqn{ \rho(u,g,h;A_n)=
\prod_{i=1}^{n+1}\frac{\prod_{m=1}^{n+1}\Gamma(t_mz_i^{-1};\tilde r,\tilde p)
\prod_{j=1}^{n+2}\Gamma(s_jz_i;\tilde r,\tilde p)\, \Gamma(St_i;\tilde r,\tilde p)}
{\Gamma(TSz_i;\tilde r,\tilde p)\prod_{j=1}^{n+2}\Gamma(t_is_j;\tilde r,\tilde p)}
}&& \\ &&
\times \prod_{1\leq i<j\leq n+1}\frac{1}{\Gamma(z_iz_j^{-1},z_jz_i^{-1};\tilde r,\tilde p)}
\frac{1}{\Gamma(T;\tilde r,\tilde p)}\prod_{j=1}^{n+2}
\frac{\Gamma(TSs_j^{-1};\tilde r,\tilde p)}{\Gamma(Ss_j^{-1};\tilde r,\tilde p)},
\label{rho-A'}\end{eqnarray}
где $z_j=e^{-2\pi i u_j/\omega_3},$ $t_m=e^{-2\pi i g_m/\omega_3},$
$s_j=e^{-2\pi i h_j/\omega_3},$ $T=\prod_{m=1}^{n+1}t_m,$
$S=\prod_{j=1}^{n+2}s_j.$
Подобно $C_n$ случаю, интегралы по $u_j$ становятся эквивалентными
интегралам по $z_j\in\T$. Применим $A_n$ формулу интегрирования (\ref{ell-int-A})
(что позволено благодаря ограничениям, наложенным на параметры
 $|t_m|<1,$ $|s_l|<1,$ и $|pq|<|TS|$) с $q$ и $p$ замененными на
$\tilde r$ и $\tilde p$. Это приводит к тому, что интеграл в (\ref{A_n-circle})
становится равным $\kappa^n (n+1)!.$
\hfill{Q.E.D.}

\smallskip

Рассмотрим теперь модификацию эллиптического бета-интеграла, предложенного Варнааром и
 автором в работе \cite{spi-war:inversions}, на случай $|q|\leq 1$.

\begin{theorem}
Возьмем переменные интегрирования $u=(u_1,\ldots,u_n)$ $\in\C^n$
и комплексные параметры $g=(g_1,\ldots,g_{n+3})$, $h=(h_1,\ldots,h_n)$.
Обозначив $u_{n+1}=-\sum_{k=1}^n u_k$, мы определяем ядро
\begin{eqnarray}\nonumber
\lefteqn{ \delta(u,g,h;A_n)=\prod_{1\leq i<j\leq n+1}\frac{G(F-u_i-u_j)}
{G(u_i-u_j,u_j-u_i)}
\prod_{j=1}^{n+3}\prod_{m=1}^n\frac{G(F+h_m-g_j)}{G(h_m+g_j)}
} && \\ && \times
\prod_{i=1}^{n+1}\frac{\prod_{m=1}^n G(h_m+u_i)
\prod_{j=1}^{n+3}G(g_j-u_i)}{\prod_{m=1}^n G(F+h_m-u_i)}
\prod_{1\leq j<k\leq n+3}
\frac{1}{G(F-g_j-g_k)},
\label{kernel-D-unit}\end{eqnarray}
где $F=\sum_{j=1}^{n+3}g_j$.
Предположим, что $\mbox{Im} (\omega_1/\omega_2)\geq 0$ и $\mbox{Im}
(\omega_3/\omega_1),$ $\mbox{Im} (\omega_3/\omega_2) >0$, а так же
$\mbox{Im}(g_j/\omega_3)<0,\;$ $\mbox{Im}(h_m/\omega_3)<0,\;$
$\mbox{Im}((F+h_m-\omega_1-\omega_2)/\omega_3)>0$
для $j=1,\ldots,n+3,\, m=1,\ldots,n,$ тогда
\begin{equation}
\int_{-\frac{\omega_3}{2}}^{\frac{\omega_3}{2}}\cdots
\int_{-\frac{\omega_3}{2}}^{\frac{\omega_3}{2}} \delta(u,g,h;A_n)
\frac{du_1}{\omega_2}\cdots \frac{du_n}{\omega_2}
= \kappa^n (n+1)!
\label{D_n-circle}\end{equation}
с тем же параметром $\kappa$ как и в предыдущих случаях.
\end{theorem}
{\bf Доказательство.}
Подставляя (\ref{gamma-tr}) в (\ref{kernel-D-unit}), получаем
$$
\delta(u,g,h;A_n)=e^{-\pi i nP(0)}\Delta(z,t,s;A_n),
$$
где
\begin{eqnarray}\label{rho-D'}
\lefteqn{ \Delta(z,t,s;A_n)=\prod_{1\leq i<j\leq n+1}
\frac{\Gamma(Dz_i^{-1}z_j^{-1};\tilde r,\tilde p)}
{\Gamma(z_iz_j^{-1},z_i^{-1}z_j;\tilde r,\tilde p)}
\prod_{j=1}^{n+3}\prod_{m=1}^n\frac{\Gamma(Ds_mt_j^{-1};\tilde r,\tilde p)}
{\Gamma(s_mt_j;\tilde r,\tilde p)}
 }&&  \\ &&
\times\prod_{i=1}^{n+1}\frac{\prod_{m=1}^{n}\Gamma(s_mz_i;\tilde r,\tilde p)
\prod_{j=1}^{n+3}\Gamma(t_jz_i^{-1};\tilde r,\tilde p)}
{\prod_{m=1}^{n}\Gamma(Ds_mz_i^{-1};\tilde r,\tilde p)}
\prod_{1\leq j<k\leq n+3}
\frac{1}{\Gamma(Dt_j^{-1}t_k^{-1};\tilde r,\tilde p)}
\nonumber\end{eqnarray}
с $z_j=e^{-2\pi i u_j/\omega_3},$ $t_j=e^{-2\pi i g_j/\omega_3},$
$s_m=e^{-2\pi i h_m/\omega_3},$ $D=\prod_{j=1}^{n+3}t_j.$
Как показано в главе 3, при $|t_j|, |s_m|<1$ и $|\tilde p\tilde r|<|Ds_m|$
(что верно благодаря взятым нами ограничениям) имеет место следующая $A_n$
(или $``D_n"$) формула интегрирования
\begin{equation}
\int_{\T^n}\Delta(z,t,s;A_n)\frac{dz}{z}
=\frac{(2\pi i)^n (n+1)!}{(\tilde r;\tilde r)_\infty^n
(\tilde p;\tilde p)_\infty^n}.
\label{ell-int-D2}\end{equation}
Применяя это равенство к нашему случаю, получаем что левая часть
 (\ref{D_n-circle}) равна $\kappa^n(n+1)!$, как и требовалось в формулировке теоремы.
\hfill{Q.E.D.}

\section{Предел $p\to 0$ и $q$-бета-интегралы c $|q|\leq 1$}

\subsection{$C_n$ интеграл типа II}

При $p=0$, интеграл \re{SintB} сводится к Густафсоновскому многократному
интегралу (однопараметрическому расширению $q$-интеграла Сельберга),
доказанному в работе \cite{gus:some2}. Построим соответствующий многомерный
аналог интеграла \re{red1}.  Необходимая формула может быть формально
получена из модифицированного интеграла \re{ds-circle} с
$\omega_1/\omega_2>0$ взятием предела $p,r\to 0$, так как это было
пояснено после теоремы \ref{circle-int:thm}. При этом появляются
$q$-бета-интегралы по бесконечному контуру использующие функции двойного синуса.
Некоторые точные одномерные формулы интегрирования такого типа были найдены в работах
\cite{fkv:strongly,pon-tes:clebsch,rui:int,sto:basic,tak:twisted}.

\begin{theorem}\label{nr:thm}
Пусть $\omega_1,\omega_2$ обозначают квазипериоды, удовлетворяющие
ограничениям $\mbox{Im}(\omega_1/\omega_2)\geq 0$ и
$\mbox{Re}(\omega_1/\omega_2)>0$. Более того, пусть $g$ и $g_n,$
$n=0,\ldots,4$, будут параметрами, удовлетворяющими ограничениям
$\mbox{Re}(g/\omega_1),$  $\mbox{Re}(g/\omega_2),$
$\mbox{Re}(g_n/\omega_2)>0$ и
$\mbox{Re}((\mathcal{B}-\omega_1)/\omega_2)<1$ (с тем же самым $\mathcal{B}$
как в теореме \ref{mcircle-int:thm}). Тогда
\begin{equation}\label{q-ds-circle}
\int_{\mathbb{L}^n} \Delta(\mathbf{u};\mathbf{g})
\frac{du_1}{\omega_2} \cdots \frac{du_n}{\omega_2}
={\cal N}(\mathbf{g}),
\end{equation}
где $\mathbb{L}=i\omega_2\mathbb{R}$,
\begin{equation}
\Delta(\mathbf{u};\mathbf{g})=\prod_{1\leq j<k\leq n} \frac{S(\pm
u_j \pm u_k;\boldsymbol{\omega})} {S(g\pm u_j \pm
u_k;\boldsymbol{\omega})} \prod_{j=1}^n\frac{S(\pm 2u_j,
\mathcal{B} \pm u_j;\boldsymbol{\omega})} {\prod_{n=0}^4 S(g_n \pm
u_j;\boldsymbol{\omega}) } \label{q-deg}\end{equation}
и
\begin{equation}
{\cal N}(\mathbf{g}) = (-2)^n n!\frac{(\tilde q;\tilde
q)_\infty^n}
{(q;q)_\infty^n}\prod_{j=1}^n\frac{S(g;\boldsymbol{\omega})}
{S(jg;\boldsymbol{\omega})}
\frac{\prod_{n=0}^4S((1-j)g+\mathcal{B}-g_n;\boldsymbol{\omega})}
{\prod_{0\leq n<m\leq 4}S((j-1)g+g_n+g_m;\boldsymbol{\omega})}.
\label{n}\end{equation}
%\end{subequations}
\end{theorem}

При специальных значениях параметров, Густафсоновские многомерные
интегралы, соответствующие многомерному обобщению интеграла Рахмана
\cite{rah:integral}, сводятся к многомерным $q$-бета интегралам типа
Аски-Вильсона. Соответствующее вырождение теоремы \ref{nr:thm}
выглядит так.

\begin{theorem}\label{aw:thm}
Пусть $\omega_1,\omega_2$ обозначают квазипериоды со свойствами
$\mbox{Im}(\omega_1/\omega_2)\geq 0$ и
$\mbox{Re}(\omega_1/\omega_2)> 0$, а также пусть $g, g_n,$
$n=0,\ldots,3$, обозначают параметры, подчиняющиеся условиям
$\mbox{Re}(g/\omega_1),$ $\mbox{Re}(g/\omega_2),$
$\mbox{Re}(g_n/\omega_2)>0$ и
$\mbox{Re}((\mathcal{B}-\omega_2)/\omega_1)<1$ с
$\mathcal{B}\equiv (2n-2)g+\sum_{n=0}^3g_n$. Тогда,
\begin{eqnarray}\nonumber
\lefteqn{ \int_{\mathbb{L}^n} \prod_{1\leq j<k\leq n} \frac{S(\pm
u_j \pm u_k;\boldsymbol{\omega})} {S(g\pm u_j \pm
u_k;\boldsymbol{\omega})} \prod_{j=1}^n\frac{S(\pm
2u_j;\boldsymbol{\omega})} {\prod_{n=0}^3 S(g_n \pm
u_j;\boldsymbol{\omega}) }\frac{du_1}{\omega_2} \cdots
\frac{du_n}{\omega_2} } && \nonumber \\ && = (-2)^n
n!\frac{(\tilde q;\tilde q)_\infty^n}
{(q;q)_\infty^n}\prod_{j=1}^n\frac{S(g;\boldsymbol{\omega})}
{S(jg;\boldsymbol{\omega})}
\frac{S((1-j)g+\mathcal{B};\boldsymbol{\omega})} {\prod_{0\leq
n<m\leq 3}S((j-1)g+g_n+g_m;\boldsymbol{\omega})}.
\label{q-ds-selberg}\end{eqnarray}
\end{theorem}

При $n=1$, теорема \ref{aw:thm} дает  одномерный $q$-бета интеграл типа Аски-Вильсона
\begin{equation}
\int_{\mathbb{L}}\frac{S(\pm 2u;\boldsymbol{\omega})}
{\prod_{n=0}^3 S(g_n\pm
u;\boldsymbol{\omega})}\frac{du}{\omega_2}= -2\frac{(\tilde
q;\tilde q)_\infty}{(q;q)_\infty}
\frac{S(g_0+g_1+g_2+g_3;\boldsymbol{\omega})} { \prod_{0\leq
n<m\leq 3}S(g_n+g_m;\boldsymbol{\omega})}, \end{equation}
доказанный Рюйсенаарсом \cite{rui:int} и Стокманом \cite{sto:hyperbolic}.

Формально, интеграл теоремы \ref{aw:thm} вытекает из теоремы \ref{nr:thm} с
$\omega_1/\omega_2>0$ после сдвига
$g_4\to g_4+i\omega_2t$ и взятия предела $t\to +\infty$.
Однако, такая процедура требует строгого обоснования, поскольку контур интегрирования
бесконечен. Сначала мы докажем теорему \ref{nr:thm}, а затем укажем
модификации, необходимые для перехода к теореме \ref{aw:thm}.

Предположим на некоторое время, что квазипериоды $\omega_1,\omega_2$
несоизмеримы над полем рациональных чисел $\mathbb{Q}$. Функция двойного
синуса $S(u;\boldsymbol{\omega})$ имеет простые нули, расположенные в
точках $u=-\omega_1\mathbb{N}-\omega_2\mathbb{N}$ и простые полюсы при
$u=\omega_1(1+\mathbb{N})+\omega_2 (1+\mathbb{N})$. Поэтому, ядро
$\Delta(\mathbf{u};\mathbf{g})$ в \re{q-ds-circle} имеет полюсы в точках
\begin{eqnarray}\nonumber
&& \pm u_j=-\mathcal{B}+\omega_1(1+\mathbb{N})+\omega_2(1+\mathbb{N}),
\; g_n+\omega_1\mathbb{N}+\omega_2\mathbb{N},\; n=0,\ldots,4,
\\ \nonumber && \makebox[2em]{}
g+ u_k+\omega_1\mathbb{N}+ \omega_2\mathbb{N},\;
g- u_k+\omega_1\mathbb{N}+ \omega_2\mathbb{N},\;
k=1,\ldots,n,\: k\neq j,
\label{weight-poles}\end{eqnarray}
где $j=1,\ldots, n$.

Комбинируя асимптотики \re{asymp1} и \re{asymp2} можно увидеть, что
отношения $S(u;\boldsymbol{\omega})/S(g+u;\boldsymbol{\omega})$ и
$S(2u,\mathcal{B}+u;\boldsymbol{\omega})/\prod_{n=0}^4S(g_n+u)$
ограничены на комплексной линии $\mathbb{L}=i\omega_2\mathbb{R}$.
Действительно, при $u=ix\omega_2$,
$x\in\mathbb{R}$ мы находимся в стороне от полюсов и имеем
\begin{equation*}
\frac{S(i\omega_2x;\boldsymbol{\omega})}{S(g+i\omega_2x;\boldsymbol{\omega})}=
\begin{cases}
O(1) &\mbox{for}\; x\to +\infty \\
O(e^{-2\pi xg/\omega_1})&\mbox{for}\;x\to -\infty
\end{cases}
\end{equation*}
и
\begin{equation*}
\frac{S(2i\omega_2x,\mathcal{B}+i\omega_2x;\boldsymbol{\omega})}
{\prod_{n=0}^4S(g_n+i\omega_2x;\boldsymbol{\omega})}=
\begin{cases}
O(1) &\mbox{for}\; x\to +\infty \\
O(e^{2\pi x
(2(n-1)g/\omega_1+1+\omega_2/\omega_1)})&\mbox{for}\;x\to -\infty
\end{cases}.
\end{equation*}
Таким образом подынтегральная функция $\Delta(\mathbf{u};\mathbf{g})$
гладкая и убывает экспоненциально на бесконечности области интегрирования
 $\mathbb{L}^n$. Следовательно, интеграл \re{q-ds-circle} сходится.

Для того, чтобы убедиться в справедливости формулы интегрирования,
будем следовать методу работы \cite{rai:trans} и
выделим три параметра $g_0,g_1,g_2$ и факторизуем ядро
$\Delta(\mathbf{u};\mathbf{g})=
\Delta_+(\mathbf{u})\Delta_-(\mathbf{u})$ с
\begin{equation}
\Delta_+(\mathbf{u})=\prod_{1\leq j<k\leq n} \frac{S(u_j \pm
u_k;\boldsymbol{\omega})} {S(g + u_j \pm u_k;\boldsymbol{\omega})}
\prod_{j=1}^n\frac{S(2u_j, \mathcal{B} - u_j,\omega_1+ \mathcal{C}
-u_j; \boldsymbol{\omega})}{ S(\omega_1+ \mathcal{C}
+u_j;\boldsymbol{\omega})\prod_{n=0}^4
S(g_n+u_j;\boldsymbol{\omega})}, \label{d-plus}\end{equation}
где $\mathcal{C}=(n-1)g+g_0+g_1+g_2$ и $\Delta_-(\mathbf{u})=
\Delta_+(-u_1,\ldots,-u_n)$. Аналогично, мы определяем функции
\begin{eqnarray}\label{tilde-d}
\lefteqn{ \tilde \Delta_+(\mathbf{u})=\prod_{1\leq j<k\leq n}
\frac{S(u_j \pm u_k;\boldsymbol{\omega})} {S(g + u_j \pm
u_k;\boldsymbol{\omega})} }&&
\\ && \times
\prod_{j=1}^n\frac{S( 2u_j,
\mathcal{B}+\frac{\omega_1}{2}-u_j,\mathcal{C}+\frac{\omega_1}{2}-u_j;
\boldsymbol{\omega})}{S(\mathcal{C}+\frac{\omega_1}{2}+u_j,g_3-\frac{\omega_1}{2}+u_j,
g_4-\frac{\omega_1}{2}+u_j;\boldsymbol{\omega}) \prod_{n=0}^2
S(g_n+\frac{\omega_1}{2}+u_j; \boldsymbol{\omega})},
\nonumber\end{eqnarray}
и $\tilde\Delta_-(\mathbf{u})=\tilde\Delta_+(-u_1,\ldots,-u_n)$,
которые дают факторизацию ядра со сдвинутыми значениями параметров:
$$
\tilde \Delta_+(\mathbf{u})\tilde \Delta_-(\mathbf{u})
=\Delta(\mathbf{u};g,
g_0+\frac{\omega_1}{2},g_1+\frac{\omega_1}{2},
g_2+\frac{\omega_1}{2}, g_3-\frac{\omega_1}{2},
g_4-\frac{\omega_1}{2}).
$$
Предположим, что наши параметры таковы, что сдвинутые параметры
(и их аналоги, получающиеся сдвигами на $\pm\omega_2/2$) также
принадлежат области определения, указанной в формулировке теоремы.
Тогда справедливо следующее равенство
\begin{eqnarray}\nonumber
&& \int_{\mathbb{L}^n} \tilde
\Delta_+(u_1+\frac{\omega_1}{2},\ldots,u_n+\frac{\omega_1}{2})
\Delta_-(\mathbf{u})\, du_1\cdots du_n
\\ && \makebox[4em]{}
=\int_{\mathbb{L}^n}\tilde \Delta_+(\mathbf{u})
\Delta_-(u_1-\frac{\omega_1}{2},\ldots, u_n-\frac{\omega_1}{2})\,
du_1\cdots du_n, \label{integral-shift}\end{eqnarray}
получающееся после сдвига контура интегрирования $\mathbb{L}$ в левой
части на $-\omega_1/2$. Заметим, что такие сдвиги разрешены благодаря
отсутствию полюсов в полосе между $\mathbb{L}$ и $\mathbb{L}-\omega_1/2$
и экспоненциальному убыванию ядра на бесконечности. Действительно, отношение
$S(u;\boldsymbol{\omega})/S(g+u;\boldsymbol{\omega})$ голоморфно в полосе
$\{ u=s\omega_1+ix\omega_2\mid 0\leq
s\leq 1,\; -\infty <x<\infty\}$, а отношение
$$
\frac{S( 2u,\mathcal{B}+u,\mathcal{C}+u; \boldsymbol{\omega})}
{S(\mathcal{C}+\frac{\omega_1}{2}+u,g_3-\frac{\omega_1}{2}+u,
g_4-\frac{\omega_1}{2}+u;\boldsymbol{\omega}) \prod_{n=0}^2
S(g_n+\frac{\omega_1}{2}+u; \boldsymbol{\omega})}
$$
голоморфно в полосе $\{ u=s\omega_1+ix\omega_2\mid 0\leq
s\leq 1/2,\; -\infty <x<\infty\}$. Меняя знаки $u_j\to -u_j$ у переменных
интегрирования и суммируя по всем  $2^n$ комбинациям таких инверсий, из
равенства \re{integral-shift} получаем
\begin{eqnarray}\nonumber
&& \int_{\mathbb{L}^n} \rho(\mathbf{u};\mathbf{g})
\Delta_+(\mathbf{u})\Delta_-(\mathbf{u})\, du_1\cdots du_n
\\ && \makebox[4em]{}
=\int_{\mathbb{L}^n} \tilde\rho(\mathbf{u};\mathbf{g}) \tilde
\Delta_+(\mathbf{u})\tilde\Delta_-(\mathbf{u})\, du_1\cdots du_n,
\label{int-sum}\end{eqnarray}
где
\begin{eqnarray*}\label{rho}
&& \rho(\mathbf{u};\mathbf{g})=\sum_{\nu_\ell=\pm 1\atop
\ell=1,\ldots,n}
\frac{\tilde\Delta_+(\nu_1u_1+\frac{\omega_1}{2},\ldots,
\nu_nu_n+\frac{\omega_1}{2})}{\Delta_+(\nu_1 u_1,\ldots,\nu_n u_n)},  \\
&& \tilde\rho(\mathbf{u};\mathbf{g})=\sum_{\nu_\ell=\pm 1\atop
\ell=1,\ldots,n}
\frac{\Delta_-(\nu_1u_1-\frac{\omega_1}{2},\ldots, \nu_n
u_n-\frac{\omega_1}{2})} {\tilde\Delta_-(\nu_1 u_1,\ldots,\nu_nu_n)}.
\label{t-rho}\end{eqnarray*}
После ряда упрощений, функция $\rho(\mathbf{u};\mathbf{g})$ принимает вид
\begin{eqnarray}
  \rho(\mathbf{u};\mathbf{g}) &=& \sum_{\nu_\ell=\pm 1\atop \ell=1,\ldots,n}
\prod_{1\leq j<k\leq n}\frac{1-tz_j^{\nu_j}z_k^{\nu_k}}
{1-z_j^{\nu_j}z_k^{\nu_k}} \prod_{j=1}^n
\frac{(1-t^{n-1}t_0t_1t_2z_j^{-\nu_j})\prod_{n=0}^2(1-t_nz_j^{\nu_j})}
{1-z_j^{2\nu_j}} \nonumber\\
&=& \prod_{j=1}^{n}(1-t^{j-1}t_0t_1)(1-t^{j-1}t_0t_2)
(1-t^{j-1}t_1t_2),
\label{rho-sum} \end{eqnarray}
с $t=e^{2\pi ig/\omega_2}$, $t_n=e^{2\pi ig_n/\omega_2}$,
$z_k=e^{2\pi iu_k/\omega_2}$. Действительно, ядро суммы в \re{rho-sum}
инвариантно относительно перестановок $z_k$ и отражений $z_k\to z_k^{-1}$.
Произведение $\rho(\mathbf{u};\mathbf{g})$ и множителя
$$
\prod_{1\leq j<k\leq n}\frac{(1-z_jz_k)(1-z_jz_k^{-1})}{z_j}
\prod_{j=1}^n\frac{1-z_j^2}{z_j}
$$
дает Лорановский полином переменных $z_j,j=1,\ldots,n,$ антисимметричный
по отношению к обоим преобразованиям (отдельно, перестановкам $z_j$ и отражениям
$z_j\to z_j^{-1}$). Любой такой полином пропорционален указанному выше множителю.
Константа пропорциональности находится после наложения ограничения $z_j=t_0t^{n-j}$,
что оставляет только один член в сумме (со всеми $\nu_j=1$) равный выражению в
правой части в \re{rho-sum}.

Точно таким же образом можно получить для
$\tilde\rho(\mathbf{u};\mathbf{g})$, заменив в соотношении \re{rho}
$t_{0,1,2}$ на $t_{3,4}q^{-1/2},$ $t^{n-1}t_0t_1t_2q^{1/2}$,
\begin{eqnarray*}
\tilde\rho(\mathbf{u};\mathbf{g})&=&  \sum_{\nu_\ell=\pm 1\atop
\ell=1,\ldots,n} \prod_{1\leq j<k\leq
n}\frac{1-tz_j^{\nu_j}z_k^{\nu_k}} {1-z_j^{\nu_j}z_k^{\nu_k}}
\prod_{j=1}^n \Biggl(
\frac{(1-t_3q^{-1/2}z_j^{\nu_j})(1-t_4q^{-1/2}z_j^{\nu_j})}
{1-z_j^{2\nu_j}} \nonumber \\ && \makebox[6em]{}\times
(1-t^{n-1}t_0t_1t_2q^{1/2}z_j^{\nu_j})(1-Bq^{-1/2}z_j^{-\nu_j})\Biggr)
 \nonumber \\
&=& \prod_{j=1}^{n}(1-t^{j-1}t_3t_4/q)
(1-t^{1-j}B/t_3)(1-t^{1-j}B/t_4),
\end{eqnarray*}
где $B=e^{2\pi i\mathcal{B}/\omega_2}$.
Полученные выражения показывают, что интересующие нас функции не
зависят от переменных интегрирования $\mathbf{u}$, то есть
$\rho(\mathbf{u};\mathbf{g})=\rho(\mathbf{g})$ и
$\tilde{\rho}(\mathbf{u};\mathbf{g})=\tilde{\rho}(\mathbf{g})$.
Следовательно, мы можем вынести соответствующие множители за
знаки интегрирования и переписать \re{int-sum} как
\begin{equation*}
\int_{\mathbb{L}^n} \Delta_+(\mathbf{u})\Delta_-(\mathbf{u})\,
du_1\cdots du_n =\frac{\tilde\rho(\mathbf{g})}{\rho(\mathbf{g})}
\int_{\mathbb{L}^n} \tilde
\Delta_+(\mathbf{u})\tilde\Delta_-(\mathbf{u})\, du_1\cdots du_n,
\end{equation*}
и, следовательно,
\begin{equation}\label{int-trans}
\frac{{\cal N}(\mathbf{g})}
{{\cal N}(g,g_0+\frac{\omega_1}{2},g_1+\frac{\omega_1}{2},
g_2+\frac{\omega_1}{2},g_3-\frac{\omega_1}{2},g_4-\frac{\omega_1}{2})}=
\frac{\tilde\rho(\mathbf{g})}{\rho(\mathbf{g})}.
\end{equation}

В результате мы видим, что отношение левой и правой частей формулы
\re{q-ds-circle} инвариантно при сдвигах
$g_{0,1,2}\to g_{0,1,2}+\omega_1/2, g_{3,4}\to
g_{3,4}-\omega_1/2$, и, по симметрии, любых других сдвигах, получающихся
перестановкой индексов у параметров. Функция двойного синуса $S(u;\boldsymbol{\omega})$,
подынтегральная функция в \re{q-ds-circle} и само значение интеграла
${\cal N}(\mathbf{g})$ симметричны по отношению к перестановкам $\omega_1$ и $\omega_2$
\cite{kls:unitary}. Контур интегрирования $\mathbb{L}$
нарушает симметрию между $\omega_1$ и $\omega_2$, но преобразования, использовавшиеся в
\re{d-plus}-\re{int-trans}, носят чисто алгебраический характер и не зависят от
контура интегрирования. Поэтому, указанное отношение инвариантно также относительно
сдвигов $g_{0,1,2}\to g_{0,1,2}+\omega_2/2,$ $g_{3,4}\to g_{3,4}-\omega_2/2$
(со всеми перестановками индексов у параметров).

Благодаря аналитичности, мы можем менять контур интегрирования $\mathbb{L}$ на любой
другой, охватывающий то же множество полюсов, без изменения значения интеграла.
Для правильно подобранного контура, можно добиться инвариантности указанного
отношения при сдвигах $g_n\to g_n+ k\omega_1/2+
m\omega_2/2$ с произвольными $k,m\in\mathbb{Z}$.
Взяв $\omega_1,\omega_2 > 0$, мы можем выбрать подмножество этих значений
параметров, содержащих предельные точки для фиксированного контура интегрирования $\mathbb{L}$.
Более того, совершая последовательные $\omega_{1,2}/2$ сдвиги в различных направлениях,
можно избежать сильных деформаций контура интегрирования
(аналогичный прием уже использовался в третьей главе при доказательстве
эллиптических бета-интегралов). Благодаря указанному свойству инвариантности,
отношение левой и правой частей в \re{q-ds-circle} есть функция, зависящая
только от $\omega_{1,2}$ и $g$ и обозначаемая далее как $f(\omega_1,\omega_2,g)$.

Для того, чтобы увидеть, что $f(\omega_1,\omega_2,g)$ на самом деле равна единице,
необходимо провести анализ вычетов ядра интеграла, аналогично процедуре
приведенной в конце третьей главы. Для этого возьмем один из полюсов,
скажем, $u=g_0$ и потребуем, чтобы он из полуплоскости
$\mbox{Re}(u/\omega_2)>0$ пересек контур интегрирования $\mathbb{L}$.
При этом полюс $u=-g_0$ из полуплоскости $\mbox{Re}(u/\omega_2)<0$  также
пересечет $\mathbb{L}$ (из-за симметрии относительно отражения $u\to-u$).
Все остальные полюсы остаются в своих полуплоскостях слева или справа от
$\mathbb{L}$. Этого можно добиться правильным подбором параметров $g, g_1,\ldots,
g_4$, и $\omega_{1,2}$. Нетрудно видеть, что вычеты пересекающих $\mathbb{L}$ полюсов
по переменной $u_n$ имеют особенности в точках
$u_k=\pm(g_0+g), k=1,\ldots,n-1,$ (вместо $u_k=\pm g_0$) и они все еще
находятся справа от $\mathbb{L}$ при $\mbox{Re}((g_0+g)/\omega_2)<0$. Такие
сдвиги полюсов происходят всякий раз когда вычисляются вычеты. Обозначим
$\rho_k=g_0+(k-1)g$ и потребуем, чтобы
$$
\mbox{Re}\left(\frac{\rho_k}{\omega_2}\right)<0,
\quad k=1, \ldots, n,\qquad
\mbox{Re}\left(\frac{g_0+\omega_1}{\omega_2}\right),
\mbox{Re}\left(\frac{g_0+\omega_2}{\omega_2}\right)>0.
$$

Благодаря ограничениям наложенным на параметры $g, \omega_{1,2}$, мы получаем
более простую формулу для суммы интегралов и их вычетов, чем формула
работы \cite{die-spi:elliptic} описанная в конце третьей главы:
\be
\int_{\mathbb{L}^n_d } \Delta(\mathbf{u};\mathbf{g})
\frac{du_1}{\omega_2}\cdots\frac{du_n}{\omega_2}
=\sum_{m=0}^n 2^mm! \binom{n}{m} \int_{\mathbb{L}^{n-m}}
\mu_m(\mathbf{u}) \frac{du_1}{\omega_2}\cdots\frac{du_{n-m}}{\omega_2},
\label{res}\ee
где $\mathbb{L}_d$ обозначает деформацию контура $\mathbb{L}$, избегающую
пересечений с какими-либо полюсами при указанных изменениях $g_0$. Множитель
$2^m$ возникает из-за того, что вычеты возникают попарно и их значения совпадают
(из-за симметрии $u_k\to -u_k$ и разных ориентаций контуров, окружающих
полюсы с разных сторон $\mathbb{L}$). Множители $m!$ и $\binom{n}{m}$
связаны с числом упорядочений $m$ циклов и числом способов выбрать эти
циклы из $n$ доступных возможностей.

Функции вычетов имеют вид $\mu_0(\mathbf{u})=\Delta(\mathbf{u};\mathbf{g})$
и при $m>0$
\begin{equation}
\mu_m(\mathbf{u}) = \kappa_m
\delta_{m,n-m}(\mathbf{u})\Delta_{n-m}(\mathbf{u};\mathbf{g}),
\label{mu}\end{equation}
где $\Delta_{n-m}(\mathbf{u};\mathbf{g})$ получаются из ядра
\re{q-deg} при замене в нем $n$ на $n-m$, но тем же самым
$\mathcal{B}=(2n-2)g+\sum_{n=0}^4 g_n$. Другие коэффициенты равны
\begin{equation*}
\kappa_m =(-1)^m\frac{(\tilde q;\tilde
q)_\infty^m}{(q;q)_\infty^m} \prod_{1\leq j<k\leq m}
\frac{S(\pm\rho_k - \rho_j;\boldsymbol{\omega})}
{S(g\pm\rho_k-\rho_j;\boldsymbol{\omega})} \prod_{l=1}^m
\frac{S(-2\rho_l,\mathcal{B}\pm \rho_l;\boldsymbol{\omega})}
       {\prod_{n=1}^4S(g_n\pm\rho_l;\boldsymbol{\omega})},
\end{equation*}
и
\begin{equation*}
\delta_{m,n-m}(\mathbf{u}) = \prod_{\begin{subarray}{c}1\leq j\leq
m\\ 1\leq k\leq n-m\end{subarray}} \frac{S(\pm \rho_j\pm
u_k;\boldsymbol{\omega})} {S(g\pm\rho_j\pm
u_k;\boldsymbol{\omega})}.
\end{equation*}
Выражения для $\mu_m(\mathbf{u})$ получаются по индукции.
Действительно, вид $\mu_1(\mathbf{u})$ легко устанавливается после
использования соотношения
$$
\lim_{u\to \pm g_0}\frac{u\mp g_0}{S(g_0\mp
u;\boldsymbol{\omega})} =\pm \frac{\omega_2}{2\pi i}\frac{(\tilde
q;\tilde q)_\infty} {(q;q)_\infty}
$$
и факта, что контуры окружающие эти полюсы ориентированы по часовой стрелке
для верхних знаков и против часовой стрелки для нижних (это дает общий отрицательный
знак для $\kappa_1$).

Предположим, что $\mu_m$ определяется формулой \re{mu} для некоторого $m>1$.
Чтобы найти $\mu_{m+1}$, необходимо вычислить вычеты полюсов в точках
 $u_{n-m}=\pm \rho_{m+1}$. Простой расчет показывает, что
$$
\int_{c_m}
\mu_m(\mathbf{u})\frac{du_{n-m}}{\omega_2}=\mu_{m+1}(\mathbf{u}),
$$
где $c_m$ обозначает маленький контур, ориентированный по часовой стрелке
и окружающий полюс при $u_{n-m}=\rho_{m+1}$.

Благодаря аналитичности, наши изменения параметров и деформации контуров интегрирования
не меняют значения интеграла и поэтому правая сторона суммы в \re{res} равна
$f(\omega_1,\omega_2,g)\mathcal{N}(\mathbf{g})$. Мы теперь делим обе стороны
этого равенства на $\mathcal{N}(\mathbf{g})$ и берем предел $g_4\to -g_0-(n-1)g$.
При $m<n$, коэффициенты $\kappa_m(\mathbf{g})$, которые могут быть переписаны в виде
\begin{eqnarray*}
\lefteqn{ \kappa_m =(-1)^m\frac{(\tilde q;\tilde q)_\infty^m}
{(q;q)_\infty^m}\prod_{l=1}^m \Biggl(
\frac{S(g,(2-m-l)g-2g_0;\boldsymbol{\omega})}{S(lg;\boldsymbol{\omega})}
}
\\ && \times
\frac{S(2g_0+\sum_{r=1}^4g_r+(2n+l-3)g,
\sum_{r=1}^4g_r+(2n-l-1)g;\boldsymbol{\omega})}
{\prod_{r=1}^4S(g_r+g_0+(l-1)g,g_r-g_0-(l-1)g;
\boldsymbol{\omega})} \Biggr),
\end{eqnarray*}
не содержат расходящихся в этом пределе множителей, а интегралы,
которые умножаются на эти коэффициенты, остаются ограниченными.
Поэтому, выживают только члены с $m=n$ и, после простых выкладок, получаем
$$
\lim_{g_4\to
-g_0-(n-1)g}\frac{2^nn!\kappa_n}{\mathcal{N}(\mathbf{g})}=1,
$$
что означает, что $f(\omega_1,\omega_2,g)=1$. После доказательства
равенства \re{q-ds-circle} в указанной ограниченной области значений
параметров (в которой параметры, сдвинутые на $\pm\omega_{1,2}/2$,
удовлетворяют указанным ограничениям и $\omega_{1,2}>0$),
мы можем аналитически его продолжить на значения квазипериодов
$\omega_{1,2}$ и параметров $g, g_n$ указанные в формулировке теоремы.
Таким образом, теорема \ref{nr:thm} доказана.

Обратимся теперь к интегралу типа Аски-Вильсона \re{q-ds-circle}.
Условия его сходимости существенно отличаются
от предыдущего случая. Действительно,
\begin{equation*}
\frac{S(2i\omega_2x;\boldsymbol{\omega})}
{\prod_{m=0}^3S(g_m+i\omega_2x;\boldsymbol{\omega})}=
\begin{cases}
O(1) &\mbox{при}\; x\to +\infty \\
O(e^{2\pi x(1+\omega_2/\omega_1-
\sum_{m=0}^3 g_m/\omega_1 )}) &\mbox{при}\;x\to -\infty
\end{cases}.
\end{equation*}
Комбинируя эти предельные соотношения с асимптотикой для отношения
$S(i\omega_2x;\boldsymbol{\omega})/S(g+i\omega_2x;\boldsymbol{\omega})$,
можно увидеть, что подынтегральная функция остается ограниченной
в области интегрирования $\mathbb{L}^n$ и убывает экспоненциально быстро
на его бесконечностях, если потребовать
$\mbox{Re}((\mathcal{B}-\omega_2)/\omega_1)<1$, где
$\mathcal{B}=(2n-2)g+\sum_{n=0}^3g_n$.

Доказательство инвариантности отношения левой и правой частей равенства
(\ref{q-ds-circle}) относительно описанного сдвига параметров
опиралась только на алгебраических манипуляциях над ядром интеграла.
Поэтому мы можем повторить их для подынтегральной функции, получающейся
в пределе $\mbox{Im}(g_4/\omega_2),$ $\mbox{Im}(g_4/\omega_1)\to
+\infty$ (или $t_4\to 0$). При этом подынтегральная функция в
(\ref{q-ds-circle}) переходит в ядро для интеграла \re{q-ds-selberg}.
Поэтому, предельные аналоги равенств
\re{d-plus}--\re{int-trans} показывают, что отношение левой и правой частей
в равенстве \re{q-ds-selberg} не зависит от сдвигов параметров
$g_{0,1,2}\to g_{0,1,2}+\omega_{1,2}/2,$
$g_3\to g_3-\omega_{1,2}/2$ с учетом перестановок всех индексов.
Пользуясь опять аналитическим продолжением и анализом вычетов,
упрощенного подходящим образом, получаем, что рассматриваемое отношение
равно единице. Таким образом, мы также установили справедливость
формулы  интегрирования (\ref{q-ds-selberg}).

\subsection{Интегралы типа I}

В этом параграфе рассматриваются  $q$-редукции модифицированных эллиптических
бета-интегралов типа I. Они соответствуют пределу
$\mbox{Im}(\omega_3)\to \infty$ взятому таким образом, что $p, r\to 0$.
При этом модифицированная эллиптическая гамма-функция редуцируется
к функции двойного синуса и возникают интегралы по бесконечным контурам.
Это формальный предел и он требует строго обоснования.

\begin{theorem}
Для корневой системы $C_n$ используем те же обозначения что и в
(\ref{kernel-C-unit}) и определим ядро
\begin{eqnarray}\nonumber
&& \rho(u,g;\omega;C_n)=\prod_{1\leq i<j\leq n}S(\pm u_i\pm u_j)
\prod_{i=1}^n\frac{S(\pm 2u_i,\mathcal{A}\pm u_i)}
{\prod_{m=1}^{2n+3}S(g_m\pm u_i)}
\\ &&\makebox[8em]{}\times
\frac{\prod_{1\leq m<s\leq 2n+3}S(g_m+g_s)}
{\prod_{m=1}^{2n+3}S(\mathcal{A}-g_m)},
\label{kernel-C-q}\end{eqnarray}
где $\mathcal{A}=\sum_{m=1}^{2n+3}g_m$ и $S(u)\equiv S(u;\boldsymbol{\omega})$.
Предположим теперь, что $\mbox{Im} (\omega_1/\omega_2)\geq 0$ и
$\mbox{Re} (\omega_1/\omega_2)>0$ совместно с
$\mbox{Re} (g_m/\omega_2) > 0, \;
\mbox{Re}((\mathcal{A}-\omega_1)/\omega_2)<1
$
при $m=1,\dots,2n+3.$ Тогда
\begin{equation}
\int_{\L^n}\rho(u,g;\omega;C_n)
\frac{du_1}{\omega_2}\cdots \frac{du_n}{\omega_2}
= (-2)^n n!\frac{(\tilde q;\tilde q)_\infty^n}{(q;q)_\infty^n}
\label{C_n-q}\end{equation}
с контуром интегрирования $\L=i\omega_2\R$.
\end{theorem}
{\bf Доказательство.}
Рассмотрим сначала сходимость взятого интеграла. Для этого мы предположим, что
область интегрирования ограничена $n$-мерным кубом и рассмотрим интеграл по
переменной $u_k=i\omega_2x_k$ для какого-нибудь фиксированого $k$.
Благодаря симметрии $u_j\to-u_j$, достаточно рассмотреть подобласть
$0\leq x_j<\infty,\, j=1,\dots,n$.  Наложим ограничение
$\mbox{Re} (\omega_1/\omega_2)>0$, так что при $x_k\to+\infty$ имеем
$u_k/i\omega_1\to+\infty$. Из асимптотического поведения $S(u)$,
представленного в (\ref{asymp1}) и (\ref{asymp2}), получаем
$$
\frac{S(\pm 2u_k,\mathcal{A}\pm u_k)}
{\prod_{m=1}^{2n+3}S(g_m\pm u_k)}=O(e^{-2(n-1)\pi ix_k^2\omega_2/\omega_1
-2\pi n x_k(1+\omega_2/\omega_1)}).
$$
Используя формулу отражения
$$
S(u)S(-u)=e^{-\pi iB_{2,2}(u)}(1-e^{-2\pi iu/\omega_2})
(1-e^{-2\pi iu/\omega_1}),
$$
при $x_k\to+\infty$ получаем оценку
$$
\prod_{1\leq i<j\leq n}S(\pm u_i\pm u_j)=
O(e^{-2\pi i(n-1)\sum_{j=1}^nB_{2,2}(-u_j)}),
$$
не зависящую от значений других переменных $x_j$.
Таким образом, на границе области интегрирования имеем
$$
\rho(u,g;\omega;C_n)=O(e^{-2\pi (1+\omega_2/\omega_1)\sum_{j=1}^nx_j}),
$$
т.е. ядро убывает экспоненциально быстро и интеграл сходится.

После взятия предела $p\to 0$ в (\ref{eqn-C}), мы получаем
разностное уравнение для подынтегральной функции в нашем случае
(умножение $t_m$ и $z_i$ на $q^{\pm1}$ эквивалентно сдвигам
$g_m$ и $u_i$ на $\pm\omega_1$):
\begin{eqnarray} \nonumber
&& \rho(u,g_1+\omega_1,g_2,\ldots,g_{2n+3};\omega;C_n)-
\rho(u,g;\omega;C_n)
\\ && \makebox[4em]{}
=\sum_{i=1}^n\left(f_i(u_1,...,u_i-\omega_1,\ldots,u_n,g)-f_i(u,g)\right),
\label{eqn-C-q}\end{eqnarray}
где
\begin{eqnarray}\nonumber
&& f_i(u,g)=\rho(u,g;\omega;C_n)\prod_{j=1,\neq i}^n\frac{(1-t_1z_j)(1-t_1/z_j)}
{(1-z_iz_j)(1-z_i/z_j)}
\\ && \makebox[4em]{} \times
\frac{\prod_{j=1}^{2n+3}(1-t_jz_i)}{\prod_{j=2}^{2n+3}(1-t_1t_j)}
\frac{(1-t_1A)t_1}{(1-z_i^2)(1-Az_i)z_i}
\label{g-C-q}\end{eqnarray}
с $z_j=e^{2\pi iu_j/\omega_2},\, t_k=e^{2\pi ig_k/\omega_2}$,
$A=e^{2\pi i\mathcal{A}/\omega_2}.$ Поскольку $S(u)$ симметрична по $\omega_{1,2}$,
мы имеем аналогичное уравнение для сдвигов на $\omega_2$.

Предположим временно, что квазипериоды $\omega_{1,2}$ действительны, несоизмеримы
и $\omega_1<\omega_2$. Наложим также ограничение $\mbox{Re}(\mathcal{A})<\omega_1$.
Интегрируя (\ref{eqn-C-q}) по $u_j$ вдоль контура $\L$, получим
$$
I(g_1+\omega_1,g_2,\dots,g_{2n+3})
=I(g)\equiv\int_\L \rho(u,g;\omega;C_n)du_1\dots du_n.
$$
Это равенство следует из того факта, что при взятых ограничениях функции
$f_i(u,g)$ не имеют полюсов в полосе $-\omega_1<u_i<0$ и мы можем сдвинуть контур
$\L-\omega_1$ к $\L$. По симметрии, мы имеем
$I(g_1+\omega_2,g_2,\dots,g_{2n+3})=I(g).$ После подходящей деформации контура
интегрирования, можно добиться равенства
$I(g_1+j\omega_1-k\omega_2,g_2,\dots,g_{2n+3})=I(g)$ при всех $j,k\in\N$,
таких что $j\omega_1-k\omega_2\in[0,\omega_2]$. Благодаря наличию предельной точки
в этом множестве, можно заключить, что функция $I(g)$ не зависит от $g_1$ и,
таким образом, от всех $g_m$. Явный вид интеграла, дающийся выражением в правой части
(\ref{C_n-q}), может быть найден с помощью анализа вычетов аналогичного описанному
в работе \cite{die-spi:modular}. Наконец, по аналитичности мы расширяем область
допустимых значений параметров до пределов, указанных в формулировке теоремы.
\hfill{Q.E.D.}

\begin{theorem}
Для простейшего интеграла на корневой системе $A_n$ (интеграла типа I)
примем обозначения как в (\ref{kernel-A-unit}) и определим ядро
\begin{eqnarray}\nonumber
&& \rho(u,g,h;\omega;A_n)=
\prod_{i=1}^{n+1}\frac{S(F+H+u_i)\prod_{j=1}^{n+2}S(g_i+h_j)}
{\prod_{m=1}^{n+1}S(g_m-u_i)\prod_{j=1}^{n+2}S(h_j+u_i)\, S(H+g_i)}
\\ && \makebox[4em]{}
\times \prod_{1\leq i<j\leq n+1}S(u_i-u_j,u_j-u_i)\;
S(F)\prod_{j=1}^{n+2}\frac{S(H-h_j)}{S(F+H-h_j)}.
\label{kernel-A-q}\end{eqnarray}
Теперь предположим, что $\mbox{Im} (\omega_1/\omega_2)\geq 0$ и
$\mbox{Re} (\omega_1/\omega_2)>0$ совместно с
$$
\mbox{Re} (g_m/\omega_2) > 0, \quad \mbox{Re} (h_j/\omega_2) > 0,\quad
\mbox{Re}((F+H-\omega_1)/\omega_2)<1.
$$
Тогда
\begin{equation}
\int_{\L^n}\rho(u,g,h;\omega;A_n)
\frac{du_1}{\omega_2}\cdots \frac{du_n}{\omega_2}
= (-1)^n (n+1)!\frac{(\tilde q;\tilde q)_\infty^n}{(q;q)_\infty^n}
\label{A_n-q}\end{equation}
с контуром интегрирования $\L=i\omega_2\R$.
\end{theorem}
{\bf Доказательство.}
Сходимость интеграла анализируется аналогично предыдущему случаю.
Необходимое разностное уравнение на ядро получается из предела
$p\to 0$ в (\ref{eqn-A}):
\begin{eqnarray}\nonumber
&& \rho(u,g_1+\omega_1,g_2,\ldots,g_{n+1},h;\omega;A_n)-\rho(u,g,h;\omega;A_n)
\\ && \makebox[2em]{}
=\sum_{i=1}^n\left(f_i(u_1,\ldots,u_i-\omega_1,\ldots,u_n,g,h)-f_i(u,g,h)\right),
\label{eqn-A-q}\end{eqnarray}
где
\begin{equation}
\frac{f_i(u,g,h)}{\rho(u,g,h;\omega;A_n)}=
\prod_{j=1,\neq i}^{n+1}\frac{1-t_1/z_j}{1-z_i/z_j}
\prod_{j=1}^{n+2}\frac{1-z_is_j}{1-t_1s_j}
\frac{(1-z_iT/t_1)(1-TSt_1)t_1}{(1-T)(1-TSz_i)z_i}
\label{g-A-q}\end{equation}
с $z_j=e^{2\pi iu_j/\omega_2},\, t_k=e^{2\pi ig_k/\omega_2}$,
$s_j=e^{2\pi ih_j/\omega_2}$, и
$A=e^{2\pi i\mathcal{A}/\omega_2}.$ По симметрии,
мы имеем также аналогичное уравнение для сдвигов на $\omega_2$.

Предположим временно, что квазипериоды $\omega_{1,2}$ положительны
(скажем, $\omega_1<\omega_2$) и несоизмеримы. Наложим также ограничения
$\mbox{Re}(F+H)<\omega_2$ и $\mbox{Re}(g_m)>\omega_1,\, m=2,\dots,n+1$.
Интегрируя (\ref{eqn-A-q}) по $u_j$ вдоль контура $\L$, получаем
$I(g_1+\omega_1,g_2,\dots,g_{2n+3},h)$
$=I(g,h)\equiv\int_\L \rho(u,g,h;\omega;A_n)du_1\dots du_n$ (функции
$f_i(u,g,h)$ не имеет полюсов в полосе $-\omega_1<u_i<0$).

Аналогично $C_n$-случаю, для подходящих деформаций контура имеем
$I(g_1+j\omega_1-k\omega_2,g_2,\dots,g_{2n+3},h)=I(g,h)$ при всех $j,k\in\N$
таких что $j\omega_1-k\omega_2\in[0,\omega_2]$. В результате, функция
$I(g,h)$ не зависит от $g_1$ и, таким образом, от всех $g_m$.
Применяя анализ вычетов, аналогичный проведенному в \cite{spi:integrals},
находим, что функция $I(g,h)$ не зависит так же от $h$ и равна выражению в правой части
(\ref{A_n-q}). Аналитическое продолжение по параметрам позволяет расширить область
применимости формулы интегрирования \re{A_n-q} до указанных пределов.
\hfill{Q.E.D.}

Чисто гипергеометрический $A_n$  интеграл Густафсона \cite{gus:some1,gus:some3}
связан с обычной квантовой цепочкой Тоды (частное сообщение от С.М. Харчева).
Поэтому, естественно ожидать, что соответствующий предельный случай интеграла
(\ref{A_n-q}) связан с условием нормировки для собственных функций гамильтониана
квантовой многочастичной $q$-цепочки Тоды, исследованной в работе \cite{kls:unitary}
(см. также \cite{ols-rog}).

\begin{theorem}
Для $q$-редукции формулы интегрирования (\ref{D_n-circle}), примем те же обозначения
как и в (\ref{kernel-D-unit}) и определим ядро
\begin{eqnarray}\nonumber
\lefteqn{ \delta(u,g,h;\omega;A_n)=
\prod_{1\leq i<j\leq n+1}\frac{S(u_i-u_j,u_j-u_i)}{S(F-u_i-u_j)}
\prod_{j=1}^{n+3}\prod_{m=1}^n\frac{S(h_m+g_j)}{S(F+h_m-g_j)}
}&& \\ && \times
\prod_{i=1}^{n+1}\frac{\prod_{m=1}^n S(F+h_m-u_i)}
{\prod_{m=1}^n S(h_m+u_i)\prod_{j=1}^{n+3}S(g_j-u_i)}
\prod_{1\leq j<k\leq n+3}S(F-g_j-g_k).
\label{kernel-D-q}\end{eqnarray}
Предположим теперь, что $\mbox{Im} (\omega_1/\omega_2)\geq 0$ и
$\mbox{Re} (\omega_1/\omega_2)>0$ совместно с
$$
\mbox{Re} (g_j/\omega_2) > 0, \quad \mbox{Re} (h_m/\omega_2) > 0, \quad
\mbox{Re}((F+h_m-\omega_1)/\omega_2)<1.
$$
Тогда
\begin{equation}
\int_{\L^n}\delta(u,g,h;\omega;A_n)
\frac{du_1}{\omega_2}\cdots \frac{du_n}{\omega_2}
= (-1)^n (n+1)!\frac{(\tilde q;\tilde q)_\infty^n}{(q;q)_\infty^n}
\label{D_n-q}\end{equation}
с контуром интегрирования $\L=i\omega_2\R$.
\end{theorem}

Сходимость интеграла (\ref{D_n-q}) устанавливается так же как и в предыдущих
случаях. Остальные шаги в доказательстве этой теоремы аналогичны описанным в
$A_n$ и $C_n$ ситуациях. Используя уравнение для ядра, появляющееся
из $p\to 0$ предела в уравнении  построенном в \cite{spi-war:inversions},
можно показать, что левая часть (\ref{D_n-q}) не меняется после сдвигов
$g_m\to g_m+\omega_{1,2}$. Такие сдвиги могут итерироваться после
подходящей деформации контура интегрирования и это приводит к
независимости интеграла от параметров $g_m$. Соответствующий
анализ вычетов ($p\to 0$ аналог процедуры, использованной в
\cite{spi-war:inversions}) приводит к (\ref{D_n-q}). Подробные детали
опускаются в виду их повторного характера.

\section{Одномерные биортогональные функции с $|q|\leq 1$}

Опишем кратко биортогональные функции с мерой, фиксированной эллиптическим
бета-интегралом \re{circle-int} с $|q|\leq 1$ . Определим функции
\ba\nonumber
&& R_n(u;\omega_1,\omega_2,\omega_3)=
{_{12}V_{11}}\bigl( e^{2\pi i \frac{g_3-g_4}{\omega_2}},
e^{2\pi i \frac{\omega_1-g_0-g_4}{\omega_2}},
e^{2\pi i \frac{\omega_1-g_1-g_4}{\omega_2}},
e^{2\pi i \frac{\omega_1-g_2-g_4}{\omega_2}},
\\ && \makebox[4em]{}
e^{2\pi i \frac{g_3+u}{\omega_2}},
e^{2\pi i \frac{g_3-u}{\omega_2}},
e^{-2\pi i n\frac{\omega_1}{\omega_2}},
e^{2\pi i \frac{{\cal A}+(n-1)\omega_1-g_4}{\omega_2}};
e^{2\pi i \frac{\omega_1}{\omega_2}},
e^{2\pi i \frac{\omega_3}{\omega_2}} \bigr).
\label{brf-2}\ea
Очевидно, что параметры $g_k$ связаны с параметрами $t_k$,
использовавшимися в предыдущей главе следующим образом:
$t_k=e^{2\pi i \frac{g_k}{\omega_2}}$, а $z=e^{2\pi i \frac{u}{\omega_2}}$.
Случай $q\leftrightarrow p$ симметричных двухиндексных биортогональных
функций предыдущей главы соответствует $\omega_1\leftrightarrow\omega_3$
симметричному произведению функций \re{brf-2}:
$$
R_{nm}(u)=R_n(u;\omega_1,\omega_2,\omega_3)R_m(u;\omega_3,\omega_2,\omega_1).
$$

Случай модифицированных биортогональных функций, которые остаются хорошо
оп\-ределенными и при условии $\omega_1/\omega_2>0$ (ведущего к $|q|=1$),
соответствует $\omega_1\leftrightarrow\omega_2$ симметричному произведению
\re{brf-2}:
\be
R_{nm}^{mod}(u)=R_n(u;\omega_1,\omega_2,\omega_3)
R_m(u;\omega_2,\omega_1,\omega_3).
\label{unit-brf}\ee
Аналогичным образом определяются и партнеры этих функций
$T_{nm}^{mod}(u)$, такие что удовлетворяются соотношения двухиндексной
биортогональности
$
\langle T_{nl}|R_{mk}\rangle=h_{nl}^{mod}\delta_{mn}\delta_{kl},
$
где свертка $\langle \psi|\chi\rangle$ определяется интегралом
\re{circle-int} с правильно подобранным контуром интегрирования
$C^{mod}_{mn,kl}$, зависящим от индексов $n,l,m,k$, а
$h_{nl}^{mod}=h_n(\omega_1,\omega_2,\omega_3)h(\omega_2,\omega_1,\omega_3)$,
где множители $h_n(\omega_1,\omega_2,\omega_3)$ совпадают с \re{norm}.

Ключевое отличие описанных биортогональных функций от построенных
в предыдущей главе состоит в том, что, помимо отсутствия сингулярностей
при $|q|= 1$, они остаются хорошо определенными и в пределе
$\omega_3\to \infty$  взятым таким образом, что и $p\to 0$ и $r\to 0$.
При этом возникают функции
\ba\nonumber
&& R_{nm}^{mod}(u;\omega_1,\omega_2)=
{_{10}W_{9}}\bigl( e^{2\pi i \frac{g_3-g_4}{\omega_2}},
e^{2\pi i \frac{\omega_1-g_0-g_4}{\omega_2}},
e^{2\pi i \frac{\omega_1-g_1-g_4}{\omega_2}},
e^{2\pi i \frac{\omega_1-g_2-g_4}{\omega_2}},
\\ \nonumber && \makebox[8em]{}
e^{2\pi i \frac{g_3+u}{\omega_2}},
e^{2\pi i \frac{g_3-u}{\omega_2}},
e^{-2\pi i n\frac{\omega_1}{\omega_2}},
e^{2\pi i \frac{{\cal A}+(n-1)\omega_1-g_4}{\omega_2}}; q,q \bigr)
\\ \nonumber && \makebox[6em]{}
\times
{_{10}W_{9}}\bigl( e^{2\pi i \frac{g_3-g_4}{\omega_1}},
e^{2\pi i \frac{\omega_2-g_0-g_4}{\omega_1}},
e^{2\pi i \frac{\omega_2-g_1-g_4}{\omega_1}},
e^{2\pi i \frac{\omega_2-g_2-g_4}{\omega_1}},
\\ && \makebox[8em]{}
e^{2\pi i \frac{g_3+u}{\omega_1}},
e^{2\pi i \frac{g_3-u}{\omega_1}},
e^{-2\pi i m\frac{\omega_2}{\omega_1}},
e^{2\pi i \frac{{\cal A}+(m-1)\omega_2-g_4}{\omega_1}};
\tilde q^{-1}, \tilde q^{-1}\bigr).
\label{q-brf}\ea
которые удовлетворяют соотношениям двухиндексной биортогональности
для меры, определяемой $n=1$ случаем $q$-бета-интеграла \re{q-ds-circle}
доказанным в \cite{sto:hyperbolic}. Поскольку контур интегрирования бесконечен,
необходимо проверить условие сходимости интегралов в условиях
биортогональности. Первый из $_{10}\varphi_9$ факторов в \re{q-brf}
является рациональной функцией от
$$
\gamma(u;\omega_2)=\frac{(1-e^{2\pi i\frac{u+a}{\omega_2}})
(1-e^{2\pi i\frac{u-a}{\omega_2}})}
{(1-e^{2\pi i\frac{u+b}{\omega_2}})
(1-e^{2\pi i\frac{u-b}{\omega_2}})},
$$
где $a$ и $b$ калибровочные параметры, а второй $_{10}\varphi_9$
фактор рационален по $\gamma(u;\omega_1)$. Поэтому на бесконечности
эти функции постоянны и проблем со сходимостью не возникает.
Однако если попробовать выродить эти функции на уровень полиномов
Аски-Вильсона, то возникают серьезные трудности со сходимостью,
упомянутые в работах \cite{rui:int,sto:hyperbolic}.  Аналогичным образом
можно построить и $|q|=1$ партнеры эллиптических многомерных
биортогональных функций работы \cite{rai:abelian} и их $q$-вырождения.

%% file: CONCL.TEX
% Conclusions

\chapter*{Заключение}
\addcontentsline{toc}{chapter}{Заключение}

В диссертации получены следующие основные результаты:

\begin{enumerate}

\item  Построена одномерная эллиптическая бета-функция, зависящая от пяти
независимых параметров и двух базисных переменных.
Она определяет принципиально новый класс точно вычисляемых интегралов,
обобщающих бета-интеграл Эйлера, $q$-бета-интегралы Аски-Вильсона, Рахмана
и их различные предельные случаи.

\item Предложено и полностью доказано несколько многомерных эллиптических
бета-интегралов различных типов, связанных с системами корней $A_n$ и $C_n$.
Один из этих интегралов представляет собой эллиптическое обобщение интеграла
Сельберга, имеющего фундаментальное значение из-за большого числа приложений.

\item  Построена общая теория рядов и интегралов гипергеометрического типа
связанных с тета-функ\-ция\-ми Якоби, названных тета-ги\-пер\-гео\-ме\-три\-чес\-ки\-ми
функциями. Прояснено происхождение условий балансировки,
вполне уравновешенности и совершенной уравновешенности для функций
гипергеометрического типа с точки зрения условий эллиптичности,
содержащихся в эллиптических гипергеометрических функциях.

\item  Построена система непрерывных биортогональных функций одной переменной,
для которых одномерный эллиптический  бета-интеграл определяет меру.
Построены трехчленное рекуррентное соотношение и разностное уравнение
для этих функций.
Открыт новый тип соотношений ортогональности, названный двухиндексной
биортогональностью.
На настоящее время --- это самая общая система специальных функций
одной переменной, обобщающая классические ортогональные полиномы
Аски-Вильсона и биортогональные рациональные функции Рахмана.

\item Найдено эллиптическое обобщение цепочки Бэйли для рядов гипергеометрического
типа и с его помощью выведен ряд тождеств для эллиптических гипергеометрических рядов.
Впервые показано, что формализм цепочек Бэйли имеет интегральный
аналог. С помощью одномерного эллиптического бета-интеграла построено двоичное
дерево тождеств для многократных эллиптических гипергеометрических интегралов.

\item  Вычислена самая общая известная на настоящий момент обрывающаяся
цепная дробь, связанная с функциями гипергеометрического типа.
Она выражается через обрывающийся совершенно уравновешенный $_{12}E_{11}$
эллиптический гипергеометрический ряд.

\item
Исходя из сумм вычетов в многомерных эллиптических бета-ин\-те\-гра\-лах
выведены формулы суммирования многократных эллиптических гипергеометрических
рядов на корневых системах $A_n$ и $C_n$.

\item  Предложена модифицированная эллиптическая гамма-функция, которая хорошо
оп\-ре\-де\-ле\-на в случае когда один из базисных параметров лежит на единичной
окружности, $|q|=1$.  Найдены модифицированные эллиптические бета-интегралы, для
которых подынтегральная функция и результат точного вычисления выражаются через
эту обобщенную гамма-функцию.

\item Получен ряд новых результатов для $q$-ги\-пер\-гео\-ме\-три\-чес\-ких функ\-ций.
С помощью модулярно преобразованной $q$-гамма-функции, введен
новый класс $q$-функций Мейера, который был пропущен в предыдущих исследованиях.
Стро\-го доказаны многомерные $q$-бета-ин\-те\-гра\-лы, выражающиеся через функцию
двойного синуса и появляющиеся в формальном пределе $p\to 0$ из модифицированных
эллиптических бета-интегралов на корневых системах $A_n$ и $C_n$.

\item Построена (1+1)-мерная интегрируемая цепочка с дискретным временем
($R_{II}$-це\-поч\-ка),
которая обобщает обычную и релятивистскую цепочку Тоды с дискретным временем.
Она получается из условия совместности двух разностных
уравнений с рациональной зависимостью от спектрального параметра.
Найдена автомодельная редукция $R_{II}$-цепочки, порождающая
широкий класс эллиптических решений и приводящая к
биортогональным функциям выражающимся через совершенно уравновешенные
эллиптические гипергеометрические ряды.

\end{enumerate}

%\vglue 1cm
\bigskip

В заключении я хотел бы выразить признательность
моим соавторам С.О. Варнаару, Л. Вине, Я.Ф. ван Диехену, А.С. Жеданову,
И.М. Луценко и С.К. Скорику за плодотворное сотрудничество.
Стимулирующие обсуждения с Р. Аски, А. Берковичем, А.А. Болибрухом,
Дж. Гаспером, Д. Загиром, М. Исмаилом, К. Кратенталером, Ю.И. Маниным, А. Окуньковым,
М. Рахманом, С. Рюйсенаарсом, Э. Рэйнсом, Х. Розенгреном, С.К. Сусловым, В. Тарасовым,
С.М. Харчевым и Дж. Эндрюсом
были важны для изучения различных проблем изложенных в настоящей диссертации.
Я глубоко благодарен руководству и сотрудникам Лаборатории теоретической физики
им. Н.Н. Боголюбова, а также В.А. Матвееву и В.А. Рубакову
за интерес к моей работе и поддержку. Некоторые важные результаты, представленные в
диссертации, были получены при визитах
в Институт математики им. Макса Планка (г. Бонн, Германия) --- я выражаю искреннюю
благодарность дирекции этого института за создание благоприятных условий
для работы, творческую атмосферу и гостеприимство. Я также благодарен
Российскому фонду фундаментальных исследований за финансовую поддержку
исследований, отраженных в настоящей диссертации, на
протяжении ряда лет в рамках грантов 97-01-00281, 97-01-01041,
00-01-00299, 03-01-00780, 03-01-00781.

%% file: Refsfull.tex
%
% References
%

%%%%%%%%%%%%%%%%%%%%%%%%%%%%%%%%%%%%%%%%%%%%%%%